\documentclass{amsart}
\usepackage{hyperref}
\usepackage{newlfont}
\usepackage{graphicx}
\usepackage{indentfirst}
\usepackage{fancyhdr}
\usepackage{amsmath}
\usepackage{amsthm}
\usepackage{latexsym}
\usepackage[psamsfonts]{amssymb}
\usepackage{mathrsfs}
\usepackage{frcursive}
\usepackage
{lineno}
\usepackage{lipsum}
\usepackage{calligra}
\usepackage[T1]{fontenc}
\usepackage{wesa}
\usepackage{extarrows}
\usepackage{tcolorbox}
\usepackage[utf8]{inputenc}
\usepackage[english]{babel}
\usepackage{mathtools}
\usepackage{adjustbox,expl3,etoolbox}
\usepackage{marginnote}
\DeclareMathOperator{\Char}{Char}

\DeclareMathOperator{\ad}{ad}

\DeclareMathOperator{\dist}{dist}

\usepackage{accents}
\newcommand{\dbtilde}[1]{\accentset{\approx}{#1}}

\usepackage{adjustbox,expl3,etoolbox}

\letcs\replicate{prg_replicate:nn}

\newcommand*\longsum[1][1]{%
	\mathop{\textnormal{%
			\clipbox{0pt 0pt {.5\width} 0pt}{$\displaystyle\sum$}%
			\replicate{#1}{\clipbox{{.5\width} 0pt {.4\width} 0pt}{$\displaystyle\sum$}}%
			\clipbox{{.6\width} 0pt 0pt 0pt}{$\displaystyle\sum$}}}%
}
\usepackage{tikz}

\makeatletter
\@namedef{subjclassname@2020}{\textup{2020} Mathematics Subject Classification}
\makeatother

\allowdisplaybreaks
\newtheorem{theorem}{Theorem}
\newtheorem{proposition}{Proposition}
\newtheorem{lemma}{Lemma}
\newtheorem{corollary}{Corollary}
\newtheorem{definition}{Definition}
\newtheorem{example}{Example}
\newcounter{obsctr}

\newtheorem{remark}{Remark}

\renewcommand{\theequation}{\thesection.\arabic{equation}}

\DeclareMathAlphabet{\mathcalligra}{T1}{calligra}{m}{n}
\DeclareFontShape{T1}{calligra}{m}{n}{<->s*[2.2]callig15}{}

\begin{document}
\title[Microlocal regularity of Gevrey Vectors]{On the Microlocal Regularity of the Gevrey Vectors for second order partial differential operators with non negative
	characteristic form of first kind}
\author{Gregorio Chinni}
\address{5, Piazza di Porta San Donato, Dipartimento di Matematica, Universit\'a di Bologna, 40126 Bologna Italy}
%
\email[G. Chinni]{gregorio.chinni@gmail.com}
\author{Makhlouf Derridj}
\address{5, Rue de la Juvini\'ere, 78350 Les Loges en Josas, France}
\email[M. Derridj]{maklouf.derridj@numericable.fr}
\date{\today}
%
%
\begin{abstract}
	  We study the microlocal regularity of the analytic/Gevrey vectors 
	  for the following class of second order partial differential equations
      \begin{linenomath}
      	\begin{align*}
      	P(x,D) = \sum_{\ell,j=1}^{n} a_{\ell,j}(x) D_{\ell} D_{j} +  \sum_{\ell=1}^{n} i b_{\ell}(x) D_{\ell} +c(x), 
      	\end{align*}
      \end{linenomath}
      where $a_{\ell,j}(x) = a_{j,\ell}(x)$, $b_{\ell}(x)$, $\ell,j \in \lbrace 1,\dots,\, n\rbrace$,
      are real valued real Gevrey functions of order $s$ and  $c(x)$ is a Gevrey function of order $s$,
      $s \geq 1$, on $\Omega$ open neighborhood of the origin in $\mathbb{R}^{n}$.\\
      Thus providing a microlocal version of a result due to M. Derridj in \cite{D_2020}.
\end{abstract}
%
\maketitle
\tableofcontents
\section{Introduction}
\renewcommand{\theequation}{\thesection.\arabic{equation}}
\setcounter{equation}{0} \setcounter{theorem}{0}
\setcounter{proposition}{0} \setcounter{lemma}{0}
\setcounter{corollary}{0} \setcounter{definition}{0}
\setcounter{remark}{1}
%
This work follows the one we did in \cite{CD-22}, in the same subject,
where we dealt with H\"ormander's operators of the first kind
(or commonly known as ``sums of squares of vector fields"),
and considered the case of analytic vectors of operators with analytic coefficients.
For that we used the method by F.B.I. transform. 
In the present paper we consider second order partial differential operators of the form
\begin{linenomath}
	\begin{multline}\label{H-O-R_Op}
	P(x,D) = \sum_{\ell,j=1}^{n} a_{\ell,j}(x) D_{\ell} D_{j} +  \sum_{\ell=1}^{n} i  b_{\ell}(x) D_{\ell} +c(x) 
	\\
	\doteq P^{0}(x,D) + X(x,D) + c(x), 
	\end{multline}
\end{linenomath}
on $\Omega$, open neighborhood of the origin in $\mathbb{R}^{n}$,
where $a_{\ell,j}(x) $, $b_{\ell}(x)$, $\ell,j \in \lbrace 1,\dots,\, n\rbrace$,
are real valued Gevrey functions of order $s$ on $\Omega$,
the matrix $A(x)= \left( a_{\ell,j}(x)\right)$ is real symmetric,
$a_{\ell,j}(x) = a_{j,\ell}(x)$, and $A(x) \geq 0$ on $\Omega$
(i.e. $\langle A(x)v, v\rangle \geq 0$ for $v \in \mathbb{C}^{n}$, $x \in \Omega$)
and  $c(x)$ is a Gevrey function of order $s$
($s \geq 1$) on $\Omega$.
We recall that $s=1$ corresponds to the analytic case.\\
This class of second order operators,
with non-negative characteristic form,
was first studied by O.A.~Ole\u inik and  E.V.~Radkevi\v c in \cite{OR1971}.
Our purpose is to investigate the microlocal regularity of the
analytic-Gevrey vectors ($s\geq 1$) of $P$,
using the method of a priori estimates,
(idea developed in a preceding paper by the second author, \cite{D_2020})
in order to get suitable estimates of what we call microlocalized functions
associated to the function under study (see details in the next sections).\\
Since the work of  T. Kotake and M. Narasimhan (\cite{KotakeNarasimhan62}, 1962)
where they proved the so called ``Kotake-Narasimhan property", or ``iterates property",
for elliptic operators with analytic coefficients, an intensive investigation of this property was
undertaken by many mathematicians, along with its generalizations in different directions
and the use of more and more modern tools.
In the case of elliptic operators, iterates property was extended to the systems and
for $s$-Gevrey vectors
(see \cite{bcr_1986}, \cite{D_2017}, for surveys on this question, where there are many references).\\
In 1978, G. M\'etivier (\cite{Metivier78}) showed that, in the case of $s$-Gevrey vectors
with $s>1$, the ellipticity property is necessary for ``iterates property" to hold
(meaning: $s$-Gevrey vectors are in $s$-Gevrey class).
In the case of analytic vectors, M.S. Baouendi and G. M\'etivier showed
Kotake-Narasimhan property for hypoelliptic partial differential operators
of principal type with analytic coefficients (\cite{BaouendiMetivier-1982}, 1982).
In the case of systems of vector fields with analytic coefficients,
satisfying H\"ormander's condition, we mention two papers appeared in 1980,
where iterates property was shown (\cite{DamlakhiHellfer-1980} in case of analytic vectors,
and \cite{HelfferMattera1980} in the case called ``reduced analytic vectors").  
In the case of systems of complex vector fields
R. Barostichi, P. Cordaro and G. Petronilho (\cite{BCP_2011})
studied analytic vectors in locally integrable structures in 2011.\\
Concerning the case of second order partial differential operators,
the H\"ormander operators were mostly studied, after the famous
article on the hypoellipticity by L. H\"ormander, \cite{H67}.
As we are interested on iterates property, we do not write in other properties
like analytic or Gevrey hypoellipticity (local or microlocal).
The first result, on Gevrey regularity of analytic vectors, we mention is in global context,
for a subclass of ``sums of squares". It appeared in 2016 (\cite{bccj_2016})
and dealt with products of two tori (see also \cite{CC-16} for similar result in a different contest). 
The local version of that result for general H\"ormander's operators
was proved by the second author in (\cite{D_2019}, see also \cite{D_2018}),
shortly after, for operators with non-negative characteristic form
an analogous result was proved by the second author in (\cite{D_2020}),
result for which we give in this paper the microlocal version. \\
Let us finish this introduction with the mention of some results using intensively
the method of F.B.I. transform (and generalization of it as in \cite{Sj-Ast},
\cite{BH_2012}, \cite{Ber_Hai_2017}, \cite{HoepfnerMedrado2018}, \cite{Furdos22} )
and studying mainly, now, operators in more and more classes
of ultra-differentiable functions.  
\section{Notations, definitions, preliminary facts and main result}
\renewcommand{\theequation}{\thesection.\arabic{equation}}
\setcounter{equation}{0} \setcounter{theorem}{0}
\setcounter{proposition}{0} \setcounter{lemma}{0}
\setcounter{corollary}{0} \setcounter{definition}{0}
\setcounter{remark}{1}
In this section we recall the local and microlocal
\textit{H\"ormander-Ole\u{\i}nik-Radkevi\v c condition}, \textit{H.O.R.-condition} for shortness,
the definition of the \textit{type with respect to} $P$, where $P$ as in \eqref{H-O-R_Op}
(or more details on the subject see \cite{D_2020} and \cite{Rad2009})
and we  state the sub-elliptic estimate obtained in \cite{D_2020} in order to gain the local regularity of the Gevrey vectors
for $P$. It will be the starting point to obtain our main result.\\
We introduce the differential operators 
\begin{linenomath}
	\begin{align}\label{P_k_ud-ops}
	P^{k}(x,D) = 2 \sum_{\ell =1}^{n} a_{\ell,k}(x) D_{\ell} \quad\text{ and }\quad P_{k}(x,D) = \sum_{\ell,j =1}^{n}a_{\ell,j}^{(k)}D_{\ell}D_{j},   
	\end{align}
\end{linenomath}
of order $1$ and $2$ respectively, where $a_{\ell,j}^{(k)}(x) = D_{k}a_{\ell,j}(x)$.
From now below we will adopt the following convention:
the Latin alphabet letters in the upper index will denote the derivatives with respect
the corresponding direction, i.e. $a^{(k)}(x) = D_{k}a(x)$ (as above),
and the Greek alphabet letters in the upper index will denote the usual multi-index derivatives,
i.e. $a^{(\alpha)}(x) = D^{\alpha} a(x)= D_{1}^{\alpha_{1}} \cdots D_{n}^{\alpha_{n}} a(x)$, $\alpha= (\alpha_{1}, \dots, \alpha_{n}) \in \mathbb{Z}_{+}^{n}$.
We denote by
\begin{linenomath}
	\begin{align*}
	&\hspace{-7em}p^{k}(x,\xi) = \partial_{_{\xi_{k}}} p^{0}(x,\xi) \quad (k=1, \, \dots, n),
	\\
	&\hspace{-7em}p^{k}(x,\xi) = (1+|\xi|^{2})^{\frac{1}{2}}D_{k-n} p^{0}(x,\xi) \quad (k =  n+1, \dots, 2n),
	\end{align*} 
\end{linenomath}
%
where $p^{0}(x,\xi) $ is the principal symbol of $P(x,D)$.\\
%
%
%
%
We  associate to the operator $P(x,D)$, (\ref{H-O-R_Op}), the family 
\begin{linenomath}
	$$
	\mathscr{P}=\lbrace p^{1}(x,\xi), \dots, p^{2n}(x,\xi)\rbrace
	$$
\end{linenomath}
of homogeneous symbols of order $1$.\\
We recall that the \textit{Poisson bracket} of two symbols $q_{1}(x,\xi)$ and $q_{2}(x,\xi)$ is defined by
\begin{linenomath}
	\begin{align*}
	\left\{q_{1},q_{2} \right\}(x,\xi) = \sum_{j=1}^{n} \left( \frac{\partial q_{1}}{\partial \xi_{j}} \frac{\partial q_{2}}{\partial x_{j}} 
	- \frac{\partial q_{1}}{\partial x_{j}} \frac{\partial q_{2}}{\partial \xi_{j}} \right)\left( x, \xi\right).
	\end{align*} 
\end{linenomath}  
Let $ I = \left( i_{1}, \dots, i_{r}\right) $ with $i_{\ell} \in \lbrace 1,\dots, 2n \rbrace$. We denote by $|I|=r $ its length;
we define
\begin{linenomath}
	\begin{align*}
	\begin{array}{ll}
	\hspace{-8.5em} 1)\,\, p^{I}(x,\xi) = p^{i}(x,\xi) \text{ if } I= i, \, |I|=1,\, i\in \lbrace 1,\dots, \, 2n \rbrace,\vspace*{0.6em} \\
	%
   \hspace{-8.5em} 	2)\,\, p^{I}(x,\xi) = \left\{p^{J}, p^{i_{r}}\right\}(x,\xi) \text{ if } J= \left(i_{1}, \dots,\, i_{r-1}\right), \, |J|=r-1,
	\end{array}
	%
	\end{align*} 
\end{linenomath}
where $p^{i_{\ell}}(x,\xi) \in \mathscr{P}$.
We remark that since $p^{i}(x,\xi)$ are homogeneous symbols of degree $1$ then $p^{I}(x,\xi)$ are homogeneous symbols of
degree $1$ for every $I$. We remark that in our convention if, for example $p$ and $q$ are homogeneous of degree 1
and $\{p,q\}= 0$, then, here, $0$ is considered as homogeneous of degree 1. 
%
\begin{definition}[\cite{D_2020}, \textit{H.O.R.-condition}]
	\hspace{4em}
	%
	\begin{itemize}
		\vspace{-2em}
		\item[(i)] Let $(x,\xi) \in \Omega \times \mathbb{R}^{n}\setminus\lbrace 0 \rbrace$ then the H.O.R.-condition is satisfied
		at $(x,\xi)$  if there exists $I= \left(i_{1}, \dots,\, i_{r} \right)$ such that $ p^{I}(x,\xi) \neq 0$. 
		\item[(ii)] Let $x\in \Omega$, then the H.O.R-condition is satisfied at $x$ if for every $\xi \in \mathbb{R}^{n}\setminus \lbrace 0 \rbrace$
		the H.O.R.-condition is satisfied at $\left( x, \xi\right)$.
	\end{itemize}	             
	%
\end{definition}
%
\begin{definition}[\cite{D_2020}, \textit{type with respect to} $P$]\label{Def_Type}\hspace{5em}
	%
	\begin{itemize}
		\vspace{-2em}
		\item[(i)] Let $(x,\xi) \in \Omega \times \mathbb{R}^{n}\setminus\lbrace 0 \rbrace$ such that the H.O.R.-condition is satisfied
		at $(x,\xi)$ then
		\begin{linenomath}
			\begin{align}
			\tau\left( \left(x,\xi\right); \, \mathscr{P} \right) = \inf \left\{ |I|\,:\, p^{I}(x,\xi)\neq 0 \right\}
			\end{align}
		\end{linenomath}
		is the type with respect to $\mathscr{P}$ at $(x,\xi)$.  Otherwise, $ \tau\left( \left(x,\xi\right); \, \mathscr{P} \right) = +\infty$.
		\item[(ii)] Let $x\in \Omega$. If the H.O.R-condition is satisfied at $x$ then
		\begin{linenomath}
			\begin{align}
			\tau\left( x; \, \mathscr{P} \right) \doteq \sup \left\{ \tau\left( \left(x,\xi\right); \, \mathscr{P} \right);\, \xi \in \mathbb{R}^{n}\setminus\lbrace 0 \rbrace \right\} < +\infty,
			\end{align}
		\end{linenomath}
		is the type of $x$ with respect to $\mathscr{P}$. Otherwise, $ \tau\left( x; \, \mathscr{P} \right) = +\infty$.
	\end{itemize}	             
	%
\end{definition}
Taking advantage from the Proposition 3.1 in \cite{D_2020} and from the Proposition 1.5 in \cite{BCN-82} 
the second author obtained the following basic estimate.   
\begin{theorem}\label{D-Est}
	Let $P(x,D)$ be as in (\ref{H-O-R_Op}). Let $\Omega_{1}$ be open relatively compact in $\Omega$, $\overline{\Omega}_{1} \Subset \Omega$.
	Assume $\tau\left(\Omega_{1}, \mathscr{P}\right) $ finite, then for
	$\sigma = \left(\tau\left(\Omega_{1}, \mathscr{P}\right)\right)^{-1}$, $\mathscr{P} = \lbrace p^{1}, \dots, p^{2n}\rbrace$,
	there exists a positive constant $C$ such that 
	\begin{linenomath}
		\begin{align}\label{Der-Est}
		\|v\|_{\sigma}^{2}+ \sum_{j=1}^{n} \left( \| P^{j}v\|^{2}_{0} + \|P_{j}v\|_{-1}^{2}\right)
		\leq C \left( \sum_{m=0}^{n} | \langle E_{m}P v, E_{m}v\rangle | + \|v\|_{0}^{2}\right), \,\,\, \forall v \in \mathscr{D}(\Omega_{1}),
		\end{align} 
	\end{linenomath}
	$\| \cdot\|_{0}$ denotes the norm in $L^{2}\left( \Omega_{1}\right)$, $\|\cdot\|_{s}$ the Sobolev norm of order $s$,
	$ P^{k} $ and $ P_{k} $ as in (\ref{P_k_ud-ops}), $E_{0}=1$ and $E_{m}=D_{m} \psi \Lambda_{-1}$, $m=1,\dots, n$.
	Here $\psi$ belongs to $\mathscr{D}\left(\Omega\right)$ and is identically one
	on $\Omega_{1}$ and $\Lambda_{-1}$ is the pseudodifferential operator associated to the symbol $\lambda(\xi)^{-1}\doteq \left(1+|\xi|^{2}\right)^{-1/2}$.
\end{theorem}
We recall the local notion of Gevrey vectors.
\begin{definition}\label{D_Gs_V}
	Let $P(x,D)$ a differential operator of order $m$ with Gevrey coefficients of order $s \geq 1$ in $\Omega$ open subset of $\mathbb{R}^{n}$.
	We denote by $G^{s}(\Omega; P)$, the space of the Gevrey vectors of order $s \geq 1$, in $\Omega$, with respect of $P$ i.e.: the set of all distributions 
	$u \in \mathscr{D}'\left( \Omega \right)$ such that for any compact subset $K$ of $\Omega$ and every $N\in \mathbb{N}$, $P^{N}u$ is in $L^{2}(K)$ and
	there is a positive constant $C_{K}$ such that
	\begin{linenomath}
		\begin{align}\label{G_s-vectors}
		\|P^{N}u\|_{L^{2}(K)} \leq C_{K}^{N+1} ((mN!))^{s}, 
		\end{align}
	\end{linenomath}
	When $s=1$ we set $G^{1}(\Omega; P)= \mathscr{A}(\Omega; P)$ the set of the analytic vectors with respect to $P$ in $\Omega$.
\end{definition}
We recall the notion of Gevrey wave front set, in the spirit of H\"ormander, \cite{Horm71}.
\begin{definition}\label{D_WF_s}
	Let $x_{0} \in \Omega \subset \mathbb{R}^{n}$, $\xi_{0} \in \mathbb{R}^{n}\setminus\{0\}$ and $u\in \mathscr{D}'(\Omega)$.
	We say that $(x_{0},\xi_{0}) \notin WF_{s}(u)$, $s\geq 1$, if and only if there are an open neighborhood $U$ of $x_{0}$,
	an open cone $\Gamma$ around $\xi_{0}$ and a bounded sequence $u_{N} \in \mathcal{E}'(\Omega)$
	which is equal to $u$ in $U$ such that
	\begin{linenomath}
		\begin{equation}\label{WF_s}
		|\widehat{ u}_{N}(\xi)| \leq C^{N+1} N^{sN} \left(1+ |\xi|\right)^{-N}, \qquad\, N=1,2,\dots,
		\end{equation}
	\end{linenomath} 
	is valid for some constant $C$, independent of $N$, when $\xi\in \Gamma$. 
\end{definition}
We state now the main result of the paper

\begin{theorem}\label{M-Th}
	Let $P(x,D)$ be as in (\ref{H-O-R_Op}) and $u \in G^{s}(\Omega; P)$. 
	Let $(x_{0}, \xi_{0})$ be a point in the characteristic variety of $P(x,D)$
	such that $ \tau\left( \left(x_{0},\xi_{0}\right); \, \mathscr{P} \right) =r$
	then  $(x_{0},\xi_{0}) \notin WF_{rs}(u)$.
\end{theorem}

\subsection{Remark on the case when $P$ is of H\"ormander type}~\par
\vskip-6mm
The H\"ormander's operators of the first kind are a subclass of the operators studied.\\
Let  $X_{1}(x,D), \, \dots, \, X_{m}(x,D)$ be vector fields with real-valued $s$-Gevrey coefficients
on $\Omega$, open neighborhood of the origin in $\mathbb{R}^{n}$.
Let $P_{H}$ denote the corresponding sum of squares operator
\begin{linenomath}
\begin{align}\label{Op_PH}
P_{H}(x,D) = \sum_{j=1}^{m}X_{j}^{2}(x,D) +X_{0}(x,D) + c(x),
\end{align} 
\end{linenomath} 
where $X_{0}$ is a linear combination with $s$-Gevrey coefficients in $\Omega$ of the vector fields
$X_{1}(x,D), \, \dots, \, X_{m}(x,D)$, and the H\"ormander's condition is satisfied by the system $\left\{ X_{1}(x,D), \, \dots, \, X_{m}(x,D)\right\}$.\\
Let $X_{j}(x,\xi)$ the symbol of the vector field $X_{j}$. 
Let $ I = \left( i_{1}, \dots, i_{\nu}\right) $ with $i_{\ell} \in \lbrace 1,\dots, m \rbrace$.
We denote by $|I|=\nu $ its length;
we define
\begin{linenomath}
	\begin{align*}
	\begin{array}{ll}
	\hspace{-8.5em} X^{I}(x,\xi) = X_{i}(x,\xi) \text{ if } I= i, \, |I|=1,\, i\in \lbrace 1,\dots, \, m \rbrace, \\[2pt]
	\hspace{-8.5em} 
	X^{I}(x,\xi) = \left\{X^{J}, X_{i_{\nu}}\right\}(x,\xi) \text{ if } J= \left(i_{1}, \dots,\, i_{\nu-1}\right), \, |J|=\nu-1.
	\end{array}
	%
	\end{align*} 
\end{linenomath}
\vspace*{0.4em}
\begin{definition}\label{Def_TypeH}
	Let $(x,\xi) \in \Omega \times \mathbb{R}^{n}\setminus\lbrace 0 \rbrace$ then the H\"ormander-condition
	is satisfied at $(x,\xi)$  if there exists $I= \left(i_{1}, \dots,\, i_{\nu} \right)$ such that $ X^{I}(x,\xi) \neq 0$.\\ 
	Let $(x,\xi) \in \Omega \times \mathbb{R}^{n}\setminus\lbrace 0 \rbrace$ such that 
	the H\"ormander-condition is satisfied at $(x,\xi)$ then
	\begin{linenomath}
		\begin{align}
		\tau\left( \left(x,\xi\right); \, X \right) = \inf \left\{ |I|\,:\, X^{I}(x,\xi)\neq 0 \right\}
		\end{align}
	\end{linenomath}
	is the type with respect to the system $\{ X_{1},\dots X_{m}\}$ of $(x,\xi)$.
\end{definition}
Setting $X_{j} (x,D)= \sum_{\ell=1}^{n} \widetilde{a}_{\ell,j}(x)D_{\ell}$ then
\begin{linenomath}
	$$
	P_{H}= \longsum[9]_{ \ell_{1},\ell=1}^{n} a_{\ell_{1},\ell}(x)D_{\ell_{1}}D_{\ell} + i \sum_{ \ell=1}^{n}b_{\ell}(x)D_{\ell} + c(x),
	$$
\end{linenomath}
where $ a_{\ell_{1},\ell}(x) = \sum_{j=1}^{m} \widetilde{a}_{\ell_{1},j}(x)\widetilde{a}_{\ell,j}(x)$.\\
Moreover, we have  
\begin{linenomath}
	\begin{align*}
	\begin{array}{ll}
	\hspace{-5.5em}
	 \left(P_{H}^{0}\right)^{k}= 2\displaystyle\sum_{\ell=1}^{n} a_{k,\ell}D_{\ell},
	 \text{ where }  a_{k,\ell}=\sum_{j=1}^{m} \widetilde{a}_{k,j}\widetilde{a}_{\ell,j},\\[4pt]
	\hspace{-5.5em}
	\left(P_{H}^{0}\right)_{k}=\displaystyle\longsum[9]_{ \ell_{1},\ell=1}^{n} a_{\ell_{1},\ell}^{(k)}D_{\ell_{1}}D_{\ell}
	\text{ where }  
	a_{\ell_{1},\ell}^{(k)} = \sum_{j=1}^{m}\left( \widetilde{a}_{\ell_{1},j}^{(k)}\widetilde{a}_{\ell,j} + \widetilde{a}_{\ell_{1},j}\widetilde{a}_{\ell,j}^{(k)}\right),
	\end{array}
	%
	\end{align*} 
\end{linenomath}
where $P_{H}^{0}$ is the leading part of $P_{H}$.\\
Furthermore, looking at the symbols $X_{j}(x,\xi)$ of the $X_{j}$'s, then 
the principal symbol $p^{0}$ of $P_{H}$ is
\begin{linenomath}
	$$
	p^{0}(x,\xi)= \sum_{j=1}^{m} X_{j}^{2}(x,\xi).
	$$
\end{linenomath}
Then
\begin{linenomath}
	\begin{align*}
	\begin{array}{ll}
	\hspace{-2.5em}
	p^{k}(x,\xi)=2 \displaystyle\sum_{j=1}^{m} X_{j}^{k}(x,\xi)X_{j}(x,\xi),
	\text{ where }  X_{j}^{k}(x,\xi)=\partial_{\xi_{k}}X_{j}(x,\xi) =\widetilde{a}_{k,j}(x),\\[4pt]
	\hspace{-2.5em}
	p_{k}=2\displaystyle\displaystyle\sum_{j=1}^{m} \left(X_{j}\right)_{k}(x,\xi)X_{j}(x,\xi)
	\text{ where }  
	\left(X_{j}\right)_{k}(x,\xi) = D_{k} X_{j}(x,\xi).
	\end{array}
	%
	\end{align*} 
\end{linenomath}
The following results hold
\begin{proposition}[\cite{D_2020}]\label{D-P1}
	For any $(x_{0},\xi_{0}) \in \Omega\times \mathbb{R}^{n}\setminus\{0\}$,
	one has 
	\begin{linenomath}
		$$
		\tau\left( \left(x,\xi\right); \, X \right) \leq 	\tau\left( \left(x,\xi\right); \, \mathscr{P}  \right).
		$$
	\end{linenomath}
\end{proposition}
\begin{proposition}[\cite{D_2020}]\label{D-P2}
	For any $(x_{0},\xi_{0}) \in \Omega\times \mathbb{R}^{n}\setminus\{0\}$,
	assume that 
	\begin{linenomath}
	$$
	\tau\left( \left(x,\xi\right); \, X \right) \leq 2\, \text{ then } 
	\tau\left( \left(x,\xi\right); \, X \right) =	\tau\left( \left(x,\xi\right); \, \mathscr{P}  \right).
	$$
	\end{linenomath}
\end{proposition}
Let us show two examples in $\mathbb{R}^{2}$ and $\mathbb{R}^{3}$, where we compare the type
$\tau\left( \left(x,\xi\right); \, X \right) $ and $\tau\left( \left(x,\xi\right); \, \mathscr{P}  \right)$.
\begin{example}
	See Proposition 4.8 in \cite{D_2020}.
\end{example}

\begin{example}
	Let us consider in $\Omega$, open neighborhood of the origin in $\mathbb{R}^{2}$,
	the vector fields $X_{1}=  D_{x} $, $X_{2} =  x^{2n+1}D_{y}$
	and $ X_{3} =  x^{n}y^{m} D_{y}$, $n, m \, \in \mathbb{N}$ and $m\geq 1$, and the associated sum of square operator
	\begin{linenomath}
	 $$
	 P_{1}= X_{1}^{2} + X_{2}^{2} + X_{3}^{2}.
	 $$
	\end{linenomath}
	Then $\tau\left(\rho,X\right) = 2n+2$, $\tau\left( \rho; \, \mathscr{P}  \right)= 2(2n+1)$, where $\rho= (0,0;0,1)$.\\
	Denote by $(\xi,\eta)$ the dual variable of $(x,y)$, then
	\begin{linenomath}
		\begin{align*}
		\begin{array}{ll}
		\hspace{-13.5em} p^{0}(x,y,\xi,\eta) = \xi^{2}+ x^{4n+2} \eta^{2} + x^{2n} y^{2m}\eta^{2}, \\[3pt]
		\hspace{-13.5em} 
		p^{1}(x,y,\xi,\eta) = 2\xi, \\[3pt]
		\hspace{-13.5em} 
	    p^{2}(x,y,\xi,\eta) = 2 x^{4n+2}\eta  + 2 x^{2n} y^{2m}\eta ,\\[3pt]
		\hspace{-13.5em} 
		p_{1}(x,y,\xi,\eta) = (4n+2) x^{4n+1}  \eta^{2} +2n x^{2n-1}  y^{2m}\eta^{2}, \\[3pt]
		\hspace{-13.5em} 
		p_{2}(x,y,\xi,\eta) = 2m x^{2n} y^{2m-1}\eta^{2}.
		\end{array}
		%
		\end{align*} 
	\end{linenomath}
	We have
	\begin{linenomath}
		$$
		\left(\ad \, p^{1}\right)^{4n+1}\left(p_{1}\right)(\rho)\neq 0, \text{ so } \tau\left( \rho; \, \mathscr{P}  \right) \leq 2(2n+1).
		$$
	\end{linenomath}
	First we observe that in $p^0$ and $p^2$, the terms in which appear only the factors $x$ and $\eta$,
	they vanish at order greater or equal to $2(2n+1)$ with respect of $x$.
	On the other hand, if we look to the terms that have $y$ as factor,
	any time that we remove a power of $y$ via the poisson bracket, the power of $x$,
	in such terms, grows at least of order $4n+1$. 
	We conclude that  $ \tau\left( \rho; \, \mathscr{P}  \right) = 2(2n+1)$.\\
	The fact that $\tau\left(\rho,X\right) = 2n+2$ is elementary.\\
	(For more details on the study of the optimal regularity of the solution of the problem $P_{1}u=0$ see \cite{chinni23-1}
	and \cite{chinni23-2}, see also \cite{BC-2022} for similar discussion in global setting.)
\end{example}
\vspace{-0.5em}
We recall the result obtained in \cite{CD-22}.
\begin{theorem}\label{Th-CD22}
	Let $P_{H}$ be as in \eqref{Op_PH}. Assume that the coefficients of the vector fields are analytic.
	Let $u$ be an analytic vector for $P_{H}$. Let $(x_{0}, \xi_{0}) $ be a point in the characteristic set
	of $P_{H}$ and $\nu= \tau\left( (x_{0}, \xi_{0}),X\right)$. Then $ (x_{0}, \xi_{0}) \notin WF_{\nu}(u)$.
\end{theorem}

\begin{remark}
	In the analytic category, due to the Proposition \ref{D-P1} the microlocal regularity obtained in Theorem \ref{Th-CD22}
	is, in general, better than the one obtained in the Theorem \ref{M-Th};
	the results match always only when $P_{H}$ vanishes ``exactly" of order 2, $ \tau\left( \rho,X\right)=2$ (Proposition \ref{D-P2}),
	see Examples 1 and 2. 
	So we have optimality of the Theorem \ref{M-Th}, by this exception,
	and the example given in \cite{bccj_2016}, where the type $\tau\left(\rho,X\right) = 2$.
	As $\tau\left(\rho,X\right) = \tau\left( \rho; \, \mathscr{P}  \right)$ in this case,
	this example shows the optimality of our result. 
\end{remark}
%
\section{Basic microlocal estimate for \textit{H\"ormander-Ole\u{\i}nik-Radkevi\v c} operators}
\renewcommand{\theequation}{\thesection.\arabic{equation}}
\setcounter{equation}{0} \setcounter{theorem}{0}
\setcounter{proposition}{0} \setcounter{lemma}{0}
\setcounter{corollary}{0} \setcounter{definition}{0}
\setcounter{remark}{1}
Due to the Proposition 3.1 in \cite{D_2020} and Proposition 1.5 in \cite{BCN-82}
we have the microlocal version of Theorem \ref{D-Est}:
%
\begin{theorem}\label{MicLocEst}
Let $(x_{0}, \xi_{0})$ be a point in the characteristic variety of $P(x,D)$, $\Char(P)$.
Let $ r \doteq \tau\left( \left(x_{0},\xi_{0}\right); \, \mathscr{P} \right)$, Definition \ref{Def_Type}, then the following
estimate holds
\begin{linenomath}	
	\begin{multline}\label{Der-Est-M}
    \| \text{\textcursive{p}}\, v\|_{\frac{1}{r}}^{2}+ \sum_{j=1}^{n} \left( \,\| P^{j}v\|^{2} + \|P_{j}v\|_{-1}^{2}\right)
    \\
	\leq C \left( \sum_{\ell=0}^{n}\left|\langle E_{\ell}P v, E_{\ell}v\rangle\right|  + \|v\|^{2}\right), \,\, \forall v \in \mathscr{D}(\Omega_{4}),
	\end{multline} 
\end{linenomath}
where $E_{0} = 1 $, $E_{m}=D_{m} \psi \Lambda_{-1}$, $m=1,\dots, n$,
$\psi$ belongs to $\mathscr{D}\left(\Omega\right)$
and it is identically one on $\Omega_{4}$\footnote{The use of the lower index 4, will be more clear later where we will introduce $\Omega_{j}$, $j=0,1,2,3$.}, $\Omega_{4} \Subset \Omega$,
$\Lambda_{-1}$ is the pseudo-differential operator associated
to the symbol $\lambda(\xi)^{-1}\doteq \left(1+|\xi|^{2}\right)^{-1/2}$
and $\text{\textcursive{p}} \,(x,D)$ is a zero order pseudo-differential operator,
elliptic at $(x_{0}, \xi_{0})$ and vanishing for $x\notin \widetilde{\Omega}_{3}$,
$\widetilde{\Omega}_{3} $ open neighborhood of $x_{0}$ with $\widetilde{\Omega}_{3} \Subset \Omega_{4}$.
\end{theorem}
\vspace*{-1.3em}
We will assume, without loss of generality, that $\text{\textcursive{p}} \,(x,\xi)$,
the symbol associated to $\text{\textcursive{p}} \,(x,D)$,
is elliptic in $\widetilde{\Omega}_{3} \times \Gamma_{4}$,
where $\Gamma_{4}$ is a conic neighborhood of $\xi_{0}$,
and such that
$\text{\textcursive{p}} \,(x,\xi)_{\vline \Omega_{3} \times \Gamma_{4}} > c_{0}>0$, $\Omega_{3} \Subset \widetilde{\Omega}_{3}$.
\smallskip

Let $u$ be a $s$-Gevrey vector, $s\geq 1$, for $P$, $u\in G^{s}(\Omega; P)$,
and $(x_{0}, \xi_{0}) \in \Char(P)$ such that $ \tau\left( \left(x_{0},\xi_{0}\right); \! \mathscr{P} \right)= r$.
Let $M$ a given fixed integer which will be determined at the end,
having the form $M = pn + q$, $p$ and $q$ suitable integers.
Let $\varphi_{N}(x)$ and $\psi_{N}(x)$ be two Ehrenpreis-H\"ormander sequences (\cite{EhrenpreisIV60}, see also \cite{H_Book-1},\cite{Treves2022})
associated to the couples $ (\Omega_{0}, \Omega_{1})$ and $(\Omega_{1}, \Omega_{2})$ respectively, $x_{0} \in \Omega_{0}$.
More precisely
$\varphi_{N}(x) \equiv 1$ on $\Omega_{0}$ and supported in $\Omega_{1}$,
$\psi_{N}(x) \equiv 1$ on $\Omega_{1}$ and supported in $\Omega_{2}$,
with $\overline{\Omega}_{0} \Subset \Omega_{1}$,
$\overline{\Omega}_{1}\Subset \Omega_{2} \Subset \widetilde{\Omega}_{3}$, 
and there are two positive constants $C_{\varphi}$ and $C_{\psi}$ such that
\begin{linenomath}
	$$
	| D^{\alpha} \varphi_{N}(x)| \leq C_{\varphi}^{|\alpha| +1} N^{(|\alpha|-M)^{+}}
	\quad \text{ and } \quad
	| D^{\alpha} \psi_{N}(x)| \leq C_{\psi}^{|\alpha| +1} N^{(|\alpha|-M)^{+}},
	$$
\end{linenomath}
for all $\alpha \in \mathbb{Z}^{n}_{+}$ such that  $|\alpha | \leq N$.\\
Let $ \Theta_{N}(D) $ be a sequence of zero order pseudo-differential
operators with symbols $\Theta_{N}(\xi)$ of Ehrenpreis-Andersson type
associated to the couple of open cones $(\Gamma_{0}, \Gamma_{1})$,
with $\overline{\Gamma}_{0} \Subset \Gamma_{1}$,
$\Gamma_{0}$ conic neighborhood of $\xi_{0}$.
More precisely $\Theta_{N}(\xi)\equiv 1$ in
$ \Gamma_{0,N} = \Gamma_{0} \cap \lbrace \xi\in \mathbb{R}^{n}\,: \, |\xi|>N\rbrace $,
supported in
$\Gamma_{1,N/2}= \Gamma_{1} \cap \lbrace \xi\in \mathbb{R}^{n}\,: \, |\xi|>N/2\rbrace $
and there is a positive constant $C_{\Theta}$ such that
\begin{linenomath}
	$$
	|\Theta_{N}^{(\alpha)}(\xi)| \leq  C^{|\alpha|+1} _{\Theta} N^{\left(|\alpha|-M\right)^{+}} \left(1+ |\xi|\right)^{-|\alpha|}, 
	$$
\end{linenomath}
for all $\alpha \in \mathbb{Z}^{n}_{+}$ such that  $|\alpha | \leq N$.
We refer to the  Appendix for a detailed construction of the symbols of Ehrenpreis-Andersson type.\\
We recall that we will use the following convention:
the Latin alphabet letters in the upper index will denote the derivatives with respect to
the corresponding direction, precisely $a^{(k)}(x) = D_{k}a(x)$,
and the Greek alphabet letters in the upper index will denote the usual multi-index derivatives,
precisely $a^{(\alpha)}(x) = D^{\alpha} a(x)= D_{1}^{\alpha_{1}} \cdots D_{n}^{\alpha_{n}} a(x)$,
$\alpha= (\alpha_{1}, \dots, \alpha_{n}) \in \mathbb{Z}_{+}^{n}$.\\
In order to make the following more readable we use the following notations:\\
$ P^{k}u =: f$, $  \varphi_{N}^{(\delta)} f =: g$, $ \Theta_{N}^{(\gamma)} D^{\alpha}g =: w$
and $\psi_{N}^{(\beta)} w =: v$.

The following result is obtained  taking advantage from the Theorem \ref{MicLocEst}.
\begin{proposition}\label{Basic-Est-D}
Let $v$ be as previously defined, i.e. $v = \psi_{N}^{(\beta)} \Theta_{N}^{(\gamma)} D^{\alpha} \varphi_{N}^{(\delta)} P^{k} u$.
Then there are positive constants $C$, $ A $ and $ B $ independent of $N$, $\alpha$, $\beta$, $\gamma$, $\delta$ and $k$ such that
\begin{linenomath}
		\begin{multline}\label{Est-1/r-1-p}
		\|  v \|_{\frac{1}{r}}^{2}
		+\sum_{j=1}^{n} \left( \| P^{j} v\|^{2} + \| P_{j} v\|_{-1}^{2} \right)
		\leq 
		C \! \left( \sum_{\ell=0}^{n}\! \left|\langle E_{\ell} P v, E_{\ell} v \rangle\right|
		+  \|v\|^{2}
		\right)
		\\
		\quad
		+ A^{2(|\gamma|+1)} B^{2(2 m + \sigma + 1)}  
		N^{2s(2 m+ |\gamma|+\sigma + 2n + 4 - M)^{+} },
		\end{multline}
\end{linenomath}
where $P^{j}$, $P_{j}$ and $E_{\ell}$, $\ell=0,\, \dots, n$,
are the same as in the Theorem \ref{MicLocEst}, $m= |\alpha| - |\gamma| $ 
and $\sigma= |\beta| + |\delta| +2k$. 
\end{proposition}
%
\subsection{Proof of the Proposition \ref{Basic-Est-D}}~\par
\vskip-6mm
Let $\widetilde{\Theta}_{\widetilde{m}}(D)$ be a sequence of zero order
pseudo-differential operators with symbol $\widetilde{\Theta}_{\widetilde{m}}(\xi)$,
of Ehrenpreis-Andersson type, associated to the couple of open cones $(\Gamma_{2},\Gamma_{3})$,
where $\overline{\Gamma}_2\Subset \Gamma_{3}\Subset\overline{\Gamma}_{3}\Subset \Gamma_{4}$.
We point out that $\widetilde{\Theta}_{\widetilde{m}}(\xi)$ is supported in
$\Gamma_{3} \cap \lbrace \xi\in \mathbb{R}^{n}\,: \, |\xi|>\widetilde{m}/2\rbrace $,
$\widetilde{\Theta}_{\widetilde{m}}(\xi) \equiv 1$ in  $ \Gamma_{2} \cap \lbrace \xi\in \mathbb{R}^{n}\,: \, |\xi|>\widetilde{m}\rbrace $,
and satisfies the estimate (\ref{ThetaN-AH}), Lemma \ref{EA-Cutoff} in the Appendix,
for all $\alpha \in \mathbb{Z}^{n}_{+}$ with $|\alpha|\leq \widetilde{m}$. 
We recall that the sequence $ \Theta_{N}(D) $ with symbols $\Theta_{N}(\xi)$
of Ehrenpreis-Andersson type is associated to the couple of open cones
$(\Gamma_{0}, \Gamma_{1})$, $\overline{\Gamma}_{1} \Subset \Gamma_{2}$. 
We assume that $\widetilde{m} \leq N$.
Let $\widetilde{\psi} \in \mathscr{D}\left( \Omega_{3}\right)$ and
such that it is identically equal to $1$ on $\Omega_{2}$,
where $\overline{\Omega}_{2}\Subset \Omega_{3} \Subset\overline{\Omega}_{3} \Subset \widetilde{\Omega}_{3}$.
We recall that $\psi_{N}$ is supported in $\Omega_{2}$.
We use the same notations introduced in the beginning of this section.\\ 
We have
\begin{linenomath}
	\begin{multline*}
	\widetilde{\psi} \widetilde{\Theta}_{\widetilde{m}} \psi_{N}^{(\beta)} w = \widetilde{\psi} \left[ \widetilde{\Theta}_{\widetilde{m}},  \psi_{N}^{(\beta)} \right] w
	+ \widetilde{\psi} \psi_{N}^{(\beta)} \widetilde{\Theta}_{\widetilde{m}} w
	\\
	=  \psi_{N}^{(\beta)} w + \widetilde{\psi} \psi_{N}^{(\beta)} ( \widetilde{\Theta}_{\widetilde{m}}-1) w 
	+ \widetilde{\psi} \left[ \widetilde{\Theta}_{\widetilde{m}},  \psi_{N}^{(\beta)} \right] w.
	\end{multline*}
\end{linenomath}
We recall that $ w= \Theta_{N}^{(\gamma)} D^{\alpha} \varphi_{N}^{(\delta)} P^{k} u$. 
We deduce that
\begin{linenomath}
	\begin{align}\label{step-0-1/r}
	v= \widetilde{\psi}  \widetilde{\Theta}_{\widetilde{m}} \psi_{N}^{(\beta)} \Theta_{N}^{(\gamma)} D^{\alpha}g
	+ \widetilde{\psi} \left[\psi_{N}^{(\beta)}, \widetilde{\Theta}_{\widetilde{m}} \right] \Theta_{N}^{(\gamma)} D^{\alpha}g
	+\psi_{N}^{(\beta)} (1- \widetilde{\Theta}_{\widetilde{m}}) w,
	\end{align}
\end{linenomath}
where $ v = \psi_{N}^{(\beta)} \Theta_{N}^{(\gamma)} D^{\alpha} \varphi_{N}^{(\delta)} P^{k} u$
and $ g= \varphi_{N}^{(\delta)} P^{k} u$.
We point out that since $u\in G^{s}(\Omega;P)$, we know a good estimate for $L^{2}$-bound of $g$.
So then we know the estimate for $w$ in $H^{-m}$, with $m = |\alpha| - |\gamma|$.\\
We take the expansion up to order $\widetilde{m}$ 
of the bracket $ [\psi_{N}^{(\beta)}, \widetilde{\Theta}_{\widetilde{m}} ]$.
We obtain
\begin{linenomath}
	\begin{multline*}
	\widetilde{\psi} \left[\psi_{N}^{(\beta)}, \widetilde{\Theta}_{\widetilde{m}} \right] 
	\Theta_{N}^{(\gamma)} D^{\alpha}g
	\\
	=-\widetilde{\psi}
	\displaystyle\longsum[25]_{1 \leq |\mu| \leq \widetilde{m} - 1}\frac{1}{\mu!} \psi_{N}^{(\beta+\mu)} \widetilde{\Theta}_{\widetilde{m}}^{(\mu)}\Theta_{N}^{(\gamma)} D^{\alpha}g
	- \widetilde{\psi}\mathscr{R}_{\widetilde{m}}\left( \left[\psi_{N}^{(\beta)}, \widetilde{\Theta}_{\widetilde{m}} \right]  \right) w,
	\end{multline*}
\end{linenomath}
where
\begin{linenomath}
	\begin{multline}\label{Rem-0}
	\mathscr{R}_{\widetilde{m}}\left( \left[\psi_{N}^{(\beta)}, \widetilde{\Theta}_{\widetilde{m}} \right]  \right) w(x)
	=
	\frac{\widetilde{m}}{\left(2\pi\right)^{2n}} 
	\\
	\times 
	\! \sum_{|\mu| =  \widetilde{m}} \frac{1}{\mu!}\!
	\int e^{ix\xi}\! \!\int\!\! \widehat{\psi}_{N}^{(\beta+\mu)}(\xi-\eta) 
	\int_{0}^{1}\!\!  \widetilde{\Theta}_{\widetilde{m}}^{(\mu)} \left( \eta +t(\xi -\eta) \right)(1 -t)^{\widetilde{m} -1} dt \widehat{w} (\eta) d\eta  d\xi.
	\end{multline}
\end{linenomath}
So equality (\ref{step-0-1/r}) gives
\begin{linenomath}
	\begin{multline}\label{Est-1/r-0}
	\|  v \|_{\frac{1}{r}}
	\leq 
	\| \widetilde{\psi}  \widetilde{\Theta}_{\widetilde{m}} \psi_{N}^{(\beta)} w\|_{\frac{1}{r}}
	+
	\displaystyle\longsum[25]_{1 \leq |\mu| \leq \widetilde{m} - 1}\frac{1}{\mu!} 
	\|  \psi_{N}^{(\beta+\mu)} \widetilde{\Theta}_{\widetilde{m}}^{(\mu)} w \|_{\frac{1}{r}}
	\\
	+
	\| \widetilde{\psi} \mathscr{R}_{\widetilde{m}}\left( \left[\psi_{N}^{(\beta)}, \widetilde{\Theta}_{\widetilde{m}} \right]  \right) w \|_{\frac{1}{r}}
	+
	\|\psi_{N}^{(\beta)} (1- \widetilde{\Theta}_{\widetilde{m}}) w \|_{\frac{1}{r}}.
	\end{multline}
\end{linenomath}
We begin to estimate the first term on the right hand side of the above inequality.\\
In view of the Theorem \ref{MicLocEst} there is a positive constant $C$ such that the following estimate holds
\begin{linenomath}	
	\begin{align}\label{Der-Est-M-1}
	\| \text{\textcursive{p}}\, v\|_{\frac{1}{r}}^{2}+ \sum_{j=1}^{n} \left( \| P^{j}v\|^{2}_{0} + \|P_{j}v\|_{-1}^{2}\right)
	\leq C \left( \sum_{\ell=0}^{n}\left|\langle E_{\ell}P v, E_{\ell}v\rangle\right|  + \|v\|_{0}^{2}\right).
	\end{align} 
\end{linenomath}
Let $Q(x,D)$ the zero order operator associated to the symbol
\begin{linenomath}
	$$
	\text{\textcursive{q}}_{\widetilde{m}}(x,\xi) = 
	\displaystyle\frac{ \widetilde{\psi}(x)  \widetilde{\Theta}_{\widetilde{m}}(\xi)}{\text{\textcursive{p}}(x,\xi)}.
	$$
\end{linenomath}
We point out that $\text{\textcursive{q}}_{\widetilde{m}}(x,\xi) $ is well defined as $ |\text{\textcursive{p}} |\geq  c_{0}> 0$
on the support of $ \widetilde{\psi}\widetilde{\Theta}_{\widetilde{m}}$.\\ 
We have
\begin{linenomath}	
	\begin{align*}
	\left(Q\circ \text{\textcursive{p}}\right)v(x)
	= \frac{1}{(2\pi)^{2n}}
	\int  e^{ix\xi}\text{\textcursive{q}}_{\widetilde{m}}(x,\xi) \widehat{\text{\textcursive{p}}}(\xi-\eta,\eta)\widehat{v}(\eta) d\eta\,d\xi,
	\end{align*}
\end{linenomath}	
where $\widehat{\text{\textcursive{p}}}(\cdot,\cdot)$ is the Fourier transform of $\text{\textcursive{p}}$ with respect to $x$.\\
Since we have
\begin{linenomath}	
	\begin{align*}
	\text{\textcursive{q}}_{\widetilde{m}}(x,\eta+\tau)= 
	\frac{ \widetilde{\psi}(x)  \widetilde{\Theta}_{\widetilde{m}}(\eta)}{\text{\textcursive{p}}(x,\eta)}
	+\sum_{j=1 }^{n}\tau_{j}\int_{0}^{1} \left(\frac{ \widetilde{\psi}  \widetilde{\Theta}_{\widetilde{m}}}{\text{\textcursive{p}}}\right)^{(j)}(x,\eta + t\tau)\, dt,
	\end{align*}
\end{linenomath}	
we obtain
\begin{linenomath}	
	\begin{multline*}
	\left(Q\circ \text{\textcursive{p}}\right)v(x)
	= \frac{1}{(2\pi)^{2n}} 
	\int  e^{ix(\eta +\tau)} \frac{ \widetilde{\psi}(x)  \widetilde{\Theta}_{\widetilde{m}}(\eta+\tau)}{\text{\textcursive{p}}(x,\eta+\tau)}\widehat{\text{\textcursive{p}}}(\tau,\eta)\widehat{v}(\eta) d\eta\,d\tau
	\\
	=\frac{1}{(2\pi)^{2n}} 
	\int  e^{ix\eta} e^{ix\tau} \frac{ \widetilde{\psi}(x)  \widetilde{\Theta}_{\widetilde{m}}(\eta)}{\text{\textcursive{p}}(x,\eta)}\widehat{\text{\textcursive{p}}}(\tau,\eta)\widehat{v}(\eta) d\eta\,d\tau
	\\
	+
	\frac{1}{(2\pi)^{2n}} \sum_{j=1}^{n}
	\iint \int_{0}^{1} e^{ix(\eta +\tau)} \left(\frac{ \widetilde{\psi}  \widetilde{\Theta}_{\widetilde{m}}}{\text{\textcursive{p}}}\right)^{(j)}(x,\eta + t\tau) \tau_{j} 
	\widehat{\text{\textcursive{p}}}(\tau,\eta)\widehat{v}(\eta)\, dt\, d\eta\,d\tau
	\\
	=
	\frac{1}{(2\pi)^{n}} 
	\int  e^{ix\eta} \widetilde{\psi}(x)  \widetilde{\Theta}_{\widetilde{m}}(\eta)\widehat{v}(\eta) d\eta\,
	\\
	+
	\frac{1}{(2\pi)^{2n}} \sum_{j=1}^{n}
	\int e^{ix\eta} \int\int_{0}^{1} e^{ix \tau} \left(\frac{ \widetilde{\psi}  \widetilde{\Theta}_{\widetilde{m}}}{\text{\textcursive{p}}}\right)^{(j)}
	\!\!\!(x,\eta + t\tau) \tau_{j} 
	\widehat{\text{\textcursive{p}}}(\tau,\eta)\, dt\,d\tau \,\,  \widehat{v}(\eta)  \, d\eta
	\\
	=\widetilde{\psi} \widetilde{\Theta}_{\widetilde{m}}v(x) +\mathscr{R}_{1}v(x),
	\end{multline*}
\end{linenomath}	
where we use that $(2\pi)^{-n}\int e^{ix\tau} \widehat{\text{\textcursive{p}}}(\tau,\eta) \, d\tau=\text{\textcursive{p}}(x,\eta)$.\\
The symbol associated to the operator $\mathscr{R}_{1}(x,D)$ is
\begin{linenomath}	
	\begin{multline*}
	\text{\textcursive{r}}_{1}(x,\xi) =
	\sum_{j=1 }^{n} \iint e^{i(y-x)(\xi-\eta)} \text{\textcursive{p}}_{(j)}(y,\xi) 
	\int_{0}^{1} \text{\textcursive{q}}_{\widetilde{m}}^{(j)}(x, \xi + t (\eta-\xi)) dt\, dy\,\frac{d\eta}{(2\pi)^{n}},
	\end{multline*}
\end{linenomath}	
where as usual the lower indexes denote the derivatives with respect to the variable and the upper indexes denote the derivatives with respect the co-variable;
moreover the following estimate holds
\begin{linenomath}	
	\begin{align*}
	\left| \text{\textcursive{r}}_{1}(x,\xi)\right|
	\leq \tilde{C} \left( 1 + |\xi|^{2}\right)^{-1/2}, 
	\end{align*}
\end{linenomath}	
where $\tilde{C} $ depends only on $n$ and on the derivatives of
$ \text{\textcursive{p}}(y,\xi) $ up to order $\left\lfloor \frac{n}{2} \right\rfloor + 2$
with respect to $y$.
We obtain
\begin{linenomath}
	\begin{align*}
	\| \widetilde{\psi}  \widetilde{\Theta}_{\widetilde{m}} v\|_{\frac{1}{r}}
	\leq	\|  Q\, \text{\textcursive{p}} v\|_{\frac{1}{r}} + \| \mathscr{R}_{1}v\|_{\frac{1}{r}}.
	\end{align*}
\end{linenomath}
By the Calderon-Vaillancourt theorem applied to the zero order operator $Q$, (see \cite{KG} end \cite{Hwang-1987}), we have
\begin{linenomath}
	\begin{align}\label{Est-1/r-2}
	\| \widetilde{\psi}  \widetilde{\Theta}_{\widetilde{m}} v\|_{\frac{1}{r}}^{2}
	\leq
	C_{1}	\| \text{\textcursive{p}} v\|_{\frac{1}{r}}^{2} + \tilde{C}_{1} \| v\|_{-1+\frac{1}{r}}^{2},
	\end{align}
\end{linenomath}
where the constants $C_{1}$ and $\tilde{C}_{1}$ do not depend on $\widetilde{m}$.
By (\ref{Der-Est-M-1}) and (\ref{Est-1/r-2}) we get
\begin{align}\label{Est-1/r-3}
\| \widetilde{\psi}  \widetilde{\Theta}_{\widetilde{m}} v\|_{\frac{1}{r}}^{2}
+\sum_{j=1}^{n} \left( \| P^{j} v\|^{2} + \| P_{j} v\|_{-1}^{2} \right)
\leq 
C_{2} 
\left( \sum_{\ell=0}^{n}\left|\langle E_{\ell}P v, E_{\ell}v\rangle\right|  + \|v\|_{0}^{2}\right),
\end{align}
where $P^{j}$,  $P_{j}$ and $E_{\ell}$ are as in Theorem \ref{MicLocEst}.

\smallskip
\textbf{Estimate of the second to last term on the right hand side of (\ref{Est-1/r-0}).}\\
We have:
\begin{lemma}\label{Rem0}
Let $ \mathscr{R}_{\widetilde{m}}\left( \left[\psi_{N}^{(\beta)}, \widetilde{\Theta}_{\widetilde{m}} \right]  \right) \! w(x)$ be as in (\ref{Rem-0}),
$w= \Theta_{N}^{(\gamma)} D^{\alpha} \varphi_{N}^{(\delta)} P^{k} u$ 
and $ \widetilde{m}= |\alpha|- |\gamma| + \left\lfloor \frac{n}{2} \right\rfloor +1$,
then there are two positive constants $\widehat{C}_{1}$ and $\widehat{C}_{2}$, independent of $N,\, \beta  \, \gamma, \, \alpha, \, \delta$ and  $k $,
such that
\begin{linenomath}
	\begin{align}\label{Rem0-Est}
		\| \mathscr{R}_{\widetilde{m}}\left( \left[\psi_{N}^{(\beta)}, \widetilde{\Theta}_{\widetilde{m}} \right]  w\right)\|_{1/r}
		\leq 
		\widehat{C}_{1}^{\sigma+1} \widehat{C}_{2}^{2 m +|\gamma| + 2n + 4} 
		N^{s(2 m +|\gamma|+ \sigma + 2n + 4 - M)^{+}}
	\end{align}
\end{linenomath}
where $m =|\alpha| -|\gamma|$ and $\sigma= |\beta| + |\delta| +2k$.
\end{lemma}
Before to proof the above Lemma we need the following technical Lemma.
\begin{lemma}\label{L-1}
For every $N \in  \mathbb{N}^{*}$, one has 
\begin{linenomath}
	\begin{align}
    	k^{j} \leq B^{j} N^{(k-M)^{+}} \text{ for } B= \sup\left( M, 3\right), \, j \leq k \leq N.
	\end{align}
\end{linenomath}
\end{lemma}
\begin{proof}
\textit{First case.} Let us assume that $M \geq 3$,
we have to prove $k^{j} \leq  M^{j} N^{(k -M)^{+}}$, $ 1\leq j \leq k \leq N$.
If $k \leq M$, it is reduced to $k^{j} \leq M^{j}$ which is clear.\\
Let $k \geq M$, we have to show
\begin{linenomath}
	\begin{align*}
		k^{j} \leq M^{j} N^{(k-M)}, \text{ for }  \, j \leq k \leq N.
	\end{align*} 
\end{linenomath}
For our purpose it is enough to show that $j \left(\log k - \log M\right)- (k-M)\log N$ is non positive.
As $M \leq k \leq N$, it is sufficient to show:
\begin{linenomath}
	\begin{align*}
		E(k) \doteq k \left(\log k - \log M\right)- (k-M)\log N \leq 0.
	\end{align*}
\end{linenomath}
Let $E(x) \doteq  x \left(\log x - \log M\right)- (x-M)\log N$,
obtained replacing $k$ by $x$ with $x \in \left[ M, N\right]$. Then $E'(x)$ is equal to
$  \log x - \log M - \log N +1$; it is negative as $M \geq 3$ and $x \in \left[ M, N\right]$.
It is sufficient that $E(M) \leq 0$.
But $E(M) =0$ and $E$ is decreasing on $\left[ M, N\right]$.
So the Lemma is proved in case $M \geq 3$.\\
\textit{Second case.}
Let us assume that $M < 3$, so $M $ is equal to $2$, $1$, or $0$.
As for $M=0$, it is trivial, we just take $M=2$ or $1$. We know that
$k^{j} \leq 3^{j} N^{(k-3)^{+}}$, $j \leq k \leq N$.
But $ (k-3)^{+} \leq (k-2)^{+}$, so $ k^{j} \leq 3^{j} N^{(k-2)^{+}}$,
and $ (k-2)^{+} \leq (k-1)^{+}$ 
so $ k^{j} \leq 3^{j} N^{(k-1)^{+}}$, $j \leq k \leq N$.\\
The Lemma is proved for $B =\sup\left( M, 3\right)$. 
\end{proof}
\begin{remark}\label{Rk-1}
We will use the following elementary fact: if $p_1,\, \dots, p_{\ell} $ and $ M $ are integers then
\begin{linenomath}
	\begin{align*}
		N^{(p_1- M)^{+}} \cdots N^{(p_{\ell}- M)^{+}} \leq N^{(p_1+ \cdots + p_{\ell}- M)^{+}} .
	\end{align*}
\end{linenomath}
\end{remark}
\begin{remark}\label{Rem0-Est-1}
Choosing $M\geq 2n+4 $, the estimate (\ref{Rem0-Est}) implies the inequality
\begin{linenomath}
	\begin{align*}
		\| \mathscr{R}_{\widetilde{m}}\left( \left[\psi_{N}^{(\beta)}, \widetilde{\Theta}_{\widetilde{m}} \right] \right) w\|_{\frac{1}{r}}
		\leq 
		\widehat{C}_{3}^{\sigma+1} \widehat{C}_{4}^{2 m +|\gamma|+ 1}  
		N^{s\left(2 m + |\gamma|+ \sigma\right) },
	\end{align*}
\end{linenomath}
where $m =|\alpha| -|\gamma|$ and $\sigma= |\beta| + |\delta| +2k$.
\end{remark}
\begin{proof}[\textbf{Proof of Lemma \ref{Rem0}}]
Since $\left( 1+ |\xi|^{2}\right)^{t/2} \leq \left( 1 + |\eta|\right)^{t} \left( 1+ | \xi -\eta|^{2}\right)^{t/2}$, $t \geq 0$,
we have
\begin{linenomath}
	\begin{multline}\label{Rem0-S1}
		\| \widetilde{\psi} \mathscr{R}_{\widetilde{m}}\left( \left[\psi_{N}^{(\beta)}, \widetilde{\Theta}_{\widetilde{m}} \right]  \right) w \|_{\frac{1}{r}}
		\\
		\leq \|\mathscr{R}_{\widetilde{m}}\left( \left[\psi_{N}^{(\beta)}, \widetilde{\Theta}_{\widetilde{m}} \right]  \right) w \|_{\frac{1}{r}}
		\int\left( 1 + |\eta|\right)^{1/r} |\widehat{\widetilde{\psi} }(\eta)| d\eta 
		\\
		\leq C_{\widetilde{\psi}}\|\mathscr{R}_{\widetilde{m}}\left( \left[\psi_{N}^{(\beta)}, \widetilde{\Theta}_{\widetilde{m}} \right]  \right) w \|_{\frac{1}{r}},
	\end{multline}
\end{linenomath}
where  $C_{\widetilde{\psi}}$ depends only on the derivatives up to order $n+2$ of $\widetilde{\psi}$.\\
We have
\begin{linenomath}	
	\begin{multline}\label{Rem0-S2}
		(2\pi)^{4n}\|\mathscr{R}_{\widetilde{m}}\left( \left[\psi_{N}^{(\beta)}, \widetilde{\Theta}_{\widetilde{m}} \right]  \right) w \|_{\frac{1}{r}}
		\\
		= (2\pi)^{4n}
		\left(\int  \left( 1+ |\xi|^{2}\right)^{\frac{1}{r}} 
		\left| \mathscr{R}_{\widetilde{m}}\left( \left[\psi_{N}^{(\beta)}, \widetilde{\Theta}_{\widetilde{m}} \right]  \right) w(\xi) \right|^{2} d\xi \right)^{1/2}
		\\
		=
		\left( \int \Big| 
		\, \int
		\underbracket{ \left( 1+ |\xi|^{2}\right)^{\frac{1}{2r}} \left(\sum_{|\mu| =  \widetilde{m}} \frac{\widetilde{m}}{\mu!}\widehat{\psi}_{N}^{(\beta+\mu)}(\xi-\eta) 
			\int_{0}^{1}  \widetilde{\Theta}_{\widetilde{m}}^{(\mu)} \left( \eta +t(\xi -\eta) \right)(1 -t)^{\widetilde{m}-1}\, dt\,\right)}_{\doteq \widetilde{g}(\xi,\eta)}
		\right.
		\\
		\left.
		\times
		\widehat{w} (\eta)\, d\eta\, \Bigg|^{2} d\xi
		\right)^{1/2} 
		\\
		\leq 
		\left(\int   \left( \int  |\widetilde{g}(\xi,\eta)|  | \Theta_{N}^{(\gamma)}(\eta) \eta^{\alpha}| |\widehat{g} (\eta) |\, d\eta\, \right)^{2} d\xi \right)^{\frac{1}{2}} 
		\qquad ( \widehat{w} (\eta)= \Theta_{N}^{(\gamma)}(\eta) \eta^{\alpha} \widehat{g}(\eta))
		\\
		\leq 
		\left(\int \left( \|\widehat{g}  \|^{2}_{L^{2}_{\eta}} \,
		\int  |\widetilde{g}(\xi,\eta)|^{2}  | \Theta_{N}^{(\gamma)}(\eta) \eta^{\alpha}|^{2} \, d\eta  \right) d\xi \right)^{\frac{1}{2}} 
		\\
		\leq 
		\|\widehat{g}  \|_{L^{2}_{\eta}}
		\left( \, \int\,  \int |\widetilde{g}(\xi,\eta)|^{2}  | \Theta_{N}^{(\gamma)}(\eta) \eta^{\alpha}|^{2} \, d\eta \, d\xi \right)^{\frac{1}{2}},
	\end{multline}
\end{linenomath}
where $\|\widehat{g}\|_{L^{2}_{\eta}} = \| \varphi_{N}^{(\delta)} P^{k} u \|_{0} $.\\
We have to estimate the second factor on the right hand side of the above inequality.
Taking advantage from the Peetre's inequality,
for every $t$ in $\left[0,1\right] $ we have
\begin{linenomath}
	\begin{multline}
		\nonumber
		\left| \widetilde{\Theta}_{\widetilde{m}}^{(\mu)} \left( \eta +t(\xi -\eta)\right)\right|
		\leq C_{\widetilde{\Theta}}^{|\mu| +1}  \widetilde{m}^{\left(|\mu|-M\right)^{+}} \left(1 +  \left| \eta +t(\xi -\eta)\right|^{2}\right)^{-|\mu|/2}
		\\
		\nonumber
		\leq C_{\widetilde{\Theta}}^{|\mu| +1}  \widetilde{m}^{\left(|\mu|-M\right)^{+}} 2^{|\mu|/2} 
		\left(1 +  | \eta |^{2}\right)^{-|\mu|/2}  \left(1 +  \left|\xi -\eta\right|^{2}\right)^{|\mu|/2}.
	\end{multline}
\end{linenomath}
So
\begin{linenomath}
	\begin{multline}\label{Est-g-til}
		\int |\widetilde{g}(\xi,\eta)|^{2} \, d\xi
		\leq
		\\
		\int
		\left( 1+ |\xi|^{2}\right)^{1/r} \!
		\left(\sum_{|\mu| =  \widetilde{m}} \frac{\widetilde{m}}{\mu!} |\widehat{\psi}_{N}^{(\beta+\mu)}(\xi-\eta)| 
		\int_{0}^{1} | \widetilde{\Theta}_{\widetilde{m}}^{(\mu)} \left( \eta +t(\xi -\eta) \right) | (1 -t)^{\widetilde{m}-1}\, dt\right)^{2} d\xi
		\\
		\hspace{-21em}
		\leq C_{\widetilde{\Theta}}^{2|\mu| + 2}  2^{|\mu| +\frac{1}{r}} 
		\left(1 +  | \eta |^{2}\right)^{-|\mu|+\frac{1}{r}}
		\\
		\times
		\int
		\left(\sum_{|\mu| =  \widetilde{m}} \frac{\widetilde{m}^{\left(|\mu|-M\right)^{+}}}{\mu!}
		\left(1 +  \left|\xi -\eta\right|^{2}\right)^{\frac{|\mu|}{2}+\frac{1}{2r}} |\widehat{\psi}_{N}^{(\beta+\mu)}(\xi-\eta)| \right)^{2} d\xi
		\\
		\hspace{-18em}
		\leq 
		C_{\widetilde{\Theta}}^{2|\mu| + 2}   2^{|\mu| +\frac{n+1}{2}+\frac{1}{r}} 
		\left(1 +  | \eta |^{2}\right)^{-|\mu|+\frac{1}{r}}\times
		\\
		\int\! \frac{1}{\left(1 +  \left|\xi -\eta\right|\right)^{n+1}}\!
		\left(\sum_{|\mu| =  \widetilde{m}} \frac{\widetilde{m}^{\left(|\mu|-M\right)^{+}}}{\mu!}
		\left(1 +  \left|\xi -\eta\right|^{2}\right)^{ \frac{1}{2}\left(  |\mu|+n+1 +\frac{1}{r} \right)} |\widehat{\psi}_{N}^{(\beta+\mu)}(\xi-\eta)| \right)^{2} 
		\!\!\!d\xi,
	\end{multline}
\end{linenomath}
where we use the fact that $ ( \widetilde{m})\int_{0}^{1}(1-t)^{ \widetilde{m}-1} dt=1$.
Since $ |\mu|! \leq n^{|\mu|} \mu!$, $ |\mu|^{|\mu|} \leq (en)^{|\mu|} |\mu|!$ and 
recalling that the number of the multi-indexes with  $|\mu| = \widetilde{m} $
is given by $\binom{\widetilde{m}+n-1}{n-1}$,
smaller than $2^{\widetilde{m}}$, the term in the square brackets
in the above formula can be handled in the following way  
\begin{linenomath}
	\begin{multline}
		\nonumber
		\sum_{|\mu| =  \widetilde{m}} \frac{\widetilde{m}^{\left(|\mu|-M\right)^{+}}}{\mu!}
		\left(1 +  \left|\xi -\eta\right|^{2}\right)^{ \frac{1}{2}\left(  |\mu|+n+1 +\frac{1}{r} \right)} |\widehat{\psi}_{N}^{(\beta+\mu)}(\xi-\eta)| 
		\\
		\leq 2^{\widetilde{m} + n + 2} 
		\sum_{|\mu| =  \widetilde{m}} \frac{\widetilde{m}^{\left(|\mu|-M\right)^{+}}}{\mu!}
		\left( 1 +  \left|\xi -\eta\right|^{|\mu|+  n+2} \right)
		|\widehat{\psi}_{N}^{(\beta+\mu)}(\xi-\eta)| 
		\\
		\leq 2^{\widetilde{m} + n + 2} (en)^{\widetilde{m}} 
		\sum_{|\mu| =  \widetilde{m} }
		\frac{\widetilde{m}^{\left(|\mu|-M\right)^{+}} }{\widetilde{m}^{|\mu|} }
		\left( 1 +  \left|\xi -\eta\right|^{|\mu|+  n+2} \right)
		|\widehat{\psi}_{N}^{(\beta+\mu)}(\xi-\eta)| 
		\\
		\leq 2^{n+3} (2en)^{ \widetilde{m}}
		C_{\psi}^{|\beta|+2\widetilde{m} +n + 2} N^{ \left( |\beta|+2 \widetilde{m} + n + 2 -M\right)^{+}}
		\sum_{|\mu| =  \widetilde{m}}  1
		\\
		\leq 2^{n+3} (4en)^{ \widetilde{m}}
		C_{\psi}^{|\beta|+2\widetilde{m} +n + 2} N^{ \left( |\beta|+2 \widetilde{m} + n + 2 -M\right)^{+}}.
	\end{multline}
\end{linenomath}
The left hand side of (\ref{Est-g-til}) can be estimated by
\begin{linenomath}
	\begin{align}\label{Est-g-til-1}
		C_{1} 4^{n+4}  (4en)^{2\widetilde{m}} C_{\widetilde{\Theta}}^{2\widetilde{m}+ 2} 
		C_{\psi}^{2(|\beta|+2\widetilde{m} +n + 2)} N^{2\left( |\beta|+2 \widetilde{m} + n + 2 -M\right)^{+}} 
		\left(1 +  | \eta |^{2}\right)^{-\widetilde{m}+\frac{1}{r}}.
	\end{align}
\end{linenomath}
By (\ref{Est-g-til}) and (\ref{Est-g-til-1}) we get
\begin{linenomath}
	\begin{multline}\label{Est-g-til-f}
		\left( \, \int\, \int |\widetilde{g}(\xi,\eta)|^{2} \, d\xi\,  | \Theta_{N}^{(\gamma)}(\eta) \eta^{\alpha}|^{2} \, d\eta \right)^{\frac{1}{2}}
		\\
		\leq
		C_{1} 4^{n+4}  (4en)^{2\widetilde{m}} C_{\widetilde{\Theta}}^{\widetilde{m}+ 1} 
		C_{\psi}^{|\beta|+2\widetilde{m} +n + 2} N^{\left( |\beta|+2 \widetilde{m} + n + 2 -M\right)^{+}} 
		\\
		\times \left( \int   \left(1 +  | \eta |^{2}\right)^{-\widetilde{m}+\frac{1}{r}}  | \Theta_{N}^{(\gamma)}(\eta) \eta^{\alpha}|^{2} d\eta\right)^{1/2}
		\\
		\leq
		C_{1} 4^{n+4}  (4en)^{2\widetilde{m}} C_{\widetilde{\Theta}}^{\widetilde{m}+ 1} 
		C_{\psi}^{|\beta|+2\widetilde{m} +n + 2}  C_{\Theta}^{|\gamma|+1} N^{\left( |\beta|+2 \widetilde{m} + n + 2 -M\right)^{+}} N^{\left( |\gamma|-M\right)^{+}}
		\\
		\left( \int   \left(1 +  | \eta |^{2}\right)^{|\alpha|- |\gamma|-  \widetilde{m} +\frac{1}{r}}  d\eta\right)^{1/2}.
	\end{multline}
\end{linenomath}
Choosing $ \widetilde{m}= |\alpha|- |\gamma| + \left\lfloor \frac{n}{2} \right\rfloor +1 = m  + \left\lfloor \frac{n}{2} \right\rfloor +1$, 
by (\ref{Rem0-S1}), (\ref{Rem0-S2}) and (\ref{Est-g-til-f}) we obtain  
\begin{linenomath}
	\begin{align}
		\nonumber
		\| \widetilde{\psi} \mathscr{R}_{\widetilde{m}}( [\psi_{N}^{(\beta)}, \widetilde{\Theta}_{\widetilde{m}} ] ) w \|_{\frac{1}{r}}
		&\leq
		C_{\widetilde{\psi}} 
		C_{1}C_{2}
		4^{n+4}  (4en)^{2\left(m+ \left\lfloor \frac{n}{2} \right\rfloor +1 \right) } 
		C_{\widetilde{\Theta}}^{m  + \left\lfloor \frac{n}{2} \right\rfloor + 2} C_{\Theta}^{|\gamma|+1}
		\\
		&
		\nonumber
		\quad
		\times
		C_{\psi}^{2m+|\beta| + 2n + 4} 
		N^{\left( 2m+ |\beta|+|\gamma|+ 2n + 4 -M\right)^{+}} 
		\| \varphi_{N}^{(\delta)} P^{k} u \|_{0}.
	\end{align}
\end{linenomath}
Let $K_{0}$ be a compact set contained in $\Omega_{2}$ and containing all the supports of $\varphi_{N}$, $\Omega_{1} \subset K_{0}$.
Since $u $ is a $G^{s}$-vector for $P$ in $\Omega$ we have $\|P^{k}u\|_{L^{2}(K_{0})} \leq C_{K_{0}}^{2k+1} k^{2sk}$,
moreover by the Lemma \ref{L-1} the following estimate holds $k^{2sk} \leq B^{2sk} N^{s(2k-M)^{+}}$.
So, we get
\begin{linenomath}
\begin{align}\label{Est_g}
	\| g \|_{0} 
	=\| \varphi_{N}^{(\delta)} P^{k} u \|_{0}
	\leq
	\|P^{k}u\|_{L^{2}(K_{0})} 
	\leq C_{\varphi}^{|\delta|+1}\tilde{B}^{2k+1} N^{s(|\delta|+2k-M)^{+}}.
\end{align}
\end{linenomath}
Using this estimate and the Remark \ref{Rk-1},
a suitable choice of $\widehat{C}_{1}$ and $\widehat{C}_{2}$
allows us to gain the estimate (\ref{Rem0-Est}).
This concludes the proof of the Lemma \ref{Rem0}.
\end{proof}
\textbf{Estimate of the terms in the sum on the right hand side of (\ref{Est-1/r-0}).}\\ 
Since $\left( 1+ |\xi|^{2}\right)^{t/2} \leq \left( 1 + |\eta|\right)^{t} \left( 1+ | \xi -\eta|^{2}\right)^{t/2}$, $t \geq 0$,
we have
\begin{linenomath}
\begin{multline}\label{Est-Sum-S1}
	\displaystyle\longsum[25]_{1 \leq |\mu| \leq \widetilde{m} - 1}\frac{1}{\mu!} 
	\|  \psi_{N}^{(\beta+\mu)} \widetilde{\Theta}_{\widetilde{m}}^{(\mu)}\Theta_{N}^{(\gamma)} D^{\alpha} g \|_{\frac{1}{r}}
	\\
	\leq
	\displaystyle\longsum[25]_{1 \leq |\mu| \leq \widetilde{m} - 1}\frac{1}{\mu!} 
	\int (1+|\xi|)^{\frac{1}{r}} |\widehat{\psi}_{N}^{(\beta+\mu)} (\xi) |\, d\xi 
	\| \widetilde{\Theta}_{\widetilde{m}}^{(\mu)}\Theta_{N}^{(\gamma)} D^{\alpha} g \|_{\frac{1}{r}}
	\\
	\leq
	C_{1} 
	\displaystyle\longsum[25]_{1 \leq |\mu| \leq \widetilde{m} - 1}\frac{1}{\mu!} 
	C_{\psi}^{|\beta|+|\mu| +n+3} N^{ \left( |\beta|+|\mu|+ n + 3 -M\right)^{+}} 
	\| \widetilde{\Theta}_{\widetilde{m}}^{(\mu)}\Theta_{N}^{(\gamma)} D^{\alpha} g \|_{\frac{1}{r}}.
\end{multline}
\end{linenomath}
Since $ \widetilde{\Theta}_{\widetilde{m}} (\xi)$ and $ \Theta_{N} (\xi)$
satisfy the estimate (\ref{ThetaN-AH}), Lemma \ref{EA-Cutoff}, in Appendix, we have
\begin{linenomath}
	\begin{multline*}
	\left(1 +|\xi|^{2}\right)^{\frac{1}{r}} |\widetilde{\Theta}_{\widetilde{m}}^{(\mu)} (\xi) | \Theta_{N}^{(\gamma)} (\xi)| |\xi^{\alpha}|
	\\
	\leq C_{\widetilde{\Theta}}^{|\mu| + 1} C_{\Theta}^{|\gamma| + 1} \widetilde{m}^{ \left(|\mu| -M\right)^{+}}  N^{ \left( |\gamma| -M\right)^{+}} 
	\left(1 +|\xi|\right)^{ |\alpha| -|\gamma| -|\mu| + \frac{2}{r}}.
	\end{multline*}
\end{linenomath}	
As before, we set
$ \widetilde{m}= |\alpha|- |\gamma| + \left\lfloor \frac{n}{2} \right\rfloor +1 = m  + \left\lfloor \frac{n}{2} \right\rfloor +1 <N$.\\
We distinguish tow cases.\\
Case $|\mu|> |\alpha| - |\gamma|$
(we remark that in the sum on the right hand side of \eqref{Est-Sum-S1}
the number of multi-index $\mu$ such that $|\alpha| - |\gamma|< |\mu| \leq |\alpha| - |\gamma| +  \left\lfloor \frac{n}{2} \right\rfloor +1$
is finite and it can be roughly estimated $2^{|\alpha| - |\gamma| +  \left\lfloor \frac{n}{2} \right\rfloor +1}$.)\\
Since $|\alpha| -|\gamma| -|\mu| + \frac{2}{r}<0$ we have
$\left(1 +|\xi|\right)^{ |\alpha| -|\gamma| -|\mu| + \frac{2}{r}} \leq 1$, 
then 
\begin{linenomath}
	\begin{multline*}
     \frac{1}{\mu!} 
     C_{\psi}^{|\beta|+|\mu| +n+3} N^{ \left( |\beta|+|\mu|+ n + 3 -M\right)^{+}} 
     \| \widetilde{\Theta}_{\widetilde{m}}^{(\mu)}\Theta_{N}^{(\gamma)} D^{\alpha} g \|_{\frac{1}{r}}
     \\
     \leq
     \frac{1}{\mu!} 
     C_{\psi}^{|\beta|+|\mu| +n+3}  C_{\widetilde{\Theta}}^{|\mu| + 1} C_{\Theta}^{|\gamma| + 1} 
     N^{ \left( |\beta|+2|\mu|+|\gamma|+ n + 3 -M\right)^{+}} \|g\|_{0}
     \\
     \leq 
     \frac{1}{\mu!} 
     C_{1}^{2(|\alpha|-|\gamma|)+|\gamma|+1}  C_{2}^{|\beta| + 1} 
     N^{ \left( 2(|\alpha|-|\gamma|)+|\gamma|+|\beta|+ n + 3 -M\right)^{+}} \|g\|_{0},
	\end{multline*}
\end{linenomath}	
where $C_{1}$ and $C_{2}$ are suitable constants independent of $N$, $\alpha$, $\beta$ and $\gamma$.\\
Case $ |\mu|\leq |\alpha| - |\gamma|$:
since $ \widetilde{\Theta}_{\widetilde{m}}^{(\mu)} (\xi)\Theta_{N}^{(\gamma)} (\xi)$
is supported in the region $\left\{ \xi\in \mathbb{R}^{n}: \frac{N}{2} \leq |\xi| \leq \widetilde{m} \right\}$
we have
\begin{linenomath}
	\begin{multline*}
	\frac{1}{\mu!} 
	C_{\psi}^{|\beta|+|\mu| +n+3} N^{ \left( |\beta|+|\mu|+ n + 3 -M\right)^{+}} 
	\| \widetilde{\Theta}_{\widetilde{m}}^{(\mu)}\Theta_{N}^{(\gamma)} D^{\alpha} g \|_{\frac{1}{r}}
	\\
	\leq
	\frac{1}{\mu!} 
	C_{\psi}^{|\beta|+|\mu| +n+3}  C_{\widetilde{\Theta}}^{|\mu| + 1} C_{\Theta}^{|\gamma| + 1} 
	 N^{ \left( |\beta|+|\mu|+|\gamma|+ n + 3 -M\right)^{+}} 
	 \widetilde{m}^{ \left(|\mu| -M\right)^{+}} 
	 \left(1 +\widetilde{m}\right)^{ |\alpha| -|\gamma| -|\mu| + 1} \|g\|_{0}
	 \\
	 \leq 
	 \frac{1}{\mu!} 
	 C_{\psi}^{|\beta|+|\mu| +n+3}  C_{\widetilde{\Theta}}^{|\mu| + 1} C_{\Theta}^{|\gamma| + 1} 
	 N^{ \left( |\beta|+|\mu|+|\gamma|+ n + 3 -M\right)^{+}} 
	 \left(1 +\widetilde{m}\right)^{ |\alpha| -|\gamma| + 1} \|g\|_{0}.
	\end{multline*}
\end{linenomath}	
 By Lemma \ref{L-1} there is a constant $C_{3}$
 such that 
 \begin{linenomath}
 $$
  (1+\widetilde{m})^{|\alpha| -|\gamma|  + 1} \leq C_{3}^{|\alpha| -|\gamma| + 1} N^{ \left( 1+\widetilde{m} -M\right)^{+}}
 = C_{3}^{|\alpha| - |\gamma| +1} N^{ \left( |\alpha|-|\gamma|+\left\lfloor \frac{n}{2} \right\rfloor +2 -M\right)^{+}},
 $$
 \end{linenomath}
then the right hand side of the above inequality can be estimated by
\begin{linenomath}
     \begin{align*} 
     \frac{1}{\mu!} 
     C_{4}^{2(|\alpha|-|\gamma|)+|\gamma|+1}  C_{5}^{|\beta| + 1} 
     N^{ \left( 2(|\alpha|-|\gamma|)+|\gamma|+|\beta|+ 2n + 5 -M\right)^{+}} \|g\|_{0}.
     \end{align*}
\end{linenomath}
Since  $\| g\|_{0} = \| \varphi_{N}^{(\delta)} P^{k} u \|_{0} $, $u$ $G^{s}$-vector for $P$ in $\Omega$,
by \eqref{Est_g} and taking advantage from Lemma \ref{L-1} and Remark \ref{Rk-1}
we conclude that there are two positive constants
$\widetilde{C}_{1}$ and $\widetilde{C}_{2}$ such that
\begin{linenomath}
	\begin{multline}
	\label{Est-Sum-S1-f}
	\displaystyle\longsum[25]_{1 \leq |\mu| \leq \widetilde{m} - 1}\frac{1}{\mu!} 
	\|  \psi_{N}^{(\beta+\mu)} \widetilde{\Theta}_{\widetilde{m}}^{(\mu)}\Theta_{N}^{(\gamma)} D^{\alpha} g \|_{\frac{1}{r}}
	\\
	\leq
	\widetilde{C}_{1}^{\sigma+ 1}  \widetilde{C}_{2}^{2 m +|\gamma|+ 2n + 4} 
	N^{ s\left( 2m+|\gamma| + \sigma + 2 n + 5 -M\right)^{+}},
	\end{multline}
\end{linenomath}
where $m= |\alpha|- |\gamma|$ and $\sigma= |\beta|+|\delta|+ 2k$.

\textbf{Estimate of the last term on the right hand side of (\ref{Est-1/r-0}).}\\
Since  $(1- \widetilde{\Theta}_{\widetilde{m}}) \Theta_{N}^{(\gamma)} $
is supported in the region  $\left\{ \xi\in \mathbb{R}^{n}: \frac{N}{2} \leq |\xi| \leq \widetilde{m} \right\}$,
using the same strategy used to handle the terms in the sum,
we conclude that there are two positive constants $C_{2}$ and $C_{3}$ such that
\begin{linenomath}
	\begin{multline}
	\label{Est-last-t}
	\|\psi_{N}^{(\beta)} (1- \widetilde{\Theta}_{\widetilde{m}}) w \|_{\frac{1}{r}}
	\\
	\leq
	C_{1}^{ |\beta|+|\delta|+ 2k+ 1}  C_{2}^{m +|\gamma|+ n + 4} 
	N^{ s\left( m + |\beta|+  |\gamma| + |\delta| +2k + 2 n + 3 -M\right)^{+}}.
	\end{multline}
\end{linenomath}

By (\ref{Est-1/r-0}), (\ref{Est-1/r-3}),  (\ref{Rem0-Est}), (\ref{Est-Sum-S1-f}) and (\ref{Est-last-t}),
we conclude that there are positive constants $C$, $ A $ and $ B $
independent of $N$, $\alpha$, $\beta$, $\gamma$, $\delta$ and $k$
such that the estimate (\ref{Est-1/r-1-p}) holds.
This concludes the proof of the Proposition \ref{Basic-Est-D}.
%
\section{Estimate for the associated microlocal sequence of a Gevrey vector of $P$}
\renewcommand{\theequation}{\thesection.\arabic{equation}}
\setcounter{equation}{0} \setcounter{theorem}{0}
\setcounter{proposition}{0} \setcounter{lemma}{0}
\setcounter{corollary}{0} \setcounter{definition}{0}
\setcounter{remark}{1}
The purpose of the present section is to obtain a suitable estimate 
of $\| v \|_{p/r}^{2}$, $p=1, 2,\dots, r-1$, where
$ v =\psi_{N}^{(\beta)} \Theta_{N}^{(\gamma)} D^{\alpha} \varphi_{N}^{(\delta)} P^{k} u$,
here $\psi_{N}$, $\Theta_{N}$ and $\varphi_{N}$
are as in the previous section.
This will allow us to obtain in the next section the microlocal regularity
of $u$ at the point $(x_{0},\xi_{0}) \in \Char(P)$.\\
We will use the same notation of the previous section:
$ P^{k}u =: f$, $  \varphi_{N}^{(\delta)} f =: g$, $ \Theta_{N}^{(\gamma)} D^{\alpha}g =: w$
and $\psi_{N}^{(\beta)} w =: v$.
\vspace{-1em}
\subsection{Estimates in $H^{1/r}$}~\par
\vskip-4mm
Our goal is to obtain a suitable estimate for 
$ \|v\|_{1/r}=\| \psi_{N}^{(\beta)} \Theta_{N}^{(\gamma)} D^{\alpha} \varphi_{N}^{(\delta)} P^{k} u\|_{1/r}$.
In order to obtain it we take advantage from Proposition \ref{Basic-Est-D}. 
To characterize the microlocal regularity of the $s$-Gevrey vector $u$
at $(x_{0},\xi_{0})$ ($ \tau\left( \left(x_{0},\xi_{0}\right); \! \mathscr{P} \right)= r$)
we have to obtain a suitable estimate of the left hand
side of (\ref{Est-1/r-1-p}) in Proposition \ref{Basic-Est-D}.
In order to make more readable the manuscript  we recall it
\begin{linenomath}
	\begin{multline}\label{Est-1/r-1}
	\|  v \|_{\frac{1}{r}}^{2}
	+\sum_{j=1}^{n} \left( \| P^{j} v\|^{2} + \| P_{j} v\|_{-1}^{2} \right)
	\leq 
	C \! \left( \sum_{\ell=0}^{n}\! \left|\langle E_{\ell} P v, E_{\ell} v \rangle\right|
	+  \|v\|^{2}
	\right)
	\\
	\quad
	+ A^{2(|\gamma|+1)} B^{2(2 m + \sigma + 1)}  
	N^{2s(2 m+ |\gamma|+\sigma) },
	\end{multline}
\end{linenomath}
where $A$, $B$ and $C$ are suitable constants
independent of $\alpha$, $\beta$, $\gamma$, $\delta$, $k$ and $N$,
$P^{j}$, $P_{j}$ and $E_{\ell}$, $\ell=0,\, \dots, n$, are the same as in the 
Theorem \ref{MicLocEst}, $m= |\alpha| - |\gamma| $ 
and $\sigma= |\beta| + |\delta| +2k$.\\
Now we give the estimate of the terms $|\langle E_{\ell}P v, E_{\ell} v\rangle|$ $(\ell =0,1,\dots,n)$ in \eqref{Est-1/r-1}.
We begin to handle the first term when $\ell =0$, i.e. $|\langle E_{0}P v, E_{0} v\rangle| = |\langle P v, v\rangle|$.\\
For technical reasons, that is in order to handle the commutators
of $P$, $P_{j}$ and $P^{j}$ with the pseudodifferential operators
$\Theta_{_{N}}(D)D^{\alpha}$ and $E_{\ell}$
we introduce a new Ehrenpreis-H\"ormander sequence $\widetilde{\psi}_{N}$
associated to the couple $(\Omega_{2},\Omega_{3})$,
where $\Omega_{1}\Subset\Omega_{2}\Subset \Omega_{3} \Subset \widetilde{\Omega}_{3}\Subset \Omega_{4}$,
that is supported in $\Omega_{3}$ and identically one on the closure of $\Omega_{2}$.
So $\widetilde{\psi}_{N}$ is identically one on the supports of
$\psi_{N}$ and  $\varphi_{N}$, in particular
$\psi_{N}\widetilde{\psi}_{N} = \psi_{N}$ and $\varphi_{N}\widetilde{\psi}_{N}= \varphi_{N}$.
We set
\begin{linenomath}
	\begin{align}\label{t-PN}
	\nonumber
	&\widetilde{P}_{N}(x,D) = \longsum[7]_{j_{1},j=1 }^{n} \tilde{a}_{N,j,j_{1}}(x)D_{j}D_{j_{1}}
	+ \sum_{ j_{1}=1}^{n} i \tilde{b}_{N,j_{1}}(x) D_{\ell} + \tilde{c}_{N}(x);
	\\
	&
	\widetilde{P}^{j}_{N}(x,D)  = 2 \sum_{j_{1} =1}^{n} \widetilde{a}_{N,j_{1},j}(x) D_{j_{1}} ;
	\\
	\nonumber
	&\widetilde{P}_{N,j}(x,D) = \longsum[7]_{j_{1},j_{2} =1}^{n}\widetilde{a}_{N,j_{1},j_{2}}^{(j)}D_{j_{1}}D_{j_{2}}; 
	\end{align}
\end{linenomath}
where $ \tilde{a}_{N,j, j_{1}}(x)=\widetilde{\psi}_{N}(x) a_{j,j_{1}}(x) $,
$\tilde{b}_{N,\ell}(x)=\widetilde{\psi}_{N}(x)b_{j}(x) $ and
$\tilde{c}_{N}(x)=\widetilde{\psi}_{N}(x) c(x)$.
We remark that $\psi_{N}\widetilde{P}_{N} =\psi_{N} P$, $\widetilde{P}_{N}\varphi_{N}= P\varphi_{N}$,
$[\psi_{N},\widetilde{P}_{N}] =[\psi_{N}, P]$ and  $[\widetilde{P}_{N},\varphi_{N}]= [P,\varphi_{N}]$,
the same holds if we replace $\widetilde{P}_{N}$ by $\widetilde{P}^{j}_{N}$ or $\widetilde{P}_{N,j}$.\\
Since
\begin{linenomath}
	\begin{multline}
	\label{Com-Pv}
	P v = 
	[ P, \psi_{N}^{(\beta)} ]  \Theta^{(\gamma)}_{N} D^{\alpha} g 
	+ \psi_{N}^{(\beta)} [ \widetilde{P}_{N},  \Theta^{(\gamma)}_{N} D^{\alpha} ] g
	+ \psi_{N}^{(\beta)}  \Theta^{(\gamma)}_{N} D^{\alpha} [ P, \varphi_{N}^{(\delta)}] f 
	\\
	+ \psi_{N}^{(\beta)}  \Theta^{(\gamma)}_{N} D^{\alpha} \varphi_{N}^{(\delta)} P^{k+1}u,
	\end{multline}
\end{linenomath}
we have
\begin{linenomath}
	\begin{multline}\label{Term-Pvv}
	\left|\langle P v, v\rangle \right|
	\leq 
	\left|\langle [ P, \psi_{N}^{(\beta)} ]  \Theta^{(\gamma)}_{N} D^{\alpha} g, v\rangle \right|
	+
	\left|\langle  \psi_{N}^{(\beta)} [ \widetilde{P}_{N} , \Theta^{(\gamma)}_{N} D^{\alpha} ] g, v\rangle \right|
	\\
	+
	\left|\langle  \psi_{N}^{(\beta)} \Theta^{(\gamma)}_{N} D^{\alpha} [ P, \varphi_{N}^{(\delta)}] f, v \rangle \right|
	+
	\left|\langle  \psi_{N}^{(\beta)} \Theta^{(\gamma)}_{N} D^{\alpha} \varphi_{N}^{(\delta)} P^{k+1} u, v \rangle \right|
	\\
	=\sum_{j=1}^{4} I_{j}.
	\end{multline}
\end{linenomath}
\textbf{Estimate of the term} $I_{1}$. We have
\begin{linenomath}
	\begin{multline}\label{Cm_Ppsi}
	[ P,  \psi_{N}^{(\beta)}]=  \sum_{\ell,j=1 }^{n} a_{j,\ell}(x)[D_{j}D_{\ell}, \psi_{N}^{(\beta)}] + i \sum_{ \ell=1}^{n} b_{\ell}(x) [D_{\ell}, \psi_{N}^{(\beta)} ] 
	\\
	=  \sum_{j=1}^{n}  P^{j} \psi_{N}^{(\beta+j)} + \sum_{\ell,j =1}^{n} a_{j,\ell}(x) \psi_{N}^{(\beta+\ell+j)}
	+  i\sum_{ \ell=1}^{n} b_{\ell}(x) \psi_{N}^{(\beta+\ell)}
	\\
	= 
	 \sum_{j=1}^{n}  P^{j} \psi_{N}^{( \beta+j)} + \sum_{|\mu| = 2 }  P^{\mu}  \psi_{N}^{(\beta+\mu)}
	+  i\sum_{ \ell=1}^{n} b_{\ell}(x) \psi_{N}^{(\beta+\ell)},
	\end{multline}
\end{linenomath}
where $ P^{\mu}= \partial_{\xi}^{\mu}p^{0}(x,\xi) =\partial_{\xi_{j}}\partial_{\xi_{\ell}}p^{0}(x,\xi)= a_{j,\ell}(x)$.\\
Due to the fact that$ \left( P^{j} \right)^{*} =  P^{j} + 2  \sum_{\ell =1}^{n} a_{j,\ell}^{(\ell)}(x) $ we obtain
\begin{linenomath}
	\begin{multline}\label{I_1}
	I_{1} 
	\leq \sum_{j=1}^{n} \left|\langle P^{j} \psi^{(\beta+j)}_{N} \Theta_{N}^{(\gamma) } D^{\alpha} g, v \rangle \right|
	+ \sum_{|\mu| =2} \left|\langle P^{\mu} \psi_{N}^{(\beta+\mu)} \Theta_{N}^{(\gamma)} D^{\alpha} g, v\rangle\right|
	\\
	+  \sum_{\ell =1}^{n} \left|\langle b_{\ell} \psi_{N}^{(\beta+\ell)}  \Theta_{N}^{(\gamma) } D^{\alpha} g, v \rangle\right|
	\\
	\leq \varepsilon \sum_{j=1}^{n} \| P^{j} v\|^{2} + C_{\varepsilon} \sum_{j=1}^{n} \| \psi_{N}^{(\beta+j)} \Theta_{N}^{(\gamma)} D^{\alpha} g\|^{2}
	+ 2\sum_{\ell,j =1}^{n}\| a_{j,\ell}^{(\ell)}\psi_{N}^{(\beta+j)} \Theta_{N}^{(\gamma)} D^{\alpha} g\| \|v\| 
	\\
	+  \sum_{ \ell=1}^{n} \| b_{\ell} \psi_{N}^{(\beta+\ell)}  \Theta_{N}^{(\gamma)} D^{\alpha} g\| \|v\| 
	+ \sum_{|\mu| =2} \| P^{\mu} \psi_{N}^{(\beta+\mu)} \Theta_{N}^{(\gamma)} D^{\alpha} g\| \| v\| 
	\\
	\leq \varepsilon \sum_{j=1}^{n} \| P^{j} v\|^{2} + C_{\varepsilon} \sum_{j=1}^{n} \| \psi_{N}^{(\beta+j)} \Theta_{N}^{(\gamma)} D^{\alpha} g\|^{2}
	+
	C_{2} \sum_{j =1}^{n}\| \psi_{N}^{(\beta+j)} \Theta_{N}^{(\gamma)} D^{\alpha} g\| \|v\|
	\\
	+C_{3}  \sum_{|\mu| =2} \|\psi_{N}^{(\beta+\mu)} \Theta_{N}^{(\gamma)} D^{\alpha} g\| \| v\| ,
	\end{multline}
\end{linenomath}
where $\varepsilon$ is a small suitable constant and $C_{\varepsilon}= \varepsilon^{-1}$.
We point out that the first term on the right
hand side can be absorbed by the left hand side of (\ref{Est-1/r-1}).

\textbf{Estimate of the term $I_{2}$ in (\ref{Term-Pvv})}. We have
\begin{linenomath}
	\begin{multline*}
	[\widetilde{P}_{N}, \Theta_{N}^{(\gamma)} D^{\alpha}]  
	= \sum_{j=1}^{n} \left(\widetilde{P}_{N}\right)^{(j)} \left(\Theta_{N}^{(\gamma)} D^{\alpha} \right)^{(j)} 
	+ \longsum[32]_{2 \leq |\mu| \leq  |\alpha| -|\gamma| + 1} \frac{1}{\mu!} \widetilde{P}_{N,\mu}\left(\Theta_{N}^{(\gamma)} D^{\alpha} \right)^{(\mu)}
	\\
	+ \mathscr{R}_{|\alpha| -|\gamma| + 2}\left( [\widetilde{P}_{N}, \Theta_{N}^{(\gamma)} D^{\alpha}]  \right), 
	\end{multline*}
\end{linenomath}
where
\begin{linenomath}
	\begin{align*}
	&\left( \widetilde{P}_{N} \right)^{(j)} \!\!\!= \!\!\sum_{\ell,j_{1}=1}^{n} \widetilde{a}_{N,\ell,j_{1}}^{(j)}(x) D_{\ell} D_{j_{1}} 
	+  i \sum_{\ell=0}^{n} \widetilde{b}_{N,\ell}^{(j)}(x) D_{\ell} +  \widetilde{c}_{N}^{(j)}(x)
	\\
	&\qquad\quad
	=\widetilde{\psi}_{N} \left( P_{j} +    i\sum_{\ell=0}^{n} b_{\ell}^{(j)}(x) D_{\ell} + c^{(j)}(x)\right)
	+ \widetilde{\psi}_{N}^{(j)}P,
	\\	    
	&\widetilde{P}_{N,\mu}=\sum_{\ell,j=1}^{n} \widetilde{a}_{N,\ell,j}^{(\mu)}(x) D_{\ell} D_{j} 
	+ i\sum_{\ell=0}^{n} \widetilde{b}_{N,\ell}^{(\mu)}(x) D_{\ell} + \widetilde{c}_{N}^{(\mu)}(x),
	\end{align*}
\end{linenomath}	
and	$\mathscr{R}_{|\alpha| -|\gamma| + 2}\left( [\widetilde{P}_{N}, \Theta_{N}^{(\gamma)} D^{\alpha}] \right)$
is a pseudodifferential operator with associated symbol
\begin{linenomath}
	\begin{multline}
	\label{Rem_PThD}
	\mathscr{R}_{|\alpha| -|\gamma| + 2}\left( [ \widetilde{P}_{N}, \Theta_{N}^{(\gamma)} D^{\alpha}]  \right)(x,\xi) 
	\\
	=
	\frac{1}{(2\pi)^{n}}  \longsum[16]_{ |\mu| = |\alpha| -|\gamma| + 2} \frac{ |\alpha| -|\gamma| + 2}{\mu!}
	\iint e^{i(x-y)(\eta-\xi)} (\eta-\xi)^{\mu} \widetilde{P}_{N}(y,\xi) 
	\\
	\qquad \times
	\int_{0}^{1}\!\!\!\!\! \left(1-t\right)^{|\alpha| -|\gamma| + 1}
	\left(\sigma\left(\Theta_{N}^{(\gamma)}D^{\alpha}\right)\right)^{(\mu)} (\xi+t(\eta-\xi))dt dy d\eta,
	\end{multline}
\end{linenomath}
here $\sigma\left(\Theta_{N}^{(\gamma)}D^{\alpha}\right) $
denotes the symbol associated to the operator $\Theta_{N}^{(\gamma)}D^{\alpha}$.\\
Since
\begin{linenomath}	
	\begin{align*}
	\left(\Theta_{N}^{(\gamma)} D^{\alpha} \right)^{(\mu)}
	=  \longsum[7]_{\substack{\nu\leq \mu \\ \nu\leq \alpha}} \binom{\mu}{\nu} \frac{\alpha!}{(\alpha -\nu)!} \Theta_{N}^{(\gamma +\mu-\nu)} D^{\alpha-\nu},
	\end{align*}
\end{linenomath}
we conclude
\begin{linenomath}	
	\begin{multline}
	\label{Com-PThetaD}
	[\widetilde{P}_{N}, \Theta_{N}^{(\gamma)} D^{\alpha}] 
	\\
	=
	\sum_{j=1}^{n}\left[ \widetilde{\psi}_{N}\left( P_{j} + i\sum_{\ell=0}^{n}  b_{\ell}^{(j)}(x) D_{\ell} + c^{(j)}(x) \right) +  \widetilde{\psi}_{N}^{(j)}P\right] 
	\left(\Theta_{N}^{(\gamma +j)} D^{\alpha} + \alpha_{j} \Theta_{N}^{(\gamma)} D^{\alpha-j}\right)
	\\
	\qquad\qquad\quad
	+ \longsum[32]_{2 \leq |\mu| \leq |\alpha| -|\gamma| + 1} \longsum[7]_{\substack{\nu\leq \mu \\ \nu\leq \alpha}}
	\frac{1}{\mu!}  \binom{\mu}{\nu} \frac{\alpha!}{(\alpha -\nu)!} \widetilde{P}_{N,\mu} \Theta_{N}^{(\gamma +\mu-\nu)} D^{\alpha-\nu}
	\\
	\qquad\qquad\qquad\qquad\qquad\qquad\qquad\qquad\quad
	+ \mathscr{R}_{|\alpha| -|\gamma| + 2}\left( [ \widetilde{P}_{N}, \Theta_{N}^{(\gamma)} D^{\alpha}]  \right).
	\end{multline}
\end{linenomath}
Since  $\psi_{N}^{(\beta)}\widetilde{\psi}_{N}^{(\mu)} = 0$
for any non zero $\mu$ and any $\beta$, we have
\begin{linenomath}	
	\begin{multline}\label{I_2-1}
	I_{2} 
	\leq \sum_{j=1}^{n} \left|\langle \psi_{N}^{(\beta)} P_{j} \left(\Theta_{N}^{(\gamma)} D^{\alpha}\right)^{(j)} \!\!\!g, v \rangle \right|
	+ 
	\sum_{j=1}^{n}
	\left|\langle  \psi_{N}^{(\beta)} \left(B_{j} + c^{(j)}(x) \right) \left(\Theta_{N}^{(\gamma)} D^{\alpha}\right)^{(j)} g , v \rangle \right|
	\\
	\qquad\quad
	+ \longsum[32]_{2 \leq |\mu| \leq |\alpha| -|\gamma| + 1} \longsum[7]_{\substack{\nu\leq \mu \\ \nu\leq \alpha}}
	\frac{1}{\mu!}  \binom{\mu}{\nu} \frac{\alpha!}{(\alpha -\nu)!}
	\left|\left(  \psi_{N}^{(\beta)} P_{\mu} \Theta_{N}^{(\gamma +\mu-\nu)} D^{\alpha-\nu} g, v \right)\right|
	\\
	\qquad\qquad\qquad\qquad\qquad\qquad\qquad
	+\left|\left(  \psi_{N}^{(\beta)} \mathscr{R}_{|\alpha| -|\gamma| + 2}\left( [ \widetilde{P}_{N}, \Theta_{N}^{(\gamma)} D^{\alpha}]  \right)  g, v \right)\right|
	\\
	\\
	= I_{2,1} + I_{2,2} + I_{2,3} + I_{2,4},
	\end{multline}
\end{linenomath}
where
$ P_{\mu} = \sum_{\ell,j_{1}=1}^{n} a_{\ell,j_{1}}^{(\mu)} D_{\ell} D_{j_{1}} + i \sum_{\ell=0}^{n} b_{\ell}^{(\mu)} D_{\ell} +c^{(\mu)}$
and $B_{j} =  \sum_{\ell=0}^{n} b_{\ell}^{(j)} D_{\ell}$.\\
We estimate each of the terms obtained separately.\\
We observe that
\begin{linenomath}
	\begin{align*}
	\left(P_{j}\right)^{*}
	= P_{j}
	+ \sum_{\ell,j_{1}=1}^{n} \left(a_{\ell,j_{1}}^{(j+\ell)} D_{j_{1}} +  a_{\ell,j_{1}}^{(j+j_{1})} D_{\ell}\right)
	+\sum_{\ell,j_{1}=1}^{n} a_{\ell,j_{1}}^{(j+j_{1}+ \ell)}, 
	\end{align*}
\end{linenomath}
and
\begin{linenomath}
	\begin{align*}
	\left[\psi_{N}^{(\beta)}, P_{j}\right] = 
	\sum_{\ell,j_{1}=1}^{n} a_{\ell,j_{1}}^{(j)} \left(\psi_{N}^{(\beta+\ell)}D_{j_{1}} + \psi_{N}^{(\beta+j_{1})}D_{\ell}- \psi_{N}^{(\beta+j_{1}+\ell)}\right).
	\end{align*}
\end{linenomath}
Then the terms $I_{2,j}$, $j= 1,2,3,4$ can be estimated as follows
\vspace{-1.7em}
\begin{description}
\item[\hspace{-2em}Term $I_{2,1}$] 
  \begin{linenomath}
		\begin{multline}\label{Est_I21}
		I_{2,1} 
		\leq 
		\sum_{j=1}^{n} \left|\langle
		\psi_{N}^{(\beta)}\left(\Theta_{N}^{(\gamma) } D^{\alpha}\right)^{(j)} g,  P_{j} v \rangle\right|
		\\
		+
		\sum_{j=1}^{n} \left|\langle  \left( P_{j} - P_{j}^{*} \right)
		\psi_{N}^{(\beta)}\left(\Theta_{N}^{(\gamma) } D^{\alpha}\right)^{(j)} g, v \rangle\right|
		+  \sum_{j=1}^{n} \left|\langle
		\left[\psi_{N}^{(\beta)}, P_{j}\right] \left(\Theta_{N}^{(\gamma) } D^{\alpha}\right)^{(j)} g, v \rangle\right|
		\\
		\leq  \varepsilon  \sum_{j=1}^{n}\| P_{j} v \|_{-1}^{2} + C_{\varepsilon} \sum_{j=1}^{n} \| \psi_{N}^{(\beta)} (\Theta_{N}^{(\gamma) } D^{\alpha} )^{(j)} g \|_{1}^{2}
		+
		C_{3}  \sum_{j=1}^{n} \| \psi_{N}^{(\beta)} (\Theta_{N}^{(\gamma) } D^{\alpha})^{(j)} g \| \| v\|
		\\
		\quad
		+ 2C_{4}  \sum_{j,\ell=1}^{n} \| \psi_{N}^{(\beta)}  (\Theta_{N}^{(\gamma) } D^{\alpha})^{(j)}D_{\ell} g \| \| v \|
		+ 2C_{5}  \longsum[10]_{j,j_{1},\ell=1}^{n} \| \psi_{_{N}}^{(\beta+j_{1})}  (\Theta_{N}^{(\gamma) } D^{\alpha})^{(j)}D_{\ell} g \| \| v\|
		\\
		\quad
		+ C_{5}  \longsum[9]_{j,j_{1},\ell=1}^{n} \| \psi_{N}^{(\beta+j_{1} +\ell)}  (\Theta_{N}^{(\gamma) } D^{\alpha})^{(j)} g \| \| v\|,
		\end{multline}
  \end{linenomath}
 where $\varepsilon$ is a small suitable constant, 
 $C_{\varepsilon}= \varepsilon^{-1}$, 
 $C_{3} = \displaystyle\sup_{\ell,j_{1}, |\mu|=3}\{|a_{\ell,j_{1}}^{(\mu)|}\}$,
 $C_{4} = \displaystyle\sup_{\ell,j_{1}, |\mu|=2}\{|a_{\ell,j_{1}}^{(\mu)}|\}$
 and $C_{5} = \displaystyle\sup_{\ell,j_{1}, |\mu|=1}\{|a_{\ell,j_{1}}^{(\mu)}|\}$. 
 We point out that  a suitable choice of $\varepsilon$ will allow to absorb the first term on the right hand side by the left hand side of (\ref{Est-1/r-1}).
\item[\hspace{-2em}Term $I_{2,2}$] 
   \begin{linenomath}
		\begin{align}\label{Est_I22}
		I_{2,2}
		\leq 
		C_{6}  \sum_{\ell,j=1}^{n} \| \psi_{N}^{(\beta)} (\Theta_{N}^{(\gamma) } D^{\alpha})^{(j)}D_{\ell} g \| \| v \|
		+ C_{7}  \sum_{j=1}^{n} \| \psi_{N}^{(\beta)}  (\Theta_{N}^{(\gamma) } D^{\alpha})^{(j)} g \| \| v \|
		\end{align}
	\end{linenomath}
  where $C_{6} = \displaystyle\sup_{j,\ell} \{|b_{\ell}^{(j)}|\}$ and $C_{7} = \displaystyle\sup_{j} \{|c^{(j)}| \}$.
\item[\hspace{-2em}Term $I_{2,3}$] 
	\begin{linenomath}
		\begin{multline}\label{Est_I23}
		I_{2,3}
		\leq
		\longsum[32]_{2 \leq |\mu| \leq |\alpha| -|\gamma| + 1} \longsum[7]_{\substack{\nu\leq \mu \\ \nu\leq \alpha}}
		\frac{1}{\left(\mu-\nu\right)!}  \binom{\alpha}{\nu}
		\| v \|
		\\
		\times\left(
		\longsum[8]_{j,\ell=1}^{n} \|  a_{\ell,j}^{(\mu)} \psi_{N}^{(\beta)}  \Theta_{N}^{(\gamma +\mu-\nu)} D^{\alpha-\nu}D_{j}D_{\ell} g \|
		+
		\sum_{ \ell=1}^{n} \|  b_{\ell}^{(\mu)} \psi_{N}^{(\beta)}  \Theta_{N}^{(\gamma +\mu-\nu)} D^{\alpha-\nu}D_{\ell} g \|
		\right.
		\\
		+
		\|  c^{(\mu)} \psi_{N}^{(\beta)}  \Theta_{N}^{(\gamma +\mu-\nu)} D^{\alpha-\nu}g \|
		\Biggr)
		.
		\end{multline}
	\end{linenomath}
	We stress that the order of $ \Theta_{_{N}}^{(\gamma +\mu-\nu)} D^{\alpha-\nu}D_{j}D_{\ell}$ is less or equal $|\alpha|-|\gamma|$.
\item[\hspace{-2em}Term $I_{2,4}$] 
	\begin{linenomath}
		\begin{align}\label{Est_I24_0}
		I_{2,4}
		\leq
		\|  \psi_{N}^{(\beta)} \mathscr{R}_{|\alpha| -|\gamma| + 2}\left( [ \widetilde{P}_{N}, \Theta_{N}^{(\gamma)} D^{\alpha}]  \right)  g \|  \| v \|.
		\end{align}
	\end{linenomath}
\end{description}
In order to estimate the first factor on the right hand side of \eqref{Est_I24_0}
we take advantage from the Theorem $18.1.11'$, page 75, in \cite{H_Book-3}.
We rewrite (\ref{Rem_PThD}) more explicitly   
\begin{linenomath}
	\begin{multline*}
	\mathscr{R}_{|\alpha| -|\gamma| + 2}\left( [ \widetilde{P}_{N}, \Theta_{N}^{(\gamma)} D^{\alpha}]  \right)(x,\xi) 
	\\
	=
	\frac{1}{(2\pi)^{n}}  \longsum[16]_{ |\mu| = |\alpha| -|\gamma| + 2} \frac{ |\alpha| -|\gamma| + 2}{\beta!}
	\iint e^{i(x-y)(\eta-\xi)} 
	\\
	\qquad \times
	\left( \sum_{\ell,j=1}^{n} \widetilde{a}_{N,\ell,j}^{(\mu )}(y)\xi_{\ell} \xi_{j} 
	+\sum_{\ell=1}^{n} \widetilde{b}_{N,\ell}^{(\mu)}(y)\xi_{\ell} + \widetilde{c}_{N}^{(\mu)}(y)
	\right)
	\\
	\qquad \times
	\int_{0}^{1}\!\!\!\!\! \left(1-t\right)^{|\alpha| -|\gamma| + 1}
	\left(\sigma\left(\Theta_{N}^{(\gamma)}D^{\alpha}\right)\right)^{(\mu)} (\xi+t(\eta-\xi))dt dy d\eta.
	\end{multline*}
\end{linenomath}
Now we want to bound the term under integral:
\begin{linenomath}
	\begin{multline*}
	\left|\left(\sigma\left(\Theta_{N}^{(\gamma)}D^{\alpha}\right)\right)^{(\mu)}(\xi+t(\eta-\xi))\right|
	\\
	\leq
	\longsum[7]_{\substack{\nu\leq \mu \\ \nu\leq \alpha}} \binom{\mu}{\nu} 
	\frac{\alpha!}{(\alpha -\nu)!} \left|\Theta_{N}^{(\gamma +\mu-\nu)}(\xi+t(\eta-\xi))\right| 
	\left|\xi+t(\eta-\xi)\right|^{|\alpha-\nu|}.
		\end{multline*}
\end{linenomath}
By \eqref{ThetaN-AH} in Appendix and since $\frac{\alpha!}{(\alpha -\nu)!} \leq \mu! 2^{|\alpha|}$,  the above term is bounded by
\begin{linenomath}
	\begin{align*}
	2^{|\alpha|} \mu!
	\longsum[6]_{\substack{\nu\leq \mu \\ \nu\leq \alpha}} \binom{\mu}{\nu} 
	C_{\Theta}^{|\gamma| +|\mu|-|\nu| +1} N^{(|\gamma| +|\mu| -|\nu|-M)^{+}} 
	\left( 1 + |\xi+t(\eta-\xi)|\right)^{|\alpha|-|\gamma| -|\mu|}.
		\end{align*}
\end{linenomath}
Using that $|\mu|= |\alpha|-|\gamma| +2$ it can be rewritten as
\begin{linenomath}
	\begin{align*}
	2^{|\alpha|} \mu!
	\longsum[6]_{\substack{\nu\leq \mu \\ \nu\leq \alpha}} \binom{\mu}{\nu} 
	C_{\Theta}^{|\alpha| + 3 -|\nu|} N^{(|\alpha| + 2 -|\nu|-M)^{+}} 
	\left( 1 + |\xi+t(\eta-\xi)|\right)^{-2}.
	\end{align*}
\end{linenomath}
By the Peetre's inequality and since $|\nu|\geq 0$ and $\sum_{\nu \leq \mu} \binom{\mu}{\nu} = 2^{|\mu|}$,
we conclude that there is a suitable positive constant $\tilde{C}$ such that
\begin{linenomath}
	\begin{multline*}
	\left|\left(\sigma\left(\Theta_{N}^{(\gamma)}D^{\alpha}\right)\right)^{(\mu)}(\xi+t(\eta-\xi))\right|
	\\
	\leq 
	\tilde{C}^{ |\alpha| + 1} \,\mu! \,N^{(|\alpha| +2-M)^{+}} \left( 1+ |\eta-\xi|^{2}\right)\left( 1+ |\xi|^{2}\right)^{-1}, \quad \forall  t \in \left[0,1\right].
	\end{multline*}
\end{linenomath}
Moreover since 
\begin{linenomath}
	\begin{multline*}
	D_{x}  e^{i(x-y)(\eta-\xi)}= -D_{y} e^{i(x-y)(\eta-\xi)}
	\text{ and } 
	\frac{(1+ \Delta_{y})^{\lfloor\frac{n}{2}\rfloor+2}}{\left(1+ |\eta -\xi|^{2}\right)^{\lfloor\frac{n}{2}\rfloor+2}} e^{-iy(\eta-\xi)}
	=e^{-iy(\eta-\xi)},
	\end{multline*}
\end{linenomath}
we obtain, using Lemma \ref{L-1} and  Remark \ref{Rk-1}, that
\begin{linenomath}
	\begin{multline*}
	\sup_{\xi\in \mathbb{R}^{n}}\longsum[12]_{\mu_{1} \leq n+1} \int \left| D_{x}^{\mu_{1}}
	\psi_{N}^{(\beta)} (x)\mathscr{R}_{|\alpha| -|\gamma| + 2}\left( [ \widetilde{P}_{N}, \Theta_{N}^{(\gamma)} D^{\alpha}]  \right) (x,\xi)\right|\, dx
	\\
	\leq 
	\tilde{C}_{1}^{ 2|\alpha|-|\gamma| +|\beta| + 1}  \,N^{s(2|\alpha|-|\gamma| + |\beta| +2n+5-M)^{+}}.
	\end{multline*}
\end{linenomath}
Then, by the Theorem $18.1.11'$ in \cite{H_Book-3}, we have 
\begin{linenomath}
	\begin{multline}\label{R_ag2}
	\|  \psi_{N}^{(\beta)} \mathscr{R}_{|\alpha| -|\gamma| + 2}\left( [ \widetilde{P}_{N}, \Theta_{N}^{(\gamma)} D^{\alpha}]  \right)  g \|
	\\
	\leq
	\tilde{C}_{2}^{ 2|\alpha|-|\gamma| +|\beta| + 1}  \,N^{s(2|\alpha|-|\gamma| + |\beta| +2n+5-M)^{+}}
	\| g \|_{0} ,
	\end{multline}
\end{linenomath}
where $\tilde{C}_{2}$ is a suitable positive constant independent
of $\alpha$, $\beta$ and $\gamma$.
Let $K_{0}$ be a compact set contained in $\Omega_{2}$ and
containing all the supports of $\varphi_{N}$, $\Omega_{1} \subset K_{0}$.
Since $u $ is a $G^{s}$-vector for $P$ we have $\|P^{k}u\|_{L^{2}(K_{0})} \leq C_{K_{0}}^{2k+1} k^{2sk}$,
moreover by the Lemma \ref{L-1} the following estimate holds $k^{2sk} \leq B^{2sk} N^{s(2k-M)^{+}}$.
So, we get
\begin{linenomath}
	\begin{multline}\label{Est_g-1}
	\| g \|_{0} 
	=\| \varphi_{N}^{(\delta)} P^{k} u \|_{0}
	\leq
	C_{\varphi}^{|\delta|+1}
	N^{(|\delta|-M)^{+}}
	\|P^{k}u\|_{L^{2}(K_{0})} 
	\\
	\leq 
	C_{\varphi}^{|\delta|+1}
	\tilde{B}^{2k+1} N^{s(|\delta|+2k-M)^{+}}.
	\end{multline}
\end{linenomath}
So by \eqref{Est_I24_0}, \eqref{R_ag2} and the above estimate
and taking advantage from Remark \ref{Rk-1}, we get
\begin{linenomath}
	\begin{align}\label{Est_I24}
	I_{2,4}
	\leq
	\tilde{C}_{3}^{ 2|\alpha| -|\gamma|+|\beta|+|\delta|+2k+ 1}\,N^{s(2|\alpha|-|\gamma| +|\beta|+|\delta|+2k+2n + 5 -M)^{+}} 
	\| v \|.
	\end{align}
\end{linenomath}
Summing up, by (\ref{Est_I21}), (\ref{Est_I22}), (\ref{Est_I23}) and (\ref{Est_I24}), we obtain
\begin{linenomath}	
	\begin{multline}\label{I_2}
	I_{2}
	\leq  \varepsilon \sum_{j=1}^{n}\| P_{j} v \|_{-1}^{2} 
	+ C_{\varepsilon}
	\sum_{j=1}^{n} \| \psi_{N}^{(\beta)} (\Theta_{N}^{(\gamma) } D^{\alpha} )^{(j)} g \|_{1}^{2}
	\\
	\left.
	+C
	\left[\sum_{j=1}^{n} \| \psi_{N}^{(\beta)} (\Theta_{N}^{(\gamma) } D^{\alpha})^{(j)} g \|  
	+ \sum_{j,\ell=1}^{n} \| \psi_{N}^{(\beta)}  (\Theta_{N}^{(\gamma) } D^{\alpha})^{(j)}D_{\ell} g \| 
	\right.
	\right.
	\\
	\left.
	\left.
	+ \longsum[10]_{j,j_{1},\ell=1}^{n} \| \psi_{_{N}}^{(\beta+j_{1})}  (\Theta_{N}^{(\gamma) } D^{\alpha})^{(j)}D_{\ell} g \| 
	+ \longsum[9]_{j,j_{1},\ell=1}^{n} \| \psi_{N}^{(\beta+j_{1} +\ell)}  (\Theta_{N}^{(\gamma) } D^{\alpha})^{(j)} g \| 
	\right.
	\right.
	\\
	\left.
	\left.
	+
	\longsum[35]_{2 \leq |\mu| \leq |\alpha| -|\gamma| + 1} \longsum[7]_{\substack{\nu\leq \mu \\ \nu\leq \alpha}}
	\frac{1}{\left(\mu-\nu\right)!}  \binom{\alpha}{\nu}
	\left( \longsum[8]_{j,\ell=1}^{n} 
	\|a_{\ell,j}^{(\mu)} \psi_{N}^{(\beta)}  \Theta_{N}^{(\gamma +\mu-\nu)} D^{\alpha-\nu}D_{j}D_{\ell} g \|
	\right.
	\right.
	\right.
	\\
	+
	\sum_{ \ell=1}^{n} \|  b_{\ell}^{(\mu)} \psi_{N}^{(\beta)}  \Theta_{N}^{(\gamma +\mu-\nu)} D^{\alpha-\nu}D_{\ell} g \|
	+
	\|  c^{(\mu)} \psi_{N}^{(\beta)}  \Theta_{N}^{(\gamma +\mu-\nu)} D^{\alpha-\nu}g \|
	\Biggr)
	\\
	+\,C^{ 2|\alpha| -|\gamma|+|\beta|+|\delta|+2k}\,N^{s(2|\alpha|-|\gamma| +|\beta|+|\delta|+2k+2n + 5 -M)^{+}} 
	\Biggr]  \| v\|
	,
	\end{multline}
\end{linenomath}
where $\varepsilon$ is a small suitable constant and $C = \sup\left\{ C_{3}, C_{4}, C_{5}, C_{6}, C_{7}, \tilde{C}_{3} \right\}$.

\textbf{Estimate of the term $I_{3}$ in (\ref{Term-Pvv}) } 
Since
\begin{linenomath}
	\begin{align}\label{Comm_Pphi}
	[ P,  \varphi_{N}^{(\delta)}]=
	i \sum_{j=1}^{n}  P^{j} \varphi_{N}^{(\delta+j)} - \sum_{j,\ell=1 }^{n} a_{\ell,j}  \varphi_{N}^{(\delta+j+\ell)}
	+  \sum_{ \ell=1}^{n} b_{\ell} \varphi_{N}^{(\delta+\ell)},
	\end{align}
\end{linenomath}
and using that $P^{j}\varphi_{N}^{(\delta+j)}  = \widetilde{P}^{j}_{N}\varphi_{N}^{(\delta+j)}$,
where $\widetilde{P}^{j}_{N}$ was introduced in \eqref{t-PN}, we have
\begin{linenomath}	
	\begin{multline}
	\label{I_3-0}
	I_{3}
	\leq 
	\sum_{j=1}^{n}
	\left|\langle  \psi_{N}^{(\beta)} \Theta^{(\gamma)}_{N} D^{\alpha}  \widetilde{P}^{j}_{N} \varphi_{N}^{(\delta +j)} f, v \rangle \right|
	\\
	+
	\sum_{j,\ell=1 }^{n}
	\left|\langle  \psi_{N}^{(\beta)} \Theta^{(\gamma)}_{N} D^{\alpha} \widetilde{a}_{N,\ell,j} \varphi_{N}^{(\delta+j+\ell)} f, v \rangle \right|
	+
	\sum_{ \ell=1}^{n}\left|\langle  \psi_{N}^{(\beta)} \Theta^{(\gamma)}_{N} D^{\alpha}  \widetilde{b}_{N,\ell} \varphi_{N}^{(\delta+\ell)}f, v \rangle \right|
	\\
	=I_{3,1} + I_{3,2} + I_{3,3}.
	\end{multline}
\end{linenomath}
We estimate each of the terms obtained separately.\\
Term $I_{3,1}$. We have
\begin{linenomath}	
	\begin{multline}\label{I_3-1}
	I_{3,1}
	\leq 
	\sum_{j=1}^{n}\left(
	\left|\langle P^{j} \psi_{N}^{(\beta)} \Theta^{(\gamma)}_{N} D^{\alpha}  \varphi_{N}^{(\delta +j)} f, v \rangle \right|
	+
	\left|\langle [\psi_{N}^{(\beta)}, P^{j} ]\Theta^{(\gamma)}_{N} D^{\alpha}  \varphi_{N}^{(\delta +j)} f, v \rangle \right|
	\right.
	\\
	\quad
	\left.
	+
	\left|\langle \psi_{N}^{(\beta)} [ \Theta^{(\gamma)}_{N} D^{\alpha},  \widetilde{P}^{j}_{N}] \varphi_{N}^{(\delta +j)} f, v \rangle \right|
	\right)
	\\
	\leq
	\varepsilon	\sum_{j=1}^{n} \| P^{j} v\|^{2} + C_{\varepsilon} \sum_{j=1}^{n} \| \psi_{N}^{(\beta)} \Theta^{(\gamma)}_{N} D^{\alpha}  \varphi_{N}^{(\delta +j)} f\|^{2}
	+  C_{10} \sum_{j=1}^{n} \| \psi_{N}^{(\beta)} \Theta^{(\gamma)}_{N} D^{\alpha}  \varphi_{N}^{(\delta +j)} f\| \| v\|
	\\
	\quad
	+ C_{11} \sum_{j,\ell=1}^{n} \| \psi_{N}^{(\beta+\ell)} \Theta^{(\gamma)}_{N} D^{\alpha}  \varphi_{N}^{(\delta +j)} f\| \| v\|
	\\
	+ \sum_{j,\ell=1}^{n} \longsum[35]_{1 \leq |\mu| \leq |\alpha| -|\gamma|}
	\frac{1}{\mu!} \| \widetilde{a}_{N,\ell,j}^{(\mu)}\psi_{N}^{(\beta)} \left(\Theta^{(\gamma)}_{N} D^{\alpha}\right)^{(\mu)} D_{\ell}  \varphi_{N}^{(\delta +j)} f\| \| v\|
	\\
	\quad
	+
	\sum_{j=1}^{n}
	\| \psi_{N}^{(\beta)} \mathscr{R}_{|\alpha| -|\gamma|+ 1}\left( [ \Theta^{(\gamma)}_{N} D^{\alpha},  \widetilde{P}^{j}_{N}] \right)\varphi_{N}^{(\delta +j)} f \| \| v\|,
	\end{multline}
\end{linenomath}
where $\varepsilon$ is a small constant, $C_{\varepsilon}= \varepsilon^{-1}$,
$C_{10} =2n \displaystyle\sup_{\ell,j}\{|a_{\ell,j}^{(\ell)}|\}$ and
$C_{11} = \displaystyle\sup_{\ell,j}\{|a_{\ell,j}|\}$. We recall that $f =P^{k}u$.
The first term on the right hand side can be absorbed
by the left hand side of (\ref{Est-1/r-1}).
The first factor in the last term can be handled
as the first factor on the right hand
side of (\ref{Est_I24_0}). So
\begin{linenomath}
	\begin{multline}\label{Est_Rq_I31}
	\| \psi_{N}^{(\beta)} \mathscr{R}_{|\alpha| -|\gamma|+1}\left( [ \Theta^{(\gamma)}_{N} D^{\alpha},  \widetilde{P}^{j}_{N}] \right)\varphi_{N}^{(\delta +j)} f \|
	\\
	\leq
	\tilde{C}_{1}^{ 2|\alpha| -|\gamma|+|\beta|+|\delta|+2k+ 1}\,N^{s(2|\alpha|-|\gamma| +|\beta|+|\delta|+2k+2n + 5 -M)^{+}}.
	\end{multline}
\end{linenomath}
where $\tilde{C}_{1}$ is a suitable positive constant
independent of $\alpha$, $\beta$, $\gamma$, $\delta$ and $k$.\\
Term $I_{3,2}$: From expression of $I_{3.2}$ , in \eqref{I_3-0}, we get :
\begin{linenomath}	
	\begin{multline}\label{I_3-2}
	I_{3,2}
	\leq 
	\sum_{j,\ell=1}^{n} \longsum[28]_{0 \leq |\mu| \leq |\alpha| -|\gamma|}
	\frac{1}{\mu!} \| a_{\ell,j}^{(\mu)}\psi_{N}^{(\beta)} \left(\Theta^{(\gamma)}_{N} D^{\alpha}\right)^{(\mu)}  \varphi_{N}^{(\delta+\ell +j)} f\| \| v\|
	\\
	+
	\sum_{j,\ell=1}^{n}
	\| \psi_{N}^{(\beta)} \mathscr{R}_{|\alpha| -|\gamma| + 1}\left( [ \Theta^{(\gamma)}_{N} D^{\alpha},  \widetilde{a}_{N,\ell,j}] \right)\varphi_{N}^{(\delta +\ell+j)} f \| \| v\|.
	\end{multline}
\end{linenomath}
The first factor in the last term can be handled as the first factor on the right hand
side of (\ref{Est_I24_0}), so
\begin{linenomath}
	\begin{multline}\label{Est_Rq_I32}
	\| \psi_{N}^{(\beta)} \mathscr{R}_{|\alpha| -|\gamma| + 1}\left( [ \Theta^{(\gamma)}_{N} D^{\alpha},  \widetilde{a}_{N,\ell,j}] \right)\varphi_{N}^{(\delta +\ell+j)} f \|
	\\
	\leq
	\tilde{C}_{2}^{ 2|\alpha| -|\gamma|+|\beta|+|\delta|+2k+ 1}\,N^{s(2|\alpha|-|\gamma| +|\beta|+|\delta|+2k+2n + 5 -M)^{+}}.
	\end{multline}
\end{linenomath}
where $\tilde{C}_{2}$ is a suitable positive constant independent
of $\alpha$, $\beta$, $\gamma$, $\delta$ and $k$.\\
The term $I_{3,3}$ is clearly bounded by the third term
of the second member in \eqref{I_3-1}.\\
Summing up, by (\ref{I_3-1}), (\ref{Est_Rq_I31}),
(\ref{I_3-2}) and (\ref{Est_Rq_I32}), we conclude that
\begin{linenomath}	
	\begin{multline}\label{I_3}
	I_{3}
	\leq 
	\varepsilon	\sum_{j=1}^{n} \| P^{j} v\|^{2} 
	+ C_{\varepsilon}
	\sum_{j=1}^{n} \| \psi_{N}^{(\beta)} \Theta^{(\gamma)}_{N} D^{\alpha}  \varphi_{N}^{(\delta +j)} f\|^{2}
	\\
	+ C_{12}\left[ \sum_{j=1}^{n} \| \psi_{N}^{(\beta)} \Theta^{(\gamma)}_{N} D^{\alpha}  \varphi_{N}^{(\delta +j)} f\|
	+\sum_{j,\ell=1}^{n} \| \psi_{N}^{(\beta+\ell)} \Theta^{(\gamma)}_{N} D^{\alpha}  \varphi_{N}^{(\delta +j)} f\|
	\right.
	\\
	+\longsum[8]_{j,\ell=1}^{n} \,\longsum[38]_{1 \leq |\mu| \leq |\alpha| -|\gamma|+ 1} \longsum[7]_{\substack{\nu\leq \mu \\ \nu\leq \alpha}}
	\frac{1}{\left(\mu-\nu\right)!}  \binom{\alpha}{\nu}
	\| a_{\ell,j}^{(\mu)}\psi_{N}^{(\beta)} \Theta^{(\gamma+\mu-\nu)}_{N} D^{\alpha-\nu} D_{\ell}  \varphi_{N}^{(\delta +j)} f\| 
	\\
	+\longsum[8]_{j,\ell=1}^{n} \,\longsum[35]_{0 \leq |\mu| \leq |\alpha| -|\gamma|} \longsum[7]_{\substack{\nu\leq \mu \\ \nu\leq \alpha}}
	\frac{1}{\left(\mu-\nu\right)!}  \binom{\alpha}{\nu}
	\| a_{\ell,j}^{(\mu)}\psi_{N}^{(\beta)} \Theta^{(\gamma+\mu-\nu)}_{N} D^{\alpha-\nu} \varphi_{N}^{(\delta+\ell +j)} f\| 
	\\
	+\sum_{\ell=1}^{n} \,\longsum[35]_{0 \leq |\mu| \leq |\alpha| -|\gamma| } \longsum[7]_{\substack{\nu\leq \mu \\ \nu\leq \alpha}}
	\frac{1}{\left(\mu-\nu\right)!}  \binom{\alpha}{\nu}
	\| b_{\ell}^{(\mu)}\psi_{N}^{(\beta)} \Theta^{(\gamma+\mu-\nu)}_{N} D^{\alpha-\nu} \varphi_{N}^{(\delta+\ell)} f\| 
	\\
	+C_{12}^{ 2|\alpha| -|\gamma|+|\beta|+|\delta|+2k}\,N^{s(2|\alpha|-|\gamma| +|\beta|+|\delta|+2k+2n + 5 -M)^{+}}
	\Biggr] \|v\|,
	\end{multline}
\end{linenomath}
where $C_{12} =\sup\{ C_{10}, C_{11},\tilde{C}_{1}, \tilde{C}_{2}\}$.

\textbf{Estimate of the term $I_{4}$ in (\ref{Term-Pvv}) }. We have
\begin{linenomath}	
	\begin{align}
	\label{I_4}
	I_{4}
	\leq 
	\| \psi_{N}^{(\beta)} \Theta^{(\gamma)}_{N} D^{\alpha} \varphi_{N}^{(\delta)} P^{k+1} u\| \| v \|.
	\end{align}
\end{linenomath}

\medskip
This concludes the estimate of the term $\ell=0$ on the right hand side of (\ref{Est-1/r-1}).
Using (\ref{I_1}), (\ref{I_2}), (\ref{I_3}) and (\ref{I_4}), we obtain
%
\begin{linenomath}	
	\begin{multline}\label{Est_Pvv}
	\left|\langle P v,  v \rangle\right|
	\leq
	\varepsilon \left(\sum_{j=1}^{n} \| P^{j} v\|^{2} + \sum_{j=1}^{n}\| P_{j} v \|_{-1}^{2} \right)
	+ C_{\varepsilon} \left\{
	\sum_{j=1}^{n} \| \psi_{N}^{(\beta+j)} \Theta_{N}^{(\gamma)} D^{\alpha} g\|^{2}
	\right.
	\\
	\left. 
	+
	\sum_{j=1}^{n} \| \psi_{N}^{(\beta)} (\Theta_{N}^{(\gamma) } D^{\alpha} )^{(j)} g \|^{2}_{1}
	+
	\sum_{j=1}^{n} \| \psi_{N}^{(\beta)} \Theta^{(\gamma)}_{N} D^{\alpha}  \varphi_{N}^{(\delta +j)} f\|^{2}
	\right\}
	+
	C\Biggr\{
	\| v \|^{2} 
	\\
	+
	\left[
	\| \psi_{N}^{(\beta)} \Theta^{(\gamma)}_{N} D^{\alpha} \varphi_{N}^{(\delta)} P^{k+1} u\|
	+
	\sum_{j =1}^{n}\| \psi_{N}^{(\beta+j)} \Theta_{N}^{(\gamma)} D^{\alpha} g\|
	+ 
	\sum_{|\mu| =2} \|\psi_{N}^{(\beta+\mu)} \Theta_{N}^{(\gamma)} D^{\alpha} g\| 
	\right. 
	\\
	\left.
	+\sum_{j=1}^{n} \| \psi_{N}^{(\beta)} (\Theta_{N}^{(\gamma) } D^{\alpha})^{(j)} g \|  
	+ \sum_{j,\ell=1}^{n} \| \psi_{N}^{(\beta)}  (\Theta_{N}^{(\gamma) } D^{\alpha})^{(j)}D_{\ell} g \| 
	\right.
	\\
	\left.
	+ \longsum[10]_{j,j_{1},\ell=1}^{n} \| \psi_{_{N}}^{(\beta+j_{1})}  (\Theta_{N}^{(\gamma) } D^{\alpha})^{(j)}D_{\ell} g \| 
	+ \longsum[9]_{j,j_{1},\ell=1}^{n} \| \psi_{N}^{(\beta+j_{1} +\ell)}  (\Theta_{N}^{(\gamma) } D^{\alpha})^{(j)} g \| 
	\right.
	\\
	\left.
	+ \sum_{j=1}^{n} \| \psi_{N}^{(\beta)} \Theta^{(\gamma)}_{N} D^{\alpha}  \varphi_{N}^{(\delta +j)} f\|
	+\sum_{j,\ell=1}^{n} \| \psi_{N}^{(\beta+\ell)} \Theta^{(\gamma)}_{N} D^{\alpha}  \varphi_{N}^{(\delta +j)} f\|
	\right.
	\\
	\left.
	+
	\longsum[35]_{2 \leq |\mu| \leq |\alpha| -|\gamma| + 1} \longsum[7]_{\substack{\nu\leq \mu \\ \nu\leq \alpha}}
	\frac{1}{\left(\mu-\nu\right)!}  \binom{\alpha}{\nu}
	\left( \longsum[8]_{j,\ell=1}^{n} 
	\|a_{\ell,j}^{(\mu)} \psi_{N}^{(\beta)}  \Theta_{N}^{(\gamma +\mu-\nu)} D^{\alpha-\nu}D_{j}D_{\ell} g \|
	\right.
	\right.
	\\
	+
	\sum_{ \ell=1}^{n} \|  b_{\ell}^{(\mu)} \psi_{N}^{(\beta)}  \Theta_{N}^{(\gamma +\mu-\nu)} D^{\alpha-\nu}D_{\ell} g \|
	+
	\|  c^{(\mu)} \psi_{N}^{(\beta)}  \Theta_{N}^{(\gamma +\mu-\nu)} D^{\alpha-\nu}g \|
	\Biggr)
	\\
	+\longsum[8]_{j,\ell=1}^{n} \,\longsum[38]_{1 \leq |\mu| \leq |\alpha| -|\gamma|+ 1} \longsum[7]_{\substack{\nu\leq \mu \\ \nu\leq \alpha}}
	\frac{1}{\left(\mu-\nu\right)!}  \binom{\alpha}{\nu}
	\| a_{\ell,j}^{(\mu)}\psi_{N}^{(\beta)} \Theta^{(\gamma+\mu-\nu)}_{N} D^{\alpha-\nu} D_{\ell}  \varphi_{N}^{(\delta +j)} f\|
	\\
	\left.
	+
	\longsum[35]_{0 \leq |\mu| \leq |\alpha| -|\gamma|} \longsum[7]_{\substack{\nu\leq \mu \\ \nu\leq \alpha}}
	\frac{1}{\left(\mu-\nu\right)!}  \binom{\alpha}{\nu} 
	\left( \longsum[8]_{j,\ell=1}^{n}
	\| a_{\ell,j}^{(\mu)}\psi_{N}^{(\beta)} \Theta^{(\gamma+\mu-\nu)}_{N} D^{\alpha-\nu} \varphi_{N}^{(\delta+\ell +j)} f\| 
	\right.
	\right.
	\\
	\left.
	+\sum_{\ell=1}^{n} 
	\| b_{\ell}^{(\mu)}\psi_{N}^{(\beta)} \Theta^{(\gamma+\mu-\nu)}_{N} D^{\alpha-\nu} \varphi_{N}^{(\delta+\ell)} f\| 
	\Biggr)
	\right.
	\\
	+\,C^{ 2|\alpha| -|\gamma|+|\beta|+|\delta|+2k}\,N^{s(2|\alpha|-|\gamma| +|\beta|+|\delta|+2k+2n + 5 -M)^{+}} 
	\Biggr]
	\| v\|\Biggr\},
	\end{multline}
\end{linenomath}	
where $\varepsilon$, $C_{\varepsilon}$ and $C$ are suitable constants independent of $\alpha$, $\beta$, $\gamma$, $\delta$, $k$ and $N$.
$\varepsilon$ is a small parameter that we will chose at the end in order to absorbed the first two terms by the left hand side of (\ref{Est-1/r-1});
we point out that the number of times that we use it is finite.~\par
\vskip-2mm

\textbf{Estimate of the terms $\left|\langle E_{\ell} P v, E_{\ell} v \rangle\right|$,  $\ell \geq 1$, in  (\ref{Est-1/r-1}).}~\par
\vskip-6mm
We recall that $E_{\ell}=D_{\ell} \psi \Lambda_{-1}$, where $\psi$ belongs to
$\mathscr{D}\left(\Omega\right)$ and is identically one on
$\Omega_{4}$, $\Omega_{4} \Subset \Omega$, and $\Lambda_{-1}$
is the pseudodifferential operator associated to the symbol
$\lambda(\xi)^{-1}\doteq \left(1+|\xi|^{2}\right)^{-1/2}$.
We point out that $E_{\ell}$ are zero order pseudodifferential operators,
$\|E_{\ell}\|_{L^{2}\rightarrow L^{2}} \leq C$.\\ 
We have
\begin{linenomath}
	\begin{multline}\label{Term-EPvEv}
	\left|\langle E_{\ell}P v, E_{\ell}v\rangle \right|
	\leq 
	\left|\langle E_{\ell} [ P, \psi_{N}^{(\beta)} ]  \Theta^{(\gamma)}_{N} D^{\alpha} g,  E_{\ell}v\rangle \right|
	\\
	+
	\left|\langle  E_{\ell}\psi_{N}^{(\beta)} [ P , \Theta^{(\gamma)}_{N} D^{\alpha} ] g, E_{\ell} v\rangle \right|
	+
	\left|\langle E_{\ell} \psi_{N}^{(\beta)} \Theta^{(\gamma)}_{N} D^{\alpha} [ P, \varphi_{N}^{(\delta)}] f,  E_{\ell} v \rangle \right|
	\\
	+
	\left|\langle E_{\ell} \psi_{N}^{(\beta)} \Theta^{(\gamma)}_{N} D^{\alpha} \varphi_{N}^{(\delta)} P^{k+1} u,  E_{\ell} v \rangle \right|
	=\sum_{j=1}^{4} \tilde{I}_{j}.
	\end{multline}
\end{linenomath}
\textbf{Estimate of the term $\tilde{I}_{1}$}.  By (\ref{Cm_Ppsi}), we have
\begin{linenomath}
	\begin{multline*}
	\tilde{I}_{1} 
	\leq \sum_{j=1}^{n} \left|\left(  E_{\ell} P^{j} \psi^{(\beta+j)}_{N} \Theta_{N}^{(\gamma) } D^{\alpha} g,  E_{\ell} v \right)\right|
	+ \sum_{|\mu| =2} \left| \left(  E_{\ell} P^{\mu} \psi_{N}^{(\beta+\mu)} \Theta_{N}^{(\gamma)} D^{\alpha} g,  E_{\ell}v\right)\right|
	\\
	+  \sum_{j =1}^{n} \left|\left(  E_{\ell} b_{j} \psi_{N}^{(\beta+j)}  \Theta_{N}^{(\gamma) } D^{\alpha} g,  E_{\ell} v \right)\right|
	\end{multline*}
\end{linenomath}
where $ P^{\mu}= \partial_{\xi}^{\mu}p^{0}(x,\xi) =\partial_{\xi_{j}}\partial_{\xi_{\ell}}p^{0}(x,\xi)= a_{j,\ell}(x)$.
In order to handle the first term we use $\widetilde{P}^{j}_{N}$, introduced in \eqref{t-PN}.
Since $ \left( \widetilde{P}^{j}_{N} \right)^{\hspace{-0.4em}*} =  \widetilde{P}^{j}_{N} - 2 i \sum_{q =1}^{n} \widetilde{a}_{N,j,q}^{(q)}(x)$,
we get
\begin{linenomath}
	\begin{multline*}
	\tilde{I}_{1} 
	\leq \sum_{j=1}^{n} \left|\left(  E_{\ell}\widetilde{P}^{j}_{N} \psi^{(\beta+j)}_{N} \Theta_{N}^{(\gamma) } D^{\alpha} g,  E_{\ell} v \right)\right|
	+  \sum_{ j=1}^{n} \| E_{\ell} b_{j} \psi_{N}^{(\beta+j)}  \Theta_{N}^{(\gamma)} D^{\alpha} g\| \| E_{\ell} v\| 
	\\
	+ \sum_{|\mu| =2} \|  E_{\ell} P^{\mu} \psi_{N}^{(\beta+\mu)} \Theta_{N}^{(\gamma)} D^{\alpha} g\| \|  E_{\ell}v\| 
	\\
	\leq
	\sum_{j=1}^{n} \left|\left(  E_{\ell} \psi^{(\beta+j)}_{N} \Theta_{N}^{(\gamma) } D^{\alpha} g,  E_{\ell} \widetilde{P}^{j}_{N} v \right)\right|
	+
	\sum_{j=1}^{n} \| [ E_{\ell},\widetilde{P}^{j}_{N}] \psi^{(\beta+j)}_{N} \Theta_{N}^{(\gamma) } D^{\alpha} g\|  \| E_{\ell} v \|
	\\
	+
	2\longsum[7]_{j,q=1}^{n}  \|  E_{\ell} \psi^{(\beta+j)}_{N} \Theta_{N}^{(\gamma) } D^{\alpha} g\| \|\widetilde{a}_{N,j,q}^{(q)} E_{\ell}  v \|
	+
	\sum_{j=1}^{n} \| E_{\ell} \psi^{(\beta+j)}_{N} \Theta_{N}^{(\gamma) } D^{\alpha} g\| \|[\widetilde{P}^{j}, E_{\ell}]  v \|
	\\
	+  \sum_{ j=1}^{n} \| E_{\ell} b_{j} \psi_{N}^{(\beta+j)}  \Theta_{N}^{(\gamma)} D^{\alpha} g\| \| E_{\ell} v\| 
	+ \sum_{|\mu| =2} \|  E_{\ell} P^{\mu} \psi_{N}^{(\beta+\mu)} \Theta_{N}^{(\gamma)} D^{\alpha} g\| \|  E_{\ell}v\|. 
	\end{multline*}
\end{linenomath}	
Now, $[E_{\ell},\widetilde{P}^{j}_{N}] $ are zero order pseudodifferential operators
and by the Theorem $18.1.11'$ in \cite{H_Book-3}
we have $\|[E_{\ell},\widetilde{P}^{j}_{N}] \|_{L^{2}\rightarrow L^{2}} \leq C$,
with $C$ independent of $N$.
We conclude that
\begin{linenomath}
\begin{multline}\label{Est_It_1}
	\tilde{I}_{1}
	\leq \varepsilon \sum_{j=1}^{n} \| P^{j} v\|^{2} + C_{\varepsilon} \sum_{j=1}^{n} \| \psi_{N}^{(\beta+j)} \Theta_{N}^{(\gamma)} D^{\alpha} g\|^{2}
	\hspace{12em}
	\\
	\hspace{7em}
	+
	C_{2} \left( \sum_{j =1}^{n}\| \psi_{N}^{(\beta+j)} \Theta_{N}^{(\gamma)} D^{\alpha} g\| \|v\|
	+
	\sum_{|\mu| =2} \|\psi_{N}^{(\beta+\mu)} \Theta_{N}^{(\gamma)} D^{\alpha} g\| \| v\| \right),
	\end{multline}
\end{linenomath}
where $\varepsilon$ is a small suitable positive constant;
the first term on the right hand side,
$\varepsilon \sum_{j} \| P^{j} v\|^{2}$, can be absorbed
by the left hand side of (\ref{Est-1/r-1}).

\textbf{Estimate of the term $\tilde{I}_{2}$ in (\ref{Term-EPvEv})}.
By (\ref{Com-PThetaD}), we have
\begin{linenomath}	
	\begin{multline}\label{I_2-t}
	\tilde{I}_{2} 
	\leq \sum_{j=1}^{n} \left|\left(  E_{\ell} \psi_{N}^{(\beta)} \widetilde{P}_{N,j} \left(\Theta_{N}^{(\gamma)} D^{\alpha}\right)^{(j)} g, E_{\ell} v \right)\right|
	\\
	+ 
	\sum_{j=1}^{n}
	\left|\left(  E_{\ell} \psi_{N}^{(\beta)} \left(B_{j} + c^{(j)}(x) \right) \left(\Theta_{N}^{(\gamma)} D^{\alpha}\right)^{(j)} g , E_{\ell} v \right)\right|
	\\
	\quad\quad
	+ \longsum[32]_{2 \leq |\mu| \leq |\alpha| -|\gamma| + 1} \longsum[7]_{\substack{\nu\leq \mu \\ \nu\leq \alpha}}
	\frac{1}{\mu!}  \binom{\mu}{\nu} \frac{\alpha!}{(\alpha -\nu)!}
	\left|\left(  E_{\ell} \psi_{N}^{(\beta)} P_{\mu} \Theta_{N}^{(\gamma +\mu-\nu)} D^{\alpha-\nu} g, E_{\ell} v \right)\right|
	\\
	\quad\qquad\qquad\qquad\qquad\qquad\qquad\qquad
	+\left|\left(  E_{\ell}\psi_{N}^{(\beta)} \mathscr{R}_{|\alpha| -|\gamma| + 2}\left( [ \widetilde{P}_{N}, \Theta_{N}^{(\gamma)} D^{\alpha}]  \right)  g,E_{\ell}  v \right)\right|
	\\
	\quad
	= \tilde{I}_{2,1} + \tilde{I}_{2,2} + \tilde{I}_{2,3} + \tilde{I}_{2,4},
	\end{multline}
\end{linenomath}
where
$ P_{\mu} = \sum_{q,j_{1}=1}^{n} a_{q,j_{1}}^{(\mu)} D_{q} D_{j_{1}} + i \sum_{q=0}^{n} b_{q}^{(\mu)} D_{\ell} +c^{(\mu)}$,
$B_{j} =  i\sum_{q=0}^{n} b_{q}^{(j)} D_{q}$ and $\widetilde{P}_{N,j} = \sum_{q,j_{1} =1}^{n}\widetilde{a}_{N,q,j}^{(j)}D_{q}D_{j_{1}}$,
introduced in \eqref{t-PN}.\\
The terms  $\tilde{I}_{2,2} $, $\tilde{I}_{2,3}$ and $\tilde{I}_{2,4}$ can be handled as the terms
$I_{2,2} $, $I_{2,3} $ and $I_{2,4}$ in (\ref{I_2-1}), see (\ref{Est_I22}), (\ref{Est_I23}) and (\ref{Est_I24}).\\
Concerning the term $\tilde{I}_{2,1} $, since
\begin{linenomath}
	\begin{align*}
	\left(\widetilde{P}_{N,j}\right)^{*}
	= \widetilde{P}_{N,j}
	+ \frac{1}{i}\sum_{q,j_{1}=1}^{n} \left(\widetilde{a}_{N,q,j_{1}}^{(j+q)} D_{j_{1}} +  \widetilde{a}_{N,q,j_{1}}^{(j+j_{1})} D_{q}\right)
	-\sum_{q,j_{1}=1}^{n} \widetilde{a}_{N,q,j_{1}}^{(j+j_{1}+ q)}, 
	\end{align*}
\end{linenomath}
and
\begin{linenomath}
	\begin{align*}
	\left[\psi_{N}^{(\beta)}, \widetilde{P}_{N,j}\right] = 
	i\longsum[8]_{q,j_{1}=1}^{n} \widetilde{a}_{N,q,j_{1}}^{(j)} \left(\psi_{N}^{(\beta+q)}D_{j_{1}} + \psi_{N}^{(\beta+j_{1})}D_{q}- i\psi_{N}^{(\beta+j_{1}+q)}\right),
	\end{align*}
\end{linenomath}
we have
\begin{linenomath}
	\begin{multline}\label{Est_I21t}
	\tilde{I}_{2,1} 
	\leq 
	\sum_{j=1}^{n} 
	\left|\left( E_{\ell} \psi_{N}^{(\beta)}\left(\Theta_{N}^{(\gamma) } D^{\alpha}\right)^{(j)} g, E_{\ell} P_{j} v\right)\right|
	\\
	+\sum_{j=1}^{n} \left[ 
	\left|\left( [E_{\ell}, \widetilde{P}_{N,j}^{*}] \psi_{N}^{(\beta)}\left(\Theta_{N}^{(\gamma) } D^{\alpha}\right)^{(j)} g, E_{\ell}  v\right)\right|
	+
	\left|\left( E_{\ell}\psi_{N}^{(\beta)}\left(\Theta_{N}^{(\gamma) } D^{\alpha}\right)^{(j)} g, [\widetilde{P}_{N,j},E_{\ell}]  v\right)\right|
	\right]
	\\
	+
	\sum_{j=1}^{n}
	\left|\left(  E_{\ell} \left( \widetilde{P}_{N,j} - \widetilde{P}_{N,j}^{*} \right)
	\psi_{N}^{(\beta)}\left(\Theta_{N}^{(\gamma) } D^{\alpha}\right)^{(j)} g, E_{\ell} v \right)\right|
	\\
	+  \sum_{j=1}^{n}
	\left|\left( E_{\ell}
	\left[\psi_{N}^{(\beta)}, \widetilde{P}_{N,j}\right] \left(\Theta_{N}^{(\gamma) } D^{\alpha}\right)^{(j)} g, E_{\ell} v \right)\right|
	\\
	\leq  \varepsilon  \sum_{j=1}^{n}\| P_{j} v \|_{-1}^{2} 
	+ C_{\varepsilon} \left(\sum_{j=1}^{n} \| \psi_{N}^{(\beta)} (\Theta_{N}^{(\gamma) } D^{\alpha} )^{(j)} g \|_{1}^{2} 
	+ \| \psi_{N}^{(\beta)} (\Theta_{N}^{(\gamma) } D^{\alpha} )^{(j)} g \|^{2}\right)
	\\
	+ C_{3}\left( 
	\sum_{j=1}^{n} \left(\| \psi_{N}^{(\beta)} (\Theta_{N}^{(\gamma) } D^{\alpha} )^{(j)} g \|_{1} 
	+\| \psi_{N}^{(\beta)} (\Theta_{N}^{(\gamma) } D^{\alpha} )^{(j)} g \|\right) \|v\|
	\right)
	\\
	+ C_{4}\left( 
	\sum_{j=1}^{n} \left(\| \psi_{N}^{(\beta+\ell)} (\Theta_{N}^{(\gamma) } D^{\alpha} )^{(j)} g \| 
	+\| \psi_{N}^{(\beta)} (\Theta_{N}^{(\gamma) } D^{\alpha} )^{(j)}D_{\ell} g \|\right) \|v\|
	\right)
	\\
	+ C_{5}\left( 
	\longsum[9]_{j,j_{1},q=1}^{n} \left(\| \psi_{N}^{(\beta+j_{1} )}  (\Theta_{N}^{(\gamma) } D^{\alpha})^{(j)} D_{q} g \|
	+
	\| \psi_{N}^{(\beta+j_{1} +q)}  (\Theta_{N}^{(\gamma) } D^{\alpha})^{(j)} g \| 
	\right) \|v\|
	\right),
	\end{multline}
\end{linenomath}
where $\varepsilon$ is a small suitable constant
so that the first term can be absorbed by the left
hand side of (\ref{Est-1/r-1}).

\textbf{Estimate of the term $\tilde{I}_{3}$ in (\ref{Term-EPvEv})}.
In order to handle this term we replace $P_{j}$ with $\widetilde{P}^{j}_{N}$,
introduced in \eqref{t-PN}, when it will be useful.
We recall that since $ \tilde{a}_{N,j, j_{1}}(x)=\widetilde{\psi}_{N}(x) a_{j,j_{1}}(x) $
we have  that $\psi_{N}\widetilde{P}^{j}_{N} =\psi_{N} P^{j}$, $\widetilde{P}^{j}_{N}\varphi_{N}= P^{j}\varphi_{N}$,
$[\psi_{N},\widetilde{P}^{j}_{N}] =[\psi_{N}, P^{j}]$ and  $[\widetilde{P}^{j}_{N},\varphi_{N}]= [P^{j},\varphi_{N}]$.
Using (\ref{Comm_Pphi}), we have
\begin{linenomath}	
	\begin{multline}\label{I_3t}
	\tilde{I}_{3}
	\leq 
	\sum_{j=1}^{n}
	\left|\langle  E_{\ell} \psi_{N}^{(\beta)} \Theta^{(\gamma)}_{N} D^{\alpha}  
	\widetilde{P}^{j}_{N} \varphi_{N}^{(\delta +j)} f,E_{\ell} v \rangle \right|
	\\
	+
	\sum_{j,q=1 }^{n}
	\left|\langle E_{\ell} \psi_{N}^{(\beta)} \Theta^{(\gamma)}_{N} D^{\alpha} \widetilde{a}_{N,q,j} \varphi_{N}^{(\delta+j+q)} f, E_{\ell}v \rangle \right|
	+
	\sum_{q=1}^{n}\left|\langle  E_{\ell} \psi_{N}^{(\beta)} \Theta^{(\gamma)}_{N} D^{\alpha}  
	\widetilde{b}_{N,q} \varphi_{N}^{(\delta+q)}f, E_{\ell}v \rangle \right|
	\\
	\leq 
	\sum_{j=1}^{n}\left(
	\left|\langle E_{\ell} \widetilde{P}^{j}_{N} \psi_{N}^{(\beta)} \Theta^{(\gamma)}_{N} D^{\alpha} 
	\varphi_{N}^{(\delta +j)} f, E_{\ell} v \rangle \right|
	+
	\left|\langle E_{\ell} [\psi_{N}^{(\beta)}, P^{j} ]\Theta^{(\gamma)}_{N} D^{\alpha}  \varphi_{N}^{(\delta +j)} f, E_{\ell}v \rangle \right|
	\right.
	\\
	\quad
	\left.
	+
	\left|\langle E_{\ell}\psi_{N}^{(\beta)}
	[ \Theta^{(\gamma)}_{N} D^{\alpha},  \widetilde{P}^{j}_{N}] \varphi_{N}^{(\delta +j)} f, E_{\ell}v \rangle \right|
	\right)
	\\
	+
	\sum_{j,q=1 }^{n}
	\left|\langle E_{\ell} \psi_{N}^{(\beta)} \Theta^{(\gamma)}_{N} D^{\alpha} \widetilde{a}_{N,q,j} \varphi_{N}^{(\delta+j+q)} f, E_{\ell}v \rangle \right|
	+
	\sum_{ q=1}^{n}\left|\langle  E_{\ell} \psi_{N}^{(\beta)} \Theta^{(\gamma)}_{N} D^{\alpha}  \widetilde{b}_{N,q} \varphi_{N}^{(\delta+q)}f, E_{\ell}v \rangle \right|
	\\
	=\tilde{I}_{3,1} + \tilde{I}_{3,2} + \tilde{I}_{3,3}+  \tilde{I}_{3,4} +  \tilde{I}_{3,5}.
	\end{multline}
\end{linenomath}
The term $ \tilde{I}_{3,1}$ can be handled
as the first term on the right hand side in the first line
of (\ref{Est_It_1}); $ \tilde{I}_{3,2}$ and $ \tilde{I}_{3,3}$
can be handled as the second and the third
terms on the right hand side of (\ref{I_3-1});
the terms  $ \tilde{I}_{3,4}$ and $ \tilde{I}_{3,5}$ can be estimated
as the terms $ I_{3,2}$ and $ I_{3,3}$, see (\ref{I_3-2}).

\textbf{Estimate of the term $\tilde{I}_{4}$ in (\ref{Term-EPvEv})}. We have
\begin{linenomath}	
	\begin{align}
	\label{I_4t}
	\tilde{I}_{4}
	\leq 
	C \left( 
	\| \psi_{N}^{(\beta)} \Theta^{(\gamma)}_{N} D^{\alpha} \varphi_{N}^{(\delta)} P^{k+1} u\| \| v \|
	\right).
	\end{align}
\end{linenomath}
This concludes the estimate of the terms
$\left|\langle E_{\ell}P v, E_{\ell}v\rangle \right|$, $\ell =1, \dots, n$.

Summing up, by (\ref{Est-1/r-1}), (\ref{Est_Pvv}), (\ref{Term-EPvEv}), (\ref{Est_It_1}), (\ref{Est_I21t}), (\ref{I_3t})
and (\ref{I_4t}), and since 
$(\Theta_{N}^{(\gamma) } D^{\alpha})^{(j)} = \Theta_{N}^{(\gamma+j) } D^{\alpha} +\alpha_{j} \Theta_{N}^{(\gamma) } D^{\alpha-j} $,
we obtain that there are two positive constants, $A$ and $B$ such that
\begin{linenomath}	
	\begin{multline}\label{Est-v-1/r}
	\| v\|_{\frac{1}{r}} ^{2}
	\leq
	B\left\{
	\sum_{j=1}^{n} \left(\| \psi_{N}^{(\beta+j)} \Theta_{N}^{(\gamma)} D^{\alpha} g\|^{2}
	+
	\| \psi_{N}^{(\beta)} \Theta^{(\gamma)}_{N} D^{\alpha}  \varphi_{N}^{(\delta +j)} f\|^{2}\right)
	+ 
	\| v \|^{2} 
	\right.
	\\
	\left.
	+
	\sum_{j=1}^{n}
	\left( \| \psi_{N}^{(\beta)} \Theta_{N}^{(\gamma+j) } D^{\alpha} g\|_{1}^{2}
	+\alpha_{j}^{2} \| \psi_{N}^{(\beta)}\Theta_{N}^{(\gamma) } D^{\alpha-j} g \|^{2}_{1}
	\right)
	\right.
	\\
	+
	\left[
	\sum_{j =1}^{n}
	\left( \| \psi_{N}^{(\beta+j)} \Theta_{N}^{(\gamma)} D^{\alpha} g\|
	+\| \psi_{N}^{(\beta)} \Theta_{N}^{(\gamma+j) } D^{\alpha} g\| 
	+\alpha_{j} \| \psi_{N}^{(\beta)}\Theta_{N}^{(\gamma) } D^{\alpha-j} g \| 
	\right)
	\right.
	\\
	+ 
	\sum_{|\mu| =2} \|\psi_{N}^{(\beta+\mu)} \Theta_{N}^{(\gamma)} D^{\alpha} g\| 
	+ \sum_{j,\ell=1}^{n} 
	\left(\| \psi_{N}^{(\beta)}  \Theta_{N}^{(\gamma+j) } D^{\alpha+\ell}  g \|
	+\alpha_{j} \| \psi_{N}^{(\beta)}\Theta_{N}^{(\gamma) } D^{\alpha-j+\ell} g \| 
	\right)
	\\
	+ \sum_{j=1}^{n} \| \psi_{N}^{(\beta)} \Theta^{(\gamma)}_{N} D^{\alpha}  \varphi_{N}^{(\delta +j)} f\|
	+ \longsum[8]_{j,j_{1},\ell=1}^{n}\left(
	\| \psi_{_{N}}^{(\beta+j_{1})} \Theta_{N}^{(\gamma+j) } D^{\alpha+\ell} g\|
	+ 
	\| \psi_{N}^{(\beta+j_{1} +\ell)} \Theta_{N}^{(\gamma+j) } D^{\alpha} g\|
	\right.
	\\
	\left.
	+ \alpha_{j} 
	\left(\| \psi_{N}^{(\beta+j_{1} +\ell)} \Theta_{N}^{(\gamma) } D^{\alpha-j} g \| +\| \psi_{_{N}}^{(\beta+j_{1})} \Theta_{N}^{(\gamma) } D^{\alpha-j + \ell} g \|
	\right)\right)
	\\
	\left.
	+\sum_{j,\ell=1}^{n} \left(
	\| \psi_{N}^{(\beta+\ell)} \Theta^{(\gamma)}_{N} D^{\alpha}  \varphi_{N}^{(\delta +j)} f\|
	+
	\| \psi_{N}^{(\beta+\ell)}  \Theta_{N}^{(\gamma+j) } D^{\alpha}  g \|
	+\alpha_{j} \| \psi_{N}^{(\beta+\ell)}\Theta_{N}^{(\gamma) } D^{\alpha-j} g \| 
	\right)
	\right.
	\\
	\left.
	+
	\longsum[35]_{2 \leq |\mu| \leq |\alpha| -|\gamma| + 1} \longsum[7]_{\substack{\nu\leq \mu \\ \nu\leq \alpha}}
	\frac{1}{\left(\mu-\nu\right)!}  \binom{\alpha}{\nu}
	\left( \longsum[8]_{j,\ell=1}^{n} 
	\|a_{\ell,j}^{(\mu)} \psi_{N}^{(\beta)}  \Theta_{N}^{(\gamma +\mu-\nu)} D^{\alpha-\nu}D_{j}D_{\ell} g \|
	\right.
	\right.
	\\
	+
	\sum_{ \ell=1}^{n} \|  b_{\ell}^{(\mu)} \psi_{N}^{(\beta)}  \Theta_{N}^{(\gamma +\mu-\nu)} D^{\alpha-\nu}D_{\ell} g \|
	+
	\|  c^{(\mu)} \psi_{N}^{(\beta)}  \Theta_{N}^{(\gamma +\mu-\nu)} D^{\alpha-\nu}g \|
	\Biggr)
	\\
	+\longsum[8]_{j,\ell=1}^{n} \,\longsum[38]_{1 \leq |\mu| \leq |\alpha| -|\gamma|+ 1} \longsum[7]_{\substack{\nu\leq \mu \\ \nu\leq \alpha}}
	\frac{1}{\left(\mu-\nu\right)!}  \binom{\alpha}{\nu}
	\| a_{\ell,j}^{(\mu)}\psi_{N}^{(\beta)} \Theta^{(\gamma+\mu-\nu)}_{N} D^{\alpha-\nu} D_{\ell}  \varphi_{N}^{(\delta +j)} f\|
	\\
	\left.
	+
	\longsum[35]_{0 \leq |\mu| \leq |\alpha| -|\gamma|} \longsum[7]_{\substack{\nu\leq \mu \\ \nu\leq \alpha}}
	\frac{1}{\left(\mu-\nu\right)!}  \binom{\alpha}{\nu} 
	\left( \longsum[8]_{j,\ell=1}^{n}
	\| a_{\ell,j}^{(\mu)}\psi_{N}^{(\beta)} \Theta^{(\gamma+\mu-\nu)}_{N} D^{\alpha-\nu} \varphi_{N}^{(\delta+\ell +j)} f\| 
	\right.
	\right.
	\\
	\left.
	+\sum_{\ell=1}^{n} 
	\| b_{\ell}^{(\mu)}\psi_{N}^{(\beta)} \Theta^{(\gamma+\mu-\nu)}_{N} D^{\alpha-\nu} \varphi_{N}^{(\delta+\ell)} f\| 
	\Biggr)
	\right.
	+
	A^{\sigma+1} B^{2 m + |\gamma| } 
	N^{s(2m +|\gamma| +\sigma+2n + 5 -M)^{+}} 
	\\
	+
	\| \psi_{N}^{(\beta)} \Theta^{(\gamma)}_{N} D^{\alpha} \varphi_{N}^{(\delta)} P^{k+1} u\|
	\Biggr]
	\| v\|\Biggr\}
	+ 
	A^{2(\sigma+1)} B^{2(2 m + |\gamma| + 1)}  
	N^{2s(2 m+ |\gamma|+\sigma) },
	\end{multline}
\end{linenomath} 
where $m= |\alpha| - |\gamma| $ and  $\sigma= |\beta| + |\delta| +2k$.\\
We remark that the strategy adopted in \eqref{I_1}, \eqref{Est_I21},
\eqref{I_3-1}, \eqref{Est_It_1} and \eqref{Est_I21t}, where we introduce
$\varepsilon$ in order to absorb a term on the left hand side of \eqref{Est-1/r-1},
is used a finite number of times, say at most $50$ times; this allows us to choose $\varepsilon$
so that $50 \varepsilon  \, C < 1/2$, where $C$ is the constant on the right hand side of \eqref{Est-1/r-1}. 
%
\subsection{Estimates in $H^{p/r}$, $p=2,\,\dots, r$}~\par
\vskip-4mm
The purpose of the present section is to obtain a suitable estimate
of $\| v \|_{p/r}^{2}$, $p= 2,\dots, r-1$, where
$ v =\psi_{N}^{(\beta)} \Theta_{N}^{(\gamma)} D^{\alpha} \varphi_{N}^{(\delta)} P^{k} u$,
here $\psi_{N}$, $\Theta_{N}$ and $\varphi_{N}$
are as in the previous section.
We denote by $\Lambda_{r}^{\ell}$ the pseudodifferential operator
with symbol $ \left( 1+ |\xi|^{2} \right)^{\ell/2r}$.
Let $\Omega_{3}$ be open neighborhood of $x_{0}$
such that $\Omega_{2} \Subset \Omega_{3} \Subset \widetilde{\Omega}_{3} \Subset \Omega_{4}$
and $\Gamma_{2}$ open cone around $\xi_{0}$
such that $\Gamma_{1} \Subset \Gamma_{2}$.
We introduce $\widetilde{\psi}(x) \in \mathscr{D}(\Omega_{3})$,
such that $\widetilde{\psi}\equiv 1$ on $\Omega_{2}$
and $\widetilde{\Theta}_{q}(D)$ a sequence of zero order pseudodifferential operators
with symbol $\widetilde{\Theta}_{q}(\xi)$ of Ehrenpreis-Andersson
type associated to the couple of open cones
$(\Gamma_{1}, \Gamma_{2})$, 
i.e. $\widetilde{\Theta}_{q}(\xi)\equiv 1$ in
$ \Gamma_{1} \cap \lbrace \xi\in \mathbb{R}^{n}\,: \, |\xi| \geq q\rbrace  $,
supported in $\Gamma_{2}\cap \lbrace \xi\in \mathbb{R}^{n}\,: \, |\xi| \geq q/2\rbrace $, 
and such that they satisfy the estimate (\ref{ThetaN-AH}), Lemma \ref{EA-Cutoff} in the Appendix,
for all $\alpha \in \mathbb{Z}^{n}_{+}$ with $|\alpha|\leq q$. 
We recall that the sequence $ \Theta_{N}(D) $ has symbols
$\Theta_{N}(\xi)$ of Ehrenpreis-Andersson type associated to the couple of open cones
$(\Gamma_{0}, \Gamma_{1})$.
We will use the same notation of the previous section:
$f \doteq P^{k}u$, $ g \doteq \varphi_{N}^{(\delta)} f =\varphi_{N}^{(\delta)} P^{k} u$,
$ w = \Theta_{N}^{(\gamma)} D^{\alpha}g =\Theta_{N}^{(\gamma)} D^{\alpha} \varphi_{N}^{(\delta)} P^{k} u$
and $v \doteq \psi_{N}^{(\beta)} w= \psi_{N}^{(\beta)} \Theta_{N}^{(\gamma)} D^{\alpha}g
=\psi_{N}^{(\beta)} \Theta_{N}^{(\gamma)} D^{\alpha} \varphi_{N}^{(\delta)}f
= \psi_{N}^{(\beta)} \Theta_{N}^{(\gamma)} D^{\alpha} \varphi_{N}^{(\delta)} P^{k} u$.\\
We have
\begin{linenomath}	
	\begin{multline}\label{Step-0-p/r}
	\| \psi_{N}^{(\beta)} w\|_{\frac{p}{r}}
	\leq 
	\| \psi_{N}^{(\beta)} (1-\widetilde{\Theta}_{q})\Theta_{N}^{(\gamma)} D^{\alpha} g \|_{\frac{p}{r}}
	+
	\|\widetilde{\psi}  \psi_{N}^{(\beta)} \widetilde{\Theta}_{q}\Theta_{N}^{(\gamma)} D^{\alpha} g \|_{\frac{p}{r}}
	\\
	\leq 
	\| \psi_{N}^{(\beta)} (1-\widetilde{\Theta}_{q})\Theta_{N}^{(\gamma)} D^{\alpha} g \|_{\frac{p}{r}}
	+
	\|\Lambda_{r}^{p-1}\widetilde{\psi}  \psi_{N}^{(\beta)}  \widetilde{\Theta}_{q}   \Theta_{N}^{(\gamma)} D^{\alpha} g\|_{\frac{1}{r}}
	\\
	\leq
	\| \psi_{N}^{(\beta)} (1-\widetilde{\Theta}_{q})\Theta_{N}^{(\gamma)} D^{\alpha} g \|_{\frac{p}{r}}
	+	
	\|[ \Lambda_{r}^{p-1}, \widetilde{\psi}] v\|_{\frac{1}{r}}
	+
	\|\widetilde{\psi}  \Lambda_{r}^{p-1}  \psi_{N}^{(\beta)}  \widetilde{\Theta}_{q}  w \|_{\frac{1}{r}}.
	\\
	\leq
	C_{0} \|v\|_{\frac{p}{r}-1} 
	+ \| \psi_{N}^{(\beta)} (1-\widetilde{\Theta}_{q})\Theta_{N}^{(\gamma)} D^{\alpha} g \|_{\frac{p}{r}}
	+
	\|\widetilde{\psi}  \Lambda_{r}^{p-1}  \psi_{N}^{(\beta)}  \widetilde{\Theta}_{q}  w \|_{\frac{1}{r}}.
	\end{multline}
\end{linenomath}	
Since $ (1-\widetilde{\Theta}_{q})\Theta_{N}^{(\gamma)}$
is supported in the region $\frac{N}{2} \leq |\xi | \leq q$,
here we are assuming that $q < N$; we have
\begin{linenomath}
	\begin{multline*}
	\left(1 +|\xi|^{2} \right)^{\frac{p}{2r}} | (1-\widetilde{\Theta}_{q})(\xi)| | \Theta_{N}^{(\gamma)}(\xi)| |\xi|^{\alpha}
	\leq C_{\Theta}^{|\gamma|+1} N^{(|\gamma|-M)^{+}} (1+|\xi|)^{|\alpha| -|\gamma| +\frac{p}{r}} 
	\\
	\leq C_{\Theta}^{|\gamma|+1} N^{(|\gamma|-M)^{+}} (1+q)^{|\alpha| -|\gamma| +\frac{p}{r}} 
	\leq C_{\Theta}^{|\gamma|+1} B^{|\alpha| -|\gamma| + 1} N^{(|\gamma| +q+1-M)^{+}},
	\end{multline*}
\end{linenomath}
where the last inequality was obtained
taking advantage from Lemma \ref{L-1} and Remark \ref{Rk-1}.
We have
\begin{linenomath}	
	\begin{multline*}
	\| \psi_{N}^{(\beta)} (1-\widetilde{\Theta}_{q})\Theta_{N}^{(\gamma)} D^{\alpha} g \|_{\frac{p}{r}}
	\leq \int \left(1 +|\xi|^{2} \right)^{\frac{p}{r}} |\widehat{\psi}_{N}^{(\beta)} (\xi)|d\xi\,\, 
	\|  (1-\widetilde{\Theta}_{q})\Theta_{N}^{(\gamma)} D^{\alpha} g \|_{\frac{p}{r}}
	\\
	\leq C_{\psi}^{|\beta| + n +4} N^{(|\beta| + n+3 -M)^{+}} C_{\Theta}^{|\gamma|+1} B^{|\alpha| -|\gamma| + 1} N^{(|\gamma| +q+1-M)^{+}} \|g\|_{0}.
	\end{multline*}	
\end{linenomath} 
By \eqref{Est_g} and taking advantage from the Remark \ref{Rk-1},
we conclude that there are two positive constants, $C_{1}$ and $C_{2}$,
such that
\begin{linenomath}	
	\begin{align}\label{Est-I-p/r}
	\| \psi_{N}^{(\beta)} (1-\widetilde{\Theta}_{q})\Theta_{N}^{(\gamma)} D^{\alpha} g \|_{\frac{p}{r}}
	\leq 
	C_{1}^{\sigma+1} 
	C_{2}^{ m +|\gamma|  + 1 } 
	N^{s(q+ |\gamma|+ \sigma + n+4 -M)^{+}}, 
	\end{align}	
\end{linenomath} 
where $m=  |\alpha| -|\gamma|$ and $\sigma= |\beta| + |\delta| + 2k $.\\
Now we have to handle the last term on the right hand side of (\ref{Step-0-p/r}), we have
\begin{linenomath}	
	\begin{multline}\label{Step--p/r}
	\|\widetilde{\psi}  \Lambda_{r}^{p-1}  \psi_{N}^{(\beta)}  \widetilde{\Theta}_{q}  w \|_{\frac{1}{r}}
	\leq 
	\|\widetilde{\psi}  \Lambda_{r}^{p-1}  [\psi_{N}^{(\beta)},  \widetilde{\Theta}_{q}]  w \|_{\frac{1}{r}}
	+
	\|\widetilde{\psi}  \widetilde{\Theta}_{q}\Lambda_{r}^{p-1}  \psi_{N}^{(\beta)}   w \|_{\frac{1}{r}}
	\\
	=I_{1} +I_{2}.
	\end{multline}
\end{linenomath}	

\textbf{Estimate of the term $I_{1}$}. We have
\begin{linenomath}
	\begin{multline*}
	\left[\psi_{N}^{(\beta)}, \widetilde{\Theta}_{q} \right] 
	\Theta_{N}^{(\gamma)} D^{\alpha}g
	=
	\!\displaystyle\longsum[25]_{1 \leq |\mu| \leq q - 1}\frac{1}{\mu!} \psi_{N}^{(\beta+\mu)} \widetilde{\Theta}_{q}^{(\mu)}\Theta_{N}^{(\gamma)} D^{\alpha}g
	+ \,\mathscr{R}_{q}\left( \left[\psi_{N}^{(\beta)}, \widetilde{\Theta}_{q} \right]  \right) w,
	\end{multline*}
\end{linenomath}
where 
\begin{linenomath}
	\begin{multline*}
	\mathscr{R}_{q}\left( \left[\psi_{N}^{(\beta)}, \widetilde{\Theta}_{q} \right]  \right) w(x)
	\\
	=
	- \frac{q}{\left(2\pi\right)^{4n}}  
	\! \sum_{|\mu| =  q} \frac{1}{\mu!}\!
	\int e^{ix\xi}\! \!\int\!\! \widehat{\psi}_{N}^{(\beta+\mu)}(\xi-\eta) 
	\int_{0}^{1}\!\!  \widetilde{\Theta}_{q}^{(\mu)} \left( \eta +t(\xi -\eta) \right)(1 -t)^{q -1} dt \widehat{w} (\eta) d\eta  d\xi,
	\end{multline*}
\end{linenomath}
we recall that $w= \Theta_{N}^{(\gamma)} D^{\alpha}g$.
Then:
\begin{linenomath}
	\begin{align}\label{Est-I-1-p/r}
	I_{1}
	\leq
	\displaystyle\longsum[25]_{1 \leq |\mu| \leq q - 1}\frac{1}{\mu!} \| \psi_{N}^{(\beta+\mu)} \widetilde{\Theta}_{q}^{(\mu)}w\|_{\frac{p}{r}}
	+ \| \mathscr{R}_{q}\left( \left[\psi_{N}^{(\beta)}, \widetilde{\Theta}_{q} \right]  \right) w\|_{\frac{p}{r}}.
	\end{align}
\end{linenomath}
We begin to estimate the last term on the right hand side.
To do it  we use the Lemma \ref{Rem0} adapted to this situation.
Let $q= |\alpha|- |\gamma| + \left\lfloor \frac{n}{2} \right\rfloor +1$,
then there are two positive constants $C_{1}$ and $C_{2}$ such that
\begin{linenomath}
	\begin{align}\label{Est-Rq-I1}
	\| \mathscr{R}_{q}\left( \left[\psi_{N}^{(\beta)}, \widetilde{\Theta}_{q} \right]  \right) w\|_{\frac{p}{r}}
	\leq 
	C_{1}^{\sigma +1} C_{2}^{2 m +|\gamma|+ 2n + 4} 
	N^{s(2 m +|\gamma|+ \sigma + 2n + 4 - M)^{+}},
	\end{align}
\end{linenomath}
where $m =|\alpha| -|\gamma|$ and $\sigma= |\beta| + |\delta| +2k$.\\
Now we focus on the terms in the sum in \eqref{Est-I-1-p/r}.
 Since $\widetilde{\Theta}_{q}^{(\mu)}(\xi) = 0$ 
in $ \Gamma_{1} \cap \lbrace \xi\in \mathbb{R}^{n}\,: \, |\xi| \geq q\rbrace  $,
then $ \widetilde{\Theta}_{q}^{(\mu)}(\xi)\Theta_{N}^{(\gamma)} (\xi)$
is supported in the region $2^{-1}N \leq | \xi| \leq q$. 
Using the same strategy used to obtain (\ref{Est-I-p/r}) we obtain
\begin{linenomath}
	\begin{multline*}
	\| \psi_{N}^{(\beta+\mu)} \widetilde{\Theta}_{q}^{(\mu)}w\|_{\frac{p}{r}}
	\leq
	C_{\psi}^{|\beta|+|\mu| + n +4} N^{(|\beta|+ |\mu| + n+3 -M)^{+}} 
	\| \widetilde{\Theta}_{q}^{(\mu)} \Theta_{N}^{(\gamma)} D^{\alpha} g\|_{\frac{p}{r}}
	\\
	\leq
	C_{\psi}^{|\beta|+|\mu| + n +4} N^{(|\beta|+ |\mu| + n+3 -M)^{+}} 
	C_{\widetilde{\Theta}}^{|\mu|+1} q^{(|\mu|-M)^{+}} C_{\Theta}^{|\gamma| + 1} N^{(|\gamma| - M)^{+}}
	\\
	\hspace*{15em}\times
	\| (1+|\xi|)^{|\alpha| - |\gamma| -|\mu| +\frac{p}{r}} \widehat{g}(\xi)\|_{L^{2}( 2^{-1}N < | \xi| < q )}.
	\end{multline*}
\end{linenomath}
We distinguish two cases.\\
Case $|\mu| \leq |\alpha|-|\gamma|$: the right hand side of above inequality 
is bounded by
\begin{linenomath}
	\begin{multline*}
	C_{\psi}^{|\beta|+|\mu| + n +4}  C_{\widetilde{\Theta}}^{|\mu|+1} C_{\Theta}^{|\gamma| + 1}
	N^{(|\beta|+ |\mu| +|\gamma|+ n+3 -M)^{+}} q^{(|\mu|-M)^{+}}
	(1+q)^{|\alpha| - |\gamma| -|\mu| +1} \|g\|
	\\
	\leq
	C_{\psi}^{|\beta|+|\mu| + n +4}  C_{\widetilde{\Theta}}^{|\mu|+1} C_{\Theta}^{|\gamma| + 1}
	N^{(|\beta|+ |\mu| +|\gamma|+ n+3 -M)^{+}} 
	(1+q)^{|\alpha| - |\gamma| +1} \|g\|.
	\end{multline*}
\end{linenomath}
Since $ |\alpha| - |\gamma| +1 \leq 1+q \leq N$, by Lemma \ref{L-1} there is a constant $C_{1}$
such that 
$ (1+q)^{|\alpha| - |\gamma| +1} \leq C_{1}^{|\alpha| - |\gamma| +1} N^{ \left( 1+q -M\right)^{+}}
= C_{1}^{|\alpha| - |\gamma| +1} N^{ \left(  |\alpha|- |\gamma| + \left\lfloor \frac{n}{2} \right\rfloor +2 -M\right)^{+}}$,
then we obtain 
\begin{linenomath}
	\begin{multline*}
	\| \psi_{N}^{(\beta+\mu)} \widetilde{\Theta}_{q}^{(\mu)}w\|_{\frac{p}{r}}
	\leq
	C_{\psi}^{|\beta|+|\mu| + n +4}  C_{\widetilde{\Theta}}^{|\mu|+1} C_{\Theta}^{|\gamma| + 1} C_{1}^{|\alpha| - |\gamma| +1}
	\\
	\times
	N^{(|\beta|+ |\mu| +  |\alpha| + \left\lfloor \frac{n}{2} \right\rfloor  +  n+5-M)^{+}} \|g\|
	\\
	\leq
	C_{2}^{|\beta|+1} C_{3}^{2m+|\gamma|+1} N^{(2m +|\gamma|+|\beta|+ \left\lfloor \frac{n}{2} \right\rfloor  +  n+5-M)^{+}} \|g\|,
	\end{multline*}
\end{linenomath}
where $C_{2}$ and $C_{3}$ are suitable constants independent of $N$ and $m=|\alpha|-|\gamma|$. \\
Case $|\alpha|-|\gamma| <  |\mu| \leq |\alpha|-|\gamma| +\left\lfloor \frac{n}{2} \right\rfloor  + 1$
(we remark, also in this case, that $\# \{\mu \in \mathbb{N}^{n} \,:\, |\alpha| - |\gamma|< |\mu| \leq |\alpha| - |\gamma| +  \left\lfloor \frac{n}{2} \right\rfloor +1\}$
is finite and it can be roughly estimated $2^{|\alpha| - |\gamma| +  \left\lfloor \frac{n}{2} \right\rfloor +1}$.)
We observe that 
$(1+|\xi|)^{|\alpha| - |\gamma| -|\mu| +\frac{p}{r}} \leq 1$, 
and there are at most $\left\lfloor \frac{n}{2} \right\rfloor  +1$ such terms,
then 
\begin{linenomath}
	\begin{multline*}
	\| \psi_{N}^{(\beta+\mu)} \widetilde{\Theta}_{q}^{(\mu)}w\|_{\frac{p}{r}}
	\leq
	C_{\psi}^{|\beta|+|\mu| + n +4}  C_{\widetilde{\Theta}}^{|\mu|+1} C_{\Theta}^{|\gamma| + 1}
	N^{(|\beta|+ 2|\mu| +|\gamma|+ n+3 -M)^{+}} q^{(|\mu|-M)^{+}}
    \|g\|
    \\
    \leq
    C_{4}^{|\beta|+1} C_{5}^{2m+|\gamma|+1} N^{(2m +|\gamma|+|\beta|+ 2\left\lfloor \frac{n}{2} \right\rfloor  +  n+6-M)^{+}} \|g\|,
	\end{multline*}
\end{linenomath}	
where $C_{4}$ and $C_{5}$ are suitable constants independent of $N$ and $m=|\alpha|-|\gamma|$. \\
By \eqref{Est_g} and the above considerations, we conclude that there are two new positive constants $C_{1}$ and $C_{2}$ such that
\begin{linenomath}
	\begin{multline*}
	\displaystyle\longsum[25]_{1 \leq |\mu| \leq q - 1}\frac{1}{\mu!} \| \psi_{N}^{(\beta+\mu)} \widetilde{\Theta}_{q}^{(\mu)}w\|_{\frac{p}{r}}
     \\
     \leq
     C_{1}^{|\beta|+|\delta|+ k +1} C_{2}^{2m+|\gamma|+1} 
     N^{(2m +|\gamma|+|\beta| +|\delta| + 2k + 2\left\lfloor \frac{n}{2} \right\rfloor  +  n+6-M)^{+}} .
 	\end{multline*}
 \end{linenomath}	
By the above estimate and  \eqref{Est-Rq-I1}, we obtain that
there are two positive constants $\widetilde{C}_{1}$ and $\widetilde{C}_{2}$
 such that
 \begin{linenomath}
 	\begin{align}\label{Est-I-1-p/r-f}
 	I_{1} \leq
 	\widetilde{C}_{1}^{\sigma+1}\widetilde{C}_{2}^{2 m + |\gamma|+1}
 	N^{s(2 m +|\gamma|+ \sigma + 2n + 6 - M)^{+}}
 	\end{align}
 \end{linenomath}
where $m =|\alpha| -|\gamma|$ and $\sigma= |\beta| + |\delta| +2k$.

\textbf{Estimate of the term $I_{2}$}, (\ref{Step--p/r}).
We use the same strategy used in the proof of the Proposition \ref{Basic-Est-D}.
We introduce the couple $\dbtilde{\psi}(x)$, $\dbtilde{\Theta}(D)$,
where $\dbtilde{\psi} \in \mathscr{D}(\widetilde{\Omega}_{3})$,
$\Omega_{3} \Subset \widetilde{\Omega}_{3} \Subset \Omega_{4} $
and such that $\dbtilde{\psi} (x) \equiv 1$ on $\overline{\Omega}_{3}$
and  $\dbtilde{\Theta}(D)$ is a zero order pseudodifferential operator
with associated symbol $\dbtilde{\Theta}(\xi)$ supported in 
$\Gamma_{3}$, $\Gamma_{3} \Subset   \Gamma_{4}$
and such that $\dbtilde{\Theta}(\xi) \equiv 1$
on $\overline{\Gamma}_{2} \cap \{ |\xi| \geq 1\}$.
We recall that $\widetilde{\psi} \in \mathscr{D}(\Omega_{3})$,
the sequence $\widetilde{\Theta}_{q}$
is associated to the couple $(\Gamma_{1},\Gamma_{2})$
and that the symbol $\text{\textcursive{p}} \,(x,\xi)$, in (\ref{Der-Est-M}), is elliptic
in $\widetilde{\Omega}_{3} \times \Gamma_{4}$ and $\text{\textcursive{p}} \,(x,\xi)\geq c_{0}>0$ in $\Omega_{3} \times \Gamma_{3}$.
We have, as $\dbtilde{\psi}\widetilde{\psi}= \widetilde{\psi}$:
\begin{linenomath}
	\begin{align*}
	\widetilde{\psi}  \widetilde{\Theta}_{q}\Lambda_{r}^{p-1}  \psi_{N}^{(\beta)}   w
	=
	\dbtilde{\psi} \widetilde{\psi}  \dbtilde{\Theta}\widetilde{\Theta}_{q} \Lambda_{r}^{p-1}  \psi_{N}^{(\beta)} w 
	+
	\dbtilde{\psi} \widetilde{\psi}  (1-\dbtilde{\Theta} ) \widetilde{\Theta}_{q} \Lambda_{r}^{p-1}  \psi_{N}^{(\beta)}  w, 
	\end{align*}
\end{linenomath}
so
\begin{linenomath}
	\begin{align}\label{Est-I2-0}
	I_{2} 
	\leq
	\| \dbtilde{\psi} \widetilde{\psi}  \dbtilde{\Theta}\widetilde{\Theta}_{q} \Lambda_{r}^{p-1}  \psi_{N}^{(\beta)} w \|_{\frac{1}{r}}
	+
	\| \dbtilde{\psi} \widetilde{\psi} (1-\dbtilde{\Theta} ) \widetilde{\Theta}_{q} \Lambda_{r}^{p-1}  \psi_{N}^{(\beta)}  w\|_{\frac{1}{r}}.
	\end{align}
\end{linenomath}
We begin to handle the second term on the right hand side.
Since $ (1-\dbtilde{\Theta} (\xi)) \widetilde{\Theta}_{q}(\xi)$ is supported 
in $\{ |\xi| \leq 1\}$,
we obtain
\begin{linenomath}
	\begin{multline}\label{Est-I2-01}
	\| \dbtilde{\psi} \widetilde{\psi} (1-\dbtilde{\Theta} ) \widetilde{\Theta}_{q} \Lambda_{r}^{p-1}  \psi_{N}^{(\beta)}  w\|_{\frac{1}{r}}
	\leq 
	C_{0} \| (1-\dbtilde{\Theta} ) \widetilde{\Theta}_{q} \Lambda_{r}^{p-1}  \psi_{N}^{(\beta)}  w\|_{\frac{1}{r}}
	\\
	=
	C_{0} \| ( 1+|\xi|)^{\frac{p}{r}}(1-\dbtilde{\Theta} (\xi)) \widetilde{\Theta}_{q} (\xi)  \widehat{\psi_{N}^{(\beta)}  w}(\xi)\|
	\leq
	C_{0}2^{p/r}\| \psi_{N}^{(\beta)}  w\|
	\\
	= C_{0}2^{p/r}\| \psi_{N}^{(\beta)}\Theta_{N}^{(\gamma)} D^{\alpha} \varphi_{N}^{(\delta)} P^{k} u  \| .
	\end{multline}
\end{linenomath}
We focus, now, on the first term on the right hand side of (\ref{Est-I2-0}). We have
\begin{linenomath}
	\begin{align}\label{Est-I2-1}
	\| \dbtilde{\psi} \widetilde{\psi}  \dbtilde{\Theta}\widetilde{\Theta}_{q} \Lambda_{r}^{p-1}  \psi_{N}^{(\beta)} w \|_{\frac{1}{r}}
	\leq
	\| \dbtilde{\psi} \dbtilde{\Theta} \widetilde{\psi}  \widetilde{\Theta}_{q} \Lambda_{r}^{p-1}  \psi_{N}^{(\beta)} w \|_{\frac{1}{r}}
	+
	\| \dbtilde{\psi} [\dbtilde{\Theta}, \widetilde{\psi}]  \widetilde{\Theta}_{q} \Lambda_{r}^{p-1}  \psi_{N}^{(\beta)} w \|_{\frac{1}{r}}.
	\end{align}	
\end{linenomath} 	
About the second term we have
\begin{linenomath}
	\begin{multline}\label{Est-I2-11}
	\| \dbtilde{\psi} [\dbtilde{\Theta}, \widetilde{\psi}]  \widetilde{\Theta}_{q} \Lambda_{r}^{p-1}  \psi_{N}^{(\beta)} w \|_{\frac{1}{r}}
	\leq 
	C_{0} \| \Lambda_{r}^{p-1}  \psi_{N}^{(\beta)} w \|_{-1+\frac{1}{r}}
	\\
	\leq 
	C_{0} \|  \psi_{N}^{(\beta)} w \|
	=
	C_{0} \| \psi_{N}^{(\beta)}\Theta_{N}^{(\gamma)} D^{\alpha} \varphi_{N}^{(\delta)} P^{k} u  \|,
	\end{multline}	
\end{linenomath} 	
where $C_{0}$ does not depend on $q$.\\
In order to estimate the first term on the right hand side of (\ref{Est-I2-1})
we introduce $Q(x,D)$ the zero order operator associated to the symbol
$\text{\textcursive{q}}(x,\xi) = \displaystyle\frac{ \dbtilde{\psi}(x)  \dbtilde{\Theta} (\xi)}{\text{\textcursive{p}}(x,\xi)}$,
where $\text{\textcursive{p}} \,(x,\xi)$ is the symbol associated to the zero order operator
$\text{\textcursive{p}}$, in (\ref{Der-Est-M}).
We point out that $\text{\textcursive{q}}(x,\xi) $ is well defined
as $ |\text{\textcursive{p}} |\geq  c_{0}> 0$
on the support of $\dbtilde{\psi}(x)  \dbtilde{\Theta} (\xi)$.
We have
\begin{linenomath}	
	\begin{align*}
	\left(Q\circ \text{\textcursive{p}}\right)v(x)
	= \frac{1}{(2\pi)^{2n}} \int  e^{ix\xi}\text{\textcursive{q}}(x,\xi) \widehat{\text{\textcursive{p}}}(\xi-\eta,\eta)\widehat{v}(\eta) d\eta\,d\xi,
	\end{align*}
\end{linenomath}	
where $\widehat{\text{\textcursive{p}}}(\cdot,\cdot)$ is the Fourier transform
of $\text{\textcursive{p}}$ with respect to $x$.\\
Using the Taylor expansion of $\text{\textcursive{q}}(x,\xi) $
with respect to $\tau$, $\tau= \xi-\eta$,
we have
\begin{linenomath}	
	\begin{align*}
	\text{\textcursive{q}}(x,\eta+\tau)= 
	\frac{ \dbtilde{\psi}(x)  \dbtilde{\Theta}(\eta)}{\text{\textcursive{p}}(x,\eta)}
	+\sum_{j=1 }^{n}\tau_{j}\int_{0}^{1} \left(\frac{ \dbtilde{\psi} \dbtilde{\Theta}}{\text{\textcursive{p}}}\right)^{(j)}(x,\eta + t\tau)\, dt .
	\end{align*}
\end{linenomath}	
We obtain
\begin{linenomath}	
	\begin{multline*}
	\left(Q\circ \text{\textcursive{p}}\right)v(x)
	=
	\dbtilde{\psi}(x)  \dbtilde{\Theta}v(x)
	\\
	+
	\frac{1}{(2\pi)^{2n}} \sum_{j=1}^{n}
	\int e^{ix\eta} \int\int_{0}^{1} e^{ix \tau} \left(\frac{ \dbtilde{\psi}(x)  \dbtilde{\Theta} }{\text{\textcursive{p}}}\right)^{(j)}(x,\eta + t\tau) \tau_{j} 
	\widehat{\text{\textcursive{p}}}(\tau,\eta)\widehat{v}(\eta)\, dt\,d\tau \,\,  \widehat{v}(\eta)  \, d\eta
	\\
	=\widetilde{\psi} \widetilde{\Theta}_{\widetilde{m}}v(x) +\mathscr{R}_{1}v(x).
	\end{multline*}
\end{linenomath}	
Where the symbol associated to the operator $\mathscr{R}_{1}(x,D)$ is
\begin{linenomath}	
	\begin{multline*}
	\text{\textcursive{r}}_{1}(x,\xi) =
	\sum_{j=1 }^{n} \iint e^{i(y-x)(\xi-\eta)} \text{\textcursive{p}}_{(j)}(y,\xi) 
	\int_{0}^{1} \text{\textcursive{q}}^{(j)}(x, \xi + t (\eta-\xi)) dt\, dy\,\frac{d\eta}{(2\pi)^{n}},
	\end{multline*}
\end{linenomath}	
moreover the following estimate holds
\begin{linenomath}	
	\begin{align*}
	\left| \text{\textcursive{r}}_{1}(x,\xi)\right|
	\leq \tilde{C} \left( 1 + |\xi|^{2}\right)^{-1/2}, 
	\end{align*}
\end{linenomath}	
where $\tilde{C} $ depends only on $n$ and on the derivatives of
$ \text{\textcursive{p}}(y,\xi) $ up to order $\left\lfloor \frac{n}{2} \right\rfloor + 2$
with respect to $y$.
By the Calderon-Vaillancourt theorem, see \cite{KG} or \cite{Hwang-1987},
we have
\begin{linenomath}
	\begin{align*}
	\|  \dbtilde{\psi} \dbtilde{\Theta} \widetilde{\psi}  \widetilde{\Theta}_{q} \Lambda_{r}^{p-1}  \psi_{N}^{(\beta)} w \|_{\frac{1}{r}}^{2}
	\leq	
	\|  Q\, \text{\textcursive{p}} \widetilde{\psi}  \widetilde{\Theta}_{q} \Lambda_{r}^{p-1}  \psi_{N}^{(\beta)} w \|_{\frac{1}{r}}^{2} 
	+ \| \mathscr{R}_{1} \widetilde{\psi}  \widetilde{\Theta}_{q} \Lambda_{r}^{p-1}  \psi_{N}^{(\beta)} w \|_{\frac{1}{r}}^{2}
	\\
	\leq
	C_{1}	\| \text{\textcursive{p}}  \widetilde{\psi}  \widetilde{\Theta}_{q} \Lambda_{r}^{p-1}  \psi_{N}^{(\beta)} w \|_{\frac{1}{r}}^{2} 
	+ \tilde{C}_{1} \| \widetilde{\psi}  \widetilde{\Theta}_{q} \Lambda_{r}^{p-1}  \psi_{N}^{(\beta)} w \|_{-1+\frac{1}{r}},
	\end{align*}
\end{linenomath}
where the constants $C_{1}$ and $\tilde{C}_{1}$ are suitable positive constants
independent of $\alpha$, $\beta$, $\gamma$, $\delta$, $k$ and $N$.\\
In view of the Theorem \ref{MicLocEst}, there is a positive constant $C$
such that the following estimate holds
\begin{linenomath}	
	\begin{multline}\label{Est-I2-2}
	\|  \dbtilde{\psi} \dbtilde{\Theta} \widetilde{\psi}  \widetilde{\Theta}_{q} \widetilde{w} \|_{\frac{1}{r}}^{2}
	+ \sum_{j=1}^{n} \left( \| P^{j} \widetilde{\psi}  \widetilde{\Theta}_{q} \widetilde{w} \|^{2}_{0} 
	+ \|P_{j} \widetilde{\psi}  \widetilde{\Theta}_{q} \widetilde{w}\|_{-1}^{2}\right)
	\\
	\leq C \left( \sum_{\ell=0}^{n}\left|\langle E_{\ell}P \widetilde{\psi}  \widetilde{\Theta}_{q} \widetilde{w} , E_{\ell} \widetilde{\psi}  \widetilde{\Theta}_{q} \widetilde{w}\rangle\right|  
	+ \| \widetilde{\psi}  \widetilde{\Theta}_{q} \widetilde{w}\|_{0}^{2} + \| \psi_{N}^{(\beta)}\Theta_{N}^{(\gamma)} D^{\alpha} \varphi_{N}^{(\delta)} P^{k} u  \|\right).
	\end{multline} 
\end{linenomath}
where $\widetilde{w} = \Lambda_{r}^{p-1}  \psi_{N}^{(\beta)} w =
\Lambda_{r}^{p-1} \psi_{N}^{(\beta)}\Theta_{N}^{(\gamma)} D^{\alpha} \varphi_{N}^{(\delta)} P^{k} u $
and $q= |\alpha|- |\gamma| + \left\lfloor \frac{n}{2} \right\rfloor +1$.
We recall that $E_{\ell}=D_{\ell} \psi \Lambda_{-1}$, where $\psi$ belongs to
$\mathscr{D}\left(\Omega\right)$ and is identically one on
$\Omega_{4}$, $\Omega_{4} \Subset \Omega$.
We point out that $E_{\ell}$ are zero order pseudodifferential operators. \\
We have to estimate the terms in the sum. We proceed as in the case $H^{1/r}$
with the difference that in this case we have to handle new ingredients,
in particular, the presence of the operator
$\Lambda_{r}^{p-1}$. We denote by $F_{\ell}$ the terms in the sum. We have
\begin{linenomath}
	\begin{multline}\label{Est-El-p/r}
	F_{\ell}
	\leq 
	\left|\langle E_{\ell} [P,  \widetilde{\psi}] \widetilde{\Theta}_{q}  \widetilde{w},
	E_{\ell} \widetilde{\psi} \widetilde{\Theta}_{q} \widetilde{w}\rangle\right| 
	+
	\left|\langle E_{\ell}  \widetilde{\psi} [\widetilde{P}_{N}, \widetilde{\Theta}_{q}]  \widetilde{w},
	E_{\ell} \widetilde{\psi} \widetilde{\Theta}_{q} \widetilde{w}\rangle\right| 
	\\
	+
	\left|\langle E_{\ell}  \widetilde{\psi} \widetilde{\Theta}_{q}   [\widetilde{P}_{N},\Lambda_{r}^{p-1}]  \psi_{N}^{(\beta)}  w,
	E_{\ell} \widetilde{\psi} \widetilde{\Theta}_{q} \widetilde{w}\rangle\right| 
	+
	\left|\langle E_{\ell}  \widetilde{\psi} \widetilde{\Theta}_{q}   \Lambda_{r}^{p-1} [\widetilde{P}_{N},\psi_{N}^{(\beta)}]  w,
	E_{\ell} \widetilde{\psi} \widetilde{\Theta}_{q} \widetilde{w}\rangle\right| 
	\\
	\hspace*{-5em}
	+
	\left|\langle E_{\ell}  \widetilde{\psi} \widetilde{\Theta}_{q}   \Lambda_{r}^{p-1} \psi_{N}^{(\beta)} 
	[\widetilde{P}_{N},\Theta_{N}^{(\gamma)} D^{\alpha}] \varphi_{N}^{(\delta)} P^{k} u,
	E_{\ell} \widetilde{\psi} \widetilde{\Theta}_{q} \widetilde{w}\rangle\right|
	\\
	+
	\left|\langle E_{\ell}  \widetilde{\psi} \widetilde{\Theta}_{q}   \Lambda_{r}^{p-1} \psi_{N}^{(\beta)} \Theta_{N}^{(\gamma)} D^{\alpha} 
	[\widetilde{P}_{N},\varphi_{N}^{(\delta)}] P^{k} u,
	E_{\ell} \widetilde{\psi} \widetilde{\Theta}_{q} \widetilde{w}\rangle\right|
	\\
	+
	\left|\langle E_{\ell}  \widetilde{\psi} \widetilde{\Theta}_{q}   \Lambda_{r}^{p-1} \psi_{N}^{(\beta)} \Theta_{N}^{(\gamma)} D^{\alpha}  \varphi_{N}^{(\delta)} P^{k+1} u,
	E_{\ell} \widetilde{\psi} \widetilde{\Theta}_{q} \widetilde{w}\rangle\right|
	=\sum_{i=1}^{7}F_{\ell,i},
	\end{multline}
\end{linenomath} 
where $\widetilde{P}_{N}$, was introduced in \eqref{t-PN}.
We point out that the sequence $\widetilde{\psi}_{N}$,
used to introduce $\widetilde{P}_{N}$, is associated to
the couple $(\Omega_{2},\Omega_{3})$,  $\Omega_{3} \Subset \Omega_{4}$, 
$\widetilde{\psi}$ is identically one on $\Omega_{4}$.
%

Looking at the proof of the case $1/r$, we give without much details a bound to $F_{\ell,i}$, $i=1,\dots,7$.

\textbf{Term $F_{\ell,1}$}.
In order to bound this term we need to handle $[E_{\ell},\widetilde{P}^{j}_{N}]$,
where $\widetilde{P}^{j}_{N}$ was introduced in \eqref{t-PN}. Explicitly it is given by
\begin{linenomath}
	\begin{multline}\label{CommEPjN}
	[E_{\ell},\widetilde{P}^{j}_{N}]
	=2\sum_{ j_{1}=1}^{n}[E_{\ell},\widetilde{a}_{N,j_{1},j} D_{j_{1}} ]
	=2\sum_{ j_{1}=1}^{n} [D_{\ell} \psi \Lambda_{-1}, \widetilde{a}_{N,j_{1},j} D_{j_{1}}]
	\\
	=2\sum_{ j_{1}=1}^{n}\left(D_{\ell}\psi [\Lambda_{-1},\widetilde{a}_{N,j_{1},j}] D_{j_{1}} 
	+ D_{\ell} \widetilde{a}_{N,j_{1},j}\psi^{(j_{1})} \Lambda_{-1} + \widetilde{a}_{N,j_{1},j}^{(\ell)} D_{j_{1}}\psi\Lambda_{-1}\right),
	\end{multline}
\end{linenomath}
where 
$	[ \Lambda_{-1},\widetilde{a}_{N,j_{1},j}]= 
\sum_{\varkappa =1}^{n} \widetilde{a}_{N,j_{1},j}^{(\varkappa)} \Lambda_{-1}^{(\varkappa)}
+ \mathscr{R}_{2}\left( [ \Lambda_{-1},\widetilde{a}_{N,j_{1},j}]\right)$,
the terms in the sum have order $-2$ and $\mathscr{R}_{2}$ has order $-3$.
All the terms in the summand on the right hand side of \eqref{CommEPjN}
are zero order operators,
so $[E_{\ell}, \widetilde{P}^{j}_{N}]$ is a zero order operator.
By the Theorem $18.1.11'$ in \cite{H_Book-3}
we conclude that $\|[E_{\ell},\widetilde{P}^{j}_{N}] \|_{L^{2}\rightarrow L^{2}} \leq C$,
where $C$ is independent of $N$.\\
Since $ \widetilde{P}^{j}_{N}\phantom{A}^{\!\!\!\!\!\!\!^{*}}
=  \widetilde{P}^{j}_{N} - 2 i \sum_{j_{1}=1}^{n} \widetilde{a}_{N,j,j_{1}}^{(j_{1})}(x) $,
we have
\begin{linenomath}
	\begin{multline}\label{Est-Fl1}
	F_{\ell,1} \leq 
	\sum_{j=1}^{n} \left| \langle E_{\ell} \widetilde{P}^{j}_{N} \widetilde{\psi}^{(j)} \widetilde{\Theta}_{q} \widetilde{w},
	E_{\ell}  \widetilde{\psi} \widetilde{\Theta}_{q} \widetilde{w} \rangle \right| 
	+
	\sum_{j,j_{1}=1}^{n} \left| \langle E_{\ell} a_{j,j_{1}}  \widetilde{\psi}^{(j+j_{1})} \widetilde{\Theta}_{q} \widetilde{w},
	E_{\ell}  \widetilde{\psi} \widetilde{\Theta}_{q} \widetilde{w} \rangle \right|
	\\
	\hspace*{18em}
	+
	\sum_{j_{1}=1}^{n} \left| \langle E_{\ell} b_{,j_{1}}  \widetilde{\psi}^{(j_{1})} \widetilde{\Theta}_{q} \widetilde{w},
	E_{\ell}  \widetilde{\psi} \widetilde{\Theta}_{q} \widetilde{w} \rangle \right|  
	\\
	\leq 
	\varepsilon \sum_{j=1}^{n}\| P^{j}  \widetilde{\psi}\widetilde{\Theta}_{q} \widetilde{w}\|^{2}
	+ C_{\varepsilon} \| \psi_{N}^{(\beta)} \Theta_{N}^{(\gamma)} D^{\alpha}  \varphi_{N}^{(\delta)} P^{k} u \|^{2}_{\frac{p-1}{r}}.
	\end{multline}
\end{linenomath}    
The first term on the right hand side can be absorbed by the left hand side of (\ref{Est-I2-2}).

\textbf{Term $F_{\ell,2}$}.
We have
\begin{linenomath}
	\begin{multline}\label{Est-Fl2}
	F_{\ell,2} \leq 
	\sum_{j =1}^{n} \left|\langle E_{\ell} \widetilde{\psi} \widetilde{P}_{N,j} \widetilde{\Theta}_{q}^{(j)} \widetilde{w},  E_{\ell} \widetilde{\psi}  \widetilde{\Theta}_{q}\widetilde{w}\rangle\right|
	+ C_{0} \|\widetilde{w}\|^{2} 
	+ C_{1} \| w\|^{2}
	\\
	\leq
	\sum_{j =1}^{n} 
	\left( \left|\langle [ E_{\ell} \widetilde{\psi}, \widetilde{P}_{N,j} ] \widetilde{\Theta}_{q}^{(j)} \widetilde{w},  E_{\ell} \widetilde{\psi}  \widetilde{\Theta}_{q}\widetilde{w}\rangle\right|
	+ \left|\langle \widetilde{P}_{N,j}  E_{\ell} \widetilde{\psi} \widetilde{\Theta}_{q}^{(j)} \widetilde{w},  E_{\ell} \widetilde{\psi}  \widetilde{\Theta}_{q}\widetilde{w}\rangle\right|\right)
	\\
	+ C_{0} \|\widetilde{w}\|^{2} 
	+ C_{1} \| w\|^{2}
	.
	\end{multline}
\end{linenomath}    
We handle separately the terms in the sum.
We begin with the second term.\\
Due to the fact that
\begin{linenomath}
	$$
	\widetilde{P}_{N,j}^{*} = \widetilde{P}_{N,j} - i\sum_{j_{2},j_{1}=1}^{n} \left(a_{j_{2},j_{1}}^{(j+j_{2})} D_{j_{1}} +  a_{j_{2},j_{1}}^{(j+j_{1})} D_{j_{2}}\right)
	-\sum_{j_{2},j_{1}=1}^{n} a_{j_{2},j_{1}}^{(j+j_{1}+j_{2})},
	$$
\end{linenomath} 
$\widetilde{P}_{N,j} $ as in \eqref{t-PN}, we have 
\begin{linenomath}
	\begin{multline}\label{Est-Fl2-1}
	\left|\langle \widetilde{P}_{N,j} E_{\ell} \widetilde{\psi} \widetilde{\Theta}_{q}^{(j)} \widetilde{w},  E_{\ell} \widetilde{\psi}  \widetilde{\Theta}_{q}\widetilde{w}\rangle\right|
	\leq
	\left|\langle \left(\widetilde{P}_{N,j} -\widetilde{P}_{N,j} ^{*}\right) E_{\ell} \widetilde{\psi} \widetilde{\Theta}_{q}^{(j)} \widetilde{w},  E_{\ell} \widetilde{\psi}  \widetilde{\Theta}_{q}\widetilde{w}\rangle\right|
	\\
	\hspace{13em}
	+
	\left|\langle \widetilde{P}_{N,j}^{*}  E_{\ell} \widetilde{\psi} \widetilde{\Theta}_{q}^{(j)} \widetilde{w},  E_{\ell} \widetilde{\psi}  \widetilde{\Theta}_{q}\widetilde{w}\rangle\right|
	\\
	\leq
	\left|\langle  E_{\ell} \widetilde{\psi} \widetilde{\Theta}_{q}^{(j)} \widetilde{w}, [\widetilde{P}_{N,j} , E_{\ell} ]\widetilde{\psi}  \widetilde{\Theta}_{q}\widetilde{w}\rangle\right|
	+
	\left|\langle E_{\ell} \widetilde{\psi} \widetilde{\Theta}_{q}^{(j)} \widetilde{w},  E_{\ell}  P_{j}\widetilde{\psi}  \widetilde{\Theta}_{q}\widetilde{w}\rangle\right|
	\\
	+
	C_{2}\left(\|\widetilde{w}\|^{2}  + \|E_{\ell} \widetilde{\psi} \widetilde{\Theta}_{q}^{(j)} \widetilde{w}\|_{1}\|\widetilde{w}\|\right).
	\end{multline}
\end{linenomath}    
We remark that if we choose $M\geq 2$, then
\begin{linenomath}
	$$
	\|E_{\ell} \widetilde{\psi} \widetilde{\Theta}_{q}^{(j)} \widetilde{w}\|_{1} \leq C_{3} \|\widetilde{\psi}\widetilde{\Theta}_{q}^{(j)} \widetilde{w}\|_{1}
	\leq C_{3}C_{4}\| \widetilde{w}\|,
	$$
\end{linenomath}
as $|\widetilde{\Theta}_{q}^{(\mu)}(\xi) |\leq C_{4} \left(1+|\xi|\right)^{-|\mu|}$ if $|\mu|\leq 2$,
for all $q$ as $M\geq 2$.
The second term on the right hand side can be estimated in the following way
\begin{linenomath}
	\begin{multline*}
	\left|\langle E_{\ell} \widetilde{\psi} \widetilde{\Theta}_{q}^{(j)} \widetilde{w},  E_{\ell}  P_{j}\widetilde{\psi}  \widetilde{\Theta}_{q}\widetilde{w}\rangle\right|
	\leq
	\tilde{C}_{0}\|\widetilde{w}\|\| P_{j}\widetilde{\psi}  \widetilde{\Theta}_{q}\widetilde{w}\|_{-1}
	\leq
	\varepsilon \| P_{j}\widetilde{\psi}  \widetilde{\Theta}_{q}\widetilde{w}\|_{-1}^{2} 
	+
	C_{\varepsilon} \|\widetilde{w}\|^{2}
	\\
	\leq
	\varepsilon \| P_{j}\widetilde{\psi}  \widetilde{\Theta}_{q}\widetilde{w}\|_{-1}^{2} 
	+
	\tilde{C}_{\varepsilon}\| \psi_{N}^{(\beta)} \Theta_{N}^{(\gamma)} D^{\alpha}  \varphi_{N}^{(\delta)} P^{k} u \|^{2}_{\frac{p-1}{r}} ,
	\end{multline*}
\end{linenomath}    
where $\varepsilon$ is small suitable positive constant.\\
In order to bound the first term on the right hand side of
\eqref{Est-Fl2-1} we analyze $[ \widetilde{P}_{N,j},E_{\ell}]$:
\begin{linenomath}
	\begin{multline}\label{CommEP_jN}
	[\widetilde{P}_{N,j}, E_{\ell}]
	=\longsum[7]_{ j_{1},j_{2} =1}^{n} [\widetilde{a}_{N,j_{1},j_{2}}^{(j)} D_{j_{1}} D_{j_{2}}, D_{\ell}\psi \Lambda_{-1} ]
	\\
	=\longsum[7]_{ j_{1},j_{2} =1}^{n} \left\{ D_{j_{1}} \left( \widetilde{a}_{N,j_{1},j_{2}}^{(j)}  D_{\ell}\psi^{(j_{2})} \Lambda_{-1}\right)
	+D_{j_{2}} \left( \widetilde{a}_{N,j_{1},j_{2}}^{(j)}  D_{\ell}\psi^{(j_{1})} \Lambda_{-1}\right)\right\}
	\\
	+
	\longsum[7]_{ j_{1},j_{2} =1}^{n} \left\{ D_{j_{1}} \left( \widetilde{a}_{N,j_{1},j_{2}}^{(j+\ell)} \psi \Lambda_{-1}D_{j_{2}}\right)
	+D_{j_{2}} \left( D_{\ell} \psi [\widetilde{a}_{N,j_{1},j_{2}}^{(j)},  \Lambda_{-1}] D_{j_{1}}\right)
	\right\}
	\\
	+\longsum[7]_{ j_{1},j_{2} =1}^{n} \left\{ \left( \widetilde{a}_{N,j_{1},j_{2}}^{(j+j_{1}+\ell)} \psi + \widetilde{a}_{N,j_{1},j_{2}}^{(j+\ell)} \psi^{(j_{1})}\right)
	\Lambda_{-1}D_{j_{2}} 
	+
	\left(\widetilde{a}_{N,j_{1},j_{2}}^{(j+j_{1})}  D_{\ell}\psi^{(j_{2})}
	+ \widetilde{a}_{N,j_{1},j_{2}}^{(j+j_{2})}  D_{\ell}\psi^{(j_{1})} \right)\Lambda_{-1}
	\right\}
	\\
	+
	\longsum[7]_{ j_{1},j_{2} =1}^{n} 
	D_{\ell} \left( \psi^{(j_{2})} [\widetilde{a}_{N,j_{1},j_{2}}^{(j)},  \Lambda_{-1}] 
	+  \psi [\widetilde{a}_{N,j_{1},j_{2}}^{(j+j_{2})},  \Lambda_{-1}] \right) D_{j_{1}}.
	\end{multline}
\end{linenomath}
We remark that the operators in the round brackets in the first two sums
as well as the operators in the last two sums are zero order operators.
Moreover we recall that $\psi^{(\mu)}\widetilde{\psi}_{N}^{(\nu)}=0$
for every $\mu,\, \nu \in \mathbb{N}^{n}$ with $|\mu|\geq 1$. 
So $[ \widetilde{P}_{N,j},E_{\ell}]$ are pseudodifferential operators
of order $1$ and $\| [ \widetilde{P}_{N,j},E_{\ell}] \|_{L^{2}\rightarrow H^{-1}}\leq C$,
where $C$ is independent of $N$.
Taking advantage from the above considerations,
the first term on the right hand side of (\ref{Est-Fl2-1})
can be estimated in the following way
\begin{linenomath}
	\begin{align*}
	\left|\langle  E_{\ell} \widetilde{\psi} \widetilde{\Theta}_{q}^{(j)} \widetilde{w}, 
	[\widetilde{P}_{N,j}, E_{\ell} ]\widetilde{\psi}  \widetilde{\Theta}_{q}\widetilde{w}\rangle\right|
	\leq
	C_{0} \| \widetilde{\psi} \widetilde{\Theta}_{q}^{(j)} \widetilde{w}\|_{1} \|\widetilde{\psi} \widetilde{\Theta}_{q} \widetilde{w}\|
	\leq
	\tilde{C}_{0}\| \widetilde{w}\|^{2}\,
	.
	\end{align*}
\end{linenomath}    
Concerning the first term in the sum on the right and side of (\ref{Est-Fl2}),
we observe that
$[ E_{\ell} \widetilde{\psi}, \widetilde{P}_{N,j}] = E_{\ell} [\widetilde{\psi}, \widetilde{P}_{N,j}]+ [ E_{\ell} , \widetilde{P}_{N,j}] \widetilde{\psi}$.
In view of \eqref{CommEP_jN} we have that $[ E_{\ell} \widetilde{\psi}, \widetilde{P}_{N,j}]$
are pseudodifferential operators
of order $1$, $\| [ E_{\ell} \widetilde{\psi}, \widetilde{P}_{N,j}] \|_{L^{2}\rightarrow H^{-1}} \leq C$,
where $C$ is independent of $N$.
We obtain
\begin{linenomath}
	\begin{multline*}
	\left|\langle [ E_{\ell} \widetilde{\psi}, \widetilde{P}_{N,j}] \widetilde{\Theta}_{q}^{(j)} \widetilde{w},  E_{\ell} \widetilde{\psi}  \widetilde{\Theta}_{q}\widetilde{w}\rangle\right|
	\leq 
	\|[ E_{\ell} \widetilde{\psi}, \widetilde{P}_{N,j}] \widetilde{\Theta}_{q}^{(j)} \widetilde{w}\| \| E_{\ell} \widetilde{\psi}  \widetilde{\Theta}_{q}\widetilde{w}\|
	\\
	\leq
	C_{1} \| \widetilde{\psi}\widetilde{\Theta}_{q}^{(j)} \widetilde{w}\|_{1} \| E_{\ell} \widetilde{\psi}  \widetilde{\Theta}_{q}\widetilde{w}\|
	\leq 
	C_{1}C_{2} \| \widetilde{w}\|^{2}
	\leq
	C_{3}  \| \psi_{N}^{(\beta)} \Theta_{N}^{(\gamma)} D^{\alpha}  \varphi_{N}^{(\delta)} P^{k} u \|^{2}_{\frac{p-1}{r}}
	.
	\end{multline*}
\end{linenomath}    
Summing up we conclude: there is a new positive constant $C_{\varepsilon}$ such that
\begin{linenomath}
	\begin{align}\label{Est-Fl2-f}
	F_{\ell,2}
	\leq 
	C_{\varepsilon} \| \psi_{N}^{(\beta)} \Theta_{N}^{(\gamma)} D^{\alpha}  \varphi_{N}^{(\delta)} P^{k} u \|^{2}_{\frac{p-1}{r}}
	+
	\varepsilon \sum_{j=1 }^{n}\| P_{j}\widetilde{\psi}  \widetilde{\Theta}_{q}\widetilde{w}\|_{-1}^{2}. 
	\end{align}
\end{linenomath}    
The last term on the right hand side can be absorbed by the left hand side
of (\ref{Est-I2-2}).

\textbf{Term $F_{\ell,3}$}, on the right hand side of (\ref{Est-El-p/r}). 
Since 
\begin{linenomath}
	\begin{align}\label{Comm_ta-Lp-1}
	[ \widetilde{a}_{N,j_{1},j_{2}}, \Lambda_{r}^{p-1}]= \sum_{j =1}^{n} \widetilde{a}_{N,j_{1},j_{2}}^{(j)} (\Lambda_{r}^{p-1})^{(j)} 
	+ \mathscr{R}_{2}\left( [ \widetilde{a}_{N,j_{1},j_{2}}, \Lambda_{r}^{p-1}]\right)
	\end{align}
\end{linenomath}    
where $\mathscr{R}_{2}\left( [ \widetilde{a}_{N,j_{1},j_{2}}, \Lambda_{r}^{p-1}]\right) $
is a pseudodifferential operator of order $\frac{p-1}{r}-2$,
we have
\begin{linenomath}
	\begin{multline}\label{Est-Fl3}
	F_{\ell,3} =
	\left|\langle E_{\ell}  \widetilde{\psi} \widetilde{\Theta}_{q}  [ \widetilde{P}_{N},\Lambda_{r}^{p-1}]  \psi_{N}^{(\beta)}  w,
	E_{\ell} \widetilde{\psi} \widetilde{\Theta}_{q} \widetilde{w}\rangle\right| 
	\\
	\hspace{-9em}
	\leq
	\sum_{j =1}^{n}
	\left|\langle E_{\ell}  \widetilde{\psi} \widetilde{\Theta}_{q}  \widetilde{P}_{N,j} \left(\Lambda_{r}^{p-1}\right)^{(j)}  \psi_{N}^{(\beta)}  w,
	E_{\ell} \widetilde{\psi} \widetilde{\Theta}_{q} \widetilde{w}\rangle\right|
	\\
	+
	C_{0}\left(
	\| \psi_{N}^{(\beta)} \Theta_{N}^{(\gamma)} D^{\alpha}  \varphi_{N}^{(\delta)} P^{k} u \|^{2}_{\frac{p-1}{r}}
	+
	\| \psi_{N}^{(\beta)} \Theta_{N}^{(\gamma)} D^{\alpha}  \varphi_{N}^{(\delta)} P^{k} u \|^{2}
	\right)
	\\
	\hspace{-9em}
	\leq
	\sum_{j =1}^{n}
	\left|\langle \widetilde{P}_{N,j} E_{\ell}  \widetilde{\psi} \widetilde{\Theta}_{q}  \left(\Lambda_{r}^{p-1}\right)^{(j)}  \psi_{N}^{(\beta)}  w,
	E_{\ell} \widetilde{\psi} \widetilde{\Theta}_{q} \widetilde{w}\rangle\right|.
	\\
	+
	\sum_{j =1}^{n}
	\left|\langle [ E_{\ell}  \widetilde{\psi} \widetilde{\Theta}_{q},  \widetilde{P}_{N,j}] \left(\Lambda_{r}^{p-1}\right)^{(j)}  \psi_{N}^{(\beta)}  w,
	E_{\ell} \widetilde{\psi} \widetilde{\Theta}_{q} \widetilde{w}\rangle\right|
	\\
	+
	C_{0}\left(
	\| \psi_{N}^{(\beta)} \Theta_{N}^{(\gamma)} D^{\alpha}  \varphi_{N}^{(\delta)} P^{k} u \|^{2}_{\frac{p-1}{r}}
	+
	\| \psi_{N}^{(\beta)} \Theta_{N}^{(\gamma)} D^{\alpha}  \varphi_{N}^{(\delta)} P^{k} u \|^{2}
	\right),
	\end{multline}
\end{linenomath}    
where $C_{0}$ is independent of $\alpha$, $\beta$, $\gamma$, $\delta$, $k$ and $N$.\\
The first term on the right hand side can be handled
as the second term on the right hand side of (\ref{Est-Fl2}),
see (\ref{Est-Fl2-1}); so
\begin{linenomath}
	\begin{multline*}
	\sum_{j =1}^{n}
	\left|\langle  \widetilde{P}_{N,j} E_{\ell}  \widetilde{\psi} \widetilde{\Theta}_{q}  \left(\Lambda_{r}^{p-1}\right)^{(j)}  \psi_{N}^{(\beta)}  w,
	E_{\ell} \widetilde{\psi} \widetilde{\Theta}_{q} \widetilde{w}\rangle\right|
	\leq
	\varepsilon \sum_{j=1 }^{n}\| P_{j}\widetilde{\psi}  \widetilde{\Theta}_{q}\widetilde{w}\|_{-1}^{2}
	+
	C_{\varepsilon}\| \widetilde{w}\|^{2}.
	\end{multline*}
\end{linenomath}    
where $\varepsilon$ is a suitable small parameter.\\ 
Now we handle the second term on the right hand side of (\ref{Est-Fl3}).
We begin to observe that
$ [ E_{\ell}  \widetilde{\psi} \widetilde{\Theta}_{q},   \widetilde{P}_{N,j}] = 
[ E_{\ell}  \widetilde{\psi},   \widetilde{P}_{N,j}] \widetilde{\Theta}_{q}
+ E_{\ell}  \widetilde{\psi} [\widetilde{\Theta}_{q},  \widetilde{P}_{N,j}]$. 
As previously seen, $[ E_{\ell}  \widetilde{\psi},   \widetilde{P}_{N,j}]$
are first order pseudodifferential operators, moreover
\begin{linenomath}
	\begin{multline*}
	[\widetilde{\Theta}_{q},  \widetilde{P}_{N,j}]
	= \longsum[14]_{j_{1},j_{2},j_{3}=1}^{n}\widetilde{a}_{N,j_{1},j_{2}}^{(j+j_{3})}\widetilde{\Theta}_{q}^{(j_{3})}D_{j_{1}}D_{j_{2}}
	+\longsum[10]_{j_{1},j_{2}=1}^{n} \mathscr{R}_{2}\left([ \widetilde{a}_{N,j_{1},j_{2}}^{(j)}D_{j_{1}}D_{j_{2}}, \widetilde{\Theta}_{q}]\right).
	\end{multline*}
\end{linenomath}
We point out that
$\mathscr{R}_{2}\left([ \widetilde{a}_{N,j_{1},j_{2}}^{(j)}D_{j_{1}}D_{j_{2}}, \widetilde{\Theta}_{q}]\right) $
are zero order operators.
We conclude that $ [ E_{\ell}  \widetilde{\psi} \widetilde{\Theta}_{q},   \widetilde{P}_{N,j}]$ 
are pseudodifferential operators of order $1$, 
moreover
$\| [ E_{\ell}  \widetilde{\psi} \widetilde{\Theta}_{q}, \widetilde{P}_{N,j}]\|_{L^{2}\rightarrow H^{-1}} \leq C$,
where $C$ is independent of $N$. We stress that 
$\psi^{(\mu)}\widetilde{\psi}_{N}^{(\nu)}=0$
for every $\mu,\, \nu \in \mathbb{N}^{n}$ with $|\mu|\geq 1$ and that
$M$ in the construction of $\widetilde{\psi}_{N}$ is taken grater than $n+3$. 
On the other side $\left(\Lambda_{r}^{p-1}\right)^{(j)} $ 
have order $ \frac{p-1}{r} -1$.
So we get
\begin{linenomath}
	\begin{multline*}
	\sum_{j =1}^{n}
	\left|\langle [ E_{\ell}  \widetilde{\psi} \widetilde{\Theta}_{q},  \widetilde{P}_{N,j}] \left(\Lambda_{r}^{p-1}\right)^{(j)}  \psi_{N}^{(\beta)}  w,
	E_{\ell} \widetilde{\psi} \widetilde{\Theta}_{q} \widetilde{w}\rangle\right|
	\\
	\leq
	\sum_{j =1}^{n}
	\| [ E_{\ell}  \widetilde{\psi} \widetilde{\Theta}_{q}, \widetilde{P}_{N,j}] \left(\Lambda_{r}^{p-1}\right)^{(j)}  \psi_{N}^{(\beta)}  w\|
	\| E_{\ell} \widetilde{\psi} \widetilde{\Theta}_{q} \widetilde{w}\|
	\leq
	\tilde{C}_{0}\| \psi_{N}^{(\beta)} w\|^{2}_{\frac{p-1}{r}},
	\end{multline*}
\end{linenomath}    
where the positive constant $\tilde{C}_{0}$ does not depend on $q$ and $N$.\\
Summing up we obtain
\begin{linenomath}
	\begin{align}\label{Est-Fl3-f}
	F_{\ell,3} 
	\leq
	\varepsilon \sum_{j=1 }^{n}\| P_{j}\widetilde{\psi}  \widetilde{\Theta}_{q}\widetilde{w}\|_{-1}^{2}
	+
	C_{\varepsilon} \| \psi_{N}^{(\beta)} \Theta_{N}^{(\gamma)} D^{\alpha}  \varphi_{N}^{(\delta)} P^{k} u \|^{2}_{\frac{p-1}{r}}
	.
	\end{align}
\end{linenomath}    
The first term on the right hand side can be absorbed by the left hand side
of (\ref{Est-I2-2}), once $\varepsilon$ and so $C_{\varepsilon}$ are suitably fixed.

\textbf{Term $F_{\ell,4}$} on the right hand side of (\ref{Est-El-p/r}). 
We have
\begin{linenomath}
	\begin{multline}\label{Est-Fl4}
	F_{\ell,4} =
	\left|\langle E_{\ell}  \widetilde{\psi} \widetilde{\Theta}_{q}   \Lambda_{r}^{p-1} [\widetilde{P}_{N},\psi_{N}^{(\beta)}]  w,
	E_{\ell} \widetilde{\psi} \widetilde{\Theta}_{q} \widetilde{w}\rangle\right| 
	\\
	\leq
	\sum_{j =1}^{n} \left|\langle E_{\ell}  \widetilde{\psi} \widetilde{\Theta}_{q}  \Lambda_{r}^{p-1} \widetilde{P}^{j}_{N} \psi_{N}^{(\beta+j)}  w,
	E_{\ell} \widetilde{\psi} \widetilde{\Theta}_{q} \widetilde{w}\rangle\right|
	\\
	+
	\longsum[11]_{j_{1},j_{2}=1 }^{n} 
	\left|\langle E_{\ell}  \widetilde{\psi} \widetilde{\Theta}_{q}   \Lambda_{r}^{p-1} \widetilde{a}_{N,j_{1},j_{2}}\psi_{N}^{(\beta+j_{1}+j_{2})}  w,
	E_{\ell} \widetilde{\psi} \widetilde{\Theta}_{q} \widetilde{w}\rangle\right|
	\\
	+
	\sum_{ j_{1}=1}^{n}
	\left|\langle E_{\ell}  \widetilde{\psi} \widetilde{\Theta}_{q}   \Lambda_{r}^{p-1} \widetilde{b}_{N,j_{1}} \psi_{N}^{(\beta+j_{1})}  w,
	E_{\ell} \widetilde{\psi} \widetilde{\Theta}_{q} \widetilde{w}\rangle\right|
	=F_{\ell,4,1} + F_{\ell,4,2} + F_{\ell,4,3}
	.
	\end{multline}
\end{linenomath}    
We have
\begin{linenomath}
	\begin{align*}
	F_{\ell,4,2}
	\leq
	\sum_{|\mu| = 2 } C_{0}
	\|  \psi_{N}^{(\beta+\mu)}  w\|_{\frac{p-1}{r}}
	\| \psi_{N}^{(\beta)} \Theta_{N}^{(\gamma)} D^{\alpha}  \varphi_{N}^{(\delta)} P^{k} u \|_{\frac{p-1}{r}},
	\end{align*}
\end{linenomath}    
and
\begin{linenomath}
	\begin{align*}
	F_{\ell,4,3}
	\leq
	\sum_{j_{1}= 1 }^{n} C_{1}
	\|  \psi_{N}^{(\beta+j_{1})}  w\|_{\frac{p-1}{r}}
	\| \psi_{N}^{(\beta)} \Theta_{N}^{(\gamma)} D^{\alpha}  \varphi_{N}^{(\delta)} P^{k} u \|_{\frac{p-1}{r}}
	.
	\end{align*}
\end{linenomath}    
Concerning the term $F_{\ell,4,1}$, we have
\begin{linenomath}
	\begin{multline*}
	F_{\ell,4,1}
	\leq
	\sum_{j =1}^{n} \left|\langle [E_{\ell}  \widetilde{\psi} \widetilde{\Theta}_{q}   \Lambda_{r}^{p-1}, \widetilde{P}^{j}_{N}]  \psi_{N}^{(\beta+j)}  w,
	E_{\ell} \widetilde{\psi} \widetilde{\Theta}_{q} \widetilde{w}\rangle\right|
	\\
	+
	\sum_{j =1}^{n} \left|\langle \widetilde{P}^{j}_{N} E_{\ell}  \widetilde{\psi} \widetilde{\Theta}_{q}   \Lambda_{r}^{p-1}  \psi_{N}^{(\beta+j)}  w,
	E_{\ell} \widetilde{\psi} \widetilde{\Theta}_{q} \widetilde{w}\rangle\right|
	.
	\end{multline*}
\end{linenomath}    
We observe that 
\begin{linenomath}
	\begin{multline*}
	[E_{\ell}  \widetilde{\psi} \widetilde{\Theta}_{q}   \Lambda_{r}^{p-1}, \widetilde{P}^{j}_{N}] =
	[E_{\ell}, \widetilde{P}^{j}_{N}] \widetilde{\psi} \widetilde{\Theta}_{q}\Lambda_{r}^{p-1} 
	+E_{\ell}  [\widetilde{\psi}, \widetilde{P}^{j}_{N}]  \widetilde{\Theta}_{q}   \Lambda_{r}^{p-1}
	\\
	+
	E_{\ell}  \widetilde{\psi} [\widetilde{\Theta}_{q}, \widetilde{P}^{j}_{N}]  \Lambda_{r}^{p-1}
	+
	E_{\ell}  \widetilde{\psi} \widetilde{\Theta}_{q}  [ \Lambda_{r}^{p-1}, \widetilde{P}^{j}_{N}] ,
	\end{multline*}
\end{linenomath}    
where, more explicitly, 
\begin{linenomath}
	\begin{align*}
	& [\widetilde{\Theta}_{q}, \widetilde{P}^{j}_{N}] =
	2\longsum[6]_{ j_{1},j_{2}=1}^{n} \widetilde{a}_{N,j_{1},j}^{(j_{2})}\widetilde{\Theta}_{q}^{(j_{2})}D_{j_{1}}
	+\sum_{ j_{1}=1}^{n} \mathscr{R}_{2,1}\left( [ \widetilde{\Theta}_{q},\widetilde{a}_{N,j_{1},j} D_{j_{1}}, ]\right);
	\\
	& [ \Lambda_{r}^{p-1}, \widetilde{P}^{j}_{N}] =
	2\longsum[6]_{ j_{1},j_{2}=1}^{n} \widetilde{a}_{N,j_{1},j}^{(j_{2})} (\Lambda_{r}^{p-1})^{(j_{2})} D_{j_{1}}
	+ \sum_{ j_{1}=1}^{n} \mathscr{R}_{2,2}\left( [ \widetilde{a}_{N,j_{1},j} D_{j_{1}}, \Lambda_{r}^{p-1}]\right).
	\end{align*}
\end{linenomath}    
We point out that $ [\widetilde{\Theta}_{q}, \widetilde{P}^{j}_{N}] $ are zero order operators
and $[ \Lambda_{r}^{p-1}, \widetilde{P}^{j}_{N}]  $ are operators of order $\frac{p-1}{r}$.
By \eqref{CommEP_jN} and the above consideration, we conclude that
$[E_{\ell}  \widetilde{\psi} \widetilde{\Theta}_{q}   \Lambda_{r}^{p-1}, \widetilde{P}^{j}_{N}] \in L\left( L^{2}, H^{-\frac{p-1}{r}}\right)$
and $\| [E_{\ell}  \widetilde{\psi} \widetilde{\Theta}_{q}   \Lambda_{r}^{p-1}, \widetilde{P}^{j}_{N}]\|_{L^{2}\rightarrow H^{-\frac{p-1}{r}}} \leq C$,
where $C$ depends on $p$ and $r$ but is independent of $N$.
We can handle the terms in the sums using the same strategy adopted
to handle the terms in the sum on the right hand side of (\ref{Est-Fl2});
so, we obtain
\begin{linenomath}
	\begin{align}\label{Est-Fl41}
	F_{\ell,4,1}
	\leq
	\varepsilon \sum_{j =1}^{n} \| P^{j}   \widetilde{\psi} \widetilde{\Theta}_{q}   \widetilde{w}\|^{2}
	+
	C_{\varepsilon}\left( \sum_{j =1}^{n} \| \psi_{N}^{(\beta+j)}  w\|^{2}_{\frac{p-1}{r}}
	+\|v \|^{2}_{\frac{p-1}{r}}\right).
	\end{align}
\end{linenomath}    
where $\varepsilon$ is a suitable small constant.\\
Summing up we have
\begin{linenomath}
	\begin{multline}\label{Est-Fl4-f}
	F_{\ell,4} 
	\leq
	\varepsilon \sum_{j =1}^{n} \| P^{j}   \widetilde{\psi} \widetilde{\Theta}_{q}   \widetilde{w}\|^{2}
	+
	C_{0}\sum_{|\mu| = 2 }
	\|  \psi_{N}^{(\beta+\mu)}  \Theta_{N}^{(\gamma)} D^{\alpha}  \varphi_{N}^{(\delta)} P^{k} u\|_{\frac{p-1}{r}}	
	\|v\|_{\frac{p-1}{r}}	
	\\
	+
	C_{\varepsilon}\left( \sum_{j =1}^{n} \| \psi_{N}^{(\beta+j)}  \Theta_{N}^{(\gamma)} D^{\alpha}  \varphi_{N}^{(\delta)} P^{k} u\|^{2}_{\frac{p-1}{r}}
	+
	\|v\|^{2}_{\frac{p-1}{r}}	\right)
	.
	\end{multline}
\end{linenomath}    
The first term on the right hand side can be absorbed,
taking $\varepsilon$ small enough, by the left hand side
of (\ref{Est-I2-2}).

\textbf{Term $F_{\ell,5}$} on the right hand side of (\ref{Est-El-p/r}).
We recall that
\begin{linenomath}	
	\begin{multline*}
	[\widetilde{P}_{N}, \Theta_{N}^{(\gamma)} D^{\alpha}] = 
	\sum_{j=1}^{n} \left( \widetilde{P}_{N,j} + i \sum_{j_{1}=0}^{n}  \widetilde{b}_{N,j_{1}}^{(j)}(x) D_{\ell} + \widetilde{c}^{(j)}_{N}(x)\right) 
	\left(\Theta_{N}^{(\gamma)} D^{\alpha} \right)^{(j)}
	\\
	\qquad
	+ \longsum[49]_{2 \leq |\mu| \leq |\alpha| -|\gamma| + \lfloor \frac{n}{2} \rfloor+2} 
	\frac{1}{\mu!}  \widetilde{P}_{N,\mu} \left(\Theta_{N}^{(\gamma)} D^{\alpha} \right)^{(\mu)}
	+ \mathscr{R}_{|\alpha| -|\gamma| + \lfloor \frac{n}{2} \rfloor+3}\left( [ \widetilde{P}_{N}, \Theta_{N}^{(\gamma)} D^{\alpha}]  \right),
	\end{multline*}
\end{linenomath}
where
$\widetilde{P}_{N,\mu}=\sum_{j_{1},j_{2}=1}^{n} \widetilde{a}_{N,j_{1},j_{2}}^{(\mu)}(x) D_{j_{1}} D_{j_{2}}
+ i \sum_{j_{1}=0}^{n} \widetilde{b}_{j_{1}}^{(\mu)}(x) D_{j_{1}} + \widetilde{c}^{(\mu)}_{N}(x)$,
see (\ref{Rem_PThD}) for the explicit form of $\mathscr{R}_{|\alpha| -|\gamma| + \lfloor \frac{n}{2} \rfloor+3}$.
We have
\begin{linenomath}
	\begin{multline}\label{Est-Fl5}
	F_{\ell,5} =
	\left|\langle E_{\ell}  \widetilde{\psi} \widetilde{\Theta}_{q}   \Lambda_{r}^{p-1} \psi_{N}^{(\beta)} 
	[\widetilde{P}_{N} ,\Theta_{N}^{(\gamma)} D^{\alpha}] \varphi_{N}^{(\delta)} P^{k} u,
	E_{\ell} \widetilde{\psi} \widetilde{\Theta}_{q} \widetilde{w}\rangle\right|
	\\
	\leq
	\sum_{j=1}^{n} 
	\left|\langle E_{\ell}  \widetilde{\psi} \widetilde{\Theta}_{q}   \Lambda_{r}^{p-1} \psi_{N}^{(\beta)} \widetilde{P}_{N,j}  \left(\Theta_{N}^{(\gamma)} D^{\alpha} \right)^{(j)} \varphi_{N}^{(\delta)} P^{k} u,
	E_{\ell} \widetilde{\psi} \widetilde{\Theta}_{q} \widetilde{w}\rangle\right|
	\\
	\hspace{-2em}
	+
	\longsum[11]_{j,j_{1}=1}^{n} 
	\left|\langle E_{\ell}  \widetilde{\psi} \widetilde{\Theta}_{q}   \Lambda_{r}^{p-1} \psi_{N}^{(\beta)} \widetilde{b}_{N,j_{1}}^{(j)} D_{j_{1}}
	\left(\Theta_{N}^{(\gamma)} D^{\alpha} \right)^{(j)}
	\varphi_{N}^{(\delta)} P^{k} u,
	E_{\ell} \widetilde{\psi} \widetilde{\Theta}_{q} \widetilde{w}\rangle\right|
	\\
	\qquad\qquad\qquad
	+
	\sum_{j=1}^{n} 
	\left|\langle E_{\ell}  \widetilde{\psi} \widetilde{\Theta}_{q}   \Lambda_{r}^{p-1} \psi_{N}^{(\beta)} \widetilde{c}_{N}^{(j)} \left(\Theta_{N}^{(\gamma)} D^{\alpha} \right)^{(j)} \varphi_{N}^{(\delta)} P^{k} u,
	E_{\ell} \widetilde{\psi} \widetilde{\Theta}_{q} \widetilde{w}\rangle\right|
	\\
	+
	\longsum[44]_{2 \leq |\mu| \leq |\alpha| -|\gamma| + \lfloor \frac{n}{2} \rfloor+2} 
	\frac{1}{\mu!}  
	\left|\langle E_{\ell}  \widetilde{\psi} \widetilde{\Theta}_{q}   \Lambda_{r}^{p-1} \psi_{N}^{(\beta)}\widetilde{P}_{N,\mu} \left(\Theta_{N}^{(\gamma)} D^{\alpha} \right)^{(\mu)} \varphi_{N}^{(\delta)} P^{k} u,
	E_{\ell} \widetilde{\psi} \widetilde{\Theta}_{q} \widetilde{w}\rangle\right|
	\\
	\qquad\qquad
	+
	\left|\langle E_{\ell}  \widetilde{\psi} \widetilde{\Theta}_{q}   \Lambda_{r}^{p-1} \psi_{N}^{(\beta)}
	\mathscr{R}_{|\alpha| -|\gamma| + \lfloor \frac{n}{2} \rfloor+ 3}\left( [ \widetilde{P}_{N} , \Theta_{N}^{(\gamma)} D^{\alpha}]  \right) \varphi_{N}^{(\delta)} P^{k} u,
	E_{\ell} \widetilde{\psi} \widetilde{\Theta}_{q} \widetilde{w}\rangle\right|
	\\
	=
	\sum_{r=1}^{5} F_{\ell,5,r} 
	.
	\end{multline}
\end{linenomath}  
We begin to focus on $F_{\ell,5,5}$.
In order to make the following more readable we will
set $h +1 = |\alpha| -|\gamma| + \lfloor \frac{n}{2} \rfloor+ 3$ and
we write $ \mathscr{R}_{h+1}$ instead of 
$ \mathscr{R}_{|\alpha| -|\gamma| + \lfloor \frac{n}{2} \rfloor+ 3}\left( [ \widetilde{P}_{N}, \Theta_{N}^{(\gamma)} D^{\alpha}]\right)$.
We have
\begin{linenomath}
	\begin{multline*}
	F_{\ell,5,5}
	\leq
	C_{0} \| \psi_{N}^{(\beta)} \mathscr{R}_{h+1}v\|_{\frac{p-1}{r}}\,\,
	\| E_{\ell} \widetilde{\psi} \widetilde{\Theta}_{q} \widetilde{w}\|
	\\
	\leq
	C_{0}C_{1}
	\int (1+|\xi|)^{\frac{p-1}{r}} |\widehat{\psi}_{N}^{(\beta)} (\xi) |\, d\xi\, \,\,
	\| \mathscr{R}_{h+1} g\|_{\frac{p-1}{r}}\,\,
	\| v \|_{\frac{p-1}{r}}
	\\
	\leq
	\tilde{C}_{0} 
	C_{\psi}^{|\beta|+n+4} N^{ \left( |\beta|+ n + 3 -M\right)^{+}} 
	\|\mathscr{R}_{h+1} g\|_{\frac{p-1}{r}}\,\,
	\| v \|_{\frac{p-1}{r}}.
	\end{multline*}
\end{linenomath}
We estimate the second to last factor; we have 
\begin{linenomath}
	\begin{multline*}
	\| \mathscr{R}_{h +1}  g \|_{\frac{p-1}{r}}  
	\leq
	\left(\int \left|(1+|\xi|^{2})^{\frac{p-1}{2r}} \widehat{\mathscr{R}_{h+1} g}(\xi) \right|^{2} d\xi \right) ^{\frac{1}{2}}
	\\
	=
	\left(\int \left|\int (1+|\xi|^{2})^{\frac{p-1}{2r}}  \text{\textcursive{r}}_{h+1} (\xi-\eta,\eta)\widehat{ g}(\eta) \, d\eta\right|^{2} d\xi \right) ^{\frac{1}{2}}
	\\
	\leq
	\left(\iint \left|(1+|\xi|^{2})^{\frac{p-1}{2r}}  \text{\textcursive{r}}_{h+1} (\xi-\eta,\eta)\right|^{2}\, d\eta\, d\xi \right)^{\frac{1}{2}}
	\| \widehat{ g}\|_{L^{2}_{\eta}},
	\end{multline*}
\end{linenomath}
where $\|\widehat{g}\|_{L^{2}_{\eta}} = \| \varphi_{N}^{(\delta)} P^{k} u \|_{0} $ and
\begin{linenomath}
	\begin{multline*}
	\text{\textcursive{r}}_{h+1} (\xi-\eta,\eta)
	\\
	=
	\longsum[9]_{|\mu|=h+1}
	\frac{ h+1}{i^{|\mu|}\mu!} \!
	\left( \sum_{j_{1},j=1}^{n} \widehat{\widetilde{a}}_{N,j_{1},j}^{(\mu )}(\xi-\eta)\eta_{j_{1}} \eta_{j} 
	+\sum_{j_{1}=1}^{n} \widehat{\widetilde{b}}_{N,j_{1}}^{(\mu)}(\xi-\eta)\eta_{j_{1}} + \widehat{\widetilde{c}}_{N}^{(\mu)}(\xi-\eta)
	\right)
	\\
	\qquad \qquad\qquad\qquad\qquad\times
	\int_{0}^{1}\left(1-t\right)^{h}
	\left(\sigma\left(\Theta_{N}^{(\gamma)}D^{\alpha}\right)\right)^{(\mu)} (\eta+t(\xi-\eta))\,dt .
	\end{multline*}
\end{linenomath}
For every $t$ in $\left[0,1\right] $, we have
\begin{linenomath}
	\begin{multline*}
	\left|\left(\sigma\left(\Theta_{N}^{(\gamma)}D^{\alpha}\right)\right)^{(\mu)} (\eta+t(\xi-\eta))\right|
	\\
	\leq
	\tilde{C}^{ |\alpha| + \lfloor \frac{n}{2} \rfloor+ 4} \,\mu! \,
	N^{(|\alpha| +\lfloor \frac{n}{2} \rfloor + 3-M)^{+}} 
	\frac{\,\,\,\,\left( 1+ |\xi-\eta|^{2}\right)^{\frac{1}{2}\left(\lfloor \frac{n}{2}\rfloor +3 \right)}}{\left( 1+ |\eta|\right)^{\lfloor \frac{n}{2} \rfloor+3}}.
	\end{multline*}
\end{linenomath}
Moreover, since $|\mu|=  |\alpha| -|\gamma| + \lfloor \frac{n}{2} \rfloor+ 3$,
by the Lemma \ref{L-1} and the Remark \ref{Rk-1}, for all $j_{1}, j$, we have
\begin{linenomath}
	\begin{align*}
	\left( 1+ |\xi-\eta|^{2}\right)^{\lfloor \frac{n}{2} \rfloor+3} |\widehat{\widetilde{a}}_{N,j_{1},j}^{(\mu )}(\xi-\eta)|
	\leq
	\tilde{C}_{1}^{ |\alpha| -|\gamma|+ 3 \left(\lfloor \frac{n}{2} \rfloor+ 3\right) +1}\,N^{s(|\alpha|-|\gamma| +3\left(\lfloor \frac{n}{2} \rfloor + 3\right)-M)^{+}}.
	\end{align*}
\end{linenomath}
The same estimates hold for $\widehat{\widetilde{b}}_{N,j_{1}}^{(\mu)}(\xi-\eta)$,
$j_{1}=1,\dots,n$, and  $\widehat{\widetilde{c}}_{N}^{(\mu)}(\xi-\eta)$.\\
We conclude that
\begin{linenomath}
	\begin{align*}
	\| \mathscr{R}_{h +1}  g \|_{\frac{p-1}{r}}  
	\leq
	\tilde{C}_{3}^{ 2|\alpha| -|\gamma|+ 1}\,N^{s(2|\alpha|-|\gamma| +2(n + 6) -M)^{+}} \| \varphi_{N}^{(\delta)} P^{k} u \|_{0} .
	\end{align*}
\end{linenomath}
By \eqref{Est_g} and taking advantage from the Remark \ref{Rk-1},
we get
\begin{linenomath}
	\begin{align}\label{Est-Fl55}
	F_{\ell,5,5}
	\leq
	\tilde{C}_{4}^{ \sigma+ 1}
	\tilde{C}_{5}^{ 2m+\sigma+ 1}\,N^{s(2m +|\gamma| + \sigma +3(n + 5) -M)^{+}} 
	\| v \|_{\frac{p-1}{r}}.
	\end{align}
\end{linenomath}
where $m= |\alpha| - |\gamma| $ and  $\sigma= |\beta| + |\delta| +2k$.

Term $F_{\ell,5,1}$. We have
\begin{linenomath}
	\begin{multline}\label{Est_Fl51-0}
	F_{\ell,5,1}
	\leq	
	\sum_{j=1}^{n} 
	\left|\langle [E_{\ell}  \widetilde{\psi} \widetilde{\Theta}_{q}   \Lambda_{r}^{p-1} \psi_{N}^{(\beta)} , \widetilde{P}_{N,j}] \left(\Theta_{N}^{(\gamma)} D^{\alpha} \right)^{(j)} \varphi_{N}^{(\delta)} P^{k} u,
	E_{\ell} \widetilde{\psi} \widetilde{\Theta}_{q} \widetilde{w}\rangle\right|
	\\
	+
	\sum_{j=1}^{n} 
	\left|\langle \widetilde{P}_{N,j}E_{\ell}  \widetilde{\psi} \widetilde{\Theta}_{q}   \Lambda_{r}^{p-1} \psi_{N}^{(\beta)} \left(\Theta_{N}^{(\gamma)} D^{\alpha} \right)^{(j)} \varphi_{N}^{(\delta)} P^{k} u,
	E_{\ell} \widetilde{\psi} \widetilde{\Theta}_{q} \widetilde{w}\rangle\right|.
	\end{multline}
\end{linenomath}
The last term on the right hand side can be handled
as the second term on the right hand side of (\ref{Est-Fl2}),
(see (\ref{Est-Fl2-1})), we get
\begin{linenomath}
	\begin{multline*}
	\sum_{j=1}^{n} 
	\left|\langle \widetilde{P}_{N,j} E_{\ell}  \widetilde{\psi} \widetilde{\Theta}_{q}   \Lambda_{r}^{p-1} \psi_{N}^{(\beta)} \left(\Theta_{N}^{(\gamma)} D^{\alpha} \right)^{(j)} \varphi_{N}^{(\delta)} P^{k} u,
	E_{\ell} \widetilde{\psi} \widetilde{\Theta}_{q} \widetilde{w}\rangle\right|
	\\
	\leq
	C_{\varepsilon} \sum_{j=1}^{n}
	\|  \psi_{N}^{(\beta)} \left(\Theta_{N}^{(\gamma)} D^{\alpha} \right)^{(j)} \varphi_{N}^{(\delta)} P^{k} u \|_{\frac{p-1}{r} +1}^{2}
	+
	\varepsilon 
	\sum_{j=1}^{n}
	\| P_{j}  \widetilde{\psi} \widetilde{\Theta}_{q} \widetilde{w}\|_{-1}^{2},
	\end{multline*}
\end{linenomath}
where $\varepsilon$ is a suitable small constant.
The second term on the right hand side can be absorbed by the left hand side
of (\ref{Est-I2-2}).\\
Concerning the terms in the first sum, we have 
\begin{linenomath}
	\begin{multline*}
	\|  \psi_{N}^{(\beta)} \left(\Theta_{N}^{(\gamma)} D^{\alpha} \right)^{(j)} \varphi_{N}^{(\delta)} P^{k} u \|_{\frac{p-1}{r} +1}^{2}
	\leq
	\sum_{j_{1}=1}^{n} 
	\| \psi_{N}^{(\beta +j_{1})} \left(\Theta_{N}^{(\gamma)} D^{\alpha} \right)^{(j)} \varphi_{N}^{(\delta)} P^{k} u \|_{\frac{p-1}{r}}^{2}
	\\
	+
	\sum_{j_{1}=1}^{n} 
	\| \psi_{N}^{(\beta)} D_{j_{1}} \left(\Theta_{N}^{(\gamma)} D^{\alpha} \right)^{(j)} \varphi_{N}^{(\delta)}P^{k} u \|_{\frac{p-1}{r}}^{2}
	\\
	\leq
	\sum_{j_{1}=1}^{n} 
	\left(
	\| \psi_{N}^{(\beta +j_{1})} \Theta_{N}^{(\gamma+j)} D^{\alpha} \varphi_{N}^{(\delta)} P^{k} u \|_{\frac{p-1}{r}}^{2}
	+
	\alpha_{j}^{2}\| \psi_{N}^{(\beta+j_{1})} \Theta_{N}^{(\gamma)} D^{\alpha-j} \varphi_{N}^{(\delta)} P^{k} u \|_{\frac{p-1}{r}}^{2}
	\right)
	\\
	+
	\sum_{j_{1}=1}^{n} 
	\left(
	\| \psi_{N}^{(\beta)} \Theta_{N}^{(\gamma+j)} D^{\alpha+j_{1}} \varphi_{N}^{(\delta)} P^{k} u \|_{\frac{p-1}{r}}^{2}
	+
	\alpha_{j}^{2}\| \psi_{N}^{(\beta)} \Theta_{N}^{(\gamma)} D^{\alpha-j+j_{1}} \varphi_{N}^{(\delta)} P^{k} u \|_{\frac{p-1}{r}}^{2}
	\right),
	\end{multline*}
\end{linenomath}
where we use that
$ \left(\Theta_{N}^{(\gamma)} D^{\alpha} \right)^{(j)} =
\Theta_{N}^{(\gamma+j)} D^{\alpha} + \alpha_{j} \Theta_{N}^{(\gamma)} D^{\alpha-j} $.\\
Now, we handle the first term on the right hand side of (\ref{Est_Fl51-0}).
We have
\begin{linenomath}
	\begin{multline*}
	\left|\langle [ E_{\ell}  \widetilde{\psi} \widetilde{\Theta}_{q}   \Lambda_{r}^{p-1} \psi_{N}^{(\beta)} , \widetilde{P}_{N,j}] \left(\Theta_{N}^{(\gamma)} D^{\alpha} \right)^{(j)} \varphi_{N}^{(\delta)} P^{k} u,
	E_{\ell} \widetilde{\psi} \widetilde{\Theta}_{q} \widetilde{w}\rangle\right|
	\\
	\leq
	\left|\langle E_{\ell}  \widetilde{\psi} \widetilde{\Theta}_{q}   \Lambda_{r}^{p-1} [\psi_{N}^{(\beta)} , \widetilde{P}_{N,j}] \left(\Theta_{N}^{(\gamma)} D^{\alpha} \right)^{(j)} \varphi_{N}^{(\delta)} P^{k} u,
	E_{\ell} \widetilde{\psi} \widetilde{\Theta}_{q} \widetilde{w}\rangle\right|
	\\
	+
	\left|\langle [E_{\ell}  \widetilde{\psi} \widetilde{\Theta}_{q}   \Lambda_{r}^{p-1} , \widetilde{P}_{N,j}] \psi_{N}^{(\beta)} \left(\Theta_{N}^{(\gamma)} D^{\alpha} \right)^{(j)} \varphi_{N}^{(\delta)} P^{k} u,
	E_{\ell} \widetilde{\psi} \widetilde{\Theta}_{q} \widetilde{w}\rangle\right|
	\\
	\leq
	C \left[
	\longsum[9]_{j_{1},j_{2}=1}^{n} \left( \| \psi_{N}^{(\beta +j_{1})} \Theta_{N}^{(\gamma+j)} D^{\alpha+j_{2}} \varphi_{N}^{(\delta)} P^{k} u \|_{\frac{p-1}{r}}
	+\alpha_{j} \| \psi_{N}^{(\beta +j_{1})} \Theta_{N}^{(\gamma)} D^{\alpha+j_{2}-j} \varphi_{N}^{(\delta)} P^{k} u \|_{\frac{p-1}{r}}
	\right)
	\right.
	\\
	\left.
	+\sum_{|\mu|=2} 
	\left( \| \psi_{N}^{(\beta +\mu)} \Theta_{N}^{(\gamma+j)} D^{\alpha} \varphi_{N}^{(\delta)} P^{k} u \|_{\frac{p-1}{r}}
	+\alpha_{j} \| \psi_{N}^{(\beta +\mu)} \Theta_{N}^{(\gamma)} D^{\alpha-j} \varphi_{N}^{(\delta)} P^{k} u \|_{\frac{p-1}{r}}
	\right)
	\right] \|v\|_{\frac{p-1}{r}}
	\\
	+
	\left|\langle [E_{\ell}  \widetilde{\psi} \widetilde{\Theta}_{q}   \Lambda_{r}^{p-1} , \widetilde{P}_{N,j}] \psi_{N}^{(\beta)} \left(\Theta_{N}^{(\gamma)} D^{\alpha} \right)^{(j)} \varphi_{N}^{(\delta)} P^{k} u,
	E_{\ell} \widetilde{\psi} \widetilde{\Theta}_{q} \widetilde{w}\rangle\right|.
	\end{multline*}
\end{linenomath}
Concerning the last term on the right hand side,
we observe that
\begin{linenomath}
$$[E_{\ell}  \widetilde{\psi} \widetilde{\Theta}_{q}   \Lambda_{r}^{p-1} , \widetilde{P}_{N,j}] =
E_{\ell}  \widetilde{\psi} \widetilde{\Theta}_{q}  [ \Lambda_{r}^{p-1} , \widetilde{P}_{N,j}]
+ [E_{\ell}  \widetilde{\psi} \widetilde{\Theta}_{q}, \widetilde{P}_{N,j}]  \Lambda_{r}^{p-1}.$$
\end{linenomath}
As seen previously, in the estimate of $F_{\ell,3}$ (\eqref{Est-Fl3}),
$ [E_{\ell}  \widetilde{\psi} \widetilde{\Theta}_{q}, \widetilde{P}_{N,j}]$
is a pseudodifferential operator of order one.
Moreover
\begin{linenomath}
	\begin{multline*}
	[\Lambda_{r}^{p-1},  \widetilde{P}_{N,j}]
	\\
	=
	\longsum[14]_{j_{1},j_{2},j_{3}=1}^{n}\widetilde{a}_{N,j_{1},j_{2}}^{(j+j_{3})} (\Lambda_{r}^{p-1})^{(j_{3})}D_{j_{1}}D_{j_{2}}
	+\longsum[10]_{j_{1},j_{2}=1}^{n}\sum_{|\mu|=2} \widetilde{a}_{N,j_{1},j_{2}}^{(\mu+j)} (\Lambda_{r}^{p-1})^{(\mu)}D_{j_{1}}D_{j_{2}}
	\\
	+\longsum[10]_{j_{1},j_{2}=1}^{n} \mathscr{R}_{3}\left([ \Lambda_{r}^{p-1}, \widetilde{a}_{N,j_{1},j_{2}}^{(j)}D_{j_{1}}D_{j_{2}}]\right),
	\end{multline*}
\end{linenomath}
where the terms in the first two sums are pseudodifferential operators of order
$1+\frac{p-1}{r}$ and $\frac{p-1}{r}$ respectively.
$\mathscr{R}_{3}$ are zero order operators; moreover, as operators in $L( H^{(p-1)/r} , L^{2})$,
they have norm uniformly bounded by $C$ , independent of $N$.\\
So $[E_{\ell}  \widetilde{\psi} \widetilde{\Theta}_{q}   \Lambda_{r}^{p-1} , \widetilde{P}_{N,j}]$ 
is a pseudodifferential operator of order $1+\frac{p-1}{r}$. We have
\begin{linenomath}
	\begin{multline*}
	\| \langle [E_{\ell}  \widetilde{\psi} \widetilde{\Theta}_{q}   \Lambda_{r}^{p-1} , \widetilde{P}_{N,j}] \psi_{N}^{(\beta)} \left(\Theta_{N}^{(\gamma)} D^{\alpha} \right)^{(j)} \varphi_{N}^{(\delta)} P^{k} u \|
	\\
	\leq
	C_{0}
	\longsum[9]_{ j_{1}=1}^{n} \| D_{j_{1}} \psi_{N}^{(\beta)} \left(\Theta_{N}^{(\gamma)} D^{\alpha} \right)^{(j)} g\|_{\frac{p-1}{r}}
	\\
	\leq
	\sum_{ j_{1}=1}^{n}
	\left( \| \psi_{N}^{(\beta+j_{1})} \Theta_{N}^{(\gamma +j)} D^{\alpha} g\|_{\frac{p-1}{r}}
	+ \alpha_{j} \|  \psi_{N}^{(\beta)} \Theta_{N}^{(\gamma)} D^{\alpha - j+j_{1}} g\|_{\frac{p-1}{r}}
	\right).
	\end{multline*}
\end{linenomath}
So we get, summarizing all the above estimates,
\begin{linenomath}
	\begin{multline}\label{Est_Fl51}
	F_{\ell,5,1}
	\leq
	\varepsilon
	\sum_{j=1}^{n}
	\| P_{j}  \widetilde{\psi} \widetilde{\Theta}_{q} \widetilde{w}\|_{-1}^{2}
	+
	C_{\varepsilon} \left( 
	\longsum[9]_{j_{1},j=1}^{n} 
	\| \psi_{N}^{(\beta +j_{1})} \Theta_{N}^{(\gamma+j)} D^{\alpha} g \|_{\frac{p-1}{r}}^{2}
	\right.
	\\
	+
	\| \psi_{N}^{(\beta)} \Theta_{N}^{(\gamma+j)} D^{\alpha+j_{1}} g \|_{\frac{p-1}{r}}^{2}
	+
	\alpha_{j}^{2}\| \psi_{N}^{(\beta+j_{1})} \Theta_{N}^{(\gamma)} D^{\alpha-j} g \|_{\frac{p-1}{r}}^{2}
	+
	\alpha_{j}^{2}\| \psi_{N}^{(\beta)} \Theta_{N}^{(\gamma)} D^{\alpha-j+j_{1}} g \|_{\frac{p-1}{r}}^{2}
	\Biggr)
	\\
	+
	C_{1}
	\left[
	\longsum[9]_{j,j_{1}=1}^{n}
	\left( \| \psi_{N}^{(\beta+j_{1})} \Theta_{N}^{(\gamma +j)} D^{\alpha} g\|_{\frac{p-1}{r}}
	\right.
	\right.
	\left.
	\left.
	+ \alpha_{j} \|  \psi_{N}^{(\beta)} \Theta_{N}^{(\gamma)} D^{\alpha - j+j_{1}} g\|_{\frac{p-1}{r}}
	\right)
	\right.
	\\
	+
	\longsum[11]_{j,j_{1},j_{2}=1}^{n} 
	\left( \| \psi_{N}^{(\beta +j_{1})} \Theta_{N}^{(\gamma+j)} D^{\alpha+j_{2}} g\|_{\frac{p-1}{r}}
	\right.
	\left.
	\left.
	+\alpha_{j} \| \psi_{N}^{(\beta +j_{1})} \Theta_{N}^{(\gamma)} D^{\alpha+j_{2}-j} g \|_{\frac{p-1}{r}}
	\right)
	\right.
	\\
	+
	\longsum[9]_{j, |\mu|=2} 
	\left( \| \psi_{N}^{(\beta +\mu)} \Theta_{N}^{(\gamma+j)} D^{\alpha} g \|_{\frac{p-1}{r}}
	\right.
	\left.
	+\alpha_{j} \| \psi_{N}^{(\beta +\mu)} \Theta_{N}^{(\gamma)} D^{\alpha-j} g \|_{\frac{p-1}{r}}
	\right)
	\Biggr] \|v\|_{\frac{p-1}{r}},
	\end{multline}
\end{linenomath}
we recall that $g=\varphi_{N}^{(\delta)} P^{k} u$.

\textbf{Term $F_{\ell,5,2}$}. We have
\begin{linenomath}
	\begin{multline}\label{Est-Fl52}
	F_{\ell,5,2}
	=
	\longsum[11]_{j,j_{1}=1}^{n} 
	\left|\langle E_{\ell}  \widetilde{\psi} \widetilde{\Theta}_{q}   \Lambda_{r}^{p-1} \psi_{N}^{(\beta)} \widetilde{b}_{N,j_{1}}^{(j)} D_{j_{1}}
	\left(\Theta_{N}^{(\gamma)} D^{\alpha} \right)^{(j)}\!\!
	\varphi_{N}^{(\delta)} P^{k} u,
	E_{\ell} \widetilde{\psi} \widetilde{\Theta}_{q} \widetilde{w}\rangle\right|
	\\
	\leq
	C_{0}
	\longsum[11]_{j,j_{1}=1}^{n} 
	\| \widetilde{b}_{N,j_{1}}^{(j)} \psi_{N}^{(\beta)}  D_{j_{1}} \left(\Theta_{N}^{(\gamma)} D^{\alpha} \right)^{(j)}g\|_{\frac{p-1}{r}}
	\|v\|_{\frac{p-1}{r}}
	\\
	\leq
	C_{0}
	\longsum[11]_{j,j_{1}=1}^{n} 
	\int (1+|\xi|)^{\frac{p-1}{r}} |\widehat{\widetilde{b}}_{N,j_{1}}^{(j)}(\xi) |\, d\xi\, \,\,
	\| \psi_{N}^{(\beta)}  D_{j_{1}} \left(\Theta_{N}^{(\gamma)} D^{\alpha} \right)^{(j)}g\|_{\frac{p-1}{r}}
	\|v\|_{\frac{p-1}{r}}
	\\
	\leq
	C_{0} C_{1}
	\longsum[11]_{j,j_{1}=1}^{n} 
	\left( 
	\| \psi_{N}^{(\beta)} \Theta_{N}^{(\gamma +j)} D^{\alpha+j_{1}} g\|_{\frac{p-1}{r}}
	+\alpha_{j} \| \psi_{N}^{(\beta)}  \Theta_{N}^{(\gamma)} D^{\alpha -j +j_{1})}g\|_{\frac{p-1}{r}}
	\right)
	\|v\|_{\frac{p-1}{r}},
	\end{multline}
\end{linenomath}
where $C_{1} =\displaystyle\sup_{j_{1},j} \left( \int (1+|\xi|)^{\frac{p-1}{r}} |\widehat{\widetilde{b}}_{N,j_{1}}^{(j)}(\xi) |\, d\xi \right)$.

\textbf{Term $F_{\ell,5,3}$}. It can be handled as the term $F_{\ell,5,2}$.

\textbf{Term $F_{\ell,5,4}$}. Recalling that
$\widetilde{P}_{N,\mu}=\sum_{j_{1},j_{2}=1}^{n} \widetilde{a}_{N,j_{1},j_{2}}^{(\mu)}(x) D_{j_{1}} D_{j_{2}}
+ i \sum_{j_{1}=0}^{n} \widetilde{b}_{j_{1}}^{(\mu)}(x) D_{j_{1}} + \widetilde{c}^{(\mu)}_{N}(x)$.
We have
\begin{linenomath}
	\begin{multline}\label{Est-Fl54}
	F_{\ell,5,4}
	\\
	=
	\longsum[46]_{2 \leq |\mu| \leq |\alpha| -|\gamma| + \lfloor \frac{n}{2} \rfloor+3} 
	\frac{1}{\mu!}  
	\left|\langle E_{\ell}  \widetilde{\psi} \widetilde{\Theta}_{q}  
	\Lambda_{r}^{p-1} \psi_{N}^{(\beta)}\widetilde{P}_{N,\mu} \left(\Theta_{N}^{(\gamma)} D^{\alpha} \right)^{(\mu)} \varphi_{N}^{(\delta)} P^{k} u,
	E_{\ell} \widetilde{\psi} \widetilde{\Theta}_{q} \widetilde{w}\rangle\right|
	\\
	\leq
	\longsum[46]_{2 \leq |\mu| \leq |\alpha| -|\gamma| + \lfloor \frac{n}{2} \rfloor+3} 
	\longsum[7]_{\substack{\nu\leq \mu \\ \nu\leq \alpha}}
	\frac{1}{\mu!}  
	\binom{\mu}{\nu} \frac{\alpha!}{(\alpha -\nu)!}
	\| \psi_{N}^{(\beta)}\widetilde{P}_{N,\mu} \Theta_{N}^{(\gamma+\mu-\nu)} D^{\alpha-\nu} g\|_{\frac{p-1}{r}}
	\| v \|_{\frac{p-1}{r}}
	\\
	\hspace*{-14em}
	\leq
	\longsum[46]_{2 \leq |\mu| \leq |\alpha| -|\gamma| + \lfloor \frac{n}{2} \rfloor+3} 
	\longsum[7]_{\substack{\nu\leq \mu \\ \nu\leq \alpha}}
	\frac{1}{\mu!}  
	\binom{\mu}{\nu} \frac{\alpha!}{(\alpha -\nu)!}
	\\
	\times
	\left[ \longsum[11]_{j_{1},j_{2}=1}^{n} 
	\left( \int (1+|\xi|)^{\frac{p-1}{r}} |\widehat{\widetilde{a}}_{N,j_{1},j_{2}}^{(\mu)}  (\xi) |\, d\xi \right)
	\| \psi_{N}^{(\beta)} \Theta_{N}^{(\gamma+\mu-\nu)} D^{\alpha-\nu+j_{1}+j_{2}} g\|_{\frac{p-1}{r}}
	\right.
	\\
	+
	\sum_{j_{1}=1}^{n} 
	\left( \int (1+|\xi|)^{\frac{p-1}{r}} |\widehat{\widetilde{b}}_{N,j_{1}}^{(\mu)}  (\xi) |\, d\xi \right)
	\| \psi_{N}^{(\beta)} \Theta_{N}^{(\gamma+\mu-\nu)} D^{\alpha-\nu+j_{1}} g\|_{\frac{p-1}{r}}
	\\
	\left.
	+
	\left( \int (1+|\xi|)^{\frac{p-1}{r}} |\widehat{\widetilde{c}}_{N}^{(\mu)}  (\xi) |\, d\xi \right)
	\| \psi_{N}^{(\beta)} \Theta_{N}^{(\gamma+\mu-\nu)} D^{\alpha-\nu} g\|_{\frac{p-1}{r}}
	\right]
	\| v \|_{\frac{p-1}{r}}
	\\
	\hspace*{-2em}
	\leq
	\longsum[46]_{2 \leq |\mu| \leq |\alpha| -|\gamma| + \lfloor \frac{n}{2} \rfloor+3} 
	\longsum[7]_{\substack{\nu\leq \mu \\ \nu\leq \alpha}}
	\frac{1}{\mu!}  
	\binom{\mu}{\nu} \frac{\alpha!}{(\alpha -\nu)!}
	C_{1}^{|\mu|+1} N^{ s\left( |\mu|+ n + 3 -M\right)^{+}} 
	\\
	\times
	\left[ \longsum[11]_{j_{1},j_{2}=1}^{n} 
	\| \psi_{N}^{(\beta)} \Theta_{N}^{(\gamma+\mu-\nu)} D^{\alpha-\nu+j_{1}+j_{2}} g\|_{\frac{p-1}{r}}
	+
	\| \psi_{N}^{(\beta)} \Theta_{N}^{(\gamma+\mu-\nu)} D^{\alpha-\nu} g\|_{\frac{p-1}{r}}
	\right.
	\\
	\left.
	+
	\sum_{j_{1}=1}^{n} 
	\| \psi_{N}^{(\beta)} \Theta_{N}^{(\gamma+\mu-\nu)} D^{\alpha-\nu+j_{1}} g\|_{\frac{p-1}{r}}
	\right]
	\| v \|_{\frac{p-1}{r}}
	.
	\end{multline}
\end{linenomath}

\textbf{Term $F_{\ell,6}$} on the right hand side of (\ref{Est-El-p/r}). We have
\begin{linenomath}
	\begin{multline*}
	F_{\ell,6} =
	\left|\langle E_{\ell}  \widetilde{\psi} \widetilde{\Theta}_{q}   \Lambda_{r}^{p-1} \psi_{N}^{(\beta)} \Theta_{N}^{(\gamma)} D^{\alpha} 
	[\widetilde{P}_{N},\varphi_{N}^{(\delta)}] P^{k} u, E_{\ell} \widetilde{\psi} \widetilde{\Theta}_{q} \widetilde{w}\rangle\right|
	\\
	\leq
	\sum_{j =1}^{n} \left|\langle E_{\ell}  \widetilde{\psi} \widetilde{\Theta}_{q}   \Lambda_{r}^{p-1} \psi_{N}^{(\beta)} \Theta_{N}^{(\gamma)} D^{\alpha}
	\widetilde{P}_{N}^{j} \varphi_{N}^{(\delta+j)}  P^{k} u,
	E_{\ell} \widetilde{\psi} \widetilde{\Theta}_{q} \widetilde{w}\rangle\right|
	\\
	\qquad\qquad
	+
	\longsum[11]_{j_{1},j_{2}=1 }^{n} 
	\left|\langle E_{\ell}  \widetilde{\psi} \widetilde{\Theta}_{q}   \Lambda_{r}^{p-1} \psi_{N}^{(\beta)} \Theta_{N}^{(\gamma)} D^{\alpha} 
	\widetilde{a}_{N,j_{1},j_{2}}\varphi_{N}^{(\delta+j_{1}+j_{2})} P^{k} u,
	E_{\ell} \widetilde{\psi} \widetilde{\Theta}_{q} \widetilde{w}\rangle\right|
	\\
	+
	\sum_{ j_{1}=1}^{n}
	\left|\langle E_{\ell}  \widetilde{\psi} \widetilde{\Theta}_{q}   \Lambda_{r}^{p-1} \psi_{N}^{(\beta)} \Theta_{N}^{(\gamma)} D^{\alpha} 
	\widetilde{b}_{N,j_{1}} \varphi_{N}^{(\beta+j_{1})}  P^{k} u,
	E_{\ell} \widetilde{\psi} \widetilde{\Theta}_{q} \widetilde{w}\rangle\right|
	=
	F_{\ell,6,1}+F_{\ell,6,2}+F_{\ell,6,3}.
	\end{multline*}
\end{linenomath}    
The terms $ F_{\ell,6,2}$ and $F_{\ell,6,3}$ can be handled in the same way.
We begin to estimate the term $F_{\ell,6,2} $.
Since
\begin{linenomath}	
	\begin{multline*}
	[\Theta_{N}^{(\gamma)} D^{\alpha}, \widetilde{a}_{N,j_{1},j_{2}}] = 
	\\
	\qquad
	\longsum[44]_{1 \leq |\mu| \leq |\alpha| -|\gamma| + \lfloor \frac{n}{2} \rfloor +1} 
	\longsum[7]_{\substack{\nu\leq \mu \\ \nu\leq \alpha}}
	\frac{1}{\mu!} 
	\binom{\mu}{\nu} \frac{\alpha!}{(\alpha -\nu)!} \widetilde{a}_{N,j_{1},j_{2}}^{(\mu)} \Theta_{N}^{(\gamma +\mu-\nu)} D^{\alpha-\nu}
	\\
	+ \mathscr{R}_{|\alpha| -|\gamma| + \lfloor \frac{n}{2} \rfloor+1}\left( [ \widetilde{a}_{N,j_{1},j_{2}}, \Theta_{N}^{(\gamma)} D^{\alpha}]  \right),
	\end{multline*}
\end{linenomath}
we have
\begin{linenomath}
	\begin{multline*}
	F_{\ell,6,2} 
	\leq
	C_{0}
	\longsum[11]_{j_{1},j_{2}=1 }^{n} \longsum[44]_{1 \leq |\mu| \leq |\alpha| -|\gamma| + \lfloor \frac{n}{2} \rfloor +1} 
	\longsum[7]_{\substack{\nu\leq \mu \\ \nu\leq \alpha}}
	\frac{1}{\mu!} 
	\binom{\mu}{\nu} \frac{\alpha!}{(\alpha -\nu)!}
	\\
	\times
	\left( \int (1+|\xi|)^{\frac{p-1}{r}} |\widehat{\widetilde{a}}_{N,j_{1},j_{2}}^{(\mu)}  (\xi) |\, d\xi \right)
	\| \psi_{N}^{(\beta)}  \Theta_{N}^{(\gamma +\mu-\nu)} D^{\alpha-\nu} \varphi_{N}^{(\delta+j_{1}+j_{2})} f\|_{\frac{p-1}{r}}
	\| v\|_{\frac{p-1}{r}}
	\\
	\hspace{-20em}
	+
	\left( \int (1+|\xi|)^{\frac{p-1}{r}} |\psi_{N}^{(\beta)}  (\xi) |\, d\xi \right) \times
	\\
	\times
	\| \mathscr{R}_{|\alpha| -|\gamma| + \lfloor \frac{n}{2} \rfloor+2}
	\left( [ \widetilde{a}_{N,j_{1},j_{2}}, \Theta_{N}^{(\gamma)} D^{\alpha}]  \right) \varphi_{N}^{(\delta+j_{1}+j_{2})} f\|_{\frac{p-1}{r}}
	\| v\|_{\frac{p-1}{r}}
	.
	\end{multline*}
\end{linenomath}    
where $f=P^{k}u$. Using the same strategy
used to handle the term $F_{\ell,5,5}$
(see (\ref{Est-Fl5}) and (\ref{Est-Fl55})),
the product of the first two factors of the last term on the right hand side
can be estimated by
\begin{linenomath}
	\begin{align*}
	\tilde{C}_{1}^{\sigma + 1}
	\tilde{C}_{2}^{ 2m+ |\gamma|+ 1}\,N^{s(2m +|\gamma| +\sigma+3(n + 4) -M)^{+}}.
	\end{align*}
\end{linenomath}
where $m= |\alpha| - |\gamma| $ and  $\sigma= |\beta| + |\delta| +2k$.
We obtain
\begin{linenomath}
	\begin{multline*}
	F_{\ell,6,2} 
	\leq
	C_{0}
	\longsum[11]_{j_{1},j_{2}=1 }^{n} \longsum[44]_{1 \leq |\mu| \leq |\alpha| -|\gamma| + \lfloor \frac{n}{2} \rfloor +1} 
	\longsum[7]_{\substack{\nu\leq \mu \\ \nu\leq \alpha}}
	\frac{1}{\mu!} 
	\binom{\mu}{\nu} \frac{\alpha!}{(\alpha -\nu)!}
	C_{3}^{|\mu|+1} \,N^{s(|\mu|+1 -M)^{+}}
	\\
	\times
	\| \psi_{N}^{(\beta)}  \Theta_{N}^{(\gamma +\mu-\nu)} D^{\alpha-\nu} \varphi_{N}^{(\delta+j_{1}+j_{2})} f\|_{\frac{p-1}{r}}
	\| v\|_{\frac{p-1}{r}}
	\\
	+
	\tilde{C}_{1}^{|\gamma| + 1}
	\tilde{C}_{2}^{ 2m+\sigma+ 1}\,N^{s(2m +|\gamma| +\sigma+3(n + 4) -M)^{+}} 
	\| v\|_{\frac{p-1}{r}}
	.
	\end{multline*}
\end{linenomath}    
Concerning the term $F_{\ell,6,3}$,
it can be handled as done above obtaining
\begin{linenomath}
	\begin{multline*}
	F_{\ell,6,3} 
	\leq
	C_{0}
	\longsum[6]_{j_{1}=1 }^{n} \longsum[44]_{1 \leq |\mu| \leq |\alpha| -|\gamma| + \lfloor \frac{n}{2} \rfloor +1} 
	\longsum[7]_{\substack{\nu\leq \mu \\ \nu\leq \alpha}}
	\frac{1}{\mu!} 
	\binom{\mu}{\nu} \frac{\alpha!}{(\alpha -\nu)!}
	C_{4}^{|\mu|+1} \,N^{s(|\mu|+1 -M)^{+}}
	\\
	\times
	\| \psi_{N}^{(\beta)}  \Theta_{N}^{(\gamma +\mu-\nu)} D^{\alpha-\nu} \varphi_{N}^{(\delta+j_{1})} f\|_{\frac{p-1}{r}}
	\| v\|_{\frac{p-1}{r}}
	\\
	+
	\tilde{C}_{1}^{|\gamma| + 1}
	\tilde{C}_{2}^{ 2m+\sigma+ 1}\,N^{s(2m +|\gamma| +\sigma+3(n + 4) -M)^{+}} 
	\| v\|_{\frac{p-1}{r}}
	.
	\end{multline*}
\end{linenomath}    
Terms $ F_{\ell,6,1}$. We have
\begin{linenomath}
	\begin{multline}\label{Est-Fl61}
	F_{\ell,6,1} 
	\leq
	\sum_{j =1}^{n} \left|\langle E_{\ell}  \widetilde{\psi} \widetilde{\Theta}_{q}   \Lambda_{r}^{p-1} \widetilde{P}_{N}^{j} \psi_{N}^{(\beta)}   \Theta_{N}^{(\gamma)} D^{\alpha}
	\varphi_{N}^{(\delta+j)}  P^{k} u,
	E_{\ell} \widetilde{\psi} \widetilde{\Theta}_{q} \widetilde{w}\rangle\right|
	\\
	+
	\sum_{j =1}^{n} \left|\langle E_{\ell}  \widetilde{\psi} \widetilde{\Theta}_{q}   \Lambda_{r}^{p-1} [\psi_{N}^{(\beta)} , \widetilde{P}_{N}^{j}] \Theta_{N}^{(\gamma)} D^{\alpha}
	\varphi_{N}^{(\delta+j)}  P^{k} u,
	E_{\ell} \widetilde{\psi} \widetilde{\Theta}_{q} \widetilde{w}\rangle\right|
	\\
	+
	\sum_{j =1}^{n} \left|\langle E_{\ell}  \widetilde{\psi} \widetilde{\Theta}_{q}   \Lambda_{r}^{p-1} \psi_{N}^{(\beta)}  [\widetilde{P}_{N}^{j}, \Theta_{N}^{(\gamma)} D^{\alpha}]
	\varphi_{N}^{(\delta+j)}  P^{k} u,
	E_{\ell} \widetilde{\psi} \widetilde{\Theta}_{q} \widetilde{w}\rangle\right|
	.
	\end{multline}
\end{linenomath}    
The first term on the right hand side can be handled as the term
$F_{\ell,4,1}$ (see (\ref{Est-Fl4}) and (\ref{Est-Fl41})).
It can be estimated by
\begin{linenomath}
	\begin{align*}
	\varepsilon \sum_{j =1}^{n} \| P^{j}   \widetilde{\psi} \widetilde{\Theta}_{q}   \widetilde{w}\|^{2}
	+
	C_{\varepsilon}\left( \sum_{j =1}^{n} \|\psi_{N}^{(\beta)}   \Theta_{N}^{(\gamma)} D^{\alpha}
	\varphi_{N}^{(\delta+j)}  P^{k} u\|^{2}_{\frac{p-1}{r}}
	+\|v \|^{2}_{\frac{p-1}{r}}\right).
	\end{align*}
\end{linenomath}    
Since $[\psi_{N}^{(\beta)} , \widetilde{P}_{N}^{j}] =
\sum_{ j_{1}=1}^{n} \widetilde{a}_{N,j_{1},j} \psi_{N}^{(\beta+j_{1})}$,
the second term on the right side of (\ref{Est-Fl61}) is bounded by
\begin{linenomath}
	\begin{align*}
	C_{1}\longsum[11]_{j_{1},j=1 }^{n}
	\| \psi_{N}^{(\beta+j_{1})} \Theta_{N}^{(\gamma)} D^{\alpha} \varphi_{N}^{(\delta+j)}  P^{k} u\|_{\frac{p-1}{r}}
	\|v \|_{\frac{p-1}{r}}.
	\end{align*}
\end{linenomath}    
Concerning the last term on the right hand side of (\ref{Est-Fl61}), since
\begin{linenomath}	
	\begin{multline*}
	[\widetilde{P}_{N}^{j}, \Theta_{N}^{(\gamma)} D^{\alpha}]=
	\\
	\sum_{ j_{1}=1}^{n}
	\longsum[44]_{1 \leq |\mu| \leq |\alpha| -|\gamma| + \lfloor \frac{n}{2} \rfloor +2} 
	\longsum[7]_{\substack{\nu\leq \mu \\ \nu\leq \alpha}}
	\frac{1}{\mu!} 
	\binom{\mu}{\nu} \frac{\alpha!}{(\alpha -\nu)!} \widetilde{a}_{N,j_{1},j_{2}}^{(\mu)} \Theta_{N}^{(\gamma +\mu-\nu)} D^{\alpha-\nu+j_{1}}
	\\
	+ \mathscr{R}_{|\alpha| -|\gamma| + \lfloor \frac{n}{2} \rfloor+2 }\left( [\widetilde{P}_{N}^{j}, \Theta_{N}^{(\gamma)} D^{\alpha}] \right),
	\end{multline*}
\end{linenomath}
we can estimate it by
\begin{linenomath}
	\begin{multline*}
	C_{0}
	\longsum[11]_{j,j_{1}=1 }^{n} \longsum[44]_{1 \leq |\mu| \leq |\alpha| -|\gamma| + \lfloor \frac{n}{2} \rfloor +2} 
	\longsum[7]_{\substack{\nu\leq \mu \\ \nu\leq \alpha}}
	\frac{1}{\mu!} 
	\binom{\mu}{\nu} \frac{\alpha!}{(\alpha -\nu)!}
	C_{3}^{|\mu|+1} \,N^{s(|\mu|+1 -M)^{+}}
	\\
	\times
	\| \psi_{N}^{(\beta)}  \Theta_{N}^{(\gamma +\mu-\nu)} D^{\alpha-\nu+j_{1}} \varphi_{N}^{(\delta+j)} f\|_{\frac{p-1}{r}}
	\| v\|_{\frac{p-1}{r}}
	\\
	+
	C_{2}^{|\gamma| + 1}
	C_{3}^{ 2m+\sigma+ 1}\,N^{s(2m +|\gamma| +\sigma+3(n + 5) -M)^{+}} 
	\| v\|_{\frac{p-1}{r}}
	.
	\end{multline*}
\end{linenomath}    
Summing up we have
\begin{linenomath}
	\begin{multline}\label{Est-Fl6}
	F_{\ell,6} 
	\leq
	\varepsilon \sum_{j =1}^{n} \| P^{j}   \widetilde{\psi} \widetilde{\Theta}_{q}   \widetilde{w}\|^{2}
	+
	C_{\varepsilon}\left(\sum_{j =1}^{n} \|\psi_{N}^{(\beta)}   \Theta_{N}^{(\gamma)} D^{\alpha}
	\varphi_{N}^{(\delta+j)}  P^{k} u\|^{2}_{\frac{p-1}{r}}
	\right.
	\\
	\left.
	+\|v \|^{2}_{\frac{p-1}{r}}
	\right)
	+
	C_{4}
	\Biggr[
	\longsum[11]_{j_{1},j=1 }^{n}
	\| \psi_{N}^{(\beta+j_{1})} \Theta_{N}^{(\gamma)} D^{\alpha} \varphi_{N}^{(\delta+j)}  P^{k} u\|_{\frac{p-1}{r}}
	\\
	+\longsum[11]_{j,j_{1}=1 }^{n} \longsum[44]_{1 \leq |\mu| \leq |\alpha| -|\gamma| + \lfloor \frac{n}{2} \rfloor +2} 
	\longsum[7]_{\substack{\nu\leq \mu \\ \nu\leq \alpha}}
	\frac{1}{\mu!} 
	\binom{\mu}{\nu} \frac{\alpha!}{(\alpha -\nu)!}
	C_{4}^{|\mu|} \,N^{s(|\mu|+1 -M)^{+}}
	\\
	\times
	\left(
	\| \psi_{N}^{(\beta)}  \Theta_{N}^{(\gamma +\mu-\nu)} D^{\alpha-\nu} \varphi_{N}^{(\delta+j+j_{1})} f\|_{\frac{p-1}{r}}
	+
	\| \psi_{N}^{(\beta)}  \Theta_{N}^{(\gamma +\mu-\nu)} D^{\alpha-\nu+j_{1}} \varphi_{N}^{(\delta+j)} f\|_{\frac{p-1}{r}}
	\right)
	\\
	+
	C_{4}^{|\gamma|}
	C_{5}^{ 2m+\sigma+ 1}\,N^{s(2m +|\gamma| +\sigma+3(n + 5) -M)^{+}} 
	\Biggr]
	\|v \|_{\frac{p-1}{r}}
	.
	\end{multline}
\end{linenomath}    

\textbf{Term $F_{\ell,7}$} on the right hand side of (\ref{Est-El-p/r}). 
We have
\begin{linenomath}
	\begin{multline}\label{Est-Fl7}
	F_{\ell,7} =
	\left|\langle E_{\ell}  \widetilde{\psi} \widetilde{\Theta}_{q}   \Lambda_{r}^{p-1} \psi_{N}^{(\beta)} \Theta_{N}^{(\gamma)} D^{\alpha}  \varphi_{N}^{(\delta)} P^{k+1} u,
	E_{\ell} \widetilde{\psi} \widetilde{\Theta}_{q} \widetilde{w}\rangle\right|
	\\
	\leq
	\| \psi_{N}^{(\beta)} \Theta_{N}^{(\gamma)} D^{\alpha}  \varphi_{N}^{(\delta)} P^{k+1} u\|_{\frac{p-1}{r}}
	\| \psi_{N}^{(\beta)} \Theta_{N}^{(\gamma)} D^{\alpha}  \varphi_{N}^{(\delta)} P^{k} u\|_{\frac{p-1}{r}}
	.
	\end{multline}
\end{linenomath}    

By (\ref{Est-I2-2}), (\ref{Est-El-p/r}), (\ref{Est-Fl1}), (\ref{Est-Fl2-f}), (\ref{Est-Fl3-f}), (\ref{Est-Fl4-f}),
(\ref{Est-Fl5}), (\ref{Est_Fl51}), (\ref{Est-Fl52}), (\ref{Est-Fl54}), (\ref{Est-Fl6}) and (\ref{Est-Fl7})
there are suitable positive constants independent of $\alpha$, $\beta$, $\gamma$, $\delta$ and $k$
such that  
\begin{linenomath}
	\begin{multline}\label{Est-1/r-tild}
	\|  \dbtilde{\psi} \dbtilde{\Theta} \widetilde{\psi}  \widetilde{\Theta}_{q} \widetilde{w} \|_{\frac{1}{r}}^{2}
	\leq 
	C \Biggl\{
	\|v\|_{\frac{p-1}{r}}^{2} 
	+ \sum_{j =1}^{n} \| \psi_{N}^{(\beta+j)}  w\|^{2}_{\frac{p-1}{r}}
	+\|v\|_{0}^{2} 
	\\
	+
	\longsum[9]_{j_{1},j=1}^{n} \left(
	\| \psi_{N}^{(\beta +j_{1})} \Theta_{N}^{(\gamma+j)} D^{\alpha} g \|_{\frac{p-1}{r}}^{2}
	+
	\| \psi_{N}^{(\beta)} \Theta_{N}^{(\gamma+j)} D^{\alpha+j_{1}} g \|_{\frac{p-1}{r}}^{2}
	\right.
	\\
	\left.
	\qquad\quad
	+
	\alpha_{j}^{2}\| \psi_{N}^{(\beta+j_{1})} \Theta_{N}^{(\gamma)} D^{\alpha-j} g \|_{\frac{p-1}{r}}^{2}
	+
	\alpha_{j}^{2}\| \psi_{N}^{(\beta)} \Theta_{N}^{(\gamma)} D^{\alpha-j+j_{1}} g \|_{\frac{p-1}{r}}^{2}
	\right)
	\\
	+\sum_{j =1}^{n} \|\psi_{N}^{(\beta)}   \Theta_{N}^{(\gamma)} D^{\alpha}
	\varphi_{N}^{(\delta+j)}  P^{k} u\|^{2}_{\frac{p-1}{r}}
	+
	\Biggl[\sum_{|\mu| = 2 } \|  \psi_{N}^{(\beta+\mu)}  w\|_{\frac{p-1}{r}}
	\\
	+
	\longsum[9]_{j,j_{1}=1}^{n}
	\left( 
	\| \psi_{N}^{(\beta+j_{1})} \Theta_{N}^{(\gamma +j)} D^{\alpha} g\|_{\frac{p-1}{r}}
	+ \alpha_{j} \|  \psi_{N}^{(\beta)} \Theta_{N}^{(\gamma)} D^{\alpha - j+j_{1}} g\|_{\frac{p-1}{r}}
	\right)
	\\
	+
	\longsum[11]_{j,j_{1},j_{2}=1}^{n} 
	\left( 
	\| \psi_{N}^{(\beta +j_{1})} \Theta_{N}^{(\gamma+j)} D^{\alpha+j_{2}} g\|_{\frac{p-1}{r}}
	+
	\alpha_{j} \| \psi_{N}^{(\beta +j_{1})} \Theta_{N}^{(\gamma)} D^{\alpha+j_{2}-j} g \|_{\frac{p-1}{r}}
	\right)
	\\
	+
	\sum_{ j=1}^{n}\longsum[9]_{|\mu|=2} 
	\left( \| \psi_{N}^{(\beta +\mu)} \Theta_{N}^{(\gamma+j)} D^{\alpha} g \|_{\frac{p-1}{r}}
	+
	\alpha_{j} \| \psi_{N}^{(\beta +\mu)} \Theta_{N}^{(\gamma)} D^{\alpha-j} g \|_{\frac{p-1}{r}}
	\right)
	\\
	+
	\longsum[11]_{j,j_{1}=1}^{n} 
	\left( 
	\| \psi_{N}^{(\beta)} \Theta_{N}^{(\gamma +j)} D^{\alpha+j_{1}} g\|_{\frac{p-1}{r}}
	+
	\| \psi_{N}^{(\beta+j_{1})} \Theta_{N}^{(\gamma)} D^{\alpha} \varphi_{N}^{(\delta+j)}  P^{k} u\|_{\frac{p-1}{r}}
	\right)
	\\
	+
	\longsum[44]_{2 \leq |\mu| \leq |\alpha| -|\gamma| + \lfloor \frac{n}{2} \rfloor+3} 
	\longsum[7]_{\substack{\nu\leq \mu \\ \nu\leq \alpha}}
	\frac{1}{\mu!}  
	\binom{\mu}{\nu} \frac{\alpha!}{(\alpha -\nu)!}
	C_{1}^{|\mu|+1} N^{ s\left( |\mu|+ n + 3 -M\right)^{+}} 
	\\
	\times
	\left( \longsum[11]_{j_{1},j_{2}=1}^{n} 
	\| \psi_{N}^{(\beta)} \Theta_{N}^{(\gamma+\mu-\nu)} D^{\alpha-\nu+j_{1}+j_{2}} g\|_{\frac{p-1}{r}}
	+
	\| \psi_{N}^{(\beta)} \Theta_{N}^{(\gamma+\mu-\nu)} D^{\alpha-\nu} g\|_{\frac{p-1}{r}}
	\right.
	\\
	\left.
	+
	\sum_{j_{1}=1}^{n} 
	\| \psi_{N}^{(\beta)} \Theta_{N}^{(\gamma+\mu-\nu)} D^{\alpha-\nu+j_{1}} g\|_{\frac{p-1}{r}}
	\right)
	\\
	+\longsum[11]_{j,j_{1}=1 }^{n} \longsum[44]_{1 \leq |\mu| \leq |\alpha| -|\gamma| + \lfloor \frac{n}{2} \rfloor +2} 
	\longsum[7]_{\substack{\nu\leq \mu \\ \nu\leq \alpha}}
	\frac{1}{\mu!} 
	\binom{\mu}{\nu} \frac{\alpha!}{(\alpha -\nu)!}
	C_{4}^{|\mu|} \,N^{s(|\mu|+1 -M)^{+}}
	\\
	\times
	\left(
	\| \psi_{N}^{(\beta)}  \Theta_{N}^{(\gamma +\mu-\nu)} D^{\alpha-\nu} \varphi_{N}^{(\delta+j+j_{1})} f\|_{\frac{p-1}{r}}
	+
	\| \psi_{N}^{(\beta)}  \Theta_{N}^{(\gamma +\mu-\nu)} D^{\alpha-\nu+j_{1}} \varphi_{N}^{(\delta+j)} f\|_{\frac{p-1}{r}}
	\right)
	\\
	+
	\| \psi_{N}^{(\beta)} \Theta_{N}^{(\gamma)} D^{\alpha}  \varphi_{N}^{(\delta)} P^{k+1} u\|_{\frac{p-1}{r}}
	+
	\tilde{C}_{4}^{ \sigma+ 1}
	\tilde{C}_{5}^{ 2m+|\gamma|+ 1}\,N^{s(2m +|\gamma| +\sigma+3(n + 6) -M)^{+}} 
	\Biggr]
	\|v\|_{\frac{p-1}{r}}
	\Biggr\}, 
	\end{multline} 
\end{linenomath}    
where, we recall, $ g =\varphi_{N}^{(\delta)} P^{k} u$,
$ w =\Theta_{N}^{(\gamma)} D^{\alpha} \varphi_{N}^{(\delta)} P^{k} u$
and $v = \psi_{N}^{(\beta)} \Theta_{N}^{(\gamma)} D^{\alpha} \varphi_{N}^{(\delta)} P^{k} u$.\\
We recall as  done in the case  $p=1$, that the strategy adopted in  \eqref{Est-Fl1},
\eqref{Est-Fl2-f}, \eqref{Est-Fl3-f} and \eqref{Est-Fl4-f}, where we introduce
$\varepsilon$ in order to absorb a term on the left hand side of \eqref{Est-I2-2},
is used a finite number of times, say at most $50$ times; this allow us to choose $\varepsilon$
so that $ 50\,\varepsilon \, C < 1/2$,  where $C$ is the constant on the right hand side of \eqref{Est-I2-2}.

By (\ref{Step-0-p/r}), (\ref{Est-I-p/r}), (\ref{Step--p/r}),
(\ref{Est-I-1-p/r-f}), (\ref{Est-I2-01}) and (\ref{Est-I2-11})
we conclude that
\begin{linenomath}	
	\begin{multline}\label{Est-p/r-f}
	\| v\|_{\frac{p}{r}}^{2}
	\leq 
	C_{2}^{2} \left( 
	\| \dbtilde{\psi} \dbtilde{\Theta} \widetilde{\psi}  \widetilde{\Theta}_{q} \widetilde{w} \|_{\frac{1}{r}}^{2}
	+
	\|v\|_{0}^{2}
	\right.
	\\
	\left.
	+
	C_{1}^{2(\sigma+1)} C_{2}^{2(2 m + |\gamma|)}
	N^{2s(2 m +|\gamma|+ \sigma + 2n + 4 - M)^{+}}
	\right),
	\end{multline}
\end{linenomath}	
where the first term on the right hand side can be estimated as in (\ref{Est-1/r-tild}). 
We remark that we can choose $M$ equal to $3(n+6)$.
%
\section{Microlocal Gevrey regularity of Gevrey vectors of $P$}
\renewcommand{\theequation}{\thesection.\arabic{equation}}
\setcounter{equation}{0} \setcounter{theorem}{0}
\setcounter{proposition}{0} \setcounter{lemma}{0}
\setcounter{corollary}{0} \setcounter{definition}{0}
\setcounter{remark}{1}
%
In this section we prove our main theorem concerning
the microlocal regularity of the Gevrey vectors of
$P(x,D)$, \eqref{H-O-R_Op}.
About that, we begin to prove a couple of results
consequence of the microlocal estimates established
in the previous section.    
\begin{proposition} \label{Pr1-1/r-p/r}
Let $\psi_{N}$, $\Theta_{N}$ and $\varphi_{N}$ be as in the
previous section and $u$ a Gevrey-vector of order $s$ for $P$.
There exist constants $A_{1}$ and $B_{1}$ such that, if:
\begin{linenomath}
\begin{align}\label{Pr1-1_0}
(1)_{0}
		\begin{cases}
		\| \psi_{N}^{(\beta)} \Theta_{N}^{(\gamma)} D^{\alpha} \varphi_{N}^{(\delta)} P^{k} u\| 
		\leq A_{1}^{|\sigma| +1} B_{1}^{2rm +|\gamma|+1} N^{s[rm +|\gamma|+ \sigma]},&\\
		\noalign{\vskip4pt}
		\text{for } 2r|\alpha|-(2r-1)|\gamma|+\sigma \leq N , \text{ where } \sigma = |\beta| + |\delta| +2k, &\\
		\noalign{\vskip4pt}
		m=|\alpha|-|\gamma| \text{ and } |\gamma|\leq |\alpha|. &
		\end{cases}
\end{align}
\end{linenomath}
Then, one has for $1\leq p \leq r$
\begin{linenomath}
\begin{align}\label{Pr1-1_p}
(1)_{p}
		\begin{cases}
		\| \psi_{N}^{(\beta)} \Theta_{N}^{(\gamma)} D^{\alpha} \varphi_{N}^{(\delta)} P^{k} u\|_{\frac{p}{r}} 
		\leq 
		A_{1}^{|\sigma| + p +1} B_{1}^{2rm+|\gamma| + p +1}  N^{s[rm+|\gamma|+ \sigma + p]},&\\
		\noalign{\vskip4pt}
		\text{for } 2r|\alpha|-(2r-1)|\gamma|+\sigma \leq N-2p \text{ and } |\gamma|\leq |\alpha|. &
		\end{cases}
\end{align}
\end{linenomath}
\end{proposition}
\begin{proof}
The result is obtained by induction on $p$.
The main tools are the basic estimates
in $H^{p/r}$, $1\leq p\leq r$ obtained in the previous section,
(\ref{Est-v-1/r}), (\ref{Est-1/r-tild}) and  (\ref{Est-p/r-f}).

\textbf{Step $p=1$.} We want to show that if \eqref{Pr1-1_0} holds,
then we have
\begin{linenomath}
\begin{align}\label{Pr1-1_1}
(1)_{1}
	\begin{cases}
		\| \psi_{N}^{(\beta)} \Theta_{N}^{(\gamma)} D^{\alpha} \varphi_{N}^{(\delta)} P^{k} u\|_{\frac{1}{r}} 
		\leq 
		A_{1}^{|\sigma| + 1 +1} B_{1}^{2rm+|\gamma| + 1 +1}  N^{s[rm+|\gamma|+ \sigma + 1]},&\\
		\noalign{\vskip4pt}
		\text{for } 2r|\alpha|-(2r-1)|\gamma|+\sigma \leq N-2 \text{ and } |\gamma|\leq |\alpha|. &
	\end{cases}
\end{align}
\end{linenomath}
We recall the estimate \eqref{Est-v-1/r}:
\begin{linenomath}	
\begin{multline}\label{Est-1/r-bis}
%
	\| v\|_{\frac{1}{r}} ^{2}
	\leq
	C_{2}\left\{ \| v \|^{2} 
	+
	\sum_{j=1}^{n} \left(\| \psi_{N}^{(\beta+j)} \Theta_{N}^{(\gamma)} D^{\alpha} g\|^{2}
	+
	\| \psi_{N}^{(\beta)} \Theta^{(\gamma)}_{N} D^{\alpha}  \varphi_{N}^{(\delta +j)} f\|^{2}\right)
	\right.
	\\
	\left.
	+
	\sum_{j=1}^{n}
	\left( \| \psi_{N}^{(\beta)} \Theta_{N}^{(\gamma+j) } D^{\alpha} g\|_{1}^{2}
	+\alpha_{j}^{2} \| \psi_{N}^{(\beta)}\Theta_{N}^{(\gamma) } D^{\alpha-j} g \|^{2}_{1}
	\right)
	\right.
	\\
	+
	\left[
	\sum_{j =1}^{n}
	\left( \| \psi_{N}^{(\beta+j)} \Theta_{N}^{(\gamma)} D^{\alpha} g\|
	+\| \psi_{N}^{(\beta)} \Theta_{N}^{(\gamma+j) } D^{\alpha} g\| 
	+\alpha_{j} \| \psi_{N}^{(\beta)}\Theta_{N}^{(\gamma) } D^{\alpha-j} g \| 
	\right)
	\right.
	\\
	+ 
	\sum_{|\mu| =2} \|\psi_{N}^{(\beta+\mu)} \Theta_{N}^{(\gamma)} D^{\alpha} g\| 
	+ \sum_{j,\ell=1}^{n} 
	\left(\| \psi_{N}^{(\beta)}  \Theta_{N}^{(\gamma+j) } D^{\alpha+\ell}  g \|
	+\alpha_{j} \| \psi_{N}^{(\beta)}\Theta_{N}^{(\gamma) } D^{\alpha-j+\ell} g \| 
	\right)
	\\
	+ \sum_{j=1}^{n} \| \psi_{N}^{(\beta)} \Theta^{(\gamma)}_{N} D^{\alpha}  \varphi_{N}^{(\delta +j)} f\|
	+\! \longsum[8]_{j,j_{1},\ell=1}^{n}\!\!\left(\!
	\| \psi_{_{N}}^{(\beta+j_{1})} \Theta_{N}^{(\gamma+j) } D^{\alpha+\ell} g\|
	\!+\!  
	\| \psi_{N}^{(\beta+j_{1} +\ell)} \Theta_{N}^{(\gamma+j) } D^{\alpha} g\|
	\right.
	\\
	\left.
	+ \alpha_{j} 
	\left(\| \psi_{N}^{(\beta+j_{1} +\ell)} \Theta_{N}^{(\gamma) } D^{\alpha-j} g \| +\| \psi_{_{N}}^{(\beta+j_{1})} \Theta_{N}^{(\gamma) } D^{\alpha-j + \ell} g \|
	\right)\right)
	\\
	\left.
	+\sum_{j,\ell=1}^{n} \left(
	\| \psi_{N}^{(\beta+\ell)} \Theta^{(\gamma)}_{N} D^{\alpha}  \varphi_{N}^{(\delta +j)} f\|
	+
	\| \psi_{N}^{(\beta+\ell)}  \Theta_{N}^{(\gamma+j) } D^{\alpha}  g \|
	+\alpha_{j} \| \psi_{N}^{(\beta+\ell)}\Theta_{N}^{(\gamma) } D^{\alpha-j} g \| 
	\right)
	\right.
	\\
	\left.
	+
	\longsum[35]_{2 \leq |\mu| \leq |\alpha| -|\gamma| + 1} \longsum[7]_{\substack{\nu\leq \mu \\ \nu\leq \alpha}}
	\frac{1}{\left(\mu-\nu\right)!}  \binom{\alpha}{\nu}
	\left( \longsum[8]_{j,\ell=1}^{n} 
	\|a_{\ell,j}^{(\mu)} \psi_{N}^{(\beta)}  \Theta_{N}^{(\gamma +\mu-\nu)} D^{\alpha-\nu}D_{j}D_{\ell} g \|
	\right.
	\right.
	\\
	+
	\sum_{ \ell=1}^{n} \|  b_{\ell}^{(\mu)} \psi_{N}^{(\beta)}  \Theta_{N}^{(\gamma +\mu-\nu)} D^{\alpha-\nu}D_{\ell} g \|
	+
	\|  c^{(\mu)} \psi_{N}^{(\beta)}  \Theta_{N}^{(\gamma +\mu-\nu)} D^{\alpha-\nu}g \|
	\Biggr)
	\\
	+\longsum[8]_{j,\ell=1}^{n} \,\longsum[38]_{1 \leq |\mu| \leq |\alpha| -|\gamma|+ 1} \longsum[7]_{\substack{\nu\leq \mu \\ \nu\leq \alpha}}
	\frac{1}{\left(\mu-\nu\right)!}  \binom{\alpha}{\nu}
	\| a_{\ell,j}^{(\mu)}\psi_{N}^{(\beta)} \Theta^{(\gamma+\mu-\nu)}_{N} D^{\alpha-\nu} D_{\ell}  \varphi_{N}^{(\delta +j)} f\|
	\\
	\left.
	+
	\longsum[35]_{0 \leq |\mu| \leq |\alpha| -|\gamma|} \longsum[7]_{\substack{\nu\leq \mu \\ \nu\leq \alpha}}
	\frac{1}{\left(\mu-\nu\right)!}  \binom{\alpha}{\nu} 
	\left( \longsum[8]_{j,\ell=1}^{n}
	\| a_{\ell,j}^{(\mu)}\psi_{N}^{(\beta)} \Theta^{(\gamma+\mu-\nu)}_{N} D^{\alpha-\nu} \varphi_{N}^{(\delta+\ell +j)} f\| 
	\right.
	\right.
	\\
	\left.
	+\sum_{\ell=1}^{n} 
	\| b_{\ell}^{(\mu)}\psi_{N}^{(\beta)} \Theta^{(\gamma+\mu-\nu)}_{N} D^{\alpha-\nu} \varphi_{N}^{(\delta+\ell)} f\| 
	\!\Biggr)
	\right.
	+
	C_{1}^{|\sigma|+1} C_{2}^{2 m +|\gamma| + 1 } 
	N^{s(2m +|\gamma| +\sigma+2n + 5 -M)^{+}} 
	\\
	+
	\| \psi_{N}^{(\beta)} \Theta^{(\gamma)}_{N} D^{\alpha} \varphi_{N}^{(\delta)} P^{k+1} u\|
	\Biggr]
	\| v\|\Biggr\}
	+ 
	C_{1}^{2(\sigma+1)} C_{2}^{2(2 m+|\gamma| + 1)}  
	N^{2s(2 m+ |\gamma|+\sigma) },
	\end{multline}
\end{linenomath} 
where $m= |\alpha| - |\gamma| $,  $\sigma= |\beta| + |\delta| +2k$
and $C_{1}$ and $C_{2}$ are positive constants
independent of $\alpha$, $\beta$, $\gamma$, $\delta$, $k$ and $N$
and
$f = P^{k}u$, $ g = \varphi_{N}^{(\delta)} P^{k} u$,
$v = \psi_{N}^{(\beta)} \Theta_{N}^{(\gamma)} D^{\alpha} \varphi_{N}^{(\delta)} P^{k} u$.\\
We will estimate each of the terms on the right hand side of the above inequality
separately following the order in which they are written.\\ 
In order to make the proof more readable and by simplification of the writing of it,
we introduce the following notation:
let us call the powers of $A_{1}$, $B_{1}$ and $N$ corresponding to
$ \|\psi_{N}^{(\beta)} \Theta_{N}^{(\gamma)} D^{\alpha} \varphi_{N}^{(\delta)} P^{k} u\|_{p/r}$, respectively by:
\begin{linenomath}
	\begin{align}\label{Not-exp}
	\nonumber
	& S_{p}\doteq S_{p}(\alpha,\beta,\gamma,\delta,k)
	= \sigma+p+1;
	\\
	& T_{p}\doteq T_{p}(\alpha,\beta,\gamma,\delta,k)
	= 2rm + |\gamma|+ p+1;
	\\
	\nonumber
	& U_{p}\doteq U_{p}(\alpha,\beta,\gamma,\delta,k)
	= s(rm+|\gamma|+\sigma + p);
	\end{align}
\end{linenomath}
where $m= |\alpha| - |\gamma| $,  $\sigma= |\beta| + |\delta| +2k$. 
We point out that $S_{p+1}= S_{p} +1$, $T_{p+1}= T_{p} +1$ and $U_{p+1} =U_{p} + s$.\\
Using this notation we rewrite \eqref{Pr1-1_1}: 
\begin{linenomath}
	\begin{align}\label{Pr1-1_1-n}
	(1)_{1}
	\begin{cases}
	\| \psi_{N}^{(\beta)} \Theta_{N}^{(\gamma)} D^{\alpha} \varphi_{N}^{(\delta)} P^{k} u\|_{\frac{1}{r}} 
	\leq 
	A_{1}^{S_{1}} B_{1}^{T_{1}}  N^{U_{1}},&\\
	\noalign{\vskip4pt}
	\text{for } T_{1} +S_{1} \leq N+2 \text{ and } |\gamma|\leq |\alpha|. &
	\end{cases}
	\end{align}
\end{linenomath}
The purpose is to show that all these terms on the right hand side
of \eqref{Est-1/r-bis} are smaller than $A_{1}^{2S_{1}} B_{1}^{2T_{1}}  N^{2U_{1}}$
times a factor depending on negative power of $A_{1}$ or $B_{1}$ or $N$.
A suitable choice of $A_{1}$ and $B_{1}$ will yield the summand of these factors less than one.
We remark that all the $4n+1-$tuples of the form $(\alpha', \beta', \gamma', \delta', k') \in \mathbb{N}^{4n+1}$
associated to each term on the right hand side of \eqref{Est-1/r-bis} satisfy the condition in \eqref{Pr1-1_0}.
Such condition, with the notation above introduced, can be rewritten as
\begin{linenomath}
	\begin{align}\label{Pr1-1_0-n}
	(1)_{0}
	\begin{cases}
	\| \psi_{N}^{(\beta)} \Theta_{N}^{(\gamma)} D^{\alpha} \varphi_{N}^{(\delta)} P^{k} u\| 
	\leq A_{1}^{S_{0}} B_{1}^{T_{0}} N^{U_{0}},&\\
	\noalign{\vskip4pt}
	\text{for } T_{0} +S_{0} \leq N+2 \text{ and } |\gamma|\leq |\alpha|. &
	\end{cases}
	\end{align}
\end{linenomath}
We have
\begin{linenomath}
	\begin{multline}\label{T_1}
	\| v \|^{2}  
	\leq A_{1}^{2S_{0}} B_{1}^{2T_{0}}  N^{2U_{0}} 
	=A_{1}^{2S_{1}-2} B_{1}^{2T_{1}-2}  N^{2U_{1}-2s} 
	\\
	=A_{1}^{2S_{1}} B_{1}^{2T_{1}}  N^{2U_{1}}
	\times
	\left( A_{1}^{-2} B_{1}^{-2} N^{-2s} \right).
	\end{multline}
\end{linenomath}
About the first sum on the right hand side of \eqref{Est-1/r-bis}, its terms
are associated to $4n+1-$tuples $(\alpha, \beta+j, \gamma, \delta, k)$ and
$(\alpha, \beta, \gamma, \delta+j, k)$ respectively,
both satisfy the condition \eqref{Pr1-1_1-n}.
By induction we have
\begin{linenomath}
	\begin{multline}\label{T_2}
	\sum_{j=1}^{n} \left(\| \psi_{N}^{(\beta+j)} \Theta_{N}^{(\gamma)} D^{\alpha} g\|^{2}
	+
	\| \psi_{N}^{(\beta)} \Theta^{(\gamma)}_{N} D^{\alpha}  \varphi_{N}^{(\delta +j)} f\|^{2}\right)
	\\
	\leq 
	2\sum_{j=1}^{n} A_{1}^{2S_{1}} B_{1}^{2T_{0}}  N^{2U_{1}}
	\leq n(n+1) A_{1}^{2S_{1}} B_{1}^{2T_{1}-2}  N^{2U_{1}}
	\\
	=
	A_{1}^{2S_{1}} B_{1}^{2T_{1}}  N^{2U_{1}}
	\times
	\left( n(n+1) B_{1}^{-2} \right).
	\end{multline}
\end{linenomath}
We focus, now, on the third term on the right hand side of \eqref{Est-1/r-bis}:
\begin{linenomath}
	\begin{multline*}
	\sum_{j=1}^{n}
	\left( \| \psi_{N}^{(\beta)} \Theta_{N}^{(\gamma+j) } D^{\alpha} g\|_{1}^{2}
	+\alpha_{j}^{2} \| \psi_{N}^{(\beta)}\Theta_{N}^{(\gamma) } D^{\alpha-j} g \|^{2}_{1}
	\right)
	\\
	\leq
	\longsum[8]_{j,j_{1}=1}^{n}
	\left( \| \psi_{N}^{(\beta+j_{1})} \Theta_{N}^{(\gamma+j) } D^{\alpha} g\|^{2}
	+ \| \psi_{N}^{(\beta)} \Theta_{N}^{(\gamma+j) } D^{\alpha+j_{1}} g\|^{2}
	\right)
	\\
	\hspace{6.5em}+
	\sum_{j=1}^{n} \alpha_{j}^{2}
	\left(\sum_{ j_{1}=1}^{n}
	\left(\| \psi_{N}^{(\beta+j_{1})}\Theta_{N}^{(\gamma) } D^{\alpha-j} g \|^{2}
	+
	\| \psi_{N}^{(\beta)}\Theta_{N}^{(\gamma) } D^{\alpha-j+j_{1}} g \|^{2}.
	\right)\right) ,
	\end{multline*}
\end{linenomath}
The terms in the two sums have as associated $4n+1-$tuples
$(\alpha, \beta+j_{1}, \gamma+j, \delta, k)$,
$(\alpha+j_{1}, \beta, \gamma+j, \delta, k)$,
$(\alpha- j, \beta+j_{1}, \gamma, \delta, k)$ and
$(\alpha- j+j_{1}, \beta, \gamma, \delta, k)$ respectively,
they satisfy the condition \eqref{Pr1-1_1-n}.
By induction, the right hand side of the above
inequality can be estimated by 
\begin{linenomath}
	\begin{multline*}
	\longsum[8]_{j,j_{1}=1}^{n}
	\left(
	A^{2S_{1}}_{1} B_{1}^{2T_{1}-4r}N^{2U_{1}-2s(r-1)}
	+ A^{2S_{0}}_{1} B_{1}^{2T_{1}}N^{2U_{1}}
	\right)
	\\
	\hspace{6em}+\sup_{j} \{\alpha_{j}^{2}\}
	\sum_{j=1}^{n} 
	\sum_{ j_{1}=1}^{n}
	\left( A^{2S_{1}}_{1} B_{1}^{2T_{0}-4r}N^{2U_{1}-2sr}
	+
	A^{2S_{0}}_{1} B_{1}^{2T_{0}}N^{2U_{0}}
	\right).
	\end{multline*}
\end{linenomath}
We conclude that
\begin{linenomath}
	\begin{multline}\label{T_3}
	\sum_{j=1}^{n}
	\left( \| \psi_{N}^{(\beta)} \Theta_{N}^{(\gamma+j) } D^{\alpha} g\|_{1}^{2}
	+\alpha_{j}^{2} \| \psi_{N}^{(\beta)}\Theta_{N}^{(\gamma) } D^{\alpha-j} g \|^{2}_{1}
	\right)
	\leq
	A_{1}^{2S_{1}} B_{1}^{2T_{1}}  N^{2U_{1}} 
	\\
	\times
	\frac{n^{2}(n+1)^{2}}{4}\left( B_{1}^{-4r} N^{-2s(r-1)} + A_{1}^{-2} + B_{1}^{-2(2r+1)} N^{-2s(r-1)} + A_{1}^{-2}B_{1}^{-2}\right),
	\end{multline}
\end{linenomath}
where we use that $\displaystyle\sup_{j} \alpha_{j}^{2} \leq N^{2s}$.
We stress that $r\geq 2$.\\
We handle the term in the third line of \eqref{Est-1/r-bis}:
\begin{linenomath}
	\begin{multline*}
	\left[\!
	\sum_{j =1}^{n}\!
	\left( \!
	\| \psi_{N}^{(\beta+j)} \Theta_{N}^{(\gamma)} D^{\alpha} g\|
	\!+\| \psi_{N}^{(\beta)} \Theta_{N}^{(\gamma+j) } D^{\alpha} g\| 
	\!+\alpha_{j} \| \psi_{N}^{(\beta)}\Theta_{N}^{(\gamma) } D^{\alpha-j} g \| 
	\!\right)
	\!\right]  \| v \|,
	\end{multline*}
\end{linenomath}
the terms in the sum are associated to $4n+1-$tuples  $(\alpha, \beta+j, \gamma, \delta, k)$,
$(\alpha, \beta, \gamma+j, \delta, k)$ and $(\alpha- j, \beta, \gamma, \delta, k)$ respectively.
By induction we get
\begin{linenomath}
	\begin{multline}\label{T_4}
	\left[\!
	\sum_{j =1}^{n}\!
	\left( \!
	\| \psi_{N}^{(\beta+j)} \Theta_{N}^{(\gamma)} D^{\alpha} g\|
	\!+\| \psi_{N}^{(\beta)} \Theta_{N}^{(\gamma+j) } D^{\alpha} g\| 
	\!+\alpha_{j} \| \psi_{N}^{(\beta)}\Theta_{N}^{(\gamma) } D^{\alpha-j} g \| 
	\!\right)
	\!\right]  \| v \|,
	\\  
	\leq 
	\sum_{j =1}^{n}\!
	\left(\!
	A_{1}^{S_{1}} B_{1}^{T_{0}}  N^{U_{1}}
	\!+\!  A_{1}^{S_{0}} B_{1}^{T_{1}-2r}  N^{U_{1}-sr}
	\!+\! \sup_{j} \{\alpha_{j}\}
	A_{1}^{S_{0}} B_{1}^{T_{0}-2r}  N^{U_{0}-sr}\!
	\right)\!
	A_{1}^{S_{0}} B_{1}^{T_{0}}  N^{U_{0}}
	\\
	\hspace*{-28.2em}
	\leq
	A_{1}^{2S_{1}} B_{1}^{2T_{1}}  N^{2U_{1}} 
	\\
	\times\!
	\frac{n(n+1)}{2}\left( A_{1}^{-1}B_{1}^{-2} N^{-s} \!+\! A_{1}^{-2} B_{1}^{-2r-1} N^{-s(r+1)}\! +\! A_{1}^{-2}B_{1}^{-2(r+1)}N^{-s(r+1)}\!\right),
	\end{multline}
\end{linenomath}
where we use that $\displaystyle\sup_{j} \alpha_{j} \leq N^{s}$.\\
The terms in the first sum in the fourth line of \eqref{Est-1/r-bis}
are associated to $4n+1-$tuples $(\alpha, \beta+\mu, \gamma, \delta, k)$, 
by induction we have
\begin{linenomath}
	\begin{multline}\label{T_5}
	\sum_{|\mu| =2} \|\psi_{N}^{(\beta+\mu)} \Theta_{N}^{(\gamma)} D^{\alpha} g\| 
	\| v \|
	\leq
	\left(\sum_{|\mu| =2} A_{1}^{S_{1}+1} B_{1}^{T_{0}} N^{U_{1}+s}\right) A_{1}^{S_{0}} B_{1}^{T_{0}}  N^{U_{0}}
	\\
	\leq
	A_{1}^{2S_{1}} B_{1}^{2T_{1}}  N^{2U_{1}} 
	\times
	\left(\frac{n(n+1)}{2} B_{1}^{-2}\right).
	\end{multline}
\end{linenomath}
Subsequent terms in lines four, five, six and seven, on the right hand side of \eqref{Est-1/r-bis},
can be bounded as follows
\begin{linenomath}
	\begin{multline}\label{T_6}
	\sum_{j,\ell=1}^{n} 
	\left(\| \psi_{N}^{(\beta)}  \Theta_{N}^{(\gamma+j) } D^{\alpha+\ell}  g \|
	+\alpha_{j} \| \psi_{N}^{(\beta)}\Theta_{N}^{(\gamma) } D^{\alpha-j+\ell} g \| 
	\right)
	\| v \|
	\\
	\leq
	\sum_{j,\ell=1}^{n} 
	\left(A_{1}^{S_{0}}B_{1}^{T_{1}} N^{U_{1}}
	+\sup_{j}\{\alpha_{j}\} A_{1}^{S_{0}} B_{1}^{T_{0}}N^{U_{0}}
	\right)
	A_{1}^{S_{0}} B_{1}^{T_{0}}  N^{U_{0}}
	\\
	\leq
	A_{1}^{2S_{1}} B_{1}^{2T_{1}}  N^{2U_{1}} 
	\times
	\left(
	\frac{n^{2}(n+1)^{2}}{4} A_{1}^{-2} B_{1}^{-1} N^{-s}\left(1+B_{1}^{-1}\right)
	\right);
	\end{multline}
\end{linenomath}
\begin{linenomath}
	\begin{multline}\label{T_7}
	\sum_{j=1}^{n} \| \psi_{N}^{(\beta)} \Theta^{(\gamma)}_{N} D^{\alpha}  \varphi_{N}^{(\delta +j)} f\|
	\| v \|
	\leq
	\left(\sum_{j=1}^{n} A_{1}^{S_{1}} B_{1}^{T_{0}} N^{U_{1}}\right) 	A_{1}^{S_{0}} B_{1}^{T_{0}}  N^{U_{0}}
	\\
	\leq
	A_{1}^{2S_{1}} B_{1}^{2T_{1}}  N^{2U_{1}} 
	\times
	\left(
	\frac{n(n+1)}{2} A_{1}^{-1} B_{1}^{-2} N^{-s}
	\right);
	\end{multline}
\end{linenomath}
\begin{linenomath}
	\begin{multline}\label{T_8}
	\longsum[8]_{j,j_{1},\ell=1}^{n}\left(
	\| \psi_{_{N}}^{(\beta+j_{1})} \Theta_{N}^{(\gamma+j) } D^{\alpha+\ell} g\|
	+ 
	\| \psi_{N}^{(\beta+j_{1} +\ell)} \Theta_{N}^{(\gamma+j) } D^{\alpha} g\|
	\right.
	\\
	\hspace{4em}\left.
	+ \alpha_{j} 
	\left(\| \psi_{N}^{(\beta+j_{1} +\ell)} \Theta_{N}^{(\gamma) } D^{\alpha-j} g \| +\| \psi_{_{N}}^{(\beta+j_{1})} \Theta_{N}^{(\gamma) } D^{\alpha-j + \ell} g \|
	\right)\right)
	\| v \|
	\\
	\hspace*{-6.5em}
	\leq
	\longsum[8]_{j,j_{1},\ell=1}^{n}\left(
	A_{1}^{S_{1}} B_{1}^{T_{1}} N^{U_{1}+s}
	+
	A_{1}^{S_{1}+1} B_{1}^{T_{1}-2r} N^{U_{1}-s(r-2)}
	\right.
	\\
	\hspace*{3em}
	\left.
	+\sup_{j}\{\alpha_{j}\}
	\left(
	A_{1}^{S_{1}+1} B_{1}^{T_{0}-2r} N^{U_{1}-s(r-1)}
	+
	A_{1}^{S_{1}} B_{1}^{T_{0}} N^{U_{1}}
	\right)
	\right) A_{1}^{S_{0}} B_{1}^{T_{0}}  N^{U_{0}}
	\\
	\hspace*{-22em}
	\leq
	A_{1}^{2S_{1}} B_{1}^{2T_{1}}  N^{2U_{1}} 
	\\
	\hspace{1em}
	\times
	\frac{n^{3}(n+1)^{3}}{8}
	B_{1}^{-1}
	\left(\! A_{1}^{-1}\!+\! B_{1}^{-2r} N^{-s(r-1)}\!\left(1 +B_{1}^{-1}\right)
	+\! A_{1}^{-1}B_{1}^{-1}
	\!\right),
	\end{multline}
\end{linenomath}
here, we recall that $r\geq 2$;
\begin{linenomath}
	\begin{multline}\label{T_9}
	\sum_{j,\ell=1}^{n}\left(
	\| \psi_{N}^{(\beta+\ell)} \Theta^{(\gamma)}_{N} D^{\alpha}  \varphi_{N}^{(\delta +j)} f\|
	+
	\| \psi_{N}^{(\beta+\ell)}  \Theta_{N}^{(\gamma+j) } D^{\alpha}  g \|
	\right.
	\\
	\hspace{15em}
	\left.
	+
	\alpha_{j} \| \psi_{N}^{(\beta+\ell)}\Theta_{N}^{(\gamma) } D^{\alpha-j} g \| 
	\right)
	\| v \|
	\\
	\leq
	\sum_{j,\ell=1}^{n}\! \!\left(\!
	A_{1}^{S_{1}+1} B_{1}^{T_{0}}  N^{U_{1}+s}
	\!+\!
	A_{1}^{S_{1}} B_{1}^{T_{1}-2r}  N^{U_{1}-sr}
	\!+\!
	\alpha_{j} 
	A_{1}^{S_{1}} B_{1}^{T_{0}-2r}  N^{U_{1}-sr}
	\!\right) 
	\\
	\hspace*{24em}
	\times
	A_{1}^{S_{0}} B_{1}^{T_{0}}  N^{U_{0}}
	\\
	\leq
	A_{1}^{2S_{1}} B_{1}^{2T_{1}}  N^{2U_{1}} 
	\times
	\frac{n^{2}(n+1)^{2}}{4}\!
	\left(\! B_{1}^{-2}\!+\! A_{1}^{-1} B_{1}^{-2r-1}  N^{-s(r+1)}
	\!+\!
	A_{1}^{-1} B_{1}^{-2(r+1)}  N^{-sr}
	\right),
	\end{multline}
\end{linenomath}
where we use that $\displaystyle\sup_{j} \alpha_{j} \leq N^{s}$.\\
Let us now consider the next term in \eqref{Est-1/r-bis}.
Let $K_{1}$ be a compact set containing $\Omega_{2}$,
where $\psi_{N}$ is supported, and contained in $\Omega_{3}$. 
We have
\begin{linenomath}
$$ |a_{\ell,j}^{(\mu)}(x)| \leq C_{a_{\ell,j},K_{1}}^{|\mu|+1} |\mu|^{s|\mu|}, \,\, |b_{\ell}^{(\mu)}(x)| \leq C_{b_{\ell},K_{1}}^{|\mu|+1} |\mu|^{s|\mu|}
\text{ and }  |c^{(\mu)}(x)| \leq C_{c,K_{1}}^{|\mu|+1} |\mu|^{s|\mu|},
$$
\end{linenomath}
for every $x \in K_{1}$ and $\mu\in \mathbb{N}^{n}$. 
We set $\tilde{C}= sup\{C_{a_{\ell,j},K_{1}}, C_{b_{\ell},K_{1}}, C_{c,K_{1}}\}$.\\
By induction hypothesis (we stress that $|\mu|\geq 2$) we get
\begin{linenomath}
	\begin{multline}\label{T_10}
	\longsum[35]_{2 \leq |\mu| \leq |\alpha| -|\gamma| + 1} \longsum[7]_{\substack{\nu\leq \mu \\ \nu\leq \alpha}}
	\frac{1}{\left(\mu-\nu\right)!}  \binom{\alpha}{\nu}
	\left( \longsum[8]_{j,\ell=1}^{n} 
	\|a_{\ell,j}^{(\mu)} \psi_{N}^{(\beta)}  \Theta_{N}^{(\gamma +\mu-\nu)} D^{\alpha-\nu +j+\ell} g \|
	\right.
	\\
	+
	\sum_{ \ell=1}^{n} \|  b_{\ell}^{(\mu)} \psi_{N}^{(\beta)}  \Theta_{N}^{(\gamma +\mu-\nu)} D^{\alpha-\nu+\ell} g \|
	+
	\|  c^{(\mu)} \psi_{N}^{(\beta)}  \Theta_{N}^{(\gamma +\mu-\nu)} D^{\alpha-\nu}g \|
	\Biggr)
	\|v\|
	\\
	\leq
	\longsum[35]_{2 \leq |\mu| \leq |\alpha| -|\gamma| + 1} \longsum[7]_{\substack{\nu\leq \mu \\ \nu\leq \alpha}}
	\frac{1}{\left(\mu-\nu\right)!}  \binom{\alpha}{\nu} \tilde{C}^{|\mu|+1} |\mu|^{s|\mu|}\!\!
	\left(\! \longsum[8]_{j,\ell=1}^{n} 
	\| \psi_{N}^{(\beta)}  \Theta_{N}^{(\gamma +\mu-\nu)} D^{\alpha-\nu +j+\ell} g \|\right.
	\\
	+
	\sum_{ \ell=1}^{n} \| \psi_{N}^{(\beta)}  \Theta_{N}^{(\gamma +\mu-\nu)} D^{\alpha-\nu+\ell} g \|
	+
	\|  \psi_{N}^{(\beta)}  \Theta_{N}^{(\gamma +\mu-\nu)} D^{\alpha-\nu}g \|
	\Biggr)
	\|v\|
	\\
	\hspace{-4em}
	\leq
	\frac{n^{2}(n+1)^{2}}{4} A_{1}^{S_{0}} B_{1}^{T_{0}}  N^{U_{0}}
	\longsum[35]_{2 \leq |\mu| \leq |\alpha| -|\gamma| + 1} \longsum[7]_{\substack{\nu\leq \mu \\ \nu\leq \alpha}}
	\frac{1}{\left(\mu-\nu\right)!}  \binom{\alpha}{\nu} \tilde{C}^{|\mu|+1} |\mu|^{s|\mu|}
	\\
	\times
	A_{1}^{S_{0}} B_{1}^{T_{0} +4r -|\mu|(2r-1) -|\nu|} N^{U_{0} +s(2r- |\mu|(r-1) -|\nu|) }
	\left( 1+ B_{1}^{-2r} N^{-sr} + B_{1}^{-4r} N^{-2sr} 
	\right)
	\\
	\\
	\hspace{-7em}
	\leq
	A_{1}^{2S_{1}} B_{1}^{2T_{1}}  N^{2U_{1}}\,\,\times\left[
	\frac{n^{2}(n+1)^{2}}{4} A_{1}^{-2} \left( 1+ B_{1}^{-2r} N^{-sr} + B_{1}^{-4r} N^{-2sr} \right)\right]
	\\
	\times
	\longsum[35]_{2 \leq |\mu| \leq |\alpha| -|\gamma| + 1} \longsum[7]_{\substack{\nu\leq \mu \\ \nu\leq \alpha}}
	\frac{1}{\left(\mu-\nu\right)!}  \binom{\alpha}{\nu} \tilde{C}^{|\mu|+1} |\mu|^{s|\mu|}
	B_{1}^{- (|\mu|-2)(2r-1) -|\nu|} N^{- s\left[(|\mu|-2)(r-1)+|\nu|\right] }.
	\end{multline}
\end{linenomath}
Since
\begin{linenomath}
$$
\frac{1}{\left(\mu-\nu\right)!}  \binom{\alpha}{\nu} = \frac{\alpha!}{\left(\mu-\nu\right)! \nu!\left(\alpha-\nu\right)!} \leq N^{|\nu|}
\,\,\text{ and }\,\, 
|\mu|^{s|\mu|} \leq 3^{s|\mu|} N^{s(|\mu|-2)},
$$ 
\end{linenomath}
by Lemma  \ref{L-1}, we can bound the terms in the double sum by
\begin{linenomath}
$$
\tilde{C}^{|\mu|+1}  3^{s|\mu|} B_{1}^{- (|\mu|-2)(2r-1) -|\nu|} N^{- s(|\mu|-2)(r-2)-(s-1)|\nu|}.
$$
\end{linenomath}
Moreover, since  $r\geq 2$, $s\geq 1$ and $|\mu|\geq 2$,
we have that
\begin{linenomath}
$$
- s(|\mu|-2)(r-2)-(s-1)|\nu|\leq 0.
$$
\end{linenomath}
The above quantity can be estimated by
\begin{linenomath}
$$
\tilde{C}^{|\mu|+1}  3^{s|\mu|} B_{1}^{- (|\mu|-2)(2r-1) -|\nu|}
\leq \tilde{C}^{2}  3^{2s} B_{1}^{-|\nu|}\left( \tilde{C}\, 3^{s} \,B_{1}^{- 2r+1 }\right)^{|\mu|-2} .
$$
\end{linenomath}
Now, taking $B_{1} $ greater than $2$ and large enough so
that $\tilde{C}\, 3^{s} \,B_{1}^{- 2r+1 } \leq 2^{-1} $,
we obtain
\begin{linenomath}
	\begin{multline*}
	\longsum[35]_{2 \leq |\mu| \leq |\alpha| -|\gamma| + 1} \longsum[7]_{\substack{\nu\leq \mu \\ \nu\leq \alpha}}
	\frac{1}{\left(\mu-\nu\right)!}  \binom{\alpha}{\nu} \tilde{C}^{|\mu|+1} |\mu|^{s|\mu|}
	B_{1}^{- (|\mu|-2)(2r-1) -|\nu|} N^{- s\left[(|\mu|-2)(r-1)+|\nu|\right] }
	\\
	\leq \tilde{C}^{2}  3^{2s}
	\longsum[6]_{\mu_{1} = 0}^{\infty} \left(\frac{1}{2}\right)^{\mu_{1}}\cdots 
	\longsum[6]_{\mu_{n} =0}^{\infty}  \left(\frac{1}{2}\right)^{\mu_{n}}
	\longsum[6]_{\nu_{1} = 0} \left(\frac{1}{2}\right)^{\nu_{1}}\cdots 
	\longsum[6]_{\nu_{n} =0 }^{\infty} \left(\frac{1}{2}\right)^{\nu_{n}} 
	\\
	\leq
	\tilde{C}^{2}  3^{2s} 2^{2n}.
	\end{multline*}
\end{linenomath}
Summing up, we conclude that the term on the left hand side of
\eqref{T_10} can be estimated by
\begin{linenomath}
	\begin{multline}\label{T10-1}
	A_{1}^{2S_{1}} B_{1}^{2T_{1}}  N^{2U_{1}}
	\\
	\times \left( 
	\frac{n^{2}(n+1)^{2}}{4} \tilde{C}^{3}  3^{2s} 2^{2n} \,
	A_{1}^{-2} 
	\right)
	\left( 1\!+\! B_{1}^{-2r} N^{-sr}\! + \!B_{1}^{-4r} N^{-2sr} \right).
	\end{multline}
\end{linenomath}
Using the same strategy, we have
\begin{linenomath}
	\begin{multline}\label{T_11}
	\longsum[6]_{j,\ell=1}^{n}\longsum[35]_{1 \leq |\mu| \leq |\alpha| -|\gamma|+ 1} \longsum[5]_{\substack{\nu\leq \mu \\ \nu\leq \alpha}}
	\frac{1}{\left(\mu-\nu\right)!}  \binom{\alpha}{\nu}\!
	\| a_{\ell,j}^{(\mu)}\psi_{N}^{(\beta)} \Theta^{(\gamma+\mu-\nu)}_{N} \!D^{\alpha-\nu+\ell}\varphi_{N}^{(\delta +j)}\! f\|
	\\[-7pt]
	\hspace{32em}\times
	\|v\|
	\\
	\leq
	A_{1}^{2S_{1}} B_{1}^{2T_{1}}  N^{2U_{1}}\,\,
	\times
	\left(
	\frac{n^{2}(n+1)^{2}}{4}   \tilde{C}^{2}  3^{s} 2^{2n}  A_{1}^{-1} B_{1}^{-1}
	\right);
	\phantom{ BBBBBBBBBBBB}
	\end{multline}
\end{linenomath}
and
\begin{linenomath}
	\begin{multline}\label{T_12}
	\longsum[30]_{0 \leq |\mu| \leq |\alpha| -|\gamma|} \longsum[6]_{\substack{\nu\leq \mu \\ \nu\leq \alpha}}
	\frac{1}{\left(\mu-\nu\right)!}  \binom{\alpha}{\nu} \!\!
	\!\left( \!\longsum[6]_{j,\ell=1}^{n}
	\| a_{\ell,j}^{(\mu)}\psi_{N}^{(\beta)} \Theta^{(\gamma+\mu-\nu)}_{N} D^{\alpha-\nu} \varphi_{N}^{(\delta+\ell +j)} f\| 
	\right.
	\\
	\hspace{13em}\left.
	+\sum_{\ell=1}^{n} 
	\| b_{\ell}^{(\mu)}\psi_{N}^{(\beta)} \Theta^{(\gamma+\mu-\nu)}_{N} D^{\alpha-\nu} \varphi_{N}^{(\delta+\ell)} f\| 
	\Biggr)
	\right. \|v\|
	\\
	\leq
	A_{1}^{2S_{1}} B_{1}^{2T_{1}}  N^{2U_{1}}\,\,
	\times
	\left[\frac{n^{2}(n+1)^{2}}{4}   \tilde{C}  2^{2n}  B_{1}^{-2}\left(1+A_{1}^{-1}N^{-s}\right)\right],
	\phantom{ BBBBBBBBB}
	\end{multline}
\end{linenomath}
in this case we use that $|\mu|^{s|\mu|} \leq N^{s|\mu|}$.\\
About the first term on the last line on the right hand side of \eqref{Est-1/r-bis}:
\begin{linenomath}
	\begin{align*}
	\| \psi_{N}^{(\beta)} \Theta^{(\gamma)}_{N} D^{\alpha} \varphi_{N}^{(\delta)} P^{k+1} u\|\|v\|,
	\end{align*}
\end{linenomath}
since the first factor is associated to the $4n+1-$tuple $(\alpha,\beta,\gamma,\delta,k+1)$,
it can be estimated by
\begin{linenomath}
	\begin{align}\label{T_13}
	A_{1}^{S_{1}+1} B_{1}^{T_{0}}  N^{U_{1}+1}
	\times
	A_{1}^{S_{0}} B_{1}^{T_{0}}  N^{U_{0}}
	\leq
	A_{1}^{2S_{1}} B_{1}^{2T_{1}}  N^{2U_{1}}
	\times
	B_{1}^{-2}.
	\end{align}
\end{linenomath}
Concerning the last two remaining terms in the right hand side of \eqref{Est-1/r-bis}, 
they can be handled as follows 
\begin{linenomath}
	\begin{multline}\label{T_14}
	C_{1}^{S_{0}} C_{2}^{2 m +|\gamma| + 1 } N^{s(2m +|\gamma| +\sigma+2n + 5 -M)^{+}} 
	\|v\|
	\\
	\leq 
	C_{1}^{S_{0}} C_{2}^{T_{0}-2m(r-1)} N^{U_{0}-sm(r-2)}
	A_{1}^{S_{0}} B_{1}^{T_{0}}  N^{U_{0}}
	\\
	\leq
	A_{1}^{2S_{1}} B_{1}^{2T_{1}}  N^{2U_{1}}
	\times
	\left(
	A_{1}^{-2}B_{1}^{-2}  N^{-s\left[m(r-2)+2\right]} C_{2}^{-2(r-1)m}
	\left( C_{1}A_{1}^{-1}\right)^{S_{0}} \left( C_{2}B_{1}^{-1}\right)^{T_{0}}
	\right), 
	\end{multline}
\end{linenomath}
here we use that $M=3(n+6)$,
and
\begin{linenomath}
	\begin{multline}\label{T_15}
	C_{1}^{2S_{0}} C_{2}^{2(2 m+|\gamma| + 1)}  
	N^{2s(2 m+ |\gamma|+\sigma) }
	\leq
	A_{1}^{2S_{1}} B_{1}^{2T_{1}}  N^{2U_{1}}
	\\
	\times
	A_{1}^{-2}B_{1}^{-2}
	C_{2}^{-4m(r-1)} N^{-2s[m(r-2) + 2]} 
	\left( C_{1}A_{1}^{-1}\right)^{2S_{0}} \left( C_{2}B_{1}^{-1}\right)^{2T_{0}}.
	\end{multline}
\end{linenomath}
Summing up, if $A_{1}$ and $B_{1}$ are chosen large enough, with $B_{1}$ large compared to $A_{1}$,
the sum of second factor on the right
hand side of \eqref{T_1}, \eqref{T_2}, \eqref{T_3}, \eqref{T_4}, \eqref{T_4}, \eqref{T_5}, \eqref{T_6},
\eqref{T_7}, \eqref{T_8}, \eqref{T_9}, \eqref{T10-1}, \eqref{T_11}, \eqref{T_12}, \eqref{T_13}, \eqref{T_14}
and \eqref{T_15} can be made smaller than $(2C_{2})^{-1}$. We conclude
\begin{linenomath}	
	\begin{align}\label{Est-p=1}
	\| v\|_{\frac{1}{r}} ^{2}
	\leq 
	A_{1}^{2S_{1}} B_{1}^{2T_{1}}  N^{2U_{1}}.
	\end{align}
\end{linenomath}
So we obtained \eqref{Pr1-1_1-n}.

\medskip
\textbf{Step  $p>1$.} Using the notation introduced in \eqref{Not-exp}, we assume that 
\begin{linenomath}
	\begin{align}\label{Pr1-p-1_1-n}
	(1)_{p-1}
	\begin{cases}
	\| \psi_{N}^{(\beta)} \Theta_{N}^{(\gamma)} D^{\alpha} \varphi_{N}^{(\delta)} P^{k} u\|_{\frac{p-1}{r}} 
	\leq 
	A_{1}^{S_{p-1}} B_{1}^{T_{p-1}}  N^{U_{p-1}},&\\
	\noalign{\vskip4pt}
	\text{for } T_{p-1}+S_{p-1} \leq N+2 \text{ and } |\gamma|\leq |\alpha|. &
	\end{cases}
	\end{align}
\end{linenomath}
We want to show, via the estimates obtained in previous section,
\eqref{Est-1/r-tild}  and \eqref{Est-p/r-f}, that 
\begin{linenomath}
	\begin{align}\label{Pr1-p_1-n}
	(1)_{p}
	\begin{cases}
	\| \psi_{N}^{(\beta)} \Theta_{N}^{(\gamma)} D^{\alpha} \varphi_{N}^{(\delta)} P^{k} u\|_{\frac{p}{r}} 
	\leq 
	A_{1}^{S_{p}} B_{1}^{T_{p}}  N^{U_{p}},&\\
	\noalign{\vskip4pt}
	\text{for } T_{p}+S_{p} \leq N+2 \text{ and } |\gamma|\leq |\alpha|. &
	\end{cases}
	\end{align}
\end{linenomath}
holds.\\
Combining \eqref{Est-1/r-tild}  and \eqref{Est-p/r-f} we have
\begin{linenomath}
	\begin{multline}\label{Est-p/r-bis}
	\| v\|_{\frac{p}{r}}^{2}
	\leq 
	C_{2}^{2} \! \Biggl\{ \!
	\|v\|_{\frac{p-1}{r}}^{2} \!
	+ \!\sum_{j =1}^{n} \! \| \psi_{N}^{(\beta+j)}  w\|^{2}_{\frac{p-1}{r}}
	\!+\! \sum_{j =1}^{n} \!\|\psi_{N}^{(\beta)}   \Theta_{N}^{(\gamma)} D^{\alpha}
	\varphi_{N}^{(\delta+j)}  P^{k} u\|^{2}_{\frac{p-1}{r}}
	\\
	+
	\longsum[9]_{j_{1},j=1}^{n} \left(
	\| \psi_{N}^{(\beta +j_{1})} \Theta_{N}^{(\gamma+j)} D^{\alpha} g \|_{\frac{p-1}{r}}^{2}
	+
	\| \psi_{N}^{(\beta)} \Theta_{N}^{(\gamma+j)} D^{\alpha+j_{1}} g \|_{\frac{p-1}{r}}^{2}
	\right.
	\\
	\left.
	\qquad\quad
	+
	\alpha_{j}^{2}\| \psi_{N}^{(\beta+j_{1})} \Theta_{N}^{(\gamma)} D^{\alpha-j} g \|_{\frac{p-1}{r}}^{2}
	+
	\alpha_{j}^{2}\| \psi_{N}^{(\beta)} \Theta_{N}^{(\gamma)} D^{\alpha-j+j_{1}} g \|_{\frac{p-1}{r}}^{2}
	\right)
	\\
	+
	\Biggl[\sum_{|\mu| = 2 } \| \! \psi_{N}^{(\beta+\mu)}  w\|_{\frac{p-1}{r}}
	\!+\!
	\longsum[9]_{j,j_{1}=1}^{n}\!
	\left( 
	\| \psi_{N}^{(\beta+j_{1})} \Theta_{N}^{(\gamma +j)} D^{\alpha} g\|_{\frac{p-1}{r}}
	\!+\!
	\alpha_{j} \|  \psi_{N}^{(\beta)} \Theta_{N}^{(\gamma)} D^{\alpha - j+j_{1}} g\|_{\frac{p-1}{r}}
	\right)
	\\
	+
	\longsum[11]_{j,j_{1},j_{2}=1}^{n} 
	\left( 
	\| \psi_{N}^{(\beta +j_{1})} \Theta_{N}^{(\gamma+j)} D^{\alpha+j_{2}} g\|_{\frac{p-1}{r}}
	+
	\alpha_{j} \| \psi_{N}^{(\beta +j_{1})} \Theta_{N}^{(\gamma)} D^{\alpha+j_{2}-j} g \|_{\frac{p-1}{r}}
	\right)
	\\
	+
	\sum_{ j=1}^{n}\longsum[9]_{|\mu|=2} 
	\left( \| \psi_{N}^{(\beta +\mu)} \Theta_{N}^{(\gamma+j)} D^{\alpha} g \|_{\frac{p-1}{r}}
	+
	\alpha_{j} \| \psi_{N}^{(\beta +\mu)} \Theta_{N}^{(\gamma)} D^{\alpha-j} g \|_{\frac{p-1}{r}}
	\right)
	\\
	+
	\longsum[11]_{j,j_{1}=1}^{n} 
	\left( 
	\| \psi_{N}^{(\beta)} \Theta_{N}^{(\gamma +j)} D^{\alpha+j_{1}} g\|_{\frac{p-1}{r}}
	+
	\| \psi_{N}^{(\beta+j_{1})} \Theta_{N}^{(\gamma)} D^{\alpha} \varphi_{N}^{(\delta+j)}  P^{k} u\|_{\frac{p-1}{r}}
	\right)
	\\
	+
	\longsum[44]_{2 \leq |\mu| \leq |\alpha| -|\gamma| + \lfloor \frac{n}{2} \rfloor+3} 
	\longsum[7]_{\substack{\nu\leq \mu \\ \nu\leq \alpha}}
	\frac{1}{\mu!}  
	\binom{\mu}{\nu} \frac{\alpha!}{(\alpha -\nu)!}
	C_{1}^{|\mu|+1} N^{ s\left( |\mu|+ n + 3 -M\right)^{+}} 
	\\
	\times
	\left( \longsum[11]_{j_{1},j_{2}=1}^{n} 
	\| \psi_{N}^{(\beta)} \Theta_{N}^{(\gamma+\mu-\nu)} D^{\alpha-\nu+j_{1}+j_{2}} g\|_{\frac{p-1}{r}}
	+
	\| \psi_{N}^{(\beta)} \Theta_{N}^{(\gamma+\mu-\nu)} D^{\alpha-\nu} g\|_{\frac{p-1}{r}}
	\right.
	\\
	\hspace*{18em}
	\left.
	+
	\sum_{j_{1}=1}^{n} 
	\| \psi_{N}^{(\beta)} \Theta_{N}^{(\gamma+\mu-\nu)} D^{\alpha-\nu+j_{1}} g\|_{\frac{p-1}{r}}
	\right)
	\\
	+\longsum[11]_{j,j_{1}=1 }^{n} \longsum[44]_{1 \leq |\mu| \leq |\alpha| -|\gamma| + \lfloor \frac{n}{2} \rfloor +2} 
	\longsum[7]_{\substack{\nu\leq \mu \\ \nu\leq \alpha}}
	\frac{1}{\mu!} 
	\binom{\mu}{\nu} \frac{\alpha!}{(\alpha -\nu)!}
	C_{4}^{|\mu|} \,N^{s(|\mu|+1 -M)^{+}}
	\\
	\times
	\left(
	\| \psi_{N}^{(\beta)}  \Theta_{N}^{(\gamma +\mu-\nu)} D^{\alpha-\nu} \varphi_{N}^{(\delta+j+j_{1})} f\|_{\frac{p-1}{r}}
	+
	\| \psi_{N}^{(\beta)}  \Theta_{N}^{(\gamma +\mu-\nu)} D^{\alpha-\nu+j_{1}} \varphi_{N}^{(\delta+j)} f\|_{\frac{p-1}{r}}
	\right)
	\\
	+
	\| \psi_{N}^{(\beta)} \Theta_{N}^{(\gamma)} D^{\alpha}  \varphi_{N}^{(\delta)} P^{k+1} u\|_{\frac{p-1}{r}}
	+
	C_{1}^{ \sigma+ 1}
	C_{2}^{ 2m+|\gamma|}\,N^{s(2m +|\gamma| +\sigma)} 
	\Biggr]
	\|v\|_{\frac{p-1}{r}}
	\Biggr\}, 
	\\
	+
	C_{2}^{2}\|v\|_{0}^{2} 
	+
	C_{1}^{2(\sigma+1)} C_{2}^{2(2 m + |\gamma|+1)}
	N^{2s(2 m +|\gamma|+ \sigma )}.
	\end{multline}
\end{linenomath}	
The purpose is to show that all the terms on the right hand side of \eqref{Est-p/r-bis}
are smaller than $A_{1}^{2S_{p}} B_{1}^{2T_{p}}  N^{2U_{p}}$ times a factor depending
on negative powers of $A_{1}$ or $B_{1}$ or $N$.
A suitable choice of $A_{1}$ and $B_{1}$ will yield the summand of these factors less than one.
We remark that all the $4n+1-$tuples of the form
$(\alpha', \beta', \gamma', \delta', k') \in \mathbb{N}^{4n+1}$ associated to each term 
on the right hand side of \eqref{Est-p/r-bis} satisfy the condition in \eqref{Pr1-p-1_1-n}.\\
We estimate each of the terms on the right hand side of the above inequality separately
following the order in which they are written.
Using the same strategy used in the case $p=1$, we have
\begin{linenomath}
	\begin{align}\label{T_1p}
	\|v\|_{\frac{p-1}{r}}^{2} 
	\leq A_{1}^{2S_{p}} B_{1}^{2T_{p}}  N^{2U_{p}}
	\times
	\left( A_{1}^{-2} B_{1}^{-2} N^{-2s} \right),
	\end{align}
\end{linenomath}
and
\begin{linenomath}
	\begin{multline}\label{T_2p}
	\sum_{j =1}^{n} \!\!
	\left(\!\| \psi_{N}^{(\beta+j)}  w\|^{2}_{\frac{p-1}{r}}\!
	+ \!\|\psi_{N}^{(\beta)}   \Theta_{N}^{(\gamma)} D^{\alpha}
	\varphi_{N}^{(\delta+j)}  P^{k} u\|^{2}_{\frac{p-1}{r}}
	\!\right)\!
	\\
	\leq \!
	A_{1}^{2S_{p}} B_{1}^{2T_{p}}  N^{2U_{p}}
	\times
	\left( n(n+1) B_{1}^{-2}\right),
	\end{multline}
\end{linenomath}
where we recall $w=  \Theta_{N}^{(\gamma)} D^{\alpha} \varphi_{N}^{(\delta)}  P^{k} u$.\\
Using that $\displaystyle\sup_{j} \alpha_{j} \leq N^{s}$, we get
\begin{linenomath}
	\begin{multline}\label{T_3p}
	\longsum[9]_{j_{1},j=1}^{n} \left(
	\| \psi_{N}^{(\beta +j_{1})} \Theta_{N}^{(\gamma+j)} D^{\alpha} g \|_{\frac{p-1}{r}}^{2}
	+
	\| \psi_{N}^{(\beta)} \Theta_{N}^{(\gamma+j)} D^{\alpha+j_{1}} g \|_{\frac{p-1}{r}}^{2}
	\right.
	\\
	\left.
	\qquad\quad
	+
	\alpha_{j}^{2}\| \psi_{N}^{(\beta+j_{1})} \Theta_{N}^{(\gamma)} D^{\alpha-j} g \|_{\frac{p-1}{r}}^{2}
	+
	\alpha_{j}^{2}\| \psi_{N}^{(\beta)} \Theta_{N}^{(\gamma)} D^{\alpha-j+j_{1}} g \|_{\frac{p-1}{r}}^{2}
	\right)
	\\
	\leq
	A_{1}^{2S_{p}} B_{1}^{2T_{p}}  N^{2U_{p}} \phantom{ BBBBBBBBBBBBBBBBBBBBBBBBBBBBBB}
	\\
	\times
	\frac{n^{2}(n+1)^{2}}{4}
	\left( B_{1}^{-4r} N^{-2s(r-1)} +A_{1}^{-2} 
	+ B_{1}^{-2(2r+1)} N^{-2s(r-1)} + A_{1}^{-2} B_{1}^{-2}\right),
	\end{multline}
\end{linenomath}
and
\begin{linenomath}
	\begin{align}\label{T_4p}
	\sum_{|\mu| = 2 } \| \! \psi_{N}^{(\beta+\mu)}  w\|_{\frac{p-1}{r}} \|v\|_{\frac{p-1}{r}}
	\leq
	A_{1}^{2S_{p}} B_{1}^{2T_{p}}  N^{2U_{p}}
	\times
	\frac{n(n+1)}{2} B_{1}^{-2},
	\end{align}
\end{linenomath}
here we use that the number of the multi-indexes $\mu= (\mu_{1},\dots, \mu_{n})$
such that $|\mu|=q$ is given by $\binom{q+n-1}{n-1}$.\\
Concerning the subsequent terms on the right hand side of \eqref{Est-p/r-bis}, we get
\begin{linenomath}
	\begin{multline}\label{T_5p}	
	\longsum[9]_{j,j_{1}=1}^{n}\!
	\left( 
	\| \psi_{N}^{(\beta+j_{1})} \Theta_{N}^{(\gamma +j)} D^{\alpha} g\|_{\frac{p-1}{r}}
	\!+\!
	\alpha_{j} \|  \psi_{N}^{(\beta)} \Theta_{N}^{(\gamma)} D^{\alpha - j+j_{1}} g\|_{\frac{p-1}{r}}
	\right) \|v\|_{\frac{p-1}{r}}
	\\
	\leq 
	A_{1}^{2S_{p}} B_{1}^{2T_{p}}  N^{2U_{p}}
	\times
	\frac{n^{2}(n+1)^{2}}{4}
	\left( A_{1}^{-1} B_{1}^{-2r-1} N^{-sr} +A_{1}^{-2} B_{1}^{-2} N^{-s}\right),
	\end{multline}
\end{linenomath}
\begin{linenomath}
	\begin{multline}\label{T_6p}	
	\longsum[11]_{j,j_{1},j_{2}=1}^{n} \!\!
	\left( \!
	\| \psi_{N}^{(\beta +j_{1})} \Theta_{N}^{(\gamma+j)} D^{\alpha+j_{2}} g\|_{\frac{p-1}{r}}
	\!+\!
	\alpha_{j} \| \psi_{N}^{(\beta +j_{1})} \Theta_{N}^{(\gamma)} D^{\alpha+j_{2}-j} g \|_{\frac{p-1}{r}}
	\!
	\right)\!
	\|v\|_{\frac{p-1}{r}}
	\\
	\leq 
	A_{1}^{2S_{p}} B_{1}^{2T_{p}}  N^{2U_{p}}
	\times
	\frac{n^{3}(n+1)^{3}}{8} A_{1}^{-1}B_{1}^{-1}
	\left(1+  B_{1}^{-1} \right),
	\end{multline}
\end{linenomath}
\begin{linenomath}
	\begin{multline}\label{T_7p}	
	\sum_{ j=1}^{n}\longsum[9]_{|\mu|=2} \!
	\left( \| \psi_{N}^{(\beta +\mu)} \Theta_{N}^{(\gamma+j)} D^{\alpha} g \|_{\frac{p-1}{r}}
	\!+\!
	\alpha_{j} \| \psi_{N}^{(\beta +\mu)} \Theta_{N}^{(\gamma)} D^{\alpha-j} g \|_{\frac{p-1}{r}}
	\!
	\right)
	\|v\|_{\frac{p-1}{r}}
	\\
	\leq
	A_{1}^{2S_{p}} B_{1}^{2T_{p}}  N^{2U_{p}}
	\times
	\frac{n^{2}(n+1)^{2}}{4} B_{1}^{-2r-1} N^{-s(r-1)}
	\left(1+  B_{1}^{-1} \right),
	\end{multline}
\end{linenomath}
where we use once again that $\displaystyle\sup_{j} \alpha_{j} \leq N^{s}$,
moreover we recall that $r\geq 2$.\\
About the subsequent term on the right hand side of \eqref{Est-p/r-bis}, we have
\begin{linenomath}
	\begin{multline}\label{T_8p}	
	\longsum[9]_{j,j_{1}=1}^{n} \!
	\left( \!
	\| \psi_{N}^{(\beta)} \Theta_{N}^{(\gamma +j)} D^{\alpha+j_{1}} g\|_{\frac{p-1}{r}}
	\!+\!
	\| \psi_{N}^{(\beta+j_{1})} \Theta_{N}^{(\gamma)} D^{\alpha} \varphi_{N}^{(\delta+j)}  P^{k} u\|_{\frac{p-1}{r}}
	\!
	\right)
	\!
	\|v\|_{\frac{p-1}{r}}
	\\
	\leq
	A_{1}^{2S_{p}} B_{1}^{2T_{p}}  N^{2U_{p}}
	\times
	\frac{n^{2}(n+1)^{2}}{4} B_{1}^{-1}
	\left(A_{1}^{-2}N^{-s}+  B_{1}^{-1} \right).
	\end{multline}
\end{linenomath}
In order to handle the last two sums in the right of \eqref{Est-p/r-bis},
we have to distinguish two case: $ |\mu|\leq |\alpha| - |\gamma|+2$ and 
$|\mu|> |\alpha| - |\gamma| +2 $
(
we remark that in the sum the number of multi-index $\mu$ such that $|\alpha| - |\gamma|< |\mu| \leq |\alpha| - |\gamma| +  \left\lfloor \frac{n}{2} \right\rfloor +3$
is finite and it can be roughly estimated $2^{|\alpha| - |\gamma| +  \left\lfloor \frac{n}{2} \right\rfloor +3}$.
) \\
We begin to handle the first sum. The two kinds of terms in the sum will be treated differently in the next pages.
We split the sum.
Using the same strategy adopted to obtain \eqref{T10-1}, 
via the induction hypothesis we get
\begin{linenomath}
	\begin{multline}%
	\label{T_9p-1}	
	\longsum[44]_{2 \leq |\mu| \leq |\alpha| -|\gamma| +2} 
	\longsum[7]_{\substack{\nu\leq \mu \\ \nu\leq \alpha}}
	\frac{1}{\mu!}  
	\binom{\mu}{\nu} \frac{\alpha!}{(\alpha -\nu)!}
	C_{1}^{|\mu|+1} N^{ s\left( |\mu|+ n + 3 -M\right)^{+}} 
	\\
	\times
	\left( \longsum[11]_{j_{1},j_{2}=1}^{n} 
	\| \psi_{N}^{(\beta)} \Theta_{N}^{(\gamma+\mu-\nu)} D^{\alpha-\nu+j_{1}+j_{2}} g\|_{\frac{p-1}{r}}
	+
	\| \psi_{N}^{(\beta)} \Theta_{N}^{(\gamma+\mu-\nu)} D^{\alpha-\nu} g\|_{\frac{p-1}{r}}
	\right.
	\\
	\hspace*{14em}
	\left.
	+
	\sum_{j_{1}=1}^{n} 
	\| \psi_{N}^{(\beta)} \Theta_{N}^{(\gamma+\mu-\nu)} D^{\alpha-\nu+j_{1}} g\|_{\frac{p-1}{r}}
	\right)
	\!
	\|v\|_{\frac{p-1}{r}}
	\\
	\leq
	A_{1}^{2S_{p}} B_{1}^{2T_{p}}  N^{2U_{p}}
	\times
	4^{n}\left(n^{2}(n+1)^{2}\right) C_{1} A_{1}^{-2}\left( C_{1}^{2} +  B_{1}^{-2} +C_{1} B_{1}^{-1}\right).
    \end{multline}
\end{linenomath}
About the remaining terms in the sum. Since $ -|\mu|+|\alpha|-|\gamma| + 2 \leq -1$
we remark that
\begin{linenomath}
	\begin{multline*}
     \| \psi_{N}^{(\beta)} \Theta_{N}^{(\gamma+\mu-\nu)} D^{\alpha-\nu+j_{1}+j_{2}} g\|_{\frac{p-1}{r}}
     \\
     \leq
     C_{\psi}^{|\beta|+n+3} N^{ \left( |\beta|+ n + 3 -M\right)^{+}} 
     \|\Theta_{N}^{(\gamma+\mu-\nu)} D^{\alpha-\nu+j_{1}+j_{2}} g\|_{\frac{p-1}{r}}
     \\
     \leq
     C_{\psi}^{|\beta|+n+3}
     \widetilde{C}^{|\gamma|+|\mu|-|\nu| + 1 } N^{(|\beta|+|\gamma|+|\mu|-|\nu|+n+3 -M)^{+}}
     \| \varphi_{N}^{(\delta)} P^{k} u \|_{0}
     \\
     \leq
     \widetilde{C}_{1}^{\sigma+1}
     \widetilde{C}_{2}^{|\gamma|+|\mu|-|\nu| + 1 }
      N^{s(\sigma +|\gamma|+|\mu|-|\nu|+n+3 - M)^{+}},
	\end{multline*}
\end{linenomath}
where, we recall, $\sigma= |\beta|+|\delta|+2k$. 
The terms $\| \psi_{N}^{(\beta)} \Theta_{N}^{(\gamma+\mu-\nu)} D^{\alpha-\nu} g\|_{\frac{p-1}{r}}$
and $\| \psi_{N}^{(\beta)} \Theta_{N}^{(\gamma+\mu-\nu)} D^{\alpha-\nu+j_{1}} g\|_{\frac{p-1}{r}}$
can be estimated in the same way. Moreover we recall that
\begin{linenomath}
	$$
	\frac{1}{\left(\mu-\nu\right)!}  \binom{\alpha}{\nu} = \frac{\alpha!}{\left(\mu-\nu\right)! \nu!\left(\alpha-\nu\right)!} \leq N^{|\nu|}
	$$ 
\end{linenomath}
and since $M>2n+3$ then 
\begin{linenomath}
$$
N^{ s\left( |\mu|+ n + 3 -M\right)^{+}}\leq  N^{s(|\alpha|-|\gamma|)} 
\text{ and }N^{(|\beta|+|\gamma|+|\mu|-|\nu|+n+3 -M)^{+}} \leq N^{s(|\alpha|+|\beta|-|\nu|)},
$$
\end{linenomath}
$|\alpha| -|\gamma| +3 \leq |\mu| \leq |\alpha| -|\gamma| + \lfloor \frac{n}{2} \rfloor+2$.\\
Using the induction hypothesis on the factor $\|v\|_{\frac{p-1}{r}}$ and the above estimates,
we get
\begin{linenomath}
\begin{multline} \label{T_9p-2}	
\longsum[69]_{|\alpha| -|\gamma| +3 \leq |\mu| \leq |\alpha| -|\gamma| + \lfloor \frac{n}{2} \rfloor+2}
\longsum[7]_{\substack{\nu\leq \mu \\ \nu\leq \alpha}}
\frac{1}{\mu!}  
\binom{\mu}{\nu} \frac{\alpha!}{(\alpha -\nu)!}
C_{1}^{|\mu|+1} N^{ s\left( |\mu|+ n + 3 -M\right)^{+}} 
\\
\times
\left( \longsum[11]_{j_{1},j_{2}=1}^{n} 
\| \psi_{N}^{(\beta)} \Theta_{N}^{(\gamma+\mu-\nu)} D^{\alpha-\nu+j_{1}+j_{2}} g\|_{\frac{p-1}{r}}
+
\| \psi_{N}^{(\beta)} \Theta_{N}^{(\gamma+\mu-\nu)} D^{\alpha-\nu} g\|_{\frac{p-1}{r}}
\right.
\\
\hspace*{14em}
\left.
+
\sum_{j_{1}=1}^{n} 
\| \psi_{N}^{(\beta)} \Theta_{N}^{(\gamma+\mu-\nu)} D^{\alpha-\nu+j_{1}} g\|_{\frac{p-1}{r}}
\right) 
\|v\|_{\frac{p-1}{r}}
\\
\leq
\frac{3}{4}n^{2}(n+1)^{2}\longsum[69]_{|\alpha| -|\gamma| +3 \leq |\mu| \leq |\alpha| -|\gamma| + \lfloor \frac{n}{2} \rfloor+3}
\longsum[7]_{\substack{\nu\leq \mu \\ \nu\leq \alpha}}
N^{|\nu|} C_{1}^{|\mu|+1} N^{s(|\alpha|-|\gamma|)} 
\\
\times
\left( \widetilde{C}_{1}^{\sigma+1}
\widetilde{C}_{2}^{|\gamma|+|\mu|-|\nu| + 1 }
N^{s(\sigma +|\alpha|-|\nu|)}\right)
A_{1}^{S_{p-1}} B_{1}^{T_{p-1}}  N^{U_{p-1}},
\end{multline}
\end{linenomath}
where $ S_{p-1}=|\sigma| + p  $, $T_{p-1}= 2rm+|\gamma| + p $ and
$ U_{p-1} =s[rm+|\gamma|+ \sigma + p-1]$.\\
Now we observe that since without loss of generality we can assume that $\widetilde{C}_{2}>2$ we have
\begin{linenomath}
	\begin{align*}
\longsum[7]_{\substack{\nu\leq \mu \\ \nu\leq \alpha}}
 \widetilde{C}_{2}^{-|\nu|}
 \leq \longsum[6]_{\nu_{1} = 0} \left(\frac{1}{2}\right)^{\nu_{1}}\cdots 
 \longsum[6]_{\nu_{n} =0 }^{\infty} \left(\frac{1}{2}\right)^{\nu_{n}} 
 \leq 2^{n}
	\end{align*}
\end{linenomath}
and that the number of multi-index $\mu$ such that $|\mu| \leq |\alpha| -|\gamma| + \lfloor \frac{n}{2} \rfloor+2$
is smaller than $2^{|\alpha| -|\gamma| + 2(n+1)}$. So summing up the right hand side of \eqref{T_9p-2} can be estimated
by
\begin{linenomath}
	\begin{multline}%
	\label{T_9p-2f}
	A_{1}^{2S_{p}} B_{1}^{2T_{p}}  N^{2U_{p}}
	\\
	\times
	\left[
	2^{3n+2}n^{2}(n+1)^{2} \widetilde{C}_{2}^{ \lfloor \frac{n}{2} \rfloor+2}
	A_{1}^{-p-1} \left( A_{1}^{-1} \widetilde{C}_{1}\right)^{\sigma+1}
	\left(2\widetilde{C}_{2} B_{1}^{-1}\right)^{2rm+|\gamma|}
	\right].
	\end{multline}
\end{linenomath}		
Concerning the second sum in the right hand side of \eqref{Est-p/r-bis}
we proceed as before, i.e. distinguishing two case: $ |\mu|\leq |\alpha| - |\gamma|+2$ and 
$|\mu|> |\alpha| - |\gamma| +2 $.
These two kinds of terms will be treated differently as done before.
Using the same strategy adopted to obtain \eqref{T10-1}, 
via the induction hypothesis we get
\begin{linenomath}
	\begin{multline}\label{T_10p}	
	\longsum[11]_{j,j_{1}=1 }^{n} \longsum[44]_{1 \leq |\mu| \leq |\alpha| -|\gamma| +2} 
	\longsum[7]_{\substack{\nu\leq \mu \\ \nu\leq \alpha}}
	\frac{1}{\mu!} 
	\binom{\mu}{\nu} \frac{\alpha!}{(\alpha -\nu)!}
	C_{4}^{|\mu|} \,N^{s(|\mu|+1 -M)^{+}}
	\\
	\times\!
	\left(\!
	\| \psi_{N}^{(\beta)}  \Theta_{N}^{(\gamma +\mu-\nu)} D^{\alpha-\nu} \varphi_{N}^{(\delta+j+j_{1})} f\|_{\frac{p-1}{r}}
	\!+\!
	\| \psi_{N}^{(\beta)}  \Theta_{N}^{(\gamma +\mu-\nu)} D^{\alpha-\nu+j_{1}} \varphi_{N}^{(\delta+j)} f\|_{\frac{p-1}{r}}
	\!
	\right)
	\!
	\|v\|_{\frac{p-1}{r}}
	\\
	\leq
	A_{1}^{2S_{p}} B_{1}^{2T_{p}}  N^{2U_{p}}
	\times
	4^{n} \left(n^{2}(n+1)^{2}\right)  B_{1}^{-1}\left( B_{1}^{-1} + A_{1}^{-1}C_{1} \right).
	\end{multline}
\end{linenomath}
%
On the other hand we have
\begin{linenomath}
	\begin{multline}\label{T_10p-1}	
	\longsum[11]_{j,j_{1}=1 }^{n} 
	\longsum[69]_{|\alpha| -|\gamma| +3 \leq |\mu| \leq |\alpha| -|\gamma| + \lfloor \frac{n}{2} \rfloor+2}
	\longsum[7]_{\substack{\nu\leq \mu \\ \nu\leq \alpha}}
	\frac{1}{\mu!} 
	\binom{\mu}{\nu} \frac{\alpha!}{(\alpha -\nu)!}
	C_{4}^{|\mu|} \,N^{s(|\mu|+1 -M)^{+}}
	\\
	\times\!
	\left(\!
	\| \psi_{N}^{(\beta)}  \Theta_{N}^{(\gamma +\mu-\nu)} D^{\alpha-\nu} \varphi_{N}^{(\delta+j+j_{1})} f\|_{\frac{p-1}{r}}
	\!+\!
	\| \psi_{N}^{(\beta)}  \Theta_{N}^{(\gamma +\mu-\nu)} D^{\alpha-\nu+j_{1}} \varphi_{N}^{(\delta+j)} f\|_{\frac{p-1}{r}}
	\!
	\right)
	\!
	\|v\|_{\frac{p-1}{r}}
	\\
	\hspace{-23.9em}\leq
	A_{1}^{2S_{p}} B_{1}^{2T_{p}}  N^{2U_{p}}
	\\
	\times
	\left[
	2^{3n+1}n^{2}(n+1)^{2} \widetilde{C}_{2}^{ \lfloor \frac{n}{2} \rfloor+2}
	A_{1}^{-p-1} \left( A_{1}^{-1} \widetilde{C}_{1}\right)^{\sigma+1}
	\left(2\widetilde{C}_{2} B_{1}^{-1}\right)^{2rm+|\gamma|}
	\right].
	\end{multline}
\end{linenomath}
Concerning the last four terms they can be bounded in the following way
\begin{linenomath}
	\begin{multline}\label{T_11p}
	\left(
	\| \psi_{N}^{(\beta)} \Theta_{N}^{(\gamma)} D^{\alpha}  \varphi_{N}^{(\delta)} P^{k+1} u\|_{\frac{p-1}{r}}
	+
	C_{1}^{ \sigma+ 1}
	C_{2}^{ 2m+|\gamma|}\,N^{s(2m +|\gamma| +\sigma)} 
	\right)
	\!
	\|v\|_{\frac{p-1}{r}}
	\\
	\leq
	A_{1}^{2S_{p}} B_{1}^{2T_{p}}  N^{2U_{p}}
	\times
	\left[ B_{1}^{-2} + \left( C_{1}A_{1}^{-1}\right)^{\sigma +1}  \left( C_{2}B_{1}^{-r}\right)^{2|\alpha|} \left( A_{1}B_{1}\right)^{-p-1}  \right] ;
	\end{multline}
\end{linenomath}
and
\begin{linenomath}
	\begin{multline}\label{T_12p}
	\|v\|_{0}^{2} 
	+
	C_{1}^{2(\sigma+1)} C_{2}^{2(2 m + |\gamma|)}
	N^{2s(2 m +|\gamma|+ \sigma )}
	\\
	\leq
	A_{1}^{2S_{p}} B_{1}^{2T_{p}}  N^{2U_{p}}
	\times
	\left[ \left( A_{1} B_{1} N^{-2s}\right)^{-2p} + \left( C_{1}A_{1}^{-1}\right)^{2(\sigma +1)}  \left( C_{2}B_{1}^{-2r}\right)^{2|\alpha|} \left( A_{1}B_{1}\right)^{-2p}  \right].
	\end{multline}
\end{linenomath}
Summing up, enlarging $A_{1}$ and $B_{1}$ if necessary, the summand of second factor on the right
hand side of \eqref{T_1p}, \eqref{T_2p}, \eqref{T_3p}, \eqref{T_4p}, \eqref{T_5p}, \eqref{T_6p},
\eqref{T_7p}, \eqref{T_8p}, \eqref{T_9p-1}, \eqref{T_9p-2f}, \eqref{T_10p}, \eqref{T_10p-1}\eqref{T_11p} and \eqref{T_12p}
can be made smaller than $(2C_{2}^{2})^{-1}$, we conclude
\begin{linenomath}	
	\begin{align}\label{Est-p}
	\| v\|_{\frac{p}{r}} ^{2}
	\leq 
	A_{1}^{2S_{p}} B_{1}^{2T_{p}}  N^{2U_{p}},
	\end{align}
\end{linenomath}
that is we have obtained \eqref{Pr1-p_1-n}.\\
This concludes the proof of the Proposition \ref{Pr1-1/r-p/r}. 
\end{proof}
\vspace{1.5em}
\begin{theorem}\label{Th:Pr1-1_0}
	There exist positive constants $A_{1}$ and $B_{1}$ such that property \eqref{Pr1-1_0} in Proposition \ref{Pr1-1/r-p/r} is true.
\end{theorem}
\vspace{-2em}
\begin{proof}[Proof of Theorem \ref{Th:Pr1-1_0}]
	We use the induction on $m= |\alpha| -|\gamma|$.\\
	Cases $m=0$, $|\alpha|=|\gamma|$. Since 
	\begin{linenomath}
		\begin{align*}
		\hspace{4.5em}|\Theta_{N}^{(\gamma)}(\xi)| \leq  C_{0}^{|\gamma|+1} N^{\left(|\gamma|-M\right)^{+}} \left(1+ |\xi|\right)^{-|\gamma|},
		\end{align*} 
	\end{linenomath}
	we have
	\begin{linenomath}
		\begin{multline*}
		\|\psi_{N}^{(\beta)} \Theta_{N}^{(\gamma)} D^{\alpha} \varphi_{N}^{(\delta)} P^{k} u\|
		\leq C_{1}^{|\beta|+1} N^{(|\beta| -M)^{+}} \|\Theta_{N}^{(\gamma)} D^{\alpha} g\|
		\\
		= C_{1}^{|\beta|+1} N^{(|\beta| -M)^{+}} \|\Theta_{N}^{(\gamma)} (\xi)\xi^{\alpha} \widehat{g}\|_{L^{2}_{\xi}}
		\leq
		C_{1}^{|\beta|+1} N^{(|\beta| -M)^{+}} C_{0}^{|\gamma|+1} N^{\left(|\gamma|-M\right)^{+}} \| g\|,
		\end{multline*}
	\end{linenomath}
	$g= \varphi_{N}^{(\delta)} P^{k} u$. 
	By \eqref{Est_g} and taking advantage from the Remark \ref{Rk-1},
	we conclude that
	\begin{linenomath}
		\begin{multline*}
		\|\psi_{N}^{(\beta)} \Theta_{N}^{(\gamma)} D^{\alpha} \varphi_{N}^{(\delta)} P^{k} u\|
		\leq C_{1}^{|\beta|+1} C_{0}^{|\gamma|+1} C_{3}^{|\delta|+1} C_{2}^{2k+1}
		N^{(|\beta| +|\gamma| + |\delta|-M)^{+}} k^{2sk}
		\\
		\leq
		C_{4}^{|\beta|+ |\delta|+ 2k+1} C_{0}^{|\gamma|+1} N^{(|\beta| +|\gamma| + |\delta|+ 2k-M)^{+}}.
		\end{multline*}
	\end{linenomath}
	Since $M$ is a fixed constant depending only on $n$ and greater than
	one we obtain that there are two positive constants $A_{1}$ and $B_{1}$
	such that
	\begin{linenomath}
		\begin{align*}
		\| \psi_{N}^{(\beta)} \Theta_{N}^{(\gamma)} D^{\alpha} \varphi_{N}^{(\delta)} P^{k} u\| 
		\leq A_{1}^{|\sigma| +1} B_{1}^{|\gamma| +1} N^{s[|\gamma|+ \sigma]},
		\end{align*}
	\end{linenomath}
	that is \eqref{Pr1-1_0}.
	
	Now, we assume that
	\begin{linenomath}
		\begin{align*}
		\begin{cases}
		\| \psi_{N}^{(\beta)} \Theta_{N}^{(\gamma)} D^{\alpha} \varphi_{N}^{(\delta)} P^{k} u\|
		\leq 
		A_{1}^{|\sigma|  +1} B_{1}^{2rm+|\gamma| +1}  N^{s[rm+|\gamma|+ \sigma + p]},&\\
		\noalign{\vskip4pt}
		\text{for } 2rm+|\gamma|+\sigma \leq N \text{ and } |\gamma|\leq |\alpha|. &
		\end{cases}
		\end{align*}
	\end{linenomath}
	We have to show that  if it is true for $m$ less or equal to $m_{0}$, 
	fixed non negative integer, then it is true for $m=m_{0} + 1$.\\ 
	By the Proposition \ref{Pr1-1/r-p/r}
	\begin{linenomath}
		\begin{align*}
		\begin{cases}
		\| \psi_{N}^{(\beta)} \Theta_{N}^{(\gamma)} D^{\alpha} \varphi_{N}^{(\delta)} P^{k} u\|_{\frac{p}{r}} 
		\leq 
		A_{1}^{|\sigma| + p +1} B_{1}^{2rm+|\gamma| + p +1}  N^{s[rm+|\gamma|+ \sigma + p]},&\\
		\noalign{\vskip4pt}
		\text{for } 2r|\alpha|-(2r-1)|\gamma|+\sigma \leq N-2p, \, m=|\alpha|-|\gamma| \text{ and } |\gamma|\leq |\alpha|. &
		\end{cases}
		\end{align*}
	\end{linenomath}
	holds $\forall m \leq m_{0}$. In particular when $p=r$, we get
	\begin{linenomath}
		\begin{align*}
		\begin{cases}
		\| \psi_{N}^{(\beta)} \Theta_{N}^{(\gamma)} D^{\alpha} \varphi_{N}^{(\delta)} P^{k} u\|_{1} 
		\leq 
		A_{1}^{|\sigma| + r +1} B_{1}^{2rm_{0}+|\gamma| + r +1}  N^{s[rm_{0}+|\gamma|+ \sigma + r]},&\\
		\noalign{\vskip4pt}
		\text{for } 2rm_{0} + |\gamma|+\sigma \leq N-2r, \, m_{0}=|\alpha|-|\gamma| \text{ and } |\gamma|\leq |\alpha|. &
		\end{cases}
		\end{align*}
	\end{linenomath}
	Let $(\alpha, \beta, \gamma, \delta, k) $ be in $\mathbb{N}^{4n+1}$,  with $|\alpha| -|\gamma|= m_{0}+1$, such that
	$ 2r|\alpha|- (2r-1) |\gamma|+\sigma \leq N$, then $ 2r(m_{0}+1) + |\gamma|+\sigma \leq N$; i.e. that
	$2rm_{0} + |\gamma|+\sigma \leq N-2r$. So $\alpha= \alpha_{0} +e_{j}$, $e_{j}= (0,\dots,0,1,0,\dots,0)$, with
	$|\alpha_{0}|-|\gamma| =m_{0}$. Since
	\begin{linenomath}
		$$
		2rm_{0} + |\gamma|+\sigma \leq 2rm+ |\gamma|+\sigma -2r\leq N-2r,
		$$
	\end{linenomath} 
	by inductive hypothesis, \eqref{Pr1-1_0}, in Proposition \ref{Pr1-1/r-p/r}, is true for  $(\alpha_{0}, \beta, \gamma, \delta, k) $
	and consequently, by the Proposition \ref{Pr1-1/r-p/r}, \eqref{Pr1-1_p}, with $p=r$, it is true for  $(\alpha_{0}, \beta, \gamma, \delta, k) $.
	Now, we have
	\begin{linenomath}
		\begin{multline*}
		\|\psi_{N}^{(\beta)} \Theta_{N}^{(\gamma)} D^{\alpha} \varphi_{N}^{(\delta)} P^{k} u\|
		=
		\|\psi_{N}^{(\beta)} \Theta_{N}^{(\gamma)} D_{j}D^{\alpha_{0}} \varphi_{N}^{(\delta)} P^{k} u\|
		\\
		\leq 
		\|\psi_{N}^{(\beta +j)} \Theta_{N}^{(\gamma)} D^{\alpha_{0}} \varphi_{N}^{(\delta)} P^{k} u\|
		+
		\|\psi_{N}^{(\beta)} \Theta_{N}^{(\gamma)}D^{\alpha_{0}} \varphi_{N}^{(\delta)} P^{k} u\|_{1}
		=
		I_{1} +I_{2}.
		\end{multline*}
	\end{linenomath}
	We remark that
	\begin{linenomath}
		$$
		2r|\alpha_{0}| - (2r-1)|\gamma| + \sigma+1 = 2r(|\alpha|-1) - (2r-1)|\gamma| + \sigma+1 \leq 2r|\alpha| - (2r-1)|\gamma| + \sigma,
		$$
	\end{linenomath}
	as $1 - 2r$ is negative, $r\geq 2$. Since $|\alpha_{0}| - |\gamma| =m_{0}$ and $(\alpha_{0}, \beta +e_{j}, \gamma, \delta, k) $ satisfies the condition
	in \eqref{Pr1-1_0}, then \eqref{Pr1-1_0} is true for $(\alpha_{0}, \beta + e_{j}, \gamma, \delta, k) $.
	So
	\begin{linenomath}
		\begin{align*}
		I_{1} \leq A_{1}^{\sigma +1+ 1} B_{1}^{2rm_{0} +|\gamma| +1} N^{s\left[ rm_{0} +|\gamma| + \sigma +1\right]},
		\end{align*}
	\end{linenomath}
	and
	\begin{linenomath}
		\begin{align*}
		I_{2} \leq A_{1}^{\sigma +r+ 1} B_{1}^{2rm_{0} +|\gamma|+r +1} N^{s\left[ rm_{0} +|\gamma| + \sigma +r\right]},
		\end{align*}
	\end{linenomath}
	here we use \eqref{Pr1-1_p} in Proposition \ref{Pr1-1/r-p/r}, with $p=r$, indeed
	$2r |\alpha_{0}| -(2r-1)|\gamma| +\sigma = 2r |\alpha| -(2r-1)|\gamma| +\sigma -2r \leq N-2r$.\\
	So finely we have
	\begin{linenomath}
		\begin{multline*}
		\|\psi_{N}^{(\beta)} \Theta_{N}^{(\gamma)} D^{\alpha} \varphi_{N}^{(\delta)} P^{k} u\|
		\leq
		A_{1}^{\sigma +1+ 1} B_{1}^{2r|\alpha_{0}| -(2r- 1)|\gamma| +1} N^{s\left[ r\alpha_{0} +(r-1)|\gamma| + \sigma +1\right]}
		\\
		\hspace{11em}+
		A_{1}^{\sigma +r+ 1} B_{1}^{2r|\alpha_{0}| -(2r- 1)|\gamma|+r +1} N^{s\left[ r\alpha_{0} +(r-1)|\gamma| + \sigma +r \right]}
		\\
		\leq 
		A_{1}^{\sigma + 1} B_{1}^{2r(|\alpha_{0}|+1) -(2r- 1)|\gamma| +1} N^{s\left[ r(\alpha_{0}+1) +(r-1)|\gamma| + \sigma  \right]}
		\times
		\left( A_{1} B_{1}^{-2r} N^{s(1-r)}+ A_{1}^{r}B_{1}^{-r}\right).
		\end{multline*}
	\end{linenomath}
	Since $r\geq 2$, $N^{s(1-r)} < 1$; moreover taking $A_{1}$ and $B_{1}$ large enough,  with $B_{1}$ large compared to $A_{1}$,
	we have
	\begin{linenomath}
		\begin{align*}
		\left( A_{1} B_{1}^{-2r} N^{s(1-r)}+ A_{1}^{r}B_{1}^{-r}\right) < 1.
		\end{align*} 
	\end{linenomath}
	We conclude that
	\begin{linenomath}
		\begin{align*}
		\|\psi_{N}^{(\beta)} \Theta_{N}^{(\gamma)} D^{\alpha} \varphi_{N}^{(\delta)} P^{k} u\|
		\leq 
		A_{1}^{\sigma + 1} B_{1}^{2rm  + |\gamma| +1} N^{s\left[ rm + |\gamma| + \sigma  \right]},
		\end{align*}
	\end{linenomath}
	where $m=m_{0} +1$. By induction we have obtained that \eqref{Pr1-1_0} is true for all $m\in \mathbb{N}$
	and $(\alpha, \beta, \gamma, \delta, k) \in \mathbb{N}^{4n+1}$ such that
	$2rm+ |\gamma|+\sigma \leq N$, where $\sigma = |\beta| + |\delta| + 2k$, $m= |\alpha| - |\gamma|$ and $|\gamma| \leq |\alpha| $.\\
	This conclude the proof of the theorem.
\end{proof}
\begin{remark}
If we take $|\gamma|=0$ and $\sigma = 0$ ($|\beta| = |\delta| = k =0$), \eqref{Pr1-1_0} gives
\begin{linenomath}
	\begin{align}\label{M_Rg_u}
	\|\psi_{N} \Theta_{N} D^{\alpha} \varphi_{N} u\|
	\leq 
	A_{1}B_{1}^{2r|\alpha| +1} N^{sr|\alpha|}.
	\end{align}
\end{linenomath}
\end{remark}
\medskip
\begin{corollary}\label{C1}
	Let $\psi_{N}$, $\Theta_{N}$ and $\varphi_{N}$ be as above.
	Then the following estimate holds
	\begin{linenomath}
		\begin{align}\label{Est_Wf}
		\| \Theta_{N} D^{\alpha} \psi_{N} \varphi_{N} u\|
		\leq 
		A^{|\alpha| +1} N^{sr|\alpha|},
		\end{align}
	\end{linenomath}
	where the constant $A$ is independent of $N$ and $\alpha$.
\end{corollary}
\vspace{-2em}
\begin{proof}[Proof of Corollary \ref{C1}]
We observe that
\begin{linenomath}
	\begin{align*}
	\Theta_{N} D^{\alpha} \psi_{N} \varphi_{N} u
	=\psi_{N} \Theta_{N} D^{\alpha} \varphi_{N} u
	+
	[\Theta_{N} D^{\alpha}, \psi_{N}] \varphi_{N} u,
   \end{align*}
\end{linenomath}
where	
\begin{linenomath}
		\begin{multline*}
		[\Theta_{N} D^{\alpha}, \psi_{N}] \varphi_{N} u
		=\longsum[23]_{1\leq |\mu|\leq |\alpha|-1} \frac{1}{i^{|\mu|}\mu!} \psi_{N}^{(\mu)} \left(\Theta_{N} D^{\alpha} \right)^{(\mu)} \varphi_{N} u
		\\
		+
		\mathscr{R}_{|\alpha|}\left( \left[\Theta_{N} D^{\alpha},\psi_{N} \right]  \right)\varphi_{N} u
		\\
		=\longsum[23]_{\substack{ 1\leq |\mu|\leq |\alpha|-1\\ \nu\leq \mu,\, \nu \leq \alpha }}
		\frac{\alpha!}{i^{|\mu|}\nu!(\mu-\nu)!(\alpha-\nu)!} \psi_{N}^{(\mu)} \Theta_{N}^{(\mu-\nu)} D^{\alpha-\nu} \varphi_{N} u
		\\
		+
		\mathscr{R}_{|\alpha|}\left( \left[\Theta_{N} D^{\alpha},\psi_{N} \right]  \right)\varphi_{N} u.
		\end{multline*}
\end{linenomath}
So
\begin{linenomath}
		\begin{multline}\label{Est_C1}
		\| \Theta_{N} D^{\alpha} \psi_{N} \varphi_{N} u\|
		\leq
		\| \psi_{N} \Theta_{N} D^{\alpha} \varphi_{N} u\|
		+
		\| [\Theta_{N} D^{\alpha}, \psi_{N}] \varphi_{N} u\|
		\\
		\leq
		\| \psi_{N} \Theta_{N} D^{\alpha} \varphi_{N} u\|
		+
		\longsum[23]_{\substack{ 1\leq |\mu|\leq |\alpha|-1\\ \nu\leq \mu,\, \nu \leq \alpha }}
		\frac{\alpha!}{\nu!(\mu-\nu)!(\alpha-\nu)!} \| \psi_{N}^{(\mu)} \Theta_{N}^{(\mu-\nu)} D^{\alpha-\nu} \varphi_{N} u\|
		\\
		+
		\| \mathscr{R}_{|\alpha|}\left( \left[\Theta_{N} D^{\alpha},\psi_{N} \right]  \right) \varphi_{N} u\|.
		\end{multline}
\end{linenomath}
By Theorem \ref{Th:Pr1-1_0} we have
\begin{linenomath}
		\begin{multline*}
		\| \psi_{N}^{(\mu)} \Theta_{N}^{(\mu-\nu)} D^{\alpha-\nu} \varphi_{N} u\|
		\leq
		A_{1}^{|\mu|+1} B_{1}^{2r(|\alpha| -|\mu|) + |\mu|- |\nu| +1} N^{sr|\alpha|} N^{-s[(r-2)|\mu|+|\nu|]} , 
		\end{multline*}
\end{linenomath}
moreover since $\frac{\alpha!}{\left(\mu-\nu\right)! \nu!\left(\alpha-\nu\right)!} \leq N^{|\nu|}$, $B_{1}$
is strictly greater than $2$, $r\geq 2$ and $s\geq 1$, we obtain
\begin{linenomath}
		\begin{multline*}
		\longsum[23]_{\substack{ 1\leq |\mu|\leq |\alpha|-1\\ \nu\leq \mu,\, \nu \leq \alpha }}
		\frac{\alpha!}{\nu!(\mu-\nu)!(\alpha-\nu)!} \| \psi_{N}^{(\mu)} \Theta_{N}^{(\mu-\nu)} D^{\alpha-\nu} \varphi_{N} u\|
		\\
		\leq
		A_{1}^{|\alpha|+1} B_{1}^{2r|\alpha| +1} N^{sr|\alpha|}
		\longsum[23]_{\substack{ 1\leq |\mu|\leq |\alpha|-1 }}
		\longsum[10]_{\substack{ \nu\leq \mu\\ \nu \leq \alpha }}
		B_{1}^{-|\mu|(2r-1)}B_{1}^{- |\nu|}  N^{-s[(r-2)|\mu|} N^{-(s-1)|\nu|} 
		\\
		\leq
		C_{1}^{|\alpha|+1}  N^{sr|\alpha|}
		\longsum[6]_{\mu_{1} = 0}^{\infty} \left(\frac{1}{2^{2r-1}}\right)^{\mu_{1}}\cdots 
		\longsum[6]_{\mu_{n} =0}^{\infty}  \left(\frac{1}{2^{2r-1}}\right)^{\mu_{n}}
		\longsum[6]_{\nu_{1} = 0} \left(\frac{1}{2}\right)^{\nu_{1}}\cdots 
		\longsum[6]_{\nu_{n} =0 }^{\infty} \left(\frac{1}{2}\right)^{\nu_{n}} 
		\\
		\leq
		C_{2}^{|\alpha|+1}  N^{sr|\alpha|}.\hspace{27em}
		\end{multline*}
\end{linenomath}
Using the same strategy adopted in the proof of Lemma \ref{Rem0},
see also the estimate of term $I_{2,4}$, \ref{Est_I24_0},
the last term on the right hand side of \eqref{Est_C1} can be estimated as follow
\begin{linenomath}
		\begin{align*}
		\| \mathscr{R}_{|\alpha|}\left( \left[\Theta_{N} D^{\alpha},\psi_{N} \right]  \right) \varphi_{N} u\|
		\leq 
		C_{3}^{|\alpha| +1} N^{s2|\alpha|},
		\end{align*}
\end{linenomath}
where $C_{3}$ is a suitable positive constant independent of $\alpha$. Here we use that  $M=3(n+6)$.
	
Since $r\geq 2$, by the above consideration and the estimate \eqref{M_Rg_u}
we obtain \eqref{Est_Wf}. This concludes the proof.
%
\end{proof}
Recall that, as pointed out in \cite{H_Book-1} page 283 (Lemma 8.4.4), the sequence $u_{N}$ in the Definition \ref{D_WF_s},
can always be chosen as a product of $u$ and a suitable cutoff functions,
that is we set $u_{N}= \psi_{N}\varphi_{N} u= \varphi_{N} u$, $\psi_{N}$ equals $1$ on support of $\varphi_{N}$.
Recalling that the sequence $\Theta_{N}$ is associated to the couple $(\Gamma_{0},\Gamma_{1})$,
by the above Corollary we conclude that taking $N=2r|\alpha|$ for every $\xi$ in $\widetilde{\Gamma}$, $\widetilde{\Gamma} \Subset \Gamma_{0}$,
\eqref{WF_s} is satisfied, i.e. $(x_{0},\xi_{0})\notin WF_{rs}(u)$. 
\begin{remark}
If $(x_{0},\xi_{0}) \notin WF_{rs}(u)$, that is \eqref{WF_s} is satisfied, then \eqref{M_Rg_u} holds.
\end{remark}
In view of Proposition \ref{Pr1-1/r-p/r}, Theorem \ref{Th:Pr1-1_0} and Corollary \ref{C1} we obtain
the Theorem \ref{M-Th}.
%
%
%
\section{Appendix}
\renewcommand{\theequation}{\thesection.\arabic{equation}}
\setcounter{equation}{0} \setcounter{theorem}{0}
\setcounter{proposition}{0} \setcounter{lemma}{0}
\setcounter{corollary}{0} \setcounter{definition}{0}
\setcounter{remark}{1}
Even if known, (see Lemma  2.2.1 in \cite{bcr_1986}), and in order to make
this paper as self contained as possible, we show a strategy in order to
construct the Ehrenpreis-Andersson cutoff symbols.
Let $\xi_{0} \in \mathbb{R}^{n}\setminus \{0\}$ and $r \in \mathbb{R}^{+}$ we denote by
\begin{linenomath}
	\begin{align}\label{Gamma-r}
	\Gamma_{\xi_{0}, r} = 
	\Big\{ \xi \in \mathbb{R}^{n}\setminus\{0 \} :  \left| \frac{\xi}{|\xi|} - \frac{\xi_{0}}{|\xi_{0}|} \right| < r \Big\}
	\end{align}
\end{linenomath}
a conic neighborhood of $\xi_{0}$ of size $r$.\\
We recall the classical construction of the  Ehrenpreis-H\"or\-man\-der cut-off functions.
Let $x_{0} \in \mathbb{R}^{n}$ and $\Sigma$ a neighborhood of $x_{0}$ 
then there is a constant $C_{0} > 0 $ depending only by the dimension of the ambient space, $n$,
such that given any positive $r $, a non null integer $M$ and any $ N \in \mathbb{Z}_{+}$,
there is a sequence  $ \varphi_{N}$ of smooth functions in $\mathbb{R}^{n}$, having the following properties:
\vspace*{-1.5em}
\begin{enumerate}
	\item [i)] $\varphi_{N} \equiv 1$ on $\Sigma$, $\varphi_{N}(x) = 0$ if $\dist\left( x; \Sigma\right) >\left(1+\frac{M}{2}\right) r$
	and $0\leq \varphi_{N}(x) \leq 1 $ for every $x$;
	\item[ii)] the following estimate holds
	\begin{linenomath}
		\begin{align}\label{E-H_cut-off}
		| D^{\alpha} \varphi_{N}| \leq \left(\frac{C_{0}}{2r}\right)^{|\alpha|} N^{(|\alpha|-M)^{+}},
		\qquad \text{for all } \alpha \in \mathbb{Z}^{n}_{+} \text{ such that } |\alpha | \leq N.
		\end{align}
	\end{linenomath}
\end{enumerate}
We choose a  function $\psi \in \mathscr{D}(\mathbb{R}^{n})$ with support in
$\mathscr{B}_{1/4}(0) \doteq \lbrace x\in \mathbb{R}^{n}\,: \, |x| \leq 1/4 \rbrace$ such that $\psi \geq 0$ and
$\int \psi \, dx =1$. For every $\delta > 0$ we write $\psi_{\delta}(x) = \delta^{-n} \psi\left(\frac{x}{\delta}\right)$.
Let $\chi$ be the characteristic function of the set $\lbrace x\in \mathbb{R}^{n}\,:\, \dist\left(x; \Sigma\right) < \frac{r}{2}\rbrace$.
We set
\begin{linenomath}
	\begin{align}\label{const_E-H_fu}
	\varphi_{N} = 
	\underbrace{\psi_{\frac{2r}{N}} * \psi_{\frac{2r}{N}} *\, \cdots \, *\,\psi_{\frac{2r}{N}} }_{N-\text{times}} *\, 
	\underbrace{\psi_{2r\phantom{\frac{\phantom{2r}}{\phantom{N}}}} \hspace*{-0.9em}*\, \cdots \, *\,\psi_{2r}}_{M-\text{times}} *\,\chi.
	\end{align}
\end{linenomath}
Since the support of a convolution is contained in the vector sum of the
supports of the factors in the convolution the sequence $\varphi_{N}$ satisfies the properties in $i)$.\\
Let $\alpha \in \mathbb{Z}_{+}^{n}$ with $M<|\alpha | \leq N$, we have
\begin{linenomath}
	\begin{multline*}
	D^{\alpha} \varphi_{N}  = D_{x_{j_{1}}} \!\!\!\cdots D_{x_{j_{M}}}D_{x_{j_{M+1}}}\!\!\!\cdots D_{x_{j_{|\alpha|}}} \varphi_{N}
	\\
	= \underbrace{\left(D_{x_{j_{M+1}}}\psi_{\frac{2r}{N}}\right) * \cdots *\left( D_{x_{j_{|\alpha|}}}\psi_{\frac{2r}{N}} \right) }_{(|\alpha|-M)-\text{times}} *
	\underbrace{\psi_{\frac{2r}{N}} * \cdots * \psi_{\frac{2r}{N}} }_{N-|\alpha|-\text{times}} 
	\\
	*\, \left(D_{x_{j_{1}}}\psi_{2r} \right)*\, \cdots \, *\,\left(D_{x_{j_{M}}}\psi_{2r}\right)*\,  \chi, 
	\end{multline*}
\end{linenomath}
where $j_{1}, \dots, j_{|\alpha|} $ belong to $\lbrace 1, \dots, n\rbrace$. Via the H\"older inequality we obtain
\begin{linenomath}
	\begin{multline*}
	\| D^{\alpha} \varphi_{N} \|_{\infty}\leq \prod_{i=M+1}^{|\alpha|}\|D_{x_{i}}\psi_{\frac{2r}{N}}\|_{L^{1}}
	\prod_{i=N-|\alpha|}^{N}\|\psi_{\frac{2r}{N}}\|_{L^{1}}
	\prod_{\ell=1}^{M}\|D_{x_{\ell}}\psi_{2r}\|_{L^{1}}\|\chi\|_{\infty}
	\\
	\leq \left(\frac{C_{0}}{2r}\right)^{|\alpha|} N^{(|\alpha|-M)},
	\end{multline*}
\end{linenomath}
where $C_{0} = \sup_{1\leq i\leq n} \|D_{x_{i}}\psi\|_{L^{1}}$.\\
\noindent
We set $\Theta_{0,N}= \varphi_{N}$ and $\Sigma$ the ball of radius $1/2$. We point out that the sequence $\Theta_{0,N}$
is such that $\Theta_{0,N}(\zeta)=1$ when $|\zeta| \leq1/2$ and $\Theta_{0,N}(\zeta) = 0$ when $|\zeta| \geq 1$, $\zeta\in \mathbb{R}^{n}$.
Let $\xi_{0} \in \mathbb{R}^{n}\setminus \{0\}$.
We set
\begin{linenomath}
	\begin{align}\label{ThetaN}
	\Theta_{N}(\xi) = \left( 1- \Theta_{0,N}\right)\left(\frac{\xi}{N}\right)\Theta_{0,N}\left(r^{-1}\left(\frac{\xi}{|\xi|} - \frac{\xi_{0}}{|\xi_{0}|}\right)\right),
	\end{align} 
\end{linenomath}
where $r\in\mathbb{R}^{+}$.
$\Theta_{N}(\xi) $ is supported in 
\begin{linenomath}
	\begin{align*}
	\Gamma_{\xi_{0}, \frac{N}{2} , r} =\Gamma_{\xi_{0},  r}\cap \Big\{ \xi \in \mathbb{R}^{n}\setminus\{0 \} :\, |\xi|\geq \frac{N}{2} \Big\},
	\end{align*} 
\end{linenomath}
and $\Theta_{N}(\xi) = 1$ in 
\begin{linenomath}
	\begin{align*}
	\Gamma_{\xi_{0}, N , \frac{r}{2}} = \Gamma_{\xi_{0},  \frac{r}{2}} \cap \Big\{ \xi \in \mathbb{R}^{n}\setminus\{0 \} :\, |\xi|\geq N \Big\},
	\end{align*} 
\end{linenomath}
where $\Gamma_{\xi_{0},  *}$ are as in (\ref{Gamma-r}); we remark that $\Gamma_{\xi_{0}, N,  r} \Subset \Gamma_{\xi_{0}, N,  \frac{r}{2}}$ in the sens of cones.\\
\noindent
We want to show that there is a positive constant $C$ such that
\begin{linenomath}
	\begin{align}\label{ThetaNDer}
	|\Theta_{N}^{(\alpha)}(\xi)| \leq  C^{|\alpha|+1} N^{\left(|\alpha|-M\right)^{+}} \left(1+ |\xi|\right)^{-|\alpha|} ,
	\end{align} 
\end{linenomath}
for every $\alpha \in \mathbb{Z}_{+}^{n}$ with $|\alpha| \leq N$.\\
We have
\begin{linenomath}
	\begin{align}\label{ThetaNA}
	\Theta_{N}^{(\alpha)}(\xi) = \sum_{ \beta \leq \alpha} \binom{\alpha}{\beta} 
	\left( 1- \Theta_{0,N}\right)^{(\alpha -\beta)}\left(\frac{\xi}{N}\right)\partial_{\xi}^{\beta}\Theta_{0,N}\left(r^{-1}\left(\frac{\xi}{|\xi|} - \frac{\xi_{0}}{|\xi_{0}|}\right)\right).
	\end{align} 
\end{linenomath}
In order to estimate the absolute value of second factor we use the multivariate Faa di Bruno formula. For completeness we recall
the result in \cite{BrunoM}, Theorem 2.1 and Remark 2.2. Let $h(\xi)= f(g_{1}(\xi),\dots,g_{n}(\xi))$,  then
\begin{linenomath}
	\begin{align}\label{MFaaBruno}
	\partial_{\xi}^{\beta}h(\xi)= \longsum[20]_{ 1 \leq |\gamma| \leq |\beta|} f^{(\gamma)}(\eta)
	\longsum[10]_{ p(\beta,\gamma)} \beta! \prod_{j=1}^{|\beta|} 
	\frac{ \left( \partial_{\xi}^{\boldsymbol{\ell}_{j}}g\right)^{\boldsymbol{k}_{j}}}{ \boldsymbol{k}_{j}! \left(\boldsymbol{\ell}_{j}!\right)^{|\boldsymbol{k}_{j}|} }
	\end{align} 
\end{linenomath}
where $\eta = (g_{1}(\xi),\dots,g_{n}(\xi))$, $\boldsymbol{k}_{j},\, \boldsymbol{\ell}_{j} \in \mathbb{N}^{n}$,
$\boldsymbol{k}_{j}= \left( k_{j,1},\dots,k_{j,n}\right)$, $\boldsymbol{\ell}_{j}= \left( \ell_{j,1},\dots,\ell_{j,n}\right)$,
$\boldsymbol{k}_{j}!=  k_{j,1}!\cdot \ldots \cdot k_{j,n}!$, $\boldsymbol{\ell}_{j}!=  \ell_{j,1}!\cdot \ldots \cdot \ell_{j,n}!$,
$\left( \partial_{\xi}^{\boldsymbol{\ell}_{j}}g\right)^{\boldsymbol{k}_{j}}=
\left( \partial_{\xi}^{\boldsymbol{\ell}_{j}}g_{1}\right)^{k_{j,1}} \cdot \ldots \cdot \left( \partial_{\xi}^{\boldsymbol{\ell}_{j}}g_{n}\right)^{k_{j,n}}$ 
and
\begin{linenomath}
	\begin{align}\label{pFaaBruno}
	p(\beta,\gamma) = 
	&\big\lbrace \left( \boldsymbol{k}_{1},\dots, \boldsymbol{k}_{|\beta|}; \boldsymbol{\ell}_{1}, \dots, \boldsymbol{\ell}_{|\beta|}\right)\,:\,
	\text{ for some } 1\leq s \leq |\beta|, 
	\\[8pt]
	\nonumber	
	&\quad \boldsymbol{k}_{i}=0 \text{ and } \boldsymbol{\ell}_{i}=0 
	\text{ for } 1\leq i \leq |\beta| -s;\, |\boldsymbol{k}_{i}|> 0 \text{ for } |\beta|-s+1 \leq i \leq |\beta|;
	\\
	\nonumber
	&\quad \text{and } 0 \prec \boldsymbol{\ell}_{|\beta|-s+1} \prec \cdots \prec \boldsymbol{\ell}_{|\beta|} \text{ are such that }
	\sum_{i=1}^{|\beta|}  \boldsymbol{k}_{i} =\gamma, \, \sum_{i=1}^{|\beta|}  |\boldsymbol{k}_{i}| \boldsymbol{\ell}_{i}= \beta
	\big\rbrace.
	\end{align} 
\end{linenomath}
Given $\mu= \left( \mu_{1}, \ldots, \mu_{n}\right)$ and $\nu= \left( \nu_{1}, \ldots, \nu_{n} \right)$ one writes $\mu \prec \nu$
if one of the following sentences holds:
\vspace*{-1.7em}
\begin{itemize}
	\item[$\bullet$] $ | \mu| < |\nu|$;
	\item[$\bullet$] $ | \mu| = |\nu|$ and $\mu_{1} < \nu_{1}$;
	\item[$\bullet$] $ | \mu| = |\nu|$, $\mu_{1} = \nu_{1}, \ldots, \mu_{r} = \nu_{r} $ and $\mu_{r+1} < \nu_{r+1}$ for some $1\leq r < n$.
\end{itemize}
We remark that taking an homogeneous function $v(\xi)$ of order $p$ and analytic outside $0$ since $|v^{(\mu)}(\eta)| \leq C_{1}^{|\mu|+1} \mu!$,
for every $\eta \in \mathbb{S}^{n-1}$, where $\mu \in \mathbb{Z}^{n}_{+}$ and $C_{1}$ is a positive constant independent of $\mu$, 
and moreover since $v^{(\mu)}$ is an  homogeneous function of order $p-|\mu|$ 
we have $|v^{(\mu)}(\xi)| \leq C_{1}^{|\mu|+1} \mu!|\xi|^{p-|\mu|}$ for $|\xi|>1$.\\
Let $f= \Theta_{0,N}$ and $g_{i}= r^{-1}(\xi_{i}/|\xi|-\xi_{0,i}/|\xi_{0}|)$ in  (\ref{MFaaBruno}),
we point out that $g_{i}$ are homogeneous functions of order zero and analytic outside $0$,  then 
\begin{linenomath}
	\begin{multline*}
	\left| \left( \partial_{\xi}^{\boldsymbol{\ell}_{j}}g\right)^{\boldsymbol{k}_{j}} \right|
	= \prod_{i=1}^{n}\left|  \left( \partial_{\xi}^{\boldsymbol{\ell}_{j}}g_{i}\right)^{k_{j,i}} \right|
	\leq r^{-|\gamma|} \prod_{i=1}^{n} C_{0}^{|\boldsymbol{\ell}_{j}| + 1} (\boldsymbol{\ell}_{j} !)^{k_{j,i}} |\xi|^{-|\boldsymbol{\ell}_{j}|k_{j,i}}
	\\
	\leq r^{-|\gamma|} C_{0}^{|\boldsymbol{\ell}_{j}| | \boldsymbol{k}_{j} |+ 1} (\boldsymbol{\ell}_{j} !)^{| \boldsymbol{k}_{j} |} 
	|\xi|^{-|\boldsymbol{\ell}_{j}| | \boldsymbol{k}_{j} |};
	\end{multline*} 
\end{linenomath}
we obtain
\begin{linenomath}
	\begin{multline}\label{EstTheta0N}
	\left| \partial_{\xi}^{\beta}\Theta_{0,N}\left(r^{-1}\left(\frac{\xi}{|\xi|} - \frac{\xi_{0}}{|\xi_{0}|}\right)\right)\right|
	\\
	\leq 
	\longsum[20]_{ 1 \leq |\gamma| \leq |\beta|} C_{1}^{|\gamma|+ 1} N^{|\gamma|}   r^{-|\gamma|} 
	\longsum[10]_{ p(\beta,\gamma)} \beta! \prod_{j=1}^{|\beta|} 
	\frac{ C_{0}^{|\boldsymbol{\ell}_{j}| | \boldsymbol{k}_{j} |+ 1}  |\xi|^{|\boldsymbol{\ell}_{j}| | \boldsymbol{k}_{j} |}  }{ \boldsymbol{k}_{j}!  }
	\\
	\leq r \left(\frac{4C_{1}C_{0}}{r}\right)^{|\beta| +1} N^{|\beta|} |\xi|^{-|\beta|}
	\longsum[20]_{ 1 \leq |\gamma| \leq |\beta|}  1
	\longsum[10]_{ p(\beta,\gamma)} 1
	\\
	\leq C_{2}^{|\beta| +1} N^{\left(|\beta|-M\right)^{+}} |\xi|^{-|\beta|},
	\end{multline} 
\end{linenomath}
where $C_{2}$ is a suitable positive constant independent of $\beta$.
In order to obtain the above inequality we use
that $ \prod_{j=1}^{|\beta|} |\xi|^{|\boldsymbol{\ell}_{j}| | \boldsymbol{k}_{j} |} = |\xi|^{|\beta|} $, $ \prod_{j=1}^{|\beta|} \frac{ \beta! }{ \boldsymbol{k}_{j}!  }
\leq 2^{|\beta|+ |\gamma| } (\beta- \gamma)! \leq 4^{|\beta|} N^{|\beta|-|\gamma|}$, the cardinality of the set of $\gamma \in \mathbb{Z}_{+}^{n}$ such that
$|\gamma| \leq |\beta|$ is $\binom{|\beta| + n}{|\beta|}$,  and that $p(\beta,\gamma) $ can be seen as the subset of 
$( \boldsymbol{k},  \boldsymbol{\ell})\in \mathbb{Z}_{+}^{2n|\beta|} $ such that $ | ( \boldsymbol{k},  \boldsymbol{\ell})|\leq 2|\beta|$.\\ 
In order to obtain (\ref{ThetaNDer}) we distinguish two cases. 
When $\beta \neq \alpha$, we have  
\begin{linenomath}
	\begin{multline*}
	\left(1+ |\xi|\right)^{|\alpha|} 	\left | \left( 1- \Theta_{0,N}\right)^{(\alpha -\beta)}\left(\frac{\xi}{N}\right)\right| 
	\left| \partial_{\xi}^{\beta}\Theta_{0,N}\left(r^{-1}\left(\frac{\xi}{|\xi|} - \frac{\xi_{0}}{|\xi_{0}|}\right)\right)\right|
	\\
	\leq
	2^{|\alpha|}C_{1}^{|\alpha|-|\beta| +1} C_{2}^{|\beta|+1} N^{\left(|\beta|-M\right)^{+} } |\xi|^{|\alpha|-|\beta|}
	\leq \left( 2C_{1}C_{2}\right)^{|\alpha|+1} N^{\left(|\alpha|-M \right)^{+}} ,
	\end{multline*} 
\end{linenomath}
where we take advantage from the fact that since $\beta$ is less than $\alpha$ then we have that $2^{-1}\leq  |\xi| N^{-1}\leq 1$.\\
When $\beta =\alpha$ by (\ref{EstTheta0N}) we obtain
\begin{linenomath}
	\begin{multline*}
	\left(1+ |\xi|\right)^{|\alpha|} 	\left | \left( 1- \Theta_{0,N}\right)\left(\frac{\xi}{N}\right)\right| 
	\left| \partial_{\xi}^{\alpha}\Theta_{0,N}\left(r^{-1}\left(\frac{\xi}{|\xi|} - \frac{\xi_{0}}{|\xi_{0}|}\right)\right)\right|
	\\
	\leq
	2^{|\alpha|}C_{1} C_{2}^{|\alpha|+1} N^{\left(|\alpha|-M \right)^{+}} .
	\end{multline*} 
\end{linenomath}
By (\ref{ThetaNA}) we conclude that there is a suitable positive constant $C$
independent of $\alpha$ such that (\ref{ThetaNDer}) holds.\\
We remark that by (\ref{MFaaBruno}) we have that if $\varphi_{N}$ is an  Ehrenpreis
sequence
and $g $ is an analytic function, then the sequence $\psi_{N}= \varphi_{N} \circ g $ is still an  Ehrenpreis sequence.\\
Summing up we have
\begin{lemma}\label{EA-Cutoff}
Let $\xi_{0} \in \mathbb{R}^{n}\setminus\{ 0 \} $, $\Gamma_{\xi_{0},\frac{r}{2}}$ a conic neighborhood of $\xi_{0}$, $r>0$,
and $M$ be positive integer. For every non zero integer $N$, there is a smooth function $\Theta_{N}(\xi)$ equal to $1$
in $\Gamma_{\xi_{0},\frac{r}{2}} \cap \left\{ |\xi| >N\right\}$ and supported in $\Gamma_{\xi_{0},r} \cap \left\{ |\xi| >N/2\right\}$,
$\Gamma_{\xi_{0},\frac{r}{2}} \Subset \Gamma_{\xi_{0},r}$, such that
\begin{linenomath}
	\begin{align}\label{ThetaN-AH}
	\hspace{4.5em}|\Theta_{N}^{(\alpha)}(\xi)| \leq  C^{|\alpha|+1} N^{\left(|\alpha|-M\right)^{+}} \left(1+ |\xi|\right)^{-|\alpha|} \qquad \text{ if }  |\alpha|\leq N.
	\end{align} 
\end{linenomath}
Where $C$ depends only on $n$ and $r$.
\end{lemma}
%
\section*{Declarations}
The authors state that there is no conflict of interest.

%
%

%

%
\end{document}